\magnification=1200

\loadmsam
\loadmsbm
\loadeufm
\loadeusm
\UseAMSsymbols

\font\BIGtitle=cmr10 scaled\magstep3
\font\bigtitle=cmr10 scaled\magstep1
\font\boldsectionfont=cmb10 scaled\magstep1
\font\section=cmsy10 scaled\magstep1

\def\scr#1{{\fam\eusmfam\relax#1}}

\def\scrA{{\scr A}}
\def\scrB{{\scr B}}
\def\scrC{{\scr C}}

\def\scrF{{\scr F}}
\def\scrG{{\scr G}}
\def\scrH{{\scr H}}
\def\scrK{{\scr K}}
\def\scrJ{{\scr J}}
\def\scrL{{\scr L}}
\def\scrM{{\scr M}}
\def\scrN{{\scr N}}
\def\scrO{{\scr O}}
\def\scrP{{\scr P}}

\def\scrS{{\scr S}}
\def\scrU{{\scr U}}
\def\scrR{{\scr R}}
\def\scrT{{\scr T}}
\def\gr#1{{\fam\eufmfam\relax#1}}

	\def\grb{{\gr b}}
	\def\grc{{\gr c}}

	\def\grg{{\gr g}}
	\def\grh{{\gr h}}

\def\grl{{\gr l}}
	\def\grm{{\gr m}}
	\def\grn{{\gr n}}
	\def\gro{{\gr o}}
	\def\grp{{\gr p}}

	\def\grs{{\gr s}}
	\def\grt{{\gr t}}
	\def\gru{{\gr u}}

\def\db#1{{\fam\msbfam\relax#1}}

\def\dbA{{\db A}} \def\dbB{{\db B}}
\def\dbC{{\db C}} \def\dbD{{\db D}}
 \def\dbF{{\db F}}
\def\dbG{{\db G}} \def\dbH{{\db H}}
 \def\dbJ{{\db J}}
 
 \def\dbN{{\db N}}
 \def\dbP{{\db P}}
\def\dbQ{{\db Q}} \def\dbR{{\db R}}
\def\dbS{{\db S}}

 \def\dbZ{{\db Z}}

\def\gsp{{\grg\grs\grp}}

\def\eps{{\varepsilon}}

\def\Abar{\bar{A}}

\def\Fbar{\bar{F}}

\def\Kbar{\bar{K}}

\def\Rbar{\bar{R}}

\def\Vbar{\bar{V}}

\def\ebar{\bar{e}}

\def\kbar{\bar{k}}

\def\qbar{\bar{q}}

\def\zbar{\bar{z}}
\def\dbQbar{\bar{\dbQ}}
\def\scrMbar{\bar{\scrM}}
\def\scrNbar{\bar{\scrN}}

\def\gtil{\tilde{g}}

\def\mtil{\tilde{m}}
\def\ptil{\tilde{p}}
\def\vtil{\tilde{v}}
\def\stil{\tilde{s}}
\def\Atil{\tilde{A}}
\def\Ftil{\widetilde{F}}
\def\Gtil{\widetilde{G}}
\def\Htil{\widetilde{H}}
\def\Ktil{\widetilde{K}}
\def\Ltil{\widetilde{L}}
\def\Mtil{\widetilde{M}}

\def\Rtil{\widetilde{R}}

\def\Ttil{\widetilde{T}}
\def\Xtil{\widetilde{X}}
\def\Wtil{\widetilde{W}}

\def\psitil{\tilde{\psi}}
\def\mutil{\tilde{\mu}}
\def\Shtil{\widetilde{Sh}}
\def\dbZhat{\widehat{\dbZ}}

\def\Yhat{\widehat{Y}}
\def\scrOhat{\widehat{\scrO}}
\def\scrMtil{\tilde{\scrM}}
\def\Gtil{\tilde{G}}

\def\Ker{\text{Ker}}
\def\der{\text{der}}
\def\Sh{\hbox{\rm Sh}}
\def\Sp{\text{Sp}}

\def\sc{\text{sc}}
\def\Res{\text{Res}}
\def\ab{\text{ab}}
\def\ad{\text{ad}}
\def\Ad{\text{Ad}}
\def\Gal{\text{Gal}}
\def\Hom{\text{Hom}}
\def\End{\text{End}}
\def\Aut{\text{Aut}}
\def\Spec{\text{Spec}}
\def\Tr{\text{Tr}}

\def\Lie{\text{Lie}}
\def\adnt{\text{adnt}}
\def\adnc{\text{adnc}}

\def\leaderfill{\leaders\hbox to 1em
     {\hss.\hss}\hfill}
\def\nspace{\lineskip=1pt\baselineskip=12pt\lineskiplimit=0pt}

\def\Proclaim#1{\medbreak\noindent{\bf#1\enspace}\it\ignorespaces}
\def\finishproclaim{\par\rm
     \ifdim\lastskip<\medskipamount\removelastskip
     \penalty55\medskip\fi}
\def\proof{\par\noindent {\it Proof:}\enspace}
\def\references#1{\par
  \centerline{\boldsectionfont References}\smallskip
     \parindent=#1pt\nspace}
\def\Ref[#1]{\par\smallskip\hang\indent\llap{\hbox to\parindent
     {[#1]\hfil\enspace}}\ignorespaces}
\def\Item#1{\par\smallskip\hang\indent\llap{\hbox to\parindent
     {#1\hfill$\,\,$}}\ignorespaces}
\def\ItemItem#1{\par\indent\hangindent2\parindent
     \hbox to \parindent{#1\hfill\enspace}\ignorespaces}

\def\simover#1{\overset\sim\to#1}

\def\Le{{\mathchoice{\,{\scriptstyle\le}\,}
  {\,{\scriptstyle\le}\,}
  {\,{\scriptscriptstyle\le}\,}{\,{\scriptscriptstyle\le}\,}}}
\def\Ge{{\mathchoice{\,{\scriptstyle\ge}\,}
  {\,{\scriptstyle\ge}\,}
  {\,{\scriptscriptstyle\ge}\,}{\,{\scriptscriptstyle\ge}\,}}}

\def\arrowsim{\,\smash{\mathop{\to}\limits^{\lower1.5pt
  \hbox{$\scriptstyle\sim$}}}\,}

\def\twiceover#1#2{\overset#1\to{\overset\sim\to#2}}
\def\lrdoublemapdown#1#2{\llap{$\vcenter{\hbox{$\scriptstyle#1$}}$}
     \Big\downarrow\!\Big\downarrow%
     \rlap{$\vcenter{\hbox{$\scriptstyle#2$}}$}}
\def\doublemaprights#1#2#3#4{\raise3pt\hbox{$\mathop{\,\,\hbox to
     #1pt{\rightarrowfill}\kern-30pt\lower3.95pt\hbox to
     #2pt{\rightarrowfill}\,\,}\limits_{#3}^{#4}$}}
\def\opoplus{\operatornamewithlimits{\oplus}\limits}
\def\rightcapdownarrow{\raise9pt\hbox{$\ssize\cap$}\kern-7.75pt
     \Big\downarrow}

\def\lcapmapdown#1{\rightcapdownarrow%
     \llap{$\vcenter{\hbox{$\scriptstyle#1\enspace$}}$}}
\def\rcapmapdown#1{\rightcapdownarrow\kern-1.0pt\vcenter{
     \hbox{$\scriptstyle#1$}}}
\def\lmapdown#1{\Big\downarrow\llap{$\vcenter{\hbox{$\scriptstyle#1
     \enspace$}}$}}
\def\rmapdown#1{\Big\downarrow\kern-1.0pt\vcenter{
     \hbox{$\scriptstyle#1$}}}
\def\rightsubsetarrow#1{{\ssize\subset}\kern-4.5pt\lower2.85pt
     \hbox to #1pt{\rightarrowfill}}
\def\longtwoheadedrightarrow#1{\raise2.2pt\hbox to #1pt{\hrulefill}
     \!\!\!\twoheadrightarrow}

\def\Gal{\operatorname{\hbox{Gal}}}
\def\Hom{\operatorname{\hbox{Hom}}}

\NoBlackBoxes
\parindent=25pt
\document
\footline={\hfil}

\null
\vskip 0.5 cm
\centerline{\BIGtitle Integral Canonical Models}
\medskip
\centerline{\BIGtitle of Shimura Varieties of Preabelian Type}
\vskip 0.2in
\centerline{\bigtitle Adrian Vasiu, Univ. of Arizona, July 2003}
\vskip 0.1in
\centerline{{\bf VERSION} which incorporates to the published one${}^1$ $\vfootnote{1}{Asian J. Math. Vol. 3, No. 2, pp. 401--518, June 1999. All copy rights reserved to International Press.}$ most of the addenda and errata of}
\centerline{math.NT/0104152. It adds an additional correction to 3.2.3.2 and so restates the parts}
\centerline{related to 6.4.1.1 2) which relied on 3.2.3.2 in a way which adds extra assumptions. Also: it} 
\centerline{cuts some commas, the's and above's, it replaces few ``make use" by ``using" and ``proves"}
\centerline{by ``proofs". It contains some replacements of some ``-" by ``--" as well as some replacements}
\centerline{of the form ``2.6-8" by ``2.6 to 2.8". It enlarges slightly Step B of 3.2.17.}
\medskip
\centerline{{\bf WARNING}: The initial order of \S 6 has been preserved; however, the corrections performed imply that}
\centerline{the right order to read \S6 is the one mentioned in its beginning paragraph. Six footnotes added.}  
\footline={\hfill}
\null
\vskip0.05in
\centerline{{\bf ABSTRACT}. We prove the existence of integral canonical models of Shimura}
\smallskip
\centerline{varieties of preabelian type with respect to primes of characteristic at least 5.}
\medskip
a\centerline{{\bf MSC 2000}: Primary 11G10, 11G18, 14B12, 14C30, 14D10, 14D22, 14F20, 14F30, 14G35 and 14K10.}
\smallskip
\centerline{{\bf Key words}: abelian and Shimura varieties, deformation theory, Hodge cycles,}
\centerline{reductive and $p$-divisible groups, schemes and integral models.}

\vskip0.2in
\centerline{\bigtitle Contents}

{\nspace{

\smallskip
\line{\item{\bf\S1.}{Introduction }\leaderfill 1}

\smallskip
\line{\item{\bf \S2.}{Preliminaries}\leaderfill 11}

\smallskip
\line{\item{\bf \S3.}{A general view of the integral models of
Shimura varieties} \leaderfill 20}

\smallskip
\line{\item{\bf \S4.}{Shimura varieties of Hodge type and special families of tensors} \leaderfill 59}

\smallskip
\line{\item{\bf \S5.}{The basic result} \leaderfill 81}

\smallskip
\line{\item{\bf \S6.}{The existence of integral canonical models}
     \leaderfill 104}

\smallskip
\line{\item{}{References}\leaderfill 133}

}}

\footline={\hss\tenrm \folio\hss}
\pageno=1

\bigskip
\noindent
{\boldsectionfont \S1. Introduction}

\bigskip
\Proclaim{1.0.} \rm
Let the pair $(G,X)$ define an arbitrary Shimura variety (cf. 2.3) and let $\Sh(G,X)$ be the canonical model of this Shimura variety defined over the reflex field $E(G,X)$ (cf. 2.6 to 2.8). Let $p$
be a rational prime such that $G$ is unramified over
$\dbQ_p$.
Let $v$ be a prime  of $E(G,X)$ dividing
$p$ and let $O_{(v)}$ be the localization of the ring of integers of $E(G,X)$ with respect to it. Let $H$ be a hyperspecial subgroup of $G(\dbQ_p)$. Let $\dbA_f^p$ be the ring of finite ad\`eles with the $p$-component omitted. A
smooth integral model of $\Sh(G,X)/H$ over $O_{(v)}$ is a faithfully flat $O_{(v)}$-scheme $\scrN$ together with a continuous right action (in the sense of [De2, 2.7.1]) of $G(\dbA_f^p)$ on it such that: 

\smallskip
-- its generic fibre $\scrN_{E(G,X)}$ with its induced $G(\dbA_f^p)$-action is $\Sh(G,X)/H$ with its canonical $G(\dbA_f^p)$-action;

-- there is a compact open subgroup $H_0$ of $G(\dbA_f^p)$ with the property that for any inclusion $H_1\subset H_2$ of open subgroups of $H_0$, the canonical morphism $\scrN/H_1\to\scrN/H_2$  induced by the action of $G(\dbA_f^p)$ on $\scrN$, is an \'etale morphism between smooth schemes of finite type over $O_{(v)}$.

\smallskip
In what follows it is irrelevant which hyperspecial subgroup $H$ of $G(\dbQ_p)$ we choose (cf. 3.2.7 2)), and so we often do not mention it.
 
Langlands [La, p. 411] expected the existence of a good smooth integral model of $\Sh(G,X)/H$ over $O_{(v)}$, without expressing what ``good" should mean. Milne (see [Mi4, p. 169] and [Mi3, footnote of p. 513]) conjectured the existence of a smooth integral  model of $\Sh(G,X)/H$ over $O_{(v)}$  having 
an extension property similar to the extension property enjoyed by the N\'eron model (over a discrete valuation ring $O$) of an abelian variety $A$ (over the field of fractions of O). Such a smooth integral  model, if exists, is
called an integral canonical model with respect to $v$ (and $H$) (or simply an integral canonical model, as the prime $v$ is determined by it) of
our Shimura variety $\Sh(G,X)$. For $p>2$, if it exists, it is unique due to the extension property it enjoys (cf. 3.2.4). If $p>2$ and if $\Sh(G,X)/H$ does have an integral canonical model, then this model, as an object of the category of all smooth integral models of $\Sh(G,X)/H$ over $O_{(v)}$, plays the same role (i.e. it is a final object) played by the N\'eron model of $A$, viewed as an object of the category of all  smooth models of $A$ over $O$ (i.e. of the category of all commutative smooth groups  over $O$ having $A$ as their generic fibres). Sections 3.2 to 3.5  present the general definitions and properties pertaining to integral models of Shimura varieties. Some important features are gathered in 3.2.3.2 and 3.2.12, while the descent of such integral models (based on 3.1.3.1) is explained in 3.2.13.

The extension property mentioned above is with respect to healthy regular schemes over $O_{(v)}$. We call a regular scheme $Y$ flat over a discrete valuation ring of mixed characteristic healthy if for any closed subscheme $Z$ of the special fibre of $Y$ of codimension (in $Y$) at least 2, every abelian scheme over the open subscheme of $Y$ defined by $Y\setminus Z$ extends to an abelian scheme over $Y$. We were forced to introduce the notion of a healthy regular scheme due to the fact that the statement 6.8 of [FC, p. 185] is not true in general (cf. [dJO]). However a regular scheme formally smooth over a discrete valuation ring of mixed characteristic, having a residue field of characteristic greater than one plus its ramification index, is healthy [Fa4]; see also  3.2.2 1). A complete proof of this fact is included in 3.2.17. The general theory of healthy normal  schemes as well as different extension properties (like the extended extension property) defined with their help are presented in 3.2.  As an independent result we get (cf. 3.2.2.1 and  3.2.3.3 2)):

\Proclaim{Proposition.}
If $O\hookrightarrow O_1$ is a formally \'etale homomorphism between two discrete valuation rings, with $O$ a henselian ring of mixed characteristic, then we have: 

\smallskip
1) A regular scheme $Y$  over $O$ is healthy iff $Y_{O_1}$ is healthy. 

2) An $O$-scheme $Y$ has the extension property iff $Y_{O_1}$ has the extension property.  
\finishproclaim

For the case when the inclusion $O\hookrightarrow O_1$ is just of index of ramification 1 see 2) of 3.2.2.3 A) and 3.2.2.4 a).

Integral canonical models of Shimura varieties of PEL type (these varieties are forming a subclass of the class --to be briefly reviewed in 1.2-- of Shimura varieties of Hodge type) were constructed in [LR] (cf. also the correction in [Ko]). To our knowledge no concrete integral canonical model of a Shimura variety which is not related to one of PEL type (in the sense that their adjoint varieties  are isomorphic) was previously constructed.

This paper is the first among a sequence of five papers devoted to the existence, the compactification, and the understanding of points with values in perfect fields and in (regular formally smooth rings over) Witt rings over perfect fields of the integral canonical models of Shimura varieties of preabelian type; examples will be provided. The other four papers will be [Va2] to [Va5].

In this paper we are concerned with the existence of integral canonical models of Shimura varieties of
preabelian type. A Shimura variety $\Sh(G_1,X_1)$ is said to be of preabelian type if there is a Shimura variety $\Sh(G_2,X_2)$ of Hodge type such that their adjoint Shimura varieties are isomorphic: $\Sh(G_1^{\ad},X_1^{\ad})\arrowsim \Sh(G_2^{\ad},X_2^{\ad})$. Along our work we will give a strong support to the general point of view that all properties enjoyed by the integral canonical models of Siegel modular varieties and by the universal abelian schemes over them, are also enjoyed (under proper formulation) by the integral canonical models of Shimura varieties of Hodge type (even of preabelian type) with respect to primes having a residue field of characteristic bigger than 2 and by the special abelian schemes over them (see 1.2.2 for the meaning of special used here). 
\finishproclaim

\smallskip
\Proclaim{1.1.} \rm
Our basic result (see 5.1) is: 

\Proclaim{Theorem 0.}
 With the above notations, if the Shimura variety $\Sh(G,X)$ is of Hodge type and if the pair $(G,X)$ satisfies a
slight condition $(*)$ with respect to the prime $p$ (assumed to be greater than 2) (cf. 5.1), then  $\Sh(G,X)$ has an integral canonical model with respect to any prime $v$ of $E(G,X)$ dividing $p$ (and any hyperspecial subgroup $H$ of $G(\dbQ_p)$).
\finishproclaim

Fixing the pair $(G,X)$ (of Hodge type), the condition $(*)$ is satisfied
with respect to any prime $p$ big enough (see 1.2.6). For the proof of our basic result we rely heavily on the crystalline machinery developed in [Fa1] to [Fa3].
\finishproclaim

\smallskip
\Proclaim{1.2.} \rm
To explain 1.1, we start with an injective map $f\colon (G,X)\hookrightarrow (G\Sp(W,\psi),S)$ (cf. 2.4). Here the pair $(GSp(W,\psi),S)$ defines a Siegel modular variety (cf. Example 2 of 2.5). The existence of such an injective  map is what defines the class of pairs $(G,X)$ (defining a Shimura variety) of Hodge type. Let $\dbZ_{(p)}$ be the localization of $\dbZ$ with respect to $p$. We assume the existence of a $\dbZ_{(p)}$-lattice $L$ of $W$ such that the alternating form $\psi\colon W\otimes W\to\dbQ$  induces a perfect form $\psi\colon L\otimes L\to\dbZ_{(p)}$ (i.e. the induced $\dbZ_{(p)}$-linear map from $L$ into its dual $L^*$ is an isomorphism) and there is a family of tensors $(v_{\alpha})_{\alpha\in\scrJ_0}$ in $\dbZ_{(p)}$-modules of the form $(L\otimes L^*)^{\otimes n}$, $n\in\dbN$, fixed by $G$ and of degree at most $2(p-2)$ (if $v_{\alpha}\in (L\otimes L^*)^{\otimes n}$ then the degree of $v_{\alpha}$ is $2n$), and which is $\dbZ_{(p)}$-well positioned with respect to $\psi$ for the group $G$ (see 4.3.4 for a precise definition of the notion of a well positioned family of tensors). Let $K:=\{g\in GSp(W,\psi)(\dbQ_p)\mid g(L\otimes\dbZ_p)=L\otimes\dbZ_p\}$. The hypotheses on $L$ imply that $K$ is a hyperspecial subgroup  of $G\Sp(W,\psi)(\dbQ_p)$ and that the Zariski closure $G_{\dbZ_{(p)}}$ of $G$ in $GSp(L,\psi)$ is a reductive group over $\dbZ_{(p)}$. So the intersection $H:=G(\dbQ_p)\cap K$ is a hyperspecial subgroup of $G(\dbQ_p)$. We choose a $\dbZ$-lattice $L_{\dbZ}$ of $W$ such that $\psi$ induces a perfect form $\psi\colon L_{\dbZ}\otimes L_{\dbZ}\to\dbZ$ and $L=L_{\dbZ}\otimes\dbZ_{(p)}$. Let $(v_{\alpha})_{\alpha\in\scr J}$ (with $\scr J_0\subset \scr J$) be an enlarged family of tensors in the tensor algebra of $W\oplus W^*$ ($W^*$ being the dual $\dbQ$--vector space of $W$) such that $G$ is the subgroup of $GSp(W,\psi)$ fixing the tensors of this family. The choice of the lattice $L_{\dbZ}$ and of the family $(v_{\alpha})_{\alpha\in\scr J}$ 
allows the interpretation of $\Sh(G,X)(\dbC)$ as the set of
isomorphism classes of principally polarized abelian
varieties over $\dbC$ of dimension $g$ (with $2g=\dim_{\dbQ}(W)$), having all level structures, carrying a family of Hodge cycles $(w_{\alpha})_{\alpha\in\scr J}$ and
satisfying some additional conditions (cf. 4.1).
\finishproclaim

\Proclaim{1.2.1.} \rm
It is well known that the $\dbZ_{(p)}$-scheme $\scrM$ parameterizing  isomorphism classes of principally polarized abelian schemes of dimension $g$ over $\dbZ_{(p)}$-schemes, having 
level-$N$ symplectic similitude structure for any $N\in\dbN$ relatively prime to $p$, together with the canonical action of $G\Sp(W,\psi)(\dbA_f^p)$ on it, is an integral canonical model of $\Sh(G\Sp(W,\psi),S)/K$ over $\dbZ_{(p)}$ (see 3.2.9).
\finishproclaim

\Proclaim{1.2.2.} \rm
The normalization $\scrN$ of the Zariski closure of $\Sh(G,X)/H$ in $\scrM_{O_{(v)}}$ is a normal integral  model of $\Sh(G,X)/H$ having the (extended) extension property (cf. 3.4.1; see def. 2) and 3) of  3.2.3). This integral model is an integral canonical model of $\Sh(G,X)/H$ iff $\scrN$ is formally smooth over $O_{(v)}$ (cf. 3.4.4). The universal principally polarized abelian scheme over $\scrM$ gives birth to a principally polarized abelian scheme
($\scrA,\scrP_{\scrA}$) over $\scrN$, which we call special. Let $\dbF$ be the algebraic closure of the residue field $k(v)$ of $v$. 

Hodge cycles are (presently) defined only in characteristic zero.
But the Hodge cycles (of degree not bigger than $2(p-2)$) of an abelian scheme over a discrete valuation ring  
which is finite flat over a ring $W(k)$ of Witt vectors of a perfect field $k$ of characteristic $p$ are well behaved (cf. [Fa3, Cor. 9]) with respect to the integral version
of Fontaine's comparison map (see [Fa3, Th. 7]). 
Using the above hypotheses on the $\dbZ_{(p)}$-lattice $L$, we first exploit (cf. 5.2.12) the good behavior of Hodge
cycles under the integral version of Fontaine's comparison
map (i.e. we first pass from a reductive group in the \'etale $\dbZ_p$-context to a reductive group in the integral crystalline, de Rham context). Then we use (cf. 5.2.10) de Rham conjecture of [Fa1] and [Fa2] to construct (cf. 5.3 and 5.4) local deformations of (principally polarized) abelian
schemes of dimension $g$ (having some level structures) over $W(\dbF)$, 
carrying a family of Hodge cycles and satisfying the required additional
conditions. With these deformations we prove the formal smoothness of $\scrN$. The main new idea besides the ones of (the original versions of) [Fa3] is the use of the ring $\Rtil e$ introduced in 5.2.1 (here $e\in\dbN$). It is a projective limit of artinian $W(\dbF)$-algebras; this fact plays a key role in 5.3.2.  

We detail the above two steps. Let $y\colon\Spec(\dbF)\hookrightarrow\scrN_{W(\dbF)}$ be a closed point, and let $V$ be a discrete valuation ring which is a finite flat extension of $W(\dbF)$ such that $y$ can be lifted to a point $z_V\colon\Spec(V)\to\scrN_{W(\dbF)}$. Let $e:=[V:W(\dbF)]$ and let $\Rtil^n e$ be the normalization of $\Rtil e$ in its field of fractions. First we show, starting from $z_V$, the existence of a morphism $\Spec(\Rtil^n e)\to\scrN_{W(k)}$ lifting $y$ (5.3.1.1). Using the natural epimorphism $\Rtil^n e\twoheadrightarrow W(\dbF)$ (cf. 5.3.4), we deduce the existence of a good lift $z_{W(\dbF)}\colon \Spec(W(\dbF))\to\scrN_{W(\dbF)}$ of $y$. Second we use directly [Fa3, Th. 10 and the remarks after] in the context provided by $z_{W(\dbF)}$ (see 5.4 to 5.5).
\finishproclaim

\Proclaim{1.2.3.} \rm
On the way of proving the formal smoothness of $\scrN$ we obtain (cf. 5.2.16) an improvement in the
Principle B of [Bl, 3.1].
\finishproclaim

\Proclaim{1.2.4.} \rm
Chapter 5 is entirely devoted to the construction of such local deformations and to the proof of the formal smoothness of $\scrN$, while the general (needed) theory of well positioned families of tensors for a reductive group is presented in 4.3. The  most useful well positioned families of tensors (of a general nature) are presented in 4.3.10 b) (see also 4.3.10.1) for the case of a semisimple group and in 4.3.13 for the case of a torus. For the case of a reductive group we use families of tensors formed by putting together well positioned families of tensors for its derived group and well positioned families of tensors for its toric part (i.e. for the maximal subtorus of its center): Lemma 3.1.6 allows us to do this (cf. 4.3.6 2)). The behavior of hyperspecial subgroups with respect to homomorphisms of reductive groups needed for this general theory is described in 3.1.2.

The proof of 4.3.10 b) is in two parts. The first part is a criterion of when a Lie algebra over a reduced ring $R$ comes from a semisimple group $\tilde G_R$ over $R$. The second part is a criterion of when a representation of $\Lie(\tilde G_R)$ comes from a representation of $\tilde G_R$.
\finishproclaim

\Proclaim{1.2.5.} \rm
 In 5.7.5 we illustrate our ideas in the case of classical Spin modular varieties of odd dimension (and rank two), while in 4.3.1 we show how, the previously known case of Shimura varieties of PEL type, is (for $p\Ge 3$) a particular case of our approach via well positioned families of tensors.
\finishproclaim

\Proclaim{1.2.6.} \rm
The condition $(*)$ means: there is an injective
map $f\colon (G,X)\hookrightarrow (G\Sp(W,\psi),S)$ for which there is a $\dbZ_{(p)}$-lattice $L$ of $W$ satisfying the  conditions mentioned in the first paragraph of 1.2. Fixing an injective map $f\colon (G,X)\hookrightarrow (G\Sp(W,\psi),S)$, for any rational prime $p$ big enough (with an effectively computable bound, just in terms of the representation $G\to GL(W)$), we can find a $\dbZ_{(p)}$-lattice $L$ of $W$ satisfying these  conditions   (cf. 5.8.7 and 5.8.1).

We use only $\dbZ_{(p)}$-well positioned  families of tensors having only tensors of degrees 2 and 4, and the condition $p\Ge 5$ (cf. the inequality $4\Le 2(p-2)$) is needed just for being allowed to use tensors of degree 4. The most useful tensors of degree  4 are presented in 4.3.2. In essence we are using only three tensors of degree 4. Fixing the injective map $f$, these tensors are endomorphisms of $\text{End}(W)$ (we are  identifying $\text{End}(\text{End}(W))$ with $(W\otimes W^*)^{\otimes 2}$ and $\text{End}(W)$ with $\text{End}(W^*)$): 

\smallskip
-- the first one (cf. 4.2.1) is the projection $\pi(\Lie(G^{\der}),W)$ of $\text{End}(W)$ on $\Lie(G^{\der})$ along the orthogonal complement of $\Lie(G^{\der})$ with respect to the trace bilinear form on $\text{End}(W)$; 

-- the other two tensors $B$ and $B^*$ are elements of $\text{End}(\text{End}(W))$ expressing that the Killing form on $\Lie(G^{\der})$ is perfect. 

\smallskip
Part b) of 4.3.10 together with a well known fact on Shimura varieties of Hodge type (expressed in the proof of 5.7.1 by $s(\Lie(G^{\der}),W)=2$) imply:

\Proclaim{Fact.}
The family of tensors formed by $\pi(\Lie(G^{\der}),W)$, $B$ and $B^*$ is $\dbZ_{(p)}$-well positioned for $G^{\der}$.
\finishproclaim

The role of $\psi$ is irrelevant here and so we do not need to mention with respect to $\psi$.

\Proclaim{1.2.6.1.} \rm
The condition $(*)$ is satisfied if there is a $\dbZ_{(p)}$-lattice $L$ of $W$ such that $\psi\colon L_{(p)}\otimes L_{(p)}\to\dbZ_{(p)}$ is perfect, the Zariski closure of $G$ in $GSp(L_{(p)},\psi)$ is a reductive group over $\dbZ_{(p)}$ and the above three tensors are integral with respect to it (cf. 5.7.1). 
This forms a simple criterion for a practical form of Theorem 0.
\finishproclaim

\Proclaim{1.2.6.2.} \rm
The use of tensors (Hodge cycles) of degree 4 allows us to have a uniform treatment of all Shimura varieties of Hodge type, with no preference for Shimura varieties of PEL type. But we would like to remark that, as it will be seen along our work (cf. [Va2]), the study of Shimura varieties of Hodge type of $A_l$, $B_l$ or $D_l^{\dbR}$ type (see [De2] for the possible types of a Shimura variety) is (somehow) easier than the study of Shimura varieties of Hodge type of $C_l$ or $D_l^{\dbH}$ type. 
\finishproclaim

\smallskip
\Proclaim{1.3.} \rm
We prove (6.5.1.1) the $\dbZ_{(p)}$-version of the main result of [De2]. In its simplified form (6.4.2):
\finishproclaim

\Proclaim{Theorem 1.} 
For any adjoint Shimura variety $\Sh(G_0,X_0)$ of abelian type and for any prime $p\Ge 5$ such that $G_0$ is unramified over $\dbQ_p$, there is a Shimura variety of Hodge type $\Sh(G,X)$ having $\Sh(G_0,X_0)$ as its adjoint variety, with $G$ unramified over $\dbQ_p$, and such that the pair $(G,X)$ satisfies the condition $(*)$ (of 5.1)
with respect to $p$. Moreover, for any Shimura variety $\Sh(G_1,X_1)$ of abelian type having $\Sh(G_0,X_0)$ as its adjoint variety, there is an isogeny $G^{\der}\to G_1^{\der}$.
\finishproclaim

\Proclaim{1.3.1.} \rm
There are three main tools needed for the proof of Theorem 1. The first two are provided by [De2, 2.3.10] and by the above Fact, via 1.2.6.1. But they are not enough: it is not always possible to find a $\dbZ_{(p)}$-lattice $L$ of $W$ as in 1.2.6.1. For instance the Killing form of the Lie algebra of a simple split adjoint group of $B_l$ Lie type over $W(\dbF)$ is not perfect  if $p$ divides $2l-1$, $l\in\dbN$. The third tool (cf. 6.5 and 6.6) is the construction of injective maps $f\colon (G,X)\hookrightarrow (GSp(W,\psi),S)$ such that there is a reductive subgroup $\tilde G$ of $GL(W)$ (we are not bothered if it is or it is not  contained in $GSp(W,\psi)$; however see 6.6.2) containing $G$, unramified over $\dbQ_p$, and such that: 

\smallskip
-- a variant of 1.2.6.1 (for instance cf. 5.7.4 and the proof of 6.5.1.1) can be applied to $\tilde G^{\der}$;

-- we can ``regain" $G$ out of $\tilde G$ by using endomorphisms of $W$ fixed by $G$ (i.e. we have a relative PEL situation, see 4.3.16).    
\finishproclaim

The construction of such injective maps is carried out in 6.5 and 6.6. It relies heavily on the classification [De2] of types of simple, adjoint Shimura varieties of abelian type: for each such type we have to proceed differently (but similarly).

\smallskip
\Proclaim{1.4.} \rm
Theorem 1 implies a
positive answer to Milne's conjecture in the case of Shimura
varieties of preabelian type for primes $p\Ge 5$ (6.4.1):
\finishproclaim

\Proclaim{Theorem 2.*}
If $(G,X)$ defines a Shimura variety of preabelian type and if $p\Ge 5$ is a rational prime such that $G$ is unramified over $\dbQ_p$, then $\Sh(G,X)$ has an integral canonical model with respect to any prime $v$ of $E(G,X)$ dividing $p$.
\finishproclaim

 As a scheme this model is a pro-\'etale cover of a quasi-projective smooth scheme over $O_{(v)}$. Chapter 6 is devoted to the proof of Theorem 2 (via Theorem 1). The passage from the existence of integral canonical models of Shimura varieties of Hodge type to the existence of integral canonical models of Shimura varieties of preabelian type is explained in  6.1 and 6.2. 

The passage from the  Hodge type case to the abelian type case is achieved by taking quotients through group actions (cf. 6.2.2). The groups involved are $M$-torsion groups for some $M\in\dbN$ (cf. the proof of 6.2.2). In all cases, except the case when we deal with Shimura varieties of whose adjoint varieties have a simple factor of $A_l$ type, and with a prime $p$ dividing $l+1$, this is straightforward, as $M$ is relatively prime to $p$ (cf. 3.4.5.1 and 6.2.2). When $p$ divides $M$ we have to express more concretely these group  actions and to prove that they are free actions. This is achieved in 6.2.2.1, based on properties of adjoint groups to which a particular study of adjoint filtered Lie $\sigma$-crystals attached (see 5.4.6) to maps $z_{W(\dbF)}$ as in 1.2.2 gets reduced (here $\sigma$ is the Frobenius automorphism of $W(\dbF)$). Regardless of how are $M$ and $p$ we need the fact (again cf. the proof of 6.2.2) that the integral canonical models obtained through Theorem 0 are moduli schemes (of abelian schemes). This passage is supported by simple variants (cf. 3.2.14 and 6.2.3) of [De1, 1.15].

The passage from the abelian type case to the preabelian type case is achieved via the normalization procedure (cf. 6.1). 

The paper contains a complete proof of Theorem 2 for the abelian case, while the last step (6.1.2) needed for the proof in the case of Shimura varieties of preabelian type which are not of abelian type will be presented in [Va3], as it requires the formalism of  smooth toroidal compactifications of  integral canonical models of Shimura varieties (of preabelian type). For a discussion and another approach, see 6.8. In 6.8 a complete proof of 6.1.2 is included for the case when $\scrN$ is a pro-\'etale cover of a proper scheme over $O_{(v)}$ and for the generic situation (i.e. for when $p$ is big enough) of the general case. It is based on 3.2.11 and [De2, 2.3.8]. For the sake of convenience, the results depending on the proof of 6.1.2 in the remaining cases (they are described in 6.8.6), are labeled (cf. 6.1.2.1) with a star. So also Theorems 2 and 3 are labeled. Warning: the labeled results are proved here entirely for the abelian type and for the generic situation.
The independent result 5.6.5 h) is not proved here: it is labeled with two stars.

\Proclaim{1.4.1.} \rm
For making some of the main results easily accessible to a larger mathematical community, we state in 6.4.10 a simple criterion of how to recognize an integral canonical model whose existence is guaranteed by Theorem 2.
\finishproclaim

\smallskip
\Proclaim{1.5.} \rm
Different quotients of integral canonical models of a Shimura variety $\Sh(G,X)$ of preabelian type  with respect to primes $v$ of $E(G,X)$ having a residue field of characteristic at least 5, can be glued together (see 6.4.3 and 6.4.4). To state this precisely let $\scrS$ be the set of primes whose elements are 2, 3 and the primes $p\Ge 5$ for which $G$ is  ramified over $\dbQ_p$. It is a finite set. We write the ring of finite ad\`eles as a product $\dbA_f=(\prod_{q\in\scrS} \dbQ_q)\times\dbA_f^\scrS$. Let $H^\scrS$ be a compact open subgroup of $G(\dbA_f^\scrS)$ of maximal volume (with respect to a Haar measure on the locally compact group $G(\dbA_f^{\scrS})$). It is a product of its $q$-components (for primes $q\not\in\scrS$), and every such $q$-component of it is a hyperspecial subgroup $H^q$ of $G(\dbQ_q)$ (regardless of the chosen Haar measure). We have:
\finishproclaim

\Proclaim{Theorem 3.*} 
We consider a compact open subgroup $H_{\scrS}$ of $G($$\prod_{q\in\scrS} \dbQ_p)$ such that for any prime $l\notin\scrS$, either $H_{\scrS}\times H^{\scrS}$ is $l$-smooth for $(G,X)$ in the sense of 2.11 or there is no element of $G^{\ad}_{V_0^l}(V_0^l)$ of order $l$; here $V_0^l$ is the completion of the maximal unramified extension of $\dbZ_l$ and $G_{\dbZ_l}$ is the reductive group over $\dbZ_l$ having $G_{\dbQ_l}$ as its generic fibre and $H^l$ as its group of $\dbZ_l$-valued points. Then there is a quasi-projective smooth scheme $\scrM(H_{\scrS})$ over the normalization $O_{(\scrS)}$ of $\dbZ\fracwithdelims[]1{\prod_{q\in\scrS} q}$ in $E(G,X)$, uniquely determined by the fact that its generic fibre is $\Sh(G,X)/H_{\scrS}\times H^{\scrS}$ and that, for any prime $v$ of $E(G,X)$ dividing a rational prime $q\notin\scrS$, the normalization of $\scrM(H_{\scrS})_{O_{(v)}}$ in the ring of fractions of $\Sh(G,X)/H^q$ is the integral canonical model of $\Sh(G,X)$ with respect to $v$ (and $H^q$).
\finishproclaim

These smooth schemes are the analogue of the schemes (attached to Siegel modular varieties) parameterizing principally polarized abelian schemes (of a given dimension) having a finite  symplectic similitude level-structure.
They enjoy a very important extension type property (cf. rm. 1) of 6.4.6). They are models over ``punctured" ring of integers (of number fields) of quotients of (some) finite disjoint unions of Hermitian symmetric domains by (some) arithmetic subgroups. In rm. 3) of 6.4.6 we explain why the notation $\scrM(H_{\scrS})$ is quite justified. For the compact case (i.e. when $\Sh(G,X)$ is a pro-\'etale cover of a projective  smooth $E(G,X)$-scheme) see 6.4.11. The proof of Theorem 3 is based on 6.2.4.1, which is a natural consequence of the ideas presented in 6.2.3, in the proof of 6.2.2 and in 6.5 and 6.6. 

\smallskip
\Proclaim{1.6.} \rm
We present now the part of [Va2] which brings more light to some parts of the present paper. All that follows in 1.6 to 1.8 could have been equally well presented as remarks at different places of \S 5 and \S6; but for the sake of convenience, we gathered all these results (referred to in \S 5 and \S6) here.

We extend the well known results (for Siegel modular varieties) concerning the existence of an ordinary type and the existence of the  canonical lift of an abelian variety of ordinary type, to any special principally polarized abelian scheme $(\scrA,\scrP_{\scrA})$ over an arbitrary  integral canonical model $\scrN$ of a Shimura variety $\Sh(G,X)$ of Hodge type with respect to a prime $v$ of $E(G,X)$ dividing a rational prime $p\Ge 5$ (cf. also [Va1]). Let $\scrN_{k(v)}$ be the special fibre of $\scrN$. Using the notations of 1.2 we obtain:

\smallskip
-- a $G$-ordinary type (with respect to the prime $v$ and the injective map $f\colon (G,X)\hookrightarrow (GSp(W,\psi),S)$), which is the formal isogeny type associated to abelian varieties (obtained from $\scrA$ by pull back) over the geometric points of a Zariski dense open subscheme of $\scrN_{k(v)}$;

-- $G$-ordinary points of $\scrN_{k(v)}$ (these are the points of $\scrN_{k(v)}$, with values in a field, with the property that the abelian varieties over them obtained from $\scrA$ by pull back, have as a formal isogeny type, the $G$-ordinary type);

-- $G$-canonical lifts of $G$-ordinary points with values in perfect fields (these $G$-canonical lifts are points of $\scrN$ with values in rings of Witt vectors of perfect fields).

\smallskip
The $G$-ordinary type we obtain is a usual ordinary type iff the field $k(v)$ has $p$ elements. If this is so then the abelian variety over $W(k)$ obtained from $\scrA$ by pull back through a $G$-canonical lift of a $G$-ordinary point (with values in the perfect field $k$) of $\scrN_{k(v)}$, is the canonical lift of an abelian variety of ordinary type.
\finishproclaim

\Proclaim{1.6.1.} \rm
The point defined by the generic fibre of the $G$-canonical lift of a $G$-ordinary point of $\scrN_{k(v)}$ with values in the algebraic closure of $k(v)$ is a special point of $\scrN_{E(G,X)}$ (see 2.10 for the definition of special points). 

We also prove another conjecture of Milne [Mi5, 0.1] (cf. 5.6.5, 5.6.6, 5.8.8 and [Va2]).
\finishproclaim

\Proclaim{1.6.2.} \rm 
To any point $y\colon \Spec(k)\to\scrN_{k(v)}$ (with $k$ an algebraically closed field) we attach (see 5.4.6 for the case $k=\dbF$) a Lie $\sigma$-crystal $(\grg,\varphi)$: $\grg$ is the Lie algebra of a reductive group over $W(k)$ whose fibre over $K_0:=W(k)\fracwithdelims[]1p$ is $G_{K_0}$, while $\varphi$ is a $\sigma$-linear automorphism of $\grg\otimes K_0$ ($\sigma$ being the Frobenius automorphism of $K_0$) such that $\varphi(p\grg)\subset\grg$. Any lift $z\colon \Spec(W(k))\to\scrN$ of $y$  produces naturally a filtration $0=F^2(\grg)\subset F^1(\grg)\subset F^0(\grg)\subset F^{-1}(\grg)=\grg$ such that $\varphi({\frac 1p}F^1(\grg)+F^0(\grg)+p\grg)=\grg$. So $(\grg,\varphi,F^0(\grg),F^1(\grg))$ is a $p$-divisible object of the category $\scrM\scrF_{[-1,1]}(W(k))$ (defined in [Fa1]). Also $F^0(\grg)$ is a parabolic Lie subalgebra of $\grg$ and $F^1(\grg)$ is the Lie algebra of its unipotent radical. The point $y$ is a $G$-ordinary point iff there is a lift $z$ of it to $W(k)$ which makes (by adding filtrations) the Lie $\sigma$-subcrystal of $(\grg,\varphi)$ corresponding to non-negative slopes to be a $p$-divisible object  of the category $\scrM\scrF_{[0,1]}(W(k))$ (of [Fa1]). Such a lift $z$, if exists (i.e. if $y$ is a $G$-ordinary point), is unique and defines the $G$-canonical lift of $y$.

 These Lie $\sigma$-crystals allow us to achieve a stratification of $\scrN_{k(v)}$ in $G(\dbA_f^p)$-invariant locally closed subschemes indexed by Newton polygons of the attached $\sigma$-crystals $\grg(1)$ (we tensor $\grg$ with $W(k)(1)$ to get only non-negative slopes) similar to the one enjoyed by special fibres of integral canonical models of Siegel modular varieties.  The $G$-ordinary points of $\scrN_{k(v)}$ are the points of the (generic)  Zariski dense open stratum. 
\finishproclaim

\smallskip
\Proclaim{1.7.} \rm
We also show how the results mentioned in 1.6.1, together with their proofs, can be used for handling  the Langlands--Rapoport conjecture ([LR]; see [Mi5] and [Pf] for the correct formulation) for an arbitrary integral canonical model $\scrN$  of a Shimura variety $\Sh(G,X)$ of preabelian type with respect to a prime $v$ of $E(G,X)$  having a residue field $k(v)$ of characteristic $p\Ge 5$. 

Let $\dbF$ be the algebraic closure of $k(v)$ and let $\Phi$ be the Frobenius automorphism of it having $k(v)$ as its fixed field. To the triple $(G,X,v)$ it is attached a set $M(G,X,v)$ on which $G(\dbA_f^p)$  and $\Phi$ act (cf. [Mi5] and [Pf]). The Langlands--Rapoport conjecture for $\scrN$ (or for the triple $(G,X,v)$) asserts the existence of a bijection of sets $f_{\scrN}\colon M(G,X,v)\arrowsim \scrN(\dbF)$, preserving the actions of $G(\dbA_f^p)$ and $\Phi$ on them. The existence of the canonical Lie stratification of the special fibre $\scrN_{k(v)}$ of $\scrN$ allows a formulation of the Langlands--Rapoport conjecture for any individual stratum of this stratification. To prove the Langlands--Rapoport conjecture for $\scrN$ is the same as proving  the Langlands--Rapoport conjecture for each individual stratum. We do prove it for the open stratum. 
\finishproclaim

\Proclaim{1.7.1.} \rm
The proof of [Mi5, 0.1] together with [Mi5, 6.4]  imply (cf. also [Mi5, 6.12]; [Mi5, 6.12] is worked out under the hypothesis of [Mi5, p. 24]: it can be removed, cf. 5.6.4) that the Langlands--Rapoport conjecture is true for $\scrN$ if $\Sh(G,X)$ is a Shimura variety (whose adjoint factors are) of $A_l$, $B_l$, $C_l$ or $D_l^{\dbR}$ type, modulo a sufficiently good theory  of reduction of Hodge cycles mod $p$  (very important progress was made in this direction by Milne, conform the presentation in [Mi5]). We explain why the use of such a theory of reduction of Hodge cycles can be avoided for all Shimura varieties of preabelian type.  We first prove that the integral canonical models of Shimura varieties of $D_l^{\dbH}$ type can be treated entirely as the other Shimura varieties: if there is a special principally polarized abelian scheme $(\scrA,\scrP_{\scrA})$ over $\scrN$ (this implies that $\Sh(G,X)$ is a Shimura variety of Hodge type; so $\scrN_{E(G,X)}$ is a moduli space of principally polarized abelian schemes of a given dimension, having a family of Hodge cycles  and some level structures,  and satisfying some extra conditions, while $\scrA_{E(G,X)}$ is the universal abelian scheme over $\scrN_{E(G,X)}$), then for any point $y\colon \Spec(k)\to\scrN$, with $k$ an  algebraically closed field of characteristic $p$, any principally polarized abelian variety over $k$ which is $G$-isogenous to $(A_y,p_{A_y})$ (i.e. it is isogenous in a sense involving the cycles) is $G$-isomorphic (i.e. it is isomorphic in a sense involving the cycles) with $(A_z,p_{A_z})$ for some $k$-valued point $z$ of $\scrN$ (in other words the $G$-isogeny classes are as expected to be). Here the principally polarized abelian varieties $(A_y,p_{A_y})$ and $(A_z,p_{A_z})$ are obtained from $(\scrA,\scrP_{\scrA})$ by pull back through $y$ and respectively $z$.  

As an application of this we show the existence in the general case of 1.7 of an injective map $f_{\scrN}\colon M(G,X,v)\hookrightarrow\scrN(\dbF)$ preserving the actions of $G(\dbA_f^p)$ and $\Phi$ on them (the $\dbF$-valued points of the open stratum of $\scrN_{k(v)}$ are in the image of $f_{\scrN}$). Moreover we prove that $f_{\scrN}$ is a bijection (and so that the Langlands--Rapoport conjecture for $\scrN$ is true) if the residue field $k(v^{\ad})$ of the prime $v^{\ad}$ of $E(G^{\ad},X^{\ad})$ divided by $v$ has precisely $p$ elements, or if $(G^{\ad},X^{\ad})$ has all the simple factors of $A_l$, $B_l$ or $D_l^{\dbR}$ type, with $l\in\dbN$.

\smallskip
\Proclaim{1.8.} \rm
In [Va3] we introduce the notion of an integral canonical model of a Kuga variety of Hodge type. Their existence is implied by 
the existence of integral canonical models of Shimura varieties of Hodge type. These models allow us to prove the existence of smooth toroidal  compactifications of integral canonical models of Shimura varieties of preabelian type (this has been cojectured by Milne [Mi4, 2.18]): Any integral canonical model $\scrN$ of a Shimura variety $\Sh(G,X)$ of preabelian type with respect to a prime $v$ of $E(G,X)$ dividing a rational prime $p\Ge 5$, admits
plenty of smooth toroidal compactifications and has a minimal (normal) compactification $\scrN^{mc}$. The smooth toroidal compactifications of $\scrN$ are obtained from $\scrN^{mc}$ through blowings up. In particular, if $\scrN_{E(G,X)}$ is a pro-\'etale cover of a projective scheme over $E(G,X)$, then $\scrN$ is a pro-\'etale cover  of a projective smooth scheme over $O_{(v)}$.

The toroidal compactifications of Shimura varieties of Hodge type are obtained through the same  procedure (as the integral canonical models are obtained) of taking the normalization of the Zariski closure of smooth toroidal compactifications (over number fields) of quotients of Shimura varieties of Hodge type in (extensions to \'etale  $\dbZ_{(p)}$-algebras of) smooth toroidal compactifications of quotients of integral canonical models of  Siegel modular varieties constructed in [FC] (cf. [Har] for the non-integral part over number fields). So we get special semi-abelian schemes over smooth toroidal compactifications of integral canonical models of Shimura varieties of Hodge type.

See [Va3] for definitions and for the proofs of the results mentioned in 1.8.
\finishproclaim

\smallskip
\Proclaim{1.9.} \rm
A part of the results presented in this paper, is a completely revised and improved version of the first part of our thesis [Va1]. For the sake of not making this paper too long, in 4.3.11, 6.3, 6.5.1.1 and 6.6 we use also the notations of other papers. 

The reader who is interested just to have a pretty good idea about what is going on in this article  can follow the route: 3.2 (to pick up whatever the reader is not familiar with), 4.1, 5.1, 5.6 to 5.8 and 6.4. 

We would like to thank Prof. Gerd Faltings for his encouragements to approach gradually the topics mentioned above, for numerous discussions we had about his recent results [Fa3] and [Fa4] (results without which this work would have had fewer fruits), for his advices and correction of the proofs of 5.1 and 3.4.5.1. We would like to express our gratitude to Prof. James Milne, whose very beautiful and deep work [Mi1] to [Mi6]  is highly inspiring to us and whose conjectures (see [Mi3] to [Mi5]) were the starting point of our work. We are also very much obliged to [De2]. We would like to thank Ben Moonen for asking us how healthy schemes behave with respect to the pull back operation through morphisms of schemes defined by homomorphisms of index of ramification 1 between two discrete valuation rings of mixed characteristic (this was the starting point for a great part of 3.2.1 to 3.2.3), and for the request of enlarging the presentations of 4.3.10 b) and 6.2.2. We would like to thank Prof. Pierre Deligne for pointing out a mistake in a preliminary version of 3.2.2 4).  

I would like to thank Princeton University, Max-Planck Institute from Bonn, FIM, ETH-Z\"urich and UC at Berkeley for providing us with excellent conditions for the writing of this paper. This work was partially supported by the NSF grant DMS 97-05376.
\finishproclaim

\bigskip
\noindent
{\boldsectionfont \S2. Preliminaries}

\bigskip
We fix our notations by mostly reviewing some well known facts (cf. [De1], [De2] and [Mi4]).                                                                                                                                                   
\smallskip
\Proclaim{2.1. Notations and conventions.}\rm
Reductive groups over fields are always assumed connected. Reductive group schemes are understood to have connected fibres. 
For a reductive group $G$ over a scheme we denote by $G^{\der}$, $Z(G)$,
$G^{\ab}$ and $G^{\ad}$, respectively, the derived group of
$G$, the center of $G$, the maximal abelian quotient of $G$
and the adjoint group of $G$. We say that a reductive group $G$ over $\dbQ$ is unramified over $\dbQ_p$ ($p$ being a rational prime) if $G_{\dbQ_p}$ is unramified over $\dbQ_p$. For $G$ an affine group scheme over a scheme $S$ we often denote by $\Lie(G)$ its Lie algebra, and in the case when $G$ is a reductive group scheme we denote by $Aut(G)$ the group scheme over $S$ defined by automorphisms of $G$.

If $X$ is a set endowed with an equivalence relation $R\subset X\times X$, we denote by $[x]\in X/R$ the equivalence class of $x\in X$. For a map $f\colon A\to B$ and for a subset $A_1$ of $A$, we denote by $f|A_1$ the restriction of $f$ to $A_1$. If $f\colon A\to B$ and $g\colon B\to C$ are morphisms in some category  we refer to $g\circ f$ as the composition of $f$  with $g$. All projective limits of schemes are assumed to be filtered, with affine transition morphisms. 

The expression $(G,X)$ always denotes a pair defining a
Shimura variety, while $E(G,X)$ denotes its attached reflex field. Also $\Sh(G,X)$ denotes the Shimura variety defined by $(G,X)$, identified in 2.3 to 2.8 (resp. in the rest of the paper) with the complex variety (resp. with the canonical model of the complex variety). For an arbitrary compact subgroup $K$ of $G(\dbA_f)$, we denote by $\Sh_{K}(G,X)$ the quotient of $\Sh(G,X)$ by $K$. Any $x\in X$ and any $a\in G(\dbA_f)$ define a complex point $[x,a]$ of $\Sh_K(G,X)$. 

If $k$ is a field we denote by $\kbar$ its algebraic 
closure.
For a perfect field $k$, $W(k)$ is the ring of Witt vectors
of $k$. Always $V_0$ denotes such a Witt ring for $k=\kbar$ and then $K_0$ automatically denotes its field of fractions. If $v$ is a prime of a global field $E$, we denote by $k(v)$ its residue field and by $O_{(v)}$ the localization of the ring of integers of $E$ with respect to it. The maximal abelian extension of $E$ is denoted by $E^{\ab}$. For a local ring $R$ we denote by $R^{\text h}$, $R^{\text{sh}}$ and $\widehat{R}$  respectively its henselization, its strict henselization and its completion with respect to its maximal ideal. 

Let $p$ be a rational prime. We usually write  $\dbZ_{(p)}$ instead of $O_{(p)}$. The ring of finite ad\`eles $\dbZhat\otimes_{\dbZ}\dbQ$ is
denoted by $\dbA_f$ and the ring of finite ad\`eles with the
$p$-component omitted is denoted by $\dbA_f^p$. We use freely different Tate-twists: $\dbQ(1)$, $\dbQ_{p}(1)$, $\dbZ_{p}(1)$, $\dbA_{f}(1)$ etc. For $G$ a linear group over $\dbQ$, $G(\dbA)$  is endowed with the coarser topology which makes all the maps $G(\dbA)\to\dbA_{\dbQ}^1(\dbA)=\dbA$, induced by morphisms $G\to\dbA_{\dbQ}^1$, continuous ($\dbA_{\dbQ}^1$ being the affine line over $\dbQ$). Similarly for $G(\dbA_f)$. If $G$ is a linear group over the field $K$ of fractions of a discrete valuation ring (abbreviated DVR), then $G(K)$ is endowed in the same manner with a topology. We denote by $\dbF_p$ the field with $p$ elements and by $\dbF$ its algebraic closure.

A continuous action of a totally discontinuous locally compact group  on a scheme $S$ is always in the sense of [De2, 2.7.1] and is a right action. The purity theorem stated in [SGA1, p. 275] is referred as the classical purity theorem. A quasi-projective or projective morphism is always understood in the sense of [Hart].

For every free module $M$ of finite rank over a commutative ring R we denote by $M^*$ its dual. For any non-negative integer $n$, we denote by $M^{\otimes n}$ the tensor product of $n$-copies of $M$. By the tensor algebra of $M$ we mean $\oplus_{n\in\dbN\cup \{0\}} M^{\otimes n}$. If $v_{\alpha}\in M^{\otimes n}\otimes M^{*\otimes m}$, with $n$ and $m$ non-negative integers, we denote by $\deg(v_{\alpha}):=n+m$ its degree. A family of tensors of the tensor algebra  of $M$ is usually denoted in the form $(v_{\alpha})_{\alpha\in\scrJ}$, with $\scrJ$ a set. A bilinear form on $M$ is called perfect if it induces an isomorphism from $M$ into its dual $M^*$. Occasionally we also denote by $K^*$ the group of invertible elements of a field $K$. A pair $(M,\psi)$ with M as above and with $\psi$ a perfect alternating form on it, is called a symplectic space over $R$. We  use the same notation for two perfect alternating forms if they are obtained one from another by extension of scalars.

For a finite surjective \'etale morphism $\Spec(R_1)\to\Spec(R_0)$ and for a
reductive group $G$ over $R_1$, $\Res_{R_1/R_0}G$ denotes the
reductive group over $R_0$ obtained from $G$ by restriction of
scalars.

For an  abelian variety $A$ over a field $k$ of characteristic zero we denote by $V_f(A)$ the free $\dbA_f$-module $(\varprojlim\,\ker(n\colon A_{\kbar}\to A_{\kbar}))\otimes_{\dbZ} \dbQ$. We use freely the terminology of Hodge cycles of $A$ used in [De3]. A polarization of an  abelian scheme $A$ over a scheme $Y$ is usually denoted by $p_A$ (or $p_Y$), and by abuse of notation we still denote by $p_A$ (resp. $p_Y$) the different maps on the cohomologies (homologies) of $A$ induced by it. A  pair of the form  $(A,p_A)$ (or $(A,p_Y)$) always denotes a polarized abelian scheme over $Y$. For an abelian scheme $A$ over $Y$, $A^t$ denotes the dual abelian scheme of $A$, while for any $N\in\dbN$, we denote by $A[N]$ the finite flat group scheme over $Y$ defined by the $N$-torsion points of $A$. By a level-N structure of an abelian scheme $A$ (over  $Y$) of dimension $d$, we mean an isomorphism $k\colon L(N)_Y\arrowsim A[N]$ of finite group schemes over $Y$; if moreover $A$ has a principal polarization $p_A$, then by a level-N symplectic similitude structure of $(A,p_A)$, we mean a similitude (symplectic) isomorphism $(L(N)_Y,\psi)\arrowsim (A[N],p_A)$. Here $(L(N),\psi)$ is a symplectic space over $\dbZ/N\dbZ$ of dimension $2d$; $L(N)$ is viewed as a finite flat group scheme over $\Spec(\dbZ)$. If $A$ is an abelian variety over $\dbC$ then $H^i(A,\dbQ)$, $H_i(A,\dbZ)$, etc., $i\in\dbN\cup\{0\}$, refer to groups of the Betti cohomology and homology of $A$.

We will have four more sections of notations at appropriate moments: 3.2.6, 4.3.3, 5.7.2 and 6.6.1.
\finishproclaim

\smallskip
\Proclaim{2.2. The torus $\dbS$.} \rm
Let $\dbS$ be $\Res_{\dbC/\dbR}\dbG_m$. We have:
$\dbS(\dbR)=\dbC^*$ and $\dbS(\dbC)=\dbC^*\times\dbC^*$.
The last identification is made in such a way that the inclusion
$\dbR\hookrightarrow \dbC$ induces $z\to(z,\zbar)$.
To $H$ an algebraic group over $\dbR$ and to a homomorphism
$x\colon\dbS\to H$, we associate two homomorphisms of algebraic groups:
$
\mu_x\colon \dbG_m\to H_\dbC$, given on complex points by $z\to x_\dbC(z,1),\,\,
z\in \dbG_m(\dbC)=\dbC^*$, and (the weight homomorphism) $w_x\colon \dbG_m \to H$, given on real points by $r\to x(r)^{-1},\,\, r\in\dbG_m(\dbR)=\dbR^*\subset
\dbC^*=\dbS(\dbR)$.
\finishproclaim

\smallskip
\Proclaim{2.3. Definition of a (complex) Shimura variety.} \rm 
A Shimura variety is defined by a pair $(G,X)$, called a Shimura pair, comprising from a
reductive group $G$ over $\dbQ$ and from a $G(\dbR)$-conjugacy
class $X$ of homomorphisms $\dbS\to G_\dbR$ satisfying the
following axioms:

\medskip
\parindent=30pt
\item{(SV1)}
for each $x\in X$, the Hodge structure on the Lie algebra
$\grg$ of $G$ defined by $\Ad\circ x\colon \dbS\to
GL(\grg_\dbR)$ is of type $\{(1,-1),(0,0),(-1,1)\}$;

\smallskip
\item{(SV2)}
for each $x\in X$, $\ad\,x(i)$ is a Cartan involution of
$G_\dbR^{\ad}$;

\smallskip
\item{(SV3)}
$G^{\ad}$ has no factor defined over $\dbQ$ whose real
points form a compact group.

\medskip
\parindent=25pt
Let $x\in X$. 
Let $K_\infty$ be the  subgroup of $G(\dbR)$ fixing $x$. It is a maximal compact subgroup of $G(\dbR)$ iff $G^{\ab}(\dbR)$ is compact (cf. SV2). We have $X=G(\dbR)/K_\infty$, with $x$ corresponding to the equivalence class of the identity element.

Axiom SV1 implies that the homomorphism $w_x$ is independent of $x\in X$.
We write it $w_X$.
It is called the weight of the Shimura variety defined by $(G,X)$.
Axiom SV1 also implies (cf. [De2, 1.1.14]) that $X$ has only one complex
structure such that, for all representations $\rho\colon
G_\dbR\to GL(W_{\dbR})$, with $W_{\dbR}$ a finite dimensional real vector space, the Hodge filtration $F(\rho\circ x)$ of
$W_\dbR\otimes\dbC$  depends holomorphically on $x\in X$.
This complex structure is $G(\dbR)$-invariant and the
connected components of $X$ are Hermitian symmetric domains (cf. [De2, 1.1.17]).

For each compact open subgroup $K$ of $G(\dbA_f)$
$$
\Sh_K(G,X):=G(\dbQ)\setminus X\times G(\dbA_f)/K
$$
is a finite disjoint union of quotients of $X$ by arithmetic
subgroups.
This complex space has a
natural structure of a quasi-projective (algebraic) variety over
$\dbC$ [BB], which is smooth if $K$ is small enough. In what follows $\Sh_K(G,X)$ is identified with this quasi-projective variety. 

For $K\subset L$ compact open subgroups of $G(\dbA_f)$, we
get a finite surjective morphism (of schemes) $f(L,K)\colon
\Sh_K(G,X)\to \Sh_L(G,X)$ defined by $[x,a]\to[x,a]$ ($x\in X$, $a\in G(\dbA_f)$).
If $K_1=gKg^{-1}$ with $g\in G(\dbA_f)$ we get an
isomorphism $f(K,g)\colon \Sh_K(G,X)\arrowsim
\Sh_{K_1}(G,X)$ defined by $[x,a]\to[x,ag^{-1}]$. The isomorphisms $f(K,g)$ with $g\in K$ are identity automorphisms. 

The (complex) Shimura variety $\Sh(G,X)$ is the projective limit of
the compatible system of morphisms $f(L,K)$  together with the (right) continuous action of
$G(\dbA_f)$ on it defined by the rule $[x,a]g=[x,ag]$.
The continuity property of this action implies that if $K$ is a normal subgroup of $L$,
 then $f(L,K)$
identifies $\Sh_L(G,X)$ with the quotient of $\Sh_K(G,X)$
by $L/K$ (this group acts on it through isomorphisms $f(K,g)$, with
$g\in L$). The  dimension of $\Sh(G,X)$ is the dimension of $X$ as a complex manifold.

We have 
$$
\hbox{\rm \Sh}(G,X)(\dbC)=G(\dbQ)\setminus X\times G(\dbA_f)/\overline{Z(\dbQ)},
$$
where $Z=Z(G)$ and $\overline{Z(\dbQ)}$ is the topological closure of $Z(\dbQ)$ in $Z(\dbA_f)$ ([De2, 2.1]).
\finishproclaim

\smallskip
\Proclaim{2.4. Definition of maps between Shimura varieties.} \rm
The maps from a Shimura pair $(G,X)$ into another Shimura pair $(G_1,X_1)$ are  group
homomorphisms $f\colon G\to G_1$ taking $X$ into $X_1$. We denote such a map by $f\colon(G,X)\to(G_1,X_1)$.
The maps from the Shimura variety defined by $(G,X)$ into the Shimura variety defined by $(G_1,X_1)$ are in one to one correspondence with the maps $f\colon (G,X)\to (G_1,X_1)$. The induced map $X\to X_1$ is holomorphic.
If $K$ is a compact open subgroup of $G(\dbA_f)$ and if $K_1$
is a compact open subgroup of $G_1(\dbA_f)$ such that
$f(K)\subset K_1$, then the map $f$ induces a morphism of schemes (cf. [BB]) 
$f(K_1,K)\colon$\break$ \hbox{\rm \Sh}_K
(G,X)\to\hbox{\rm \Sh}_{K_1}(G_1,X_1)$ by the rule
$[x,a]\to[f(x),f(a)]$.
Passing to the limit we get the (map or) morphism between Shimura varieties (associated to $f$ and still denoted by $f$) $f\colon
\hbox{\rm \Sh}(G,X)\to\hbox{\rm \Sh}(G_1,X_1)$. Sometimes we work with
a map between Shimura pairs  $f\colon (G,X)\to (G_1,X_1)$ and sometimes we work with the (map or) morphism between Shimura varieties associated to it $f\colon\Sh(G,X)\to\Sh(G_1,X_1)$. The map $f$ is called injective (or an embedding) if it is injective as a
group homomorphism; is called finite if the induced homomorphism at the level of derived groups is an isogeny; is called a cover if it is finite and  as a group homomorphism is surjective, having as kernel a torus $T$ satisfying $H^1(\Gal(\kbar/k),T(\kbar))=0$, for any field $k$ of characteristic zero. If $f\colon (G,X)\to (G_1,X_1)$ is a finite map, then we identify  $X$ with a disjoint union of connected components of $X_1$.

Warning: If $f\colon (G,X)\to (G_1,X_1)$ is a finite map, then we sometimes refer to $f\colon\Sh(G,X)\to\Sh(G_1,X_1)$ as a morphism (of schemes), and sometimes as a finite map (of Shimura varieties), though as a morphism it is not finite, being just pro-finite. 
\finishproclaim

\Proclaim{2.4.0. Products.} \rm
The category $Sh$ whose objects are Shimura varieties and whose morphisms are morphisms between them has finite products: If $\Sh(G_i,X_i)$, $i=\overline{1,2}$, are two Shimura varieties, then their product $\Sh(G_1,X_1)\times \Sh(G_2,X_2)$ is the Shimura variety defined by $G=G_1\times G_2$ and $X=X_1\times X_2$ (together with the logical projections defined by the projections of $G$ onto its factors $G_1$ and $G_2$).

Let $f_i\colon(G_i,X_i)\to(G,X)$ be finite maps, $i=\overline{1,2}$. So $X_1$ and $X_2$ are disjoint unions of connected components of $X$. Let $X_3$ be their intersection. It can happen that $X_3$  is empty (for an example see 2.5.1). We assume now that $X_3$ is not empty.  
Let $G_3$ be the connected component of the origin of $G_1\times_G G_2$. So $X_3$ is a set of homomorphism $\dbS\to G_{3\dbR}$ satisfying the axioms SV1 and SV2. The group $G_3(\dbR)$ acts on $X_3$ by conjugation. 
 Let $X_3=\cup_{j\in I} X_3^j$ be the disjoint union decomposition of $X_3$ into $G_3(\dbR)$-orbits. For any $j\in I$ we get a Shimura variety $\Sh(G_3,X_3^j)$ and a commutative diagram 
$$
\spreadmatrixlines{1\jot}
\CD
\Sh(G_3,X_3^j) @>{p^j_2}>> \Sh(G_2,X_2)\\
@V{p_1^j}VV @VV{f_2}V\\
\Sh(G_1,X_1) @>{f_1}>> \Sh(G,X).
\endCD
$$
The morphisms $p_i^j$, $j\in I$, are defined by the natural projections of $G_3$ on $G_i$, $i=\overline{1,2}$. 
 
We have a universal type property: for any pair $(p_1,p_2)$ of finite maps $p_i\colon (G_0,X_0)\to (G_i,X_i)$ such that $f_2\circ p_2=f_1\circ p_1$, there is a unique $j\in I$ for which there is a map $p_0\colon (G_0,X_0)\to (G_3,X_3^j)$ such that $p_i=p_i^j\circ p_0$; moreover the map $p_0$ is uniquely determined. We express this property by: the category f-$Sh$ whose objects are Shimura varieties (or pairs) and whose morphisms are the finite maps between them, has quasi fibre products. Any commutative diagram (or pair $(p_1^j,p_2^j)$)  as above (formed by finite maps) is called a quasi fibre product of the finite maps $f_1$ and $f_2$. 

So if $X_3$ is empty then $I$ is the empty set. There are examples (cf. 2.5.1) when $I$ has more than one element. However if $f_1$ or $f_2$ is a cover, then $I$ has precisely one element (cf. [Mi4, 4.11]): if this is the case we speak about the fibre product of $f_1$ and $f_2$. 
\finishproclaim

\Proclaim{2.4.1. The adjoint and toric part varieties of a Shimura variety.} \rm
Let $(G,X)$ define an arbitrary Shimura variety. Then $\hbox{\rm \Sh}(G^{\ad},X^{\ad})$ ($X^{\ad}$ being the ${G^{\ad}}(\dbR)$-conjugacy class of homomorphisms $\dbS\to G^{\ad}_{\dbR}$ containing the ones induced by $X$) is called the adjoint variety of $\hbox{\rm \Sh}(G,X)$, and $\Sh(G^{\ab},X^{\ab})$ ($X^{\ab}$ being the set with just one element defined by the homomorphism $\dbS\to G_{\dbR}^{\ab}$ induced by $X$) is called the toric part variety of $\Sh(G,X)$. We have natural maps from every Shimura variety into its adjoint variety and into its toric part variety. 
\finishproclaim

\Proclaim{2.4.2. Special pairs.} \rm
An injective map $(T,\{h\})\hookrightarrow (G,X)$ with $T$ a torus is called a special pair in $(G,X)$. 
\finishproclaim

\Proclaim{2.4.3. Automorphisms.} \rm
The group $\Aut(Sh(G,X))$ (of automorphisms of the Shimura variety $\Sh(G,X)$) is the subgroup of $Aut(G)(\dbQ)$ (it is of finite index if $G$ is an adjoint group) leaving $X$ invariant. If $G$ is adjoint and all simple factors of $(G,X)$ are such that [De2, 1.2.8 (ii)] applies, then we have $\Aut(Sh(G,X))=Aut(G)(\dbQ)$. 

\smallskip
\Proclaim{2.5. Examples of types of Shimura varieties.} \rm
\finishproclaim

\Proclaim{\it Example 1.}\rm
Let $T$ be a torus over $\dbQ$.
For any homomorphism $h\colon\dbS\to T_\dbR$, the pair $(T,\{h\})$
satisfies the axioms SV1 to SV3, and so defines a Shimura
variety of dimension 0.
We have $\hbox{\rm \Sh}(T,\{h\})(\dbC)=T(\dbA_f)/\overline{T(\dbQ)}$. Any Shimura variety of dimension 0 is obtained in this way.
\finishproclaim

\Proclaim{\it Example 2.}\rm
Let $(W,\psi)$ be a symplectic space over $\dbQ$. 
Let $GSp:=G\,Sp(W,\psi)$ be the group of its symplectic 
similitudes.
The Siegel double space $S$ consists of all rational Hodge
structure on $W$ of type $\{(-1,0),(0,-1)\}$ for which
either $2\pi i\psi$ or $-2\pi i\psi$ is a polarization.
It is a $G\Sp(\dbR)$-conjugacy class of homomorphism $\dbS\to GSp_{\dbR}$. The pair $(GSp,S)$ defines a Shimura variety.
The Shimura varieties of the form $\hbox{\rm \Sh}(GSp,S)$ 
are called Siegel modular
varieties.
\finishproclaim

\Proclaim{\it Definition 1.}\rm
A Shimura variety $\hbox{\rm \Sh}(G,X)$ 
is said to be of Hodge type if
there is an injective map from it into a Siegel modular variety. We have $\Sh(G,X)(\dbC)=G(\dbQ)\setminus X\times G(\dbA_f)$ [De2, 2.1.1]. 
\finishproclaim

The extra conditions needed to be satisfied by a Shimura
variety for being of Hodge type are:

\smallskip
\parindent=30pt
\item{(SVH4)}
the weight is defined over $\dbQ$;

\smallskip
\item{(SVH5)}
$w_X(\dbG_m)$ is the only split subtorus of $Z(G)_\dbR$;

\smallskip
\item{(SVH6)}
there is a faithful representation $\rho\colon G\hookrightarrow GL(W)$
such that the Hodge $\dbQ$--structure on $W$ defined by $\rho_{\dbR}\circ x$ is of type
$\{(-1,0),(0,-1)\}$, $\forall x\in X$.

\smallskip
This is just a reformulation of [De2, 2.3.2]: obviously SVH4 to SVH6 are satisfied by a Shimura variety of Hodge type, while  SV2 and SVH5 put together imply that for any $x\in X$, the inner automorphism of $G_{\dbR}/w_X(\dbG_m)$ defined by $x(i)$ is a Cartan involution. 

\Proclaim{\it Example 3.}\rm
The product of two Shimura varieties $\Sh(G_1,X_1)$ and $\Sh(G_2,X_2)$ of Hodge type is not of Hodge type. But the Shimura variety $\Sh(G_3,X_3)$ defined by the subgroup $G_3$ of $G_1\times G_2$ generated by $G_1^0\times G_2^0$ (with $G_i^0$ the connected subgroup of $G_i$ having the property that the quotient  homomorphism $G_i^0\to G_i/w_{X_i}(\dbG_m)$ is an isogeny, $i=\overline{1,2}$) and $w_{X_1\times X_2}(\dbG_m)$, and by an adequate union $X_3$ of some of the connected components of $X_1\times X_2$, is a Shimura variety of Hodge type: it is enough to see this for the case when $(G_1,X_1)=(GSp(W_1,\psi_1),S_1)$ and $(G_2,X_2)=(GSp(W_2,\psi_2),S_2)$; but then we have an injective map 
$$
i_3\colon (G_3,X_3)\hookrightarrow (GSp(W_1\oplus W_2, \psi_1\oplus\psi_2),S^0)$$ 
defined by the natural inclusions of $Sp(W_1,\psi_1)$ and $Sp(W_2,\psi_2)$ in $Sp(W_1\oplus W_2,\psi_1\oplus\psi_2)$ ($X_3$ in this case has two connected components, while $X_1\times X_2$ has four). 

We refer to the map $i_3$ as a Segre embedding, and to any pair $(G_3,X_3)$ as above (we do not have a unique choice for $X_3$; this is the same as the case of quasi fibre products discussed in 2.4.0 --see also 2.5.1 below--) as the Hodge quasi product of the two pairs $(G_1,X_1)$ and $(G_2,X_2)$ of Hodge type. Similarly we speak about a Hodge quasi product of $n$ Shimura varieties of Hodge type and the Segre embedding defined by the product of $n$ Siegel modular varieties, $n\in\dbN$.
\finishproclaim

\Proclaim{\it Definitions 2.}\rm
A Shimura variety defined by a pair $(G,X)$ with $G$ an adjoint group is said to be an adjoint Shimura variety or of adjoint type. If $G$ is a simple $\dbQ$--group, then $(G,X)$ is of one of the types: $A_l$, $B_l$, $C_l$,  $D_l^{\dbR}$,  $D_l^{\dbH}$,  $D_l^{\text mixed}$, $E_6$ or $E_7$ (cf. the classification [De2] of Shimura varieties of adjoint type); $\Sh(G,X)$ is called a simple adjoint Shimura variety (of $A_l$, or $B_l$, etc. type). Any adjoint Shimura variety is a product of a finite number of simple adjoint Shimura varieties. A Shimura variety is said to be of special type if its adjoint Shimura variety is a product of simple adjoint Shimura varieties of $E_6$, $E_7$ or $D_l^{\text mixed}$ type.
\finishproclaim 

\Proclaim{\it Definitions 3.}\rm
A Shimura variety $\Sh(G,X)$ is called of preabelian type if there is a Shimura variety $\Sh(G_1,X_1)$ of Hodge type such that their adjoint varieties are isomorphic. If we can choose $\Sh(G_1,X_1)$ such that $G_1^{\der}$ is a cover of $G^{\der}$, then $\Sh(G,X)$ is called of abelian type. The simple adjoint Shimura varieties of abelian type are those of  $A_l$, $B_l$, $C_l$, $D_l^{\dbR}$ or $D_l^{\dbH}$ type [De2, 2.3.8]. The product of two Shimura varieties of abelian (preabelian) type is of abelian (resp. of preabelian) type, cf. Example 3.
\finishproclaim

So for any pair $(G,X)$ which is neither of preabelian nor of special type, there is a finite map $(G,X)\to (G_1,X_1)$, with $\Sh(G_1,X_1)$ an adjoint variety which is the product of a Shimura variety of preabelian type and of a Shimura variety of special type. 
The category f-$Sh$ is a disjoint sum of categories indexed by isomorphism classes of Shimura varieties of adjoint type.

\Proclaim{\it Example 4.}\rm
A Shimura variety of dimension 1 is called a Shimura curve and a Shimura variety of dimension 2 is called a Shimura surface. For instance the Example 2, gives birth to a Shimura curve if $W$ is a vector space over $\dbQ$ of dimension 2, called the elliptic modular curve. 
\finishproclaim

\Proclaim{\it Definition 4.}\rm
A Shimura pair $(G,X)$ (resp. variety $\Sh(G,X)$) is said to be of compact type if $\Sh(G,X)$ is a pro-\'etale cover of a smooth projective $E(G,X)$-scheme.
\finishproclaim

In [BHC] it is proved: $(G,X)$ is of compact type iff the $\dbQ$--rank of $G^{\ad}$ is zero.  

\Proclaim{\it Example 5.}\rm
$(G,X)$ is of compact type if $G^{\ad}$ is a simple $\dbQ$--group such that $G^{\ad}_{\dbR}$ has compact factors.
\finishproclaim

\Proclaim{2.5.1. Extra example.} \rm
Let $(G,X)$ be such that the semisimple group $G^{\der}$ is simply connected, $G^{\ab}=\dbG_m$, and $X$ has  precisely two connected components (for instance if $\Sh(G,X)$ is a Siegel modular variety; see also 5.7.5). Let $(G_1,X_1)$ be the product of three copies of $(G,X)$. So $G_1^{\der}$ is simply connected, and $G_1^{\ab}=\dbG_m\times\dbG_m\times\dbG_m$. We consider reductive  subgroups $G_i$ of $G_1$, $i=\overline{2,4}$, containing $G_1^{\der}$. So to give such a $G_i$ is the same as to give a subtorus of $G_1^{\ab}$. We choose $G^{\ab}_i\subset G_1^{\ab}$, $i=\overline{2,4}$, to be  the diagonal embedding of $\dbG_m$, the subtorus generated by $G_2^{\ab}$ and the embedding $\dbG_m\hookrightarrow G_1^{\ab}$ corresponding to the triple of characters $(1,1,0)$ of $\dbG_m$, and respectively the subtorus generated by $G_2^{\ab}$ and the embedding $\dbG_m\hookrightarrow G_1^{\ab}$ corresponding to the triple of characters $(1,3,0)$.

We get injective finite maps $f_i\colon (G_i,X_i)\hookrightarrow (G_1,X_1)$, $i=\overline{2,4}$. Here $X_1$ has eight connected components, $X_2$ has two, while $X_3$ and $X_4$ have four. Moreover we can assume that $X_2\subset X_3=X_4$. So the maps $f_3$ and $f_4$ do not have a fibre product: they have two quasi fibre products. 

Moreover, as $G^{\ad}(\dbQ)$ is dense in $G^{\ad}(\dbR)$ (cf. [De1, 0.4]), composing the natural map $p_2\colon (G_2,X_2)\to (G_2^{\ad},X_2^{\ad})$ with a suitable automorphism  of $(G_2^{\ad},X_2^{\ad})$, we get a map $p_3\colon (G_2,X_2)\to (G_2^{\ad},X_2^{\ad})$ such that the images of $X_2$ in $X_2^{\ad}=X_1$ through the maps $p_2$ and $p_3$  have an empty intersection.  
\finishproclaim

\smallskip
\Proclaim{2.6. The reflex field.}\rm
Let $(G,X)$ be an arbitrary Shimura pair.
For any field $k$ of characteristic zero we have a right action (via conjugation) of $G(k)$ on the set $\Hom(\dbG_m,G_k)$. Let
$C(k):=\Hom(\dbG_m,G_k)/G(k)$.
An inclusion $\dbQbar\hookrightarrow\dbC$ induces a
bijection $C(\dbQbar)=C(\dbC)$.
So the element $[\mu_X]\in C(\dbC)$, corresponding
to $\mu_x$ for any $x\in X$, defines an element
$c(X)$ of $C(\dbQbar)$.
The group $\Gal(\dbQbar/\dbQ)$ acts on
$C(\dbQbar)$.
The reflex field $E(G,X)$ of the Shimura variety $\Sh(G,X)$ is the
subfield of $\dbQbar$ corresponding to the stabilizer of
$c(X)$ in $\Gal(\dbQbar/\dbQ)$.
It is a finite extension of $\dbQ$.
\finishproclaim

\smallskip
\Proclaim{2.7. The reciprocity map.} \rm
Let $(T,\{h\})$ be as in Example 1 of 2.5.
Its reflex field $E:=E(T,\{h\})$ is the field of definition
of the cocharacter $\mu_h$ of $T$.
From the homomorphism $\mu_h\colon\dbG_{m_E}\to T_E$ we get a new one
$$
N_h\colon \Res_{E/\dbQ}\dbG_{m_E} @>{\Res_{E/\dbQ}(\mu_h)}>>
\Res_{E/\dbQ}T_E @>{\text{Norm}\,E/\dbQ}>> T.
$$
So, for any $\dbQ$--algebra $A$ we get a homomorphism
$N_{h}(A)\colon\dbG_m(E\otimes A)\to T(A)$.

The reciprocity map 
$$
r(T,\{h\})\colon\Gal(E^{\ab}/E)\to
T(\dbA_f)/\overline{T(\dbQ)}
$$ 
is defined as follows:
let $\tau\in\Gal(E^{\ab}/E)$, and let $s\in\dbJ_E$ be an id\`ele (of $E$)
such that $\text{rec}_E(s)=\tau$; then
$r(T,\{h\})(\tau)=N_h(\dbA_f)(s_f)$, where $s_f$ is the finite part of $s$.
Here the Artin reciprocity map $\text{rec}_E$ is such that a
uniformizing parameter is mapped into the geometric Frobenius element.
\finishproclaim

\smallskip
\Proclaim{2.8. The canonical model of $\hbox{\rm \Sh}(G,X)$ over
$E(G,X)$.}\rm
By a model of $\hbox{\rm \Sh}(G,X)$ 
over a subfield $k$ of $\dbC$, we
mean a scheme $S$ over $k$ endowed with a continuous
action of $G(\dbA_f)$ (defined over $k$), such that there is
a $G(\dbA_f)$-equivariant isomorphism
$$
\hbox{\rm \Sh}(G,X)\arrowsim S_{\dbC}.
$$

The canonical model of $\hbox{\rm \Sh}(G,X)$ 
is the model $\scrS$ of $\hbox{\rm \Sh}(G,X)$
over $E(G,X)$ which satisfies the following property:
if $(T,\{h\})$ is a special pair in $(G,X)$ then for any
$a\in G(\dbA_f)$, the point $[h,a]$ of $\scrS(\dbC)=\Sh(G,X)(\dbC)$ is rational over
$E(T,\{h\})^{\ab}$, and every element $\tau$ of $\Gal(E(T,\{h\})^{\ab}/
E(T,\{h\}))$ acts on $[h,a]$ according to the rule
$$
\tau[h,a]=[h,ar(\tau)],
$$ 
where $r:=r(T,\{h\})$.
It exists and is uniquely determined by the above property up to a
unique isomorphism (see [De1], [De2] and [Mi2]). 

Warning: from now on by $\Sh(G,X)$ we mean $\scrS$. 
\finishproclaim

\smallskip
\Proclaim{2.9.} \rm
If $f\colon (G,X)\to (G_1,X_1)$ is a map between two Shimura
pairs, then $E(G_1,X_1)$ is a subfield of $E(G,X)$, and there
is a unique $G(\dbA_f)$-equivariant morphism (still denoted by $f$) $f\colon \Sh(G,X)\to \Sh(G_1,X_1)_{E(G,X)}$ which at the level of complex points is the map
$[x,a]\to[f(x),f(a)]$ ([De1, 5.4]). We get a $G(\dbA_f)$-equivariant morphism (still denoted by $f$) 
$$f\colon \Sh(G,X)\to \Sh(G_1,X_1)
$$ 
of $E(G_1,X_1)$-schemes.
\finishproclaim

\smallskip
\Proclaim{2.10. Definition of special points.} \rm
Let $\Sh(G,X)$ be an arbitrary Shimura variety and let $H$ be a compact subgroup of $G(\dbA_f)$. A point $w$ of $\Sh_H(G,X)$  with values in a field $k$ of characteristic zero is called special if there is a special pair $(T,\{h\})$ in $(G,X)$, such that the intersection of the $G(\dbA_f)$-orbit of $w$ in $\Sh_H(G,X)(\overline{k})$ with the image of $\Sh(T,\{h\})(\overline{k})$ in $\Sh_H(G,X)(\overline{k})$  is non-empty.
\finishproclaim

\smallskip
\Proclaim{2.11. Definition of smooth  subgroups.} \rm
Let $(G,X)$ be a Shimura pair. A  subgroup $H$ of $G(\dbA_f)$ is called smooth for $(G,X)$ if it is compact and if $\Sh(G,X)$ is a pro-\'etale cover of $\Sh_H(G,X)$. A subgroup of a $G(\dbA_f)$-conjugate of a subgroup of $G(\dbA_f)$ smooth for $(G,X)$, is itself smooth for $(G,X)$. For instance, any neat compact subgroup of $G(\dbA_f)$ is smooth for $(G,X)$.  We do not know when the converse is true. We include two examples which point out that the converse is not always true.

\medskip\noindent
{\bf Example 1.} If $G$ is a torus, then any compact subgroup of $G(\dbA_f)$ is smooth for $(G,X)$ (see 3.2.8).

\noindent
{\bf Example 2.} We assume $G$ is a $\dbQ$--simple, adjoint group. We assume $G_{\dbQ_p}$ has two non-trivial factors $G_1$ and $G_2$ and there is a non-trivial, finite subgroup $H$ of $G_1(\dbQ_p)$. Then $H$ is smooth for $(G,X)$. This is a consequence of the structure of the set $\Sh(G,X)(\dbC)$ (see [De2, 2.1.1]): as $G(\dbQ)$ has trivial intersection with any $G(\dbA_f)$-conjugate of $H$, any $g\in H$ acts freely on this set and so on $\Sh(G,X)$. We get: there are compact open subgroup of $G(\dbA_f)$ which contain $H$, are smooth for $(G,X)$ but are not neat.

\smallskip
{\bf Definition.} {\bf a)} $H$ is called S-smooth for $(G,X)$ (here S stands for strongly) if it is smooth and for each connected component $\scrC$ of ${\Sh}(G,X)_{\dbC}$ there is a compact, open subgroup $H_1$ of $G(\dbA_f)$ containing the subgroup $H(\scrC)$ of $H$ leaving invariant $\scrC$ and such that the set of complex points of the connected component of ${\Sh}_{H_1}(G,X)_{\dbC}$ dominated naturally by $\scrC$ can be identified (see [De2, p. 266--267]) with $\Sigma\setminus X^0$, with $X^0$ a connected component of $X$ and with $\Sigma$ an arithmetic subgroup of $G^{\ad}(\dbQ)$ which does not have $2$-torsion. 
\smallskip
{\bf b)} Let $p$ be a rational prime. We say $H$ is $p$-smooth if it is smooth for $(G,X)$ and if for each connected component $\scrC$ of ${\Sh}(G,X)_{\dbC}$ there is a rational prime $l(\scrC)$ different from $p$ and such that the image of any pro-$p$ subgroup $H_p(\scrC)$ of $H(\scrC)$ (with $H(\scrC)$ as in a)) in $G^{\ad}(\dbQ_{l(\scrC)})$ is trivial. If $\scrS$ is a set of rational primes having at least two elements, we say $H$ is $\scrS$-smooth, if it $p$-smooth for any $p\in\scrS$.
\finishproclaim

\smallskip
\Proclaim{2.11.1. Torsion numbers.} \rm
For any reductive group $\tilde G$ over $\dbQ$ we denote by $\scrU(\tilde G)$ the set of primes $l$ such that $\tilde G$ is unramified over $\dbQ_l$. 
If $\tilde G$ is a reductive group over a DVR of mixed characteristic $(0,p)$ or over the field of fractions of such a DVR, let $t(\tilde G)\in\dbN$ be the biggest non-negative integral power of $p$ dividing the order of an element of $\tilde G(R)$ of finite order. If $G$ is a reductive group over $\dbQ$, then its torsion number $t(G)\in\dbN$ is defined by
$$t(G):=\prod_{p\in\scrU(G)} t(G_{W(\dbF)}) \prod_{p\notin\scrU(G),\,p\,{\text{a prime}}} t(G_{W(\dbF)[{1\over p}]}).$$
\finishproclaim

\smallskip
\Proclaim{2.12. Remarks. 1)} \rm
For any Shimura pair $(G,X)$ there are finite maps $f\colon (G_1,X_1)\to (G,X)$ and $f_1\colon (G_1,X_1)\to (G_2,X_2)$ such that:

\smallskip
-- $(G_2,X_2)$ is a product of Shimura pairs $(G_i,X_i)$, $i$ running through the elements of a finite set $I$, such that $G_i^{\ad}$ is a simple $\dbQ$--group, $\forall i\in I$; 

-- they define a quasi fibre product of the natural maps $f_0\colon  (G,X)\to (G^{\ad},X^{\ad})$ and $f_2\colon (G_2,X_2)\to (G_2^{\ad},X_2^{\ad})=(G^{\ad}, X^{\ad})$;  

-- there are injective maps $(G_i,X_i)\hookrightarrow (G,X)$, $i\in I$, producing (naturally) an isogeny $\prod_{i\in I} G_i^{\der}\to G^{\der}$.

\smallskip
To see this let $G^{\ad}=\prod_{i\in I} G^{\ad}_i$ be the factorization of $G^{\ad}$ in $\dbQ$--simple factors. Let $G_i^{\der}$ be the semisimple subgroup of $G$ isogenous to $G^{\ad}_i$ and contained in the kernel of the canonical quotient homomorphism $G\to\prod_{j\in I\setminus\{i\}} G_j^{\ad}$. As $G_i$ we take the subgroup of $G$ generated by $G_i^{\der}$ and by a maximal torus of the centralizer of $G_i^{\der}$ in $G$ having the property that there is a homomorphism $\dbS\to G_{\dbR}$, corresponding to a point $x\in X$, factoring through $G_{i\dbR}$. As $X_i$ we take the $G_i(\dbR)$-conjugacy class of homomorphisms $\dbS\to G_{i\dbR}$ generated by such a factorization. 
Now we can take the maps $f$ and $f_1$ to define a quasi fibre product of the maps $f_0$ and $f_2$ (cf. 2.4.0).

2) There are Shimura varieties $\Sh(G,X)$ with $G$ a semisimple group which is not of adjoint type (plenty of examples can be constructed starting from [De2, 2.3.12]).

3) Let $(A,p_A)$ be a polarized  abelian scheme defined over an integral ring $R$ of characteristic zero. It is defined over a subring $R_1$ of $R$ admitting an embedding in $\dbC$. We get an abelian variety over $\dbC$. Passing to an isogeny we can assume that we have a principally polarized abelian variety $(A^{\prime},p_{A^{\prime}})$ over $\dbC$. Let $(W,\psi)$ be the symplectic space over $\dbQ$ defined by it, with $W:=H_1(A^{\prime},\dbQ)$. Let $G$ be the Mumford--Tate group of $A^{\prime}$. We get an injective map $(G,X)\hookrightarrow (GSp(W,\psi),S)$ of Shimura pairs, with $X$ the Hermitian symmetric domain defined by the $G(\dbR)$-conjugacy class of the homomorphism $\dbS\to G_{\dbR}$ defining the Hodge $\dbQ$--structure on $W$. We call it an injective map defined by the polarized abelian scheme $(A,p_A)$. If $A$ does not have a priori a polarization then we can pick one for its model we got over $\dbC$. 

If moreover $A$ has a family $(v_{\alpha})_{\alpha\in\scrJ}$ of Hodge cycles which is reductive with respect to the polarization $p_A$ (to be explained below), then we can choose $R_1$ such that these Hodge cycles are defined over $R_1$ (cf. [De3]). So we get a family of tensors $(t_{\alpha})_{\alpha\in\scrJ}$ of the tensor algebra of $W$ (we can assume that no Tate-twist shows up; for instance cf. 4.1). By reductive family with respect to $p_A$  we mean: the connected component of the origin of the subgroup of $GSp(W,\psi)$ fixing $t_{\alpha}$, $\forall\alpha\in\scrJ$, is a reductive group $G_1$ over $\dbQ$. The group $G_1$ together with the $G_1(\dbR)$-conjugacy class of homomorphisms $\dbS\to G_{1\dbR}$ defined by $X$ might not define a Shimura variety: axiom $SV3$ might not be satisfied. However, discarding from $G_1^{\ad}$ the factors which over $\dbR$ are compact, we get a reductive subgroup $G_2$ of $G_1$, which together with the $G_2(\dbR)$-conjugacy class $X_2$ of homomorphisms $\dbS\to G_{2\dbR}$ defined by $X$, define a Shimura variety; we have $G_2^{\ab}=G_1^{\ab}$ and $G_2^{\ad}$ is the product of the simple factors of $G_1^{\ad}$ which over $\dbR$ are non-compact. The resulting injective map $(G_2,X_2)\hookrightarrow (GSp(W,\psi))$ is call an injective map defined by $(A,p_A)$ and the reductive family of tensors $(v_{\alpha})_{\alpha\in\scrJ}$ with respect to $p_A$.  
\finishproclaim

\bigskip
\noindent
{\boldsectionfont \S3. A general view of the integral models of Shimura varieties}

\bigskip
We start by presenting in 3.1 some elements of the theory of reductive groups (and of hyperspecial subgroups) needed for applications to Shimura varieties. Then in 3.2 we introduce generalities of the theory of healthy normal schemes and of the theory of integral models of Shimura varieties. Some special features of these theories are presented in 3.3 to 3.5.

\smallskip
\Proclaim{3.1. Hyperspecial subgroups.}\rm
We restrict ourselves to what we need.
Let $V$ be a complete DVR with a perfect residue field $k$ and let $K$ be its field of fractions. Let $\pi$ be a uniformizer of $V$.
Let $G_K$ be a reductive group over $K$.
A subgroup $H$ of $G_{K}(K)$ is called hyperspecial  if there is a reductive group scheme $G$ over $V$, whose generic fibre  is $G_K$ and whose group of $V$-valued
points is $H$. It is a maximal bounded (compact if the residue field of $V$ is finite) subgroup of $G_{K}(K)$ [Ti, 3.2]. Let $\dbB$ be the building of $G_K$ over $K$ (cf. [Ti, 2.1 and 2.2]). The group $G_{K}(K)$ acts on $\dbB$. A subgroup $H$ of $G_K(K)$ is hyperspecial iff there is a hyperspecial point $x_H\in\dbB$ (see [Ti, 1.10.2 and 2.4] for the definition of such a point) such that $H=\{h\in G_{K}(K)|h(x_{H})=x_{H}\}$ (cf. [Ti, 3.8.1]).
Hyperspecial subgroups of $G_K(K)$ do exist if  $G_K$
is unramified over $K$, i.e. if $G_K$ is quasi-split and
splits over an unramified extension of $K$ ([Ti, 1.10.2]). The converse of this last statement is true if $k$ has the property that every reductive group over $k$ is quasi-split (for instance if $k$ is an algebraic extension of a finite field or if $k$ is algebraically closed).
\finishproclaim

\Proclaim{3.1.1. Remark.}\rm
Let $G_{K}\hookrightarrow {G_1}_K$ be an inclusion of reductive
groups over $K$ with $G_K^{\der}={G_1}_K^{\der}$.
We assume that ${G_1}_K$ is unramified over $K$.
Then $G_K$ is unramified over $K$ and for any hyperspecial subgroup $H_1$ of ${G_1}_{K}(K)$,
the intersection $H:=H_1\cap G(K)$ is a hyperspecial subgroup of $G_{K}(K)$ (if $G_1$ is a reductive group scheme over $V$, whose generic fibre is ${G_1}_K$ and whose group of $V$-valued points is $H_1$, then the Zariski closure of $G_K$ in $G_1$ is a reductive group scheme over $V$, whose group of $V$-valued points is $H$; the ideas needed for proving this are presented in 4.3.9).  
\finishproclaim

\Proclaim{3.1.2. The behavior of  hyperspecial subgroups with respect to homomorphisms of groups.} \rm
We digress a little bit on this subject as it is not covered in [Ti] or [Ja]. 
In this section we consider only affine group schemes of finite type over $V$ or $K$ and which have connected
fibres over $K$. 
\finishproclaim

\Proclaim{3.1.2.1. Proposition.}
a) Let $G_1$ and $G_2$ be two  smooth affine groups over $V$. Let ${f_K}\colon {G_1}_{K}\to {G_2}_{K}$ be a homomorphism such that it takes $G_1(V)$ into $G_2(V)$. If the field $k$ is infinite, then the homomorphism $f_K$ extends to a homomorphism $f\colon G_1\to G_2$.

b) a) is not true if the field $k$ is finite.

c) Let $f\colon G_1\to G_2$ be a homomorphism of smooth affine group schemes over $V$. If $G_1$ is a reductive group and if $f_K\colon {G_1}_K\to {G_2}_K$ is a closed embedding then $f$ is a closed embedding.

d) Let ${G_1}_K$ and ${G_2}_K$ be two reductive groups over $K$ such that ${G_1}_K$ is a subgroup of ${G_2}_K$ and such that they are unramified over $K$. We assume that ${G_1}_K$ is a torus which splits over an unramified extension of $K$ of odd degree, and that ${G_2}_K$ is a group of symplectic similitudes. Then there is a hyperspecial subgroup of ${G_1}_{K}(K)$ included in a hyperspecial subgroup of ${G_2}_{K}(K)$.

e) We consider a) in the case when $k$ is finite, when $G_1$ and $G_2$ are reductive groups over $V$, and when $f_K$ is a closed embedding. Then a) remains true if any one of the following two conditions is satisfied: 

\smallskip
\item{(i)} $f_K$ is an isomorphism; 

\item{(ii)} $G_1$ is a split group with a maximal split torus $T_1$ which is a closed subgroup of $G_2$, and there is a faithful representation $\rho\colon G_2\hookrightarrow GL(M)$ with $M$ a free $V$-module of dimension not bigger than the characteristic of $k$.
\finishproclaim

\proof
a) Let $G$ be the subgroup of $G_1\times G_2$ obtained by taking the Zariski closure of the graph of $f_K$. We get a homomorphism $h\colon G\to G_1$ inducing an epimorphism $G(V)\twoheadrightarrow G_1(V)$ and a homomorphism $G\to G_2$. They are defined by the projections of $G_1\times G_2$ on its factors. Let $G=\Spec(R)$ and $G_{1}=\Spec(R_1)$. To $h$ corresponds an inclusion $R_1\subset R$ which becomes an isomorphism by inverting $\pi$. 

So a) is equivalent to $R_1=R$. If $R\neq R_1$ then there is $y\in R\setminus R_1$ such that $\pi$$y\in R_1\setminus \pi$$R_1$. For any  $x\in R_1\setminus \pi$$R_1$  there is a ring homomorphism $R_1\to V$ such that $x$ goes to an invertible element of $V$ ($G_1$ is smooth over $V$ and the $k$-valued points of $G_1$ are Zariski dense in the special fibre of it, as $k$ is an infinite field [Bo, p. 218]). So, such a homomorphism $R_1\to V$, corresponding to $x=\pi$$y$, does not come from the restriction to $R_1$ of a ring homomorphism $R\to V$. This contradicts the surjectivity property of the homomorphism $G(V)\to G_1(V)$. We get $R=R_1$.

b) Example: Let $p$ be a rational prime and let $q\in\dbN\setminus\{1\}$ be  congruent to 1  mod $p-1$. Let $M$ be a free module of dimension 2 over $\dbZ_p$, and let $\{v_1,v_2\}$ be a basis of it over $\dbZ_p$. Let $\dbG_m$ be the subgroup of $GL(M)$ such that $\alpha\in\dbG_m(\dbQ_p)$ acts by multiplication with $\alpha^q$ on $v_1$ and by multiplication with $\alpha$ on $v_2$. We have $\dbG_m(\dbZ_p)\subset GL(M_1)(\dbZ_p)$, with $M_1$ the $\dbZ_p$-lattice of $M\otimes\dbQ_p$ generated by $v_1$ and $\frac {v_1+v_2}p$. But $\dbG_m$ is not a subgroup of $GL(M_1)$.

c) We can assume that the field $k$ is algebraically closed. Let $G$ be the Zariski closure of ${G_1}_K$ in $G_2$. It is a group scheme over $V$. Let ${f_1}\colon G_1\to G$ be the homomorphism induced by $f$. The homomorphism  $G_1(V)\to G(V)$ is an isomorphism, as it is a monomorphism and as $G_1(V)$ is a maximal bounded subgroup of $G_1(K)=G(K)$. The dilatation procedure [BLR, Prop. 2, p. 64] allows us to write $f_1$ as a composition $G_1\to G_n\to G$, with $n\Ge 3$ an integer, and with $G_n$ a smooth (affine) group scheme over $V$, obtained from $G$ through a finite sequence of dilatations. We have ${G_1}_K={G_n}_K=G_K$ and $G_n(V)=G(V)$. For instance, in the first step we replace $G$ by the dilatation $G_3$ of $f_1({G_1}_k)$ (the reduced group subscheme of $G_k$ defined by the image through $f_1$ of the special fibre of $G_1$) on $G$ (part $d)$ of loc. cit. shows that $G_3$ is an affine group scheme). Using [BLR, Prop. 1, p. 63], we get group homomorphisms $G_1\to G_3\to G$. Repeating the process,  we reach the case of group homomorphisms $G_1\to G_n\to G$, with $G_n$ a smooth scheme obtained from $G$ through a sequence of $n-2$ dilatations. This results from the general theory [BLR, section 3.3] of the N\'eron defect of smoothness $\rho$ of $V$-valued points of $G$. More precisely, there is a positive integer $c$, such that $\rho(x)\Le c$ for any point $x\in G(V)$ (cf. [BLR, Prop. 3, p. 66]). 
From [BLR, Lemma 4, p. 174] we get that we can take $n=c+3$.

Coming back to the situation $G_1\to G_n\to G$, we get (cf. part $a)$ above) $G_1=G_n$. But this implies that $f_1$ is an isomorphism as any dilatation of our sequence produces a commutative unipotent kernel of the special fibre. From [BLR, Prop. 2, p. 64] we see easily that it is enough to check this for the dilatation of a general linear group over $V$ with respect to the trivial subgroup of the special fibre. But this case is obvious.

We invite the reader to give another (simpler) proof of c), by just copying the proof of a) above and by using the fact that for any finite field extension $K_1$ of $K$, $G_1(V_1)$ is a maximal bounded subgroup of $G_1(K_1)$ (with $V_1$ the normalization of $V$ in $K_1$).

d) Let $T_K:=G_{1K}$ and let $G_{2K}=GSp(W_K,\psi)$, with $(W_K,\psi)$ a symplectic space over $K$. Let $K_1$ be an unramified extension  of $K$ over which $T_K$ splits. Let $V_1$ be its ring of integers. We can assume that $K_1$ is the smallest extension over which $T_K$ splits. So $K_1$ is a Galois extension of $K$. Let $C_T$ be the subset of the (additive) group of characters of $T_{K_1}$ through which $T_{K_1}$ acts on $W_K\otimes K_1$, i.e. $W_K\otimes K_1=\oplus_{\gamma\in C_T} W_{\gamma}$, with $t\in T_K(K_1)$ acting as the multiplication with $\gamma(t)$ on the non-zero $K_1$-vector space $W_{\gamma}$. We can assume that $T_K$ is a subgroup of $Sp(W_K,\psi)$ (as $GSp(W_K,\psi)$ is the extension of $Sp(W_K,\psi)$ through a one-dimensional split torus). 

As the alternating form $\psi$ is $T_K$-invariant, we deduce that if  $\alpha\in C_T$ then $-{\alpha}\in C_T$, and that $\psi(x,y)=0$ for any $x\in W_{\alpha}$ and every $y\in\oplus_{\gamma\in C_T^{\alpha}} W_{\gamma}$, where $C_T^{\alpha}:=C_T\setminus \{\alpha\}$. Moreover $\Gal(K_1/K)$ acts on $C_T$ as $T_K$ and $G_{2K}$ are  defined over $K$.
For any $\alpha\in C_T$, let $o(\alpha)$ be its orbit under the action of $\Gal(K_1/K)$ on $C_T$. 

The key fact is the following assumption (which is always satisfied if $[K_1:K]$ is an odd number): 
$$\forall\,\alpha\in C_T\setminus\{0\}, -\alpha\notin o(\alpha).$$

We can assume that $0\notin C_T$. Argument: we have a direct sum decomposition of symplectic spaces over $K$ 
$$
(W_K,\psi)=(W^0,\psi^0)\oplus (W^1,\psi^1),
$$ 
with $W^0$ the maximal subspace of $W_K$ on which $T_K$ acts trivially, with $W^1$ the subspace of $W_K$ such that $W^1\otimes K_1=\oplus_{\gamma\in C_T\setminus\{0\}} W_{\gamma}$, and with $\psi^i$ the restriction of $\psi$ to $W^i$, $i=\overline{1,2}$; so if needed, we can replace $(W_K,\psi)$ by $(W^1,\psi^1)$. 
  
Let $\alpha\in C_T$. Obviously $o(-\alpha)=\{-{\gamma}|\gamma\in o(\alpha)\}$. Let $M_{\alpha}$ be a $V_1$-lattice of $W_{\alpha}$ left invariant by the subgroup of $\Gal(K_1/K)$ fixing $\alpha$. For any $\gamma\in o(\alpha)$ let $M_{\gamma}:=\rho(M_{\alpha})$, with $\rho\in\Gal(K_1/K)$ such that $\rho(\alpha)=\gamma$, and let  $M_{-\gamma}$ be a $V_1$-lattice of $W_{-\gamma}$ such that $\psi\colon N_{\gamma}\otimes N_{\gamma}\to V_1$ is a perfect form. Here $N_{\gamma}:=M_{\gamma}\oplus M_{-\gamma}$. Let 
$M_{O(\alpha)}:=\oplus_{\gamma\in o(\alpha)\cup o(-\alpha)} M_{\gamma}$.
For any other pair $o(\alpha_1)$ and $o(-\alpha_1)$ of orbits of the action of $\Gal(K_1/K)$ on $C_T$, we define similarly a free $V_1$-submodule $M_{O(\alpha_1)}$ of $W_K\otimes K_1$. Let $M_{V_1}$ be the $V_1$-lattice of $W_K\otimes K_1$ defined by the direct sum of these $M_{O(\alpha)}$; so 
$$M_{V_1}:=\oplus_{\gamma\in C_T} M_{\gamma}.$$ 
We get that $M_{V_1}$ is stable under the action of $\Gal(K_1/K)$ and  $\psi\colon M_{V_1}\otimes M_{V_1}\to V_1$ is a perfect form. So $T_{K_1}$ extends to a subtorus $T_{V_1}$ of $Sp(M_{V_1},\psi)$. Let $M_V$ be the $V$-lattice of $W_K$ formed by elements of $M_{V_1}$ fixed by $\Gal(K_1/K)$. We have $M_{V_1}=M_V\otimes V_1$. We deduce that $T_K$ extends to a subtorus $T$ of $Sp(M_V,\psi)$. So the (unique) hyperspecial subgroup $T(V)$ of $T_K(K)$ is included in the hyperspecial subgroup $Sp(M_V,\psi)(V)$ of $Sp(W_K,\psi)(K)$.

e)  To see the first part of e), it is enough to work separately the case of a torus (which is obvious) and the case of a semisimple group. For this last case, it is enough to work  with an inner automorphism and then everything results from [Ti, first paragraph of 2.5] applied to the adjoint group. The proof of the second part of e) is similar to the proof of  4.3.10 b) (which can be used to obtain a  refinement of this part of e)), and so it is left as an exercise.

\Proclaim{3.1.2.2. Remarks.} \rm
{\bf 1)}  Let ${G_1}_K$ and ${G_2}_K$ be two reductive groups over $K$ such that ${G_1}_K$ is a subgroup of ${G_2}_K$ and such that they are unramified over $K$. It is not always true that there is a hyperspecial subgroup of ${G_1}_{K}(K)$ included in a hyperspecial subgroup of ${G_2}_{K}(K)$.

We leave to the reader to find examples for this, with ${G_2}_K=GSp(W_K,\psi)$  a group of symplectic similitudes and with ${G_1}_K$ a torus splitting over an unramified extension  of $K$ of degree 2 and for which the key fact of the proof of 3.1.2.1 d) is not true, starting from the fact that any hyperspecial subgroup of ${G_2}_K(K)$ is of the form $GSp(M_V,\psi)(V)$, with $M_V$ a $V$-lattice of $W_K$ for which there is $\eps\in\dbG_m(K)$ such that $\eps\psi\colon M_V\otimes M_V\to V$ is a perfect form. We get such examples even for $\dim_K(W_K)=4$.  
 
{\bf 2)} 1) above is always true if ${G_2}_K$ is a general or special linear group over $K$. In fact in this case any bounded subgroup $H_1$ of ${G_1}_K(K)$ is included in a hyperspecial subgroup of ${G_2}_K(K)$. To see this we can assume that $H_1$ is a maximal bounded subgroup of ${G_1}_K(K)$. Now everything results from [Ja, 10.4] and [Ti, 3.4]. 

{\bf 3)} A third proof of 3.1.2.1 c) can be obtained using maximal tori. Its advantage: it remains valid for an arbitrary DVR.
\finishproclaim

\Proclaim{3.1.3. Lemma.} 
Every reductive group $G$ over $\dbQ$ unramified over $\dbQ_p$ extends to a reductive group over $\dbZ_{(p)}$. 
\finishproclaim

\proof
It is enough to treat separately the case when $G$ is a torus and the case when $G$ is semisimple. If $G$ is a torus, then it splits over a Galois extension of $\dbQ$ which is unramified
above $p$, and so it extends to a torus over $\dbZ_{(p)}$. 

Let now $G$ be a semisimple group, and let $G_{\dbZ_p}$ be a semisimple group over $\dbZ_p$ having as generic fibre $G_{\dbQ_p}$. Let $\grg_{\dbZ_p}:=\Lie(G_{\dbZ_p})$ and let $\grg_{\dbZ_{(p)}}$ be its intersection with $\Lie(G)$. This intersection is taken inside $\Lie(G_{\dbQ_p})$. So $\grg_{\dbZ_p}=\grg_{\dbZ_{(p)}}\otimes\dbZ_p$. We get that the Zariski closure of $G^{\ad}$ in $Aut(\grg_{\dbZ_{(p)}})$ is an adjoint group $G_{\dbZ_{(p)}}^{\ad}$ over $\dbZ_{(p)}$  (this is so due to the fact that we get this over $\dbZ_p$). As $G$ is a cover of $G^{\ad}$, we conclude that $G$ extends to a semisimple group over $\dbZ_{(p)}$, obtained by taking the normalization of $G_{\dbZ_{(p)}}^{\ad}$ in the field of fractions of $G$. This ends the proof of the Lemma.

Another proof can be obtained  using the following general result of descent:
\finishproclaim

\Proclaim{3.1.3.1. Claim.}
Let $O$ be a DVR of mixed characteristic and let $L$ be its field of fractions. Let $Y_L$ be a scheme of finite type over $L$. Let $\widehat{Y}$ be a  faithfully flat scheme of finite type over the completion $\widehat{O}$ of $O$ such that its generic fibre is isomorphic to the extension of $Y_L$ to the field of fractions of $\widehat{O}$. We assume that either $\widehat{Y}$ is an affine scheme or that $O$ is a henselian ring. Then there is a unique scheme $Y$ over $O$, having $Y_L$ as its generic fibre, and such that its extension to $\widehat{O}$ is $\widehat{Y}$ (so the special fibres of $Y$ and $\widehat{Y}$ are the same).
\finishproclaim

This is not necessarily true if we do not assume that $\widehat{Y}$ is an affine scheme or that $O$ is a henselian ring: [BLR, 6.7] can be adapted to the mixed characteristic situation.
 A similar thing can be stated for morphisms (and so in terms of equivalences of some categories). 
We include a proof of the above Claim which we think is useful in the study of different integral models of (quotients of) Shimura varieties.

\Proclaim{Step -1.} \rm
We can assume that the special fibre $Y_{\scrS}$ of $\widehat{Y}$ is non-empty. It is well known that the Claim is true if $\widehat{Y}$ is affine (simple argument at the level of lattices). So we can assume that $\widehat{Y}$ is reduced. 

First we point out how the local rings of points of the special fibre of $Y$ (assumed to exist) can be recovered: for any closed point $y\colon \Spec(k)\hookrightarrow Y_{\scrS}$ (with $k$ a field), the local ring of $y$ in $Y$ is the intersection of the local ring of $y$ in $\widehat{Y}$ with the ring of fractions of $Y_L$.
Secondly we point out that the topological space underlying any scheme $Z$ over $O$ is fully determined by $Z_L$ and $Z_{\widehat{O}}$. 
These two remarks take care of the uniqueness part.
\finishproclaim

\Proclaim{Step 0.} \rm
We can work around a point $y$ as above. In particular we can assume that $Y_{\scrS}$ is affine, and so separated. Further on we can remove from $Y_L$ a closed subscheme whose Zariski closure in $\widehat{Y}$ does not contain the point $y$, in such a way that its complement in $Y_L$ and respectively the complement of this Zariski closure in $\widehat{Y}$ are separated schemes. This is possible due to the fact that $\widehat{Y}$ and $Y_L$ are noetherian schemes and due to the fact that $Y_{\scrS}$ is separated. 

Conclusion: we can assume that $\widehat{Y}$ (and so also that $Y_L$) is a separated scheme.  Now, in essence, everything results from the ideas presented in [BLR, 6.5]. We present the details.
\finishproclaim

\Proclaim{Step 1.} \rm
We can assume that $Y_L$ and $\widehat{Y}$ are normal schemes (as $O$ is an excellent ring, the argument is the same as the one needed to assume that $Y$ is reduced). 
\finishproclaim

\Proclaim{Step 2.} \rm
A. From now on we assume that $O$ is a henselian ring. 

\smallskip
B. Even better, we can assume that $O$ is a strictly henselian ring.

This admits an argument using Galois-descent (cf. [BLR, 6.2]). In other words, if $\Spec(O^1)$ is a finite Galois cover of $\Spec(O)$, with $O^1$ a DVR, and if we know that there is a scheme $Y^1$ over $O^1$ whose extension to $\widehat{O^1}$ is the extension of $\widehat{Y}$ to $\widehat{O^1}$, then the fact that the special fibre of $Y^1$ is definable over the residue field of $O$ is the extra ingredient needed to make the Galois-descent (defined by the natural action, due to the uniqueness property mentioned in  Step -1, of the Galois group $\Gal(O^1/O)$ on $Y^1$) effective.

\smallskip
C. We can assume that the residue field of $O$ is an algebraically closed field (i.e. we can replace the strictly henselian ring $O$ by a pro-finite flat DVR extension of it having the same index of ramification). 

\smallskip
D. We can replace $O$ by any finite flat DVR extension of it; so we can assume that the strictly henselian DVR $O$ is as ramified as desired.

\smallskip
E. We can replace the strictly henselian DVR $O$ by the local ring of a generic point of the special fibre of a smooth scheme $Z$ over $O$.

\smallskip
The parts C and D admit the same argument involving descent as in part B.

The part E is trivial (we can assume  that we have a scheme $Y_Z$ over $Z$ whose extension to $Z_{\widehat{O}}$ is the extension of $\widehat{Y}$ to $Z_{\widehat{O}}$; now we can take $O$-sections of $Z$ to get $Y$).

\smallskip
Replacing $O$ by another DVR $O_1$ which is obtained from $O$ by the rules described in B to E above, the nilpotent ``elements" of $Y_{\scrS}$ can be ``absorbed": the normalization $Y_1^n$ of $\widehat{Y}_{O_1}$ has a reduced special fibre. Moreover $Y_1^n$ is a projective scheme if $\widehat{Y}$ is. 
As a conclusion we can assume that $Y_{\scrS}$ is a reduced scheme and that $O$ is strictly henselian (cf. E for this last part).  
\finishproclaim

\Proclaim{Step 3.} \rm
We can assume that there is a scheme $U$ over $O$ such that its  generic fibre is $Y_L$ and its extension to $\widehat{O}$ is an open subscheme of $\widehat{Y}$ having a complement in $\widehat{Y}$ of codimension at least 2.

For this we can assume that $Y_L$ is affine (even smooth over $L$). But this is the context in which [BLR,  6, p. 161] works (if $\widehat{Y}$ is an affine scheme then the Claim is trivial). The arguments presented in loc. cit. work in the case when $\widehat{Y}$ is a normal scheme having a non-empty reduced special fibre.  
\finishproclaim

\Proclaim{Step 4.} \rm
From the Artin's approximation theorem (this is standard --see [BLR, Th. 12, p. 83]--: $O$ is an excellent ring, as $L$ has characteristic zero) we deduce easily that there is a normal scheme $Y^{\prime}$ of finite type over $O$ (we recall that $O$ is a strictly henselian ring) having $U$ as an open subscheme containing the generic fibre and all the points of codimension 1, and having $Y_{\scrS}$ as its special fibre. 
\finishproclaim

\Proclaim{Step 5.} \rm
Morally  $Y^{\prime}_{\widehat{O}}$ should be $\widehat{Y}$. The failure of being so might happen if the topology on the underlying set of $Y^{\prime}_{\widehat{O}}$ is not the expected one. If $Y^{\prime}_{\widehat{O}}$ is not $\widehat{Y}$  we have to proceed as follows.  

We can assume that $Y_L$ and $\widehat{Y}$ are normal complete schemes (cf. Nagata's embedding theorem; see [Na] and [Vo]), and that $Y_{\scrS}$ is reduced (cf. Step 2). Now we consider the normalization $\widehat{Y}_2$ of the Zariski closure of the diagonal embedding of $U_{\widehat{O}}$ in $\widehat{Y}\times Y^{\prime}_{\widehat{O}}$. We also consider the natural projections of it on the two factors. 

We first assume the existence of a scheme $Y_2$ over $O$ whose extension to $\widehat{O}$ is $\widehat{Y}_2$. This is the extra ingredient needed to be able to repeat the above arguments on the application of Artin's approximation theorem to get similarly (to $Y^{\prime}$) a scheme $Y^{''}$ over $O$ whose topology of its underlying set is as expected to be (the topology of the underlying set of $\widehat{Y}$ is a quotient topology of the topology underlying the set of $\widehat{Y}_2$). We deduce that the extension of $Y^{''}$ to $\widehat{O}$ is $\widehat{Y}$. So we can take $Y=Y^{''}$.

In fact we can replace $\widehat{Y}_2$ by any other proper scheme $\widehat{Y}_3$ over $\widehat{O}$ whose generic fibre is defined over $L$, and which admits a surjection onto $\widehat{Y}_2$. Even more, it is enough to get such a good scheme $\widehat{Y}_3$ only after we replace $O$ by an arbitrary DVR $O_1$ which is a faithfully flat $O$-algebra of the type allowed in Step 2. 
\finishproclaim

\Proclaim{Step 6.} \rm
To end the proof we are left with the proof of the existence of $Y_3$ for a suitable choice of $\widehat{Y}_3$.
A well known application of  Chow's   (cf. also [Vo, 2.5]) shows that we can assume that we are dealing with an $\widehat{Y}_2$ which is a normal faithfully flat projective $\widehat{O}$-scheme. In other words there is a surjective proper morphism $\widehat{Y}_3\to\widehat{Y}_2$, with $\widehat{Y}_3$ a normal projective scheme over $\widehat{O}$, whose generic fibre is defined over $L$.

We can assume (cf. Step 2 and the last part of Step 5) that $\widehat{Y}_2$ has a reduced special fibre. From Step 3 we deduce easily the existence of a normal projective scheme $Y_2^{\prime}$ over $O$ such that there are open subschemes $U_2^{\prime}$ and $U_2$ of $Y_2^{\prime}$ and respectively of $\widehat{Y}_2$, containing the generic fibres and the points of codimension 1, and satisfying $U_{2\widehat{O}}^{\prime}=U_2$. We can view this last identity as a rational map from $Y^{\prime}_{2\widehat{O}}$ to $\widehat{Y}_2$. But this rational map extends to a surjective morphism $Y_{3\widehat{O}}\to\widehat{Y}_2$, where $Y_3$ is a projective scheme over $O$ obtained from $Y_2^{\prime}$ through  a blowing up centered on the special fibre (cf. [Hart, 7.17.3, p. 169]; we can view $\widehat{Y}_2$ as embedded in a projective space $\dbP^n_{\widehat{O}}$, for some $n\in\dbN$). The scheme $Y_3$ is the searched for scheme over $O$. This ends the proof of the Claim.
\finishproclaim
 
\Proclaim{3.1.3.1.1. Remarks.} \rm
1) We preferred to include the above proof of 3.1.3 (it also works when $\dbZ_{(p)}$ is replaced by an arbitrary DVR $O$) as it illustrates how descent can be performed also at the level of Lie algebras. Moreover it is constructive.

2) The proof of 3.1.3.1 can be modified so that it works for an arbitrary henselian DVR $O$: the use of Artin's approximation theorem has to be replaced (in the case when $O$ is of equal positive characteristic) by the use of C and D of Step 2 (for Step 1 cf. [Ma, Cor. 2 of p. 234 and 31.H]).

3) We do expect that in 3.1.3.1 we can replace finite type  by  locally of finite type. 
\finishproclaim 

\Proclaim{3.1.3.2. Remark.} \rm
Let now $O$ be an arbitrary DVR  having a perfect residue field. Let $G_L$ be a reductive group over the field $L$ of fractions of $O$, such that its extension to the field of fractions $K$ of the completion $V$ of $O$ is unramified
over $K$. Let $H$ be a hyperspecial subgroup of $G(K)$. Then any automorphism of $G_L$ taking $H$ onto itself, extends to an automorphism of $G_O$, with $G_O$ a reductive group over $O$ such that $G_O(\widehat{O})=H$ (such a $G_O$ does exist cf. 3.1.3.1): from 3.1.2.1 a) and e) we get an automorphism of $G_V$; obviously it comes from an automorphism of $G_O$.
\finishproclaim

\Proclaim{3.1.4. Remark.} \rm
Let $(G,X)\hookrightarrow (G_1,X_1)$ be an injective map of Shimura pairs and let $p$ be a rational prime. If $G^{\der}$ and the connected component of the origin $C$ of its centralizer in $G_1$ are unramified over $\dbQ_p$, then there is an injective map $(G_0,X_0)\hookrightarrow (G_1,X_1)$ such that $G_0^{\der}=G^{\der}$, $(G_0^{\ad},X_0^{\ad})=(G^{\ad},X^{\ad})$ and $G_0$ is unramified over $\dbQ_p$. To see this it is enough to remark that there is (cf. [Ha, 5.5.3]) a maximal torus $T$ of $C$ such that:

\smallskip
--  a conjugate of some $x\in X$ by an element of $C(\dbR)$ factors through $G_{0\dbR}$, where $G_0$ is the subgroup of $G_1$ generated by $G$ and $T$;

-- $T_{\dbQ_p}$ is $C(\dbQ_p)$-conjugate to a maximal torus of $C_{\dbQ_p}$ unramified over $\dbQ_p$ (there is such a maximal torus as $C_{\dbQ_p}$ is quasi-split, cf. [Ti, 1.10]).

\Proclaim{3.1.5. Remark.} \rm
Let $G_K=G_{1K}\times G_{2K}$ be a product of reductive groups over $K$. Then $G_K$ is unramified over $K$ iff $G_{1K}$ and $G_{2K}$ are unramified over $K$, and then any hyperspecial subgroup $H$ of $G_K(K)$ is a direct product $H=H_1\times H_2$, with $H_i$ a hyperspecial subgroup of $G_{iK}(K)$, $i=\overline{1,2}$.

\Proclaim{3.1.6. Lemma.} 
 Let $R$ be an integral domain, faithfully flat over $\dbZ_{(p)}$. Let $M$ be a free $R$-module of finite rank. Let $G^0$ be a reductive subgroup of the generic fibre of $GL(M)$ such  that the Zariski closures of $G^{0\der}$ and of the connected component  $T$ of the origin of $Z(G^0)$ in $GL(M)$, are both reductive group schemes over $R$. Then the Zariski closure of $G^0$ in $GL(M)$ is a reductive group scheme over $R$.
\finishproclaim

\proof
Let $G_R^0$, $G_R^{0\der}$ and $T_R$ be respectively the Zariski closures  of $G^0$, $G^{0\der}$ and $T$ in $GL(M)$. Let $C:=G^{0\der}\cap T$. So $C$ is a finite flat group scheme over the generic fibre of $\Spec(R)$. 

We claim that the Zariski closure $C_R$ of $C$ in $GL(M)$ is a finite flat subgroup of $T_R$ and of $G^{0\der}_R$. This is a local statement in the \'etale topology of $\Spec(R)$. Of course, if $\Spec(R_1)\to\Spec(R)$ is an \'etale map, $R_1$ might not be an integral ring, and so we need to take the Zariski closure of $C\times_R \Spec(R_1)$ in $GL(M)_{R_1}$. So we can assume that $T_R$ is split and that $G^{0\der}_R$ has a maximal split torus $T^1_R$, but we no longer assume that $R$ is integral: just that it is reduced. Moreover we can assume that $R$ is a local ring. It is enough to show that the intersection $T_R\cap T_R^1$ defines a finite flat scheme over $R$. We consider the torus  $T^2_R:=T_R\times T^1_R$, and its representation $\rho$ on $M$ defined by the product of the inclusion $T_R\hookrightarrow GL(M)$ with the inverse of the inclusion $T_R^1\hookrightarrow GL(M)$. This is well defined as $T_R$ and $T_R^1$, as subtori of $GL(M)$, commute. As $R$ is local and $T_R^2$ is split we deduce that $\rho$ is a direct sum of one-dimensional representations (associated to characters of $T_R^2$). So its kernel is the extension of a finite flat commutative group scheme over $R$ by a torus; but this last torus is trivial, as this is so over the points of $\Spec(R)$ of codimension 0. So ${\text Ker}(\rho)$ is a finite flat commutative group scheme over $R$. But this kernel is $T_R\cap T_R^1$. So $C_R$ is a finite flat group scheme over $R$.     

We come back to our situation: $R$ is integral.
 Let $G_R^1$ be the quotient of $T_R\times G_R^{0\der}$ by $C_R$, where $C_R$ acts as inclusion on $T_R$ and as the inverse of the inclusion on $G_R^{0\der}$. The group $G_R^1$ is a reductive group scheme over $R$ (as we have a short exact sequence $0\to G^{0\der}_R\to G_R^1\to T_R/C_R\to 0$). We have a canonical homomorphism $q\colon G_R^1\to G_R^0$, which is an isomorphism over the generic fibre of $R$. But 3.1.2.1 c) and 3.1.2.2 3) guarantee that each fibre of $q$ is a closed embedding, and that $q$ is a proper morphism. This implies that $q$  is a closed embedding: we can assume that $R$ is finitely generated over $\dbZ_{(p)}$, and so that it is noetherian; first we deduce that $q$ is a finite morphism, and then everything results from Nakayama's Lemma. From the definition of $G_R^0$ we deduce that $q$ is an isomorphism. This implies that  $G_R^0$ is a reductive group scheme over $R$, ending the proof of the Lemma.

\Proclaim{3.1.6.1. Remark.} \rm
The above Lemma remains true if $\dbZ_{(p)}$ is replaced by an arbitrary DVR. Even more generally, it remains true if $R$ is an integral scheme, and instead of its generic fibre (over some integral scheme), we work with its generic point, cf. 3.1.2.2 3).
\finishproclaim  

\smallskip
\Proclaim{3.2. Healthy normal schemes and integral models of Shimura varieties.} \rm
In 3.2.1 and 3.2.2 we introduce the general theory of healthy normal schemes. The need of such a theory was felt when it has been discovered that the statement 6.8 of [FC, p. 185] is not true in general (for details see [dJO]). Then in 3.2.3 to 3.2.16 we present the general theory of integral models of Shimura varieties. For this theory, in essence (i.e. except [Fa4]), we need from the theory of healthy normal schemes only some definitions and remarks. However we felt that it is important to include in 3.2.1 and 3.2.2 more then just definitions (cf. the philosophy of 3.2.7 6) and rm. 3) of 3.2.3.2.1; they nourish our expectation that the theory of healthy normal schemes will blossom very much in the near future). In 3.2.17 we single aside the proof of a result of Faltings [Fa4] which plays an essential role in the theory of integral models of Shimura varieties. It introduces some of the main tools used in the study of healthy normal schemes. As these tools are referred to in 3.2.1 and 3.2.2 we suggest that after 1), 2) and 8) of 3.2.1 and 1) and 3) of 3.2.2, 3.2.17 should be studied, before the rest of 3.2.1 and 3.2.2. 

Let $p$ be a rational prime.
\finishproclaim

\Proclaim{3.2.1. Definitions.} \rm
1) A pair $(Y,U)$, with $Y$ a flat scheme over $\Spec(\dbZ)$ and with $U$ an open subscheme of $Y$ containing the generic fibre $Y_{\dbQ}$ of $Y$ and such that the complement of $U$ in $Y$ is of codimension in $Y$ at least 2, is called an extensible pair. A pair $(Y,U_Y)$, with $Y$ as before and with $U_Y$ a subset of the underlying set of $Y$ which is an intersection of the underlying sets of open subschemes $U_i$ of $Y$, $i\in I$, such that $(Y,U_i)$, $i\in I$, are extensible pairs, is called a quasi-extensible pair. Here $I$ is an arbitrary set, often infinite.

\smallskip
A normal  scheme $Y$ flat over $\Spec(\dbZ)$
is called: 

\smallskip
2) healthy if for any extensible pair $(Y,U)$, every abelian scheme over $U$ extends to an abelian scheme over $Y$; 

3) quasi healthy if for any extensible pair $(Y,U)$, every abelian scheme over $U$ extends to an abelian scheme over the normalization of $Y$ --not assumed to be finite over $Y$-- in a finite \'etale extension of the ring of fractions of $Y$;

4) almost healthy if any abelian scheme $A_{\dbQ}$ over $Y_{\dbQ}$ having level-$l^N$ structures for any $N\in\dbN$, with $l$ a rational prime which is invertible in any point of $Y$, extends to an abelian scheme over the normalization of $Y$ --not assumed to be finite over $Y$-- in a finite \'etale extension of the ring of fractions of $Y$;

5) $n$ healthy if for any extensible pair $(Y,U)$, every abelian scheme over $U$ of dimension at most $n$ extends to an abelian scheme over $Y$ (here $n\in\dbN$); 

6) locally healthy if any open subscheme of it is healthy.

\smallskip
7) Similarly we define the following types of normal schemes (flat over $\Spec(\dbZ)$): $n$ quasi healthy, $n$ almost healthy, locally quasi healthy, locally almost healthy, and locally $n$ (quasi or almost) healthy. 

\smallskip
8) Let $D$ be a Dedekind ring  which is flat over  $\dbZ\fracwithdelims[]12$. A regular scheme $Y$ flat over $D$  is called very healthy if:

i) for any prime $w$ of $D$ having a residue field $k(w)$ of positive characteristic $p$, the only open subscheme of the fiber $Y_w$ of $Y$ over $w$ containing all points of $Y_w$ having as a residue field an algebraic extension of $k(w)$, is $Y_w$ itself;
 
ii) for any geometric point $y:\Spec(\overline{k(w)})\hookrightarrow  Y_{\dbD_w}$ (with $\dbD_w$ a complete DVR faithfully flat over  the localization $D_{(w)}$ of $D$ with respect to $w$, having $\overline{k(w)}$ as its residue field and having the same ramification index as $D_{(w)}$), the completion of the local ring of $y$ is of the form $R_y=V[[x_1,x_2,...,x_m]]$, with $V$ a DVR  containing $W(\overline {k(w)})$, and such that the degree $[V: W(\overline {k(w)})]$
is less than $p-1$.

\smallskip
A normal scheme $Y$ flat over $\dbZ_{(p)}$ is called:

\smallskip
9) p-healthy if for any extensible pair $(Y,U)$, every $p$-divisible group over $U$ extends uniquely to a $p$-divisible group over $Y$; warning: here we use a hyphen (p-healthy), while in 5) we do not;

10) p-f-healthy if  for any extensible pair $(Y,U)$, every finite flat group scheme over $U$ of $p$-power order extends uniquely to a finite flat group scheme over $Y$;

11) strongly p-healthy if any $p$-divisible group over $Y_{\dbQ}$ extends uniquely to a $p$-divisible group over $Y$;

12) strongly p-f-healthy if it is p-f-healthy, and if any finite flat group scheme over $Y_{\dbQ}$ of $p$-power order extends in at most one way to a finite flat group scheme over $Y$.

\smallskip
13) As in 6) and 7) we speak about locally p-healthy, locally p-f-healthy, locally strongly p-healthy and locally strongly p-f-healthy normal schemes (flat over $\dbZ_{(p)}$).
\finishproclaim

Let now $O$ be a DVR which is a faithfully flat $\dbZ_{(p)}$-algebra, and let $e$ be its index of ramification. Let $\pi_O$ be a uniformizer of $O$.
\smallskip

\Proclaim{3.2.1.1. Remarks.} \rm
1) If $D$ is a  DVR faithfully flat over $\dbZ_{(p)}$, in order that there are very healthy regular schemes over $D$ with a non-empty special fibre,  the ramification index of $D$ has to be smaller than $p-1$. If this is so, then any projective limit of smooth schemes over $D$ with \'etale transition morphisms, for  which i) of 3.2.1 8) is true, (in particular any smooth scheme over $D$) is a very healthy regular scheme over $D$. 

2) According to a theorem of Raynaud (cf.  [Ra, 3.3.3]), if $e<p-1$, then $\Spec(O)$ is a strongly p-f-healthy regular scheme.

3) In 3.2.1 4), actually $A_{\dbQ}$ does extend to an abelian scheme over $Y$. This can be seen using the ideas of the Step $A$ of 3.2.17.

4) From the N\'eron--Ogg--Shafarevich criterion we get directly that any locally noetherian healthy normal scheme (and so any healthy regular scheme) is an almost healthy normal scheme.

5) The quotients (assumed to exist) of healthy normal schemes through finite group actions are quasi healthy normal schemes. This motivates def. 3) of 3.2.1. 

6) The regular quotients of healthy normal schemes through finite flat group actions are healthy regular schemes. This can be checked starting from Step A of 3.2.17. Similarly, the  regular quotients of locally healthy normal schemes through finite flat group actions are locally healthy regular schemes.

7) The quotients of almost healthy normal schemes through finite group actions are almost healthy normal schemes. 

8) There are plenty of examples of noetherian almost healthy normal schemes which are not regular (this results from 4) and 7) and from  3.2.2 1)),
and there are plenty of examples of healthy regular schemes which are not very healthy (cf. 3.2.2 5)).

9) Any  regular scheme $Y$ of dimension 2 flat over $\dbZ_{(p)}$ is p-f-healthy (this is a consequence of [FC, 6.2, p. 181]).

10) There are plenty of  p-f-healthy regular schemes which are not strongly p-f-healthy (like affine lines over the spectrum of a complete DVR of index of ramification $p-1$).

11) Any  p-healthy regular scheme flat over $\dbZ_{(p)}$ is healthy. The proof of this is similar to Step B of 3.2.17.

12) Any smooth scheme over a local henselian p-healthy regular scheme  $\Spec(R)$ of dimension at least two, having the property that the only open subscheme of its  special fibre (defined by $\pi_O=0$) containing its fiber over the maximal point of $\Spec(R)$, is the special fibre itself, is p-healthy. The proof of this is entirely the same as Steps B, C and D of 3.2.17 (to be compared with 3.2.2.2; an argument similar to the one of 3.2.2 4), involving Weil restriction of $p$-divisible groups, allows us to replace $R$ by a finite \'etale $R$-algebra). 

13) The class of very healthy regular schemes over a Dedekind ring $D$ flat over $\dbZ\fracwithdelims[]12$ is stable under localizations, completions (of local schemes) and passages to smooth schemes for which condition i) of 3.2.1 8) is still satisfied.

14) We could have worked out 3.2.1 for locally integral schemes instead of normal schemes. But in such a generality we have basically no results. To study any type of healthy normal schemes we can restrict our attention to integral normal schemes.

15) We consider a projective limit $Z$ of quasi-compact healthy normal schemes with transition morphisms  such that their fibres over any point of $\Spec(\dbZ)$ are dominant morphisms. We assume that one of the following two conditions is satisfied:

\smallskip
-- every fibre of $Z$ over a point of $\Spec(\dbZ)$ of finite characteristic has a finite number of points of codimension (in this fibre) zero;

-- the transition morphisms are pro-\'etale morphisms between schemes regular in points of positive characteristic and of codimension 1.

\smallskip
Then $Z$ is a healthy normal scheme (if the first condition is satisfied, the argument is standard; if the second condition is satisfied, we have to use as well [BLR, Cor. 2, p. 177]). The similar statement for almost healthy normal schemes is not true.

16) Step A of 3.2.17 explains why for checking that a normal scheme $Y$ is (quasi or almost or locally) healthy it is enough to deal with principally polarized abelian schemes. This is very useful as the moduli stack over $\Spec(\dbZ)$ parameterizing principally polarized abelian schemes of a given dimension is algebraic (and so quasi-compact and quasi-separated) [FC, 4.11, p. 23]. This means that in many situations (like in the last part of the proof of 3.2.2.1) we can work out things as in the case when we have a quasi-compact and quasi-separate
 situation.   
\finishproclaim

\Proclaim{3.2.1.2. Questions.} \rm
1) Is it true that any local healthy regular scheme over $\dbZ_{(p)}$, of dimension at least 2, is p-healthy?

2) Is the completion of a local healthy regular scheme, a healthy regular scheme itself?
\finishproclaim

We expect a positive answer to these questions. In many cases it is known that the answer to 2) is yes (cf. 3.2.2.3 B)).

\Proclaim{3.2.1.3. Problem.} \rm
Characterize the healthy regular schemes independently of the use of abelian schemes or of $p$-divisible groups.

\Proclaim{3.2.1.4. Expectations.} \rm
1)  If $Y$ is a local healthy normal scheme, then we do expect that its (strict) henselization is also a healthy normal scheme (to be compared with rm. 4) of 3.2.2). Similarly for other types of healthy local schemes. 

2)  We do expect the existence of noetherian almost healthy and of quasi healthy normal schemes which are not healthy. It should be possible to construct such examples starting from the fact that the classical purity theorem for regular rings is not true for normal noetherian rings.

3) We do expect that for any $N\in\dbN$ there are $N$ healthy normal schemes which are not $N+1$ healthy normal schemes. 

4) In 2), 3), 5), 9) and 10) of 3.2.1 we could have worked with quasi-extensible pairs instead of extensible pairs. This would have made no difference for 2), 3), 5) and 10) of 3.2.1, but we do think (we do not have an example to prove this) it would have made a difference for 3.2.1 9). The advantage of working with quasi-extensible pairs (instead of extensible pairs) consists in the fact that given a flat $\Spec(\dbZ)$-scheme $Y$ it is enough to work with only one quasi-extensible pair $(Y,U_Y)$, with $U_Y$ the subset of $Y$ defined as the intersection of the underlying sets of all open subschemes $U$ of $Y$ such that $(Y,U)$ is an extensible pair. 

5) Though we defined 2) to 7) of 3.2.1 for schemes over $\Spec(\dbZ)$, we have no understanding of the types of healthy schemes over $\dbZ_{(2)}$. In particular we do not know if $\Spec(\dbZ_2[[T]])$ is a healthy scheme; however we do expect this to be true (cf. [Va2]).  

6) We do not know even one example of a healthy normal scheme over $\Spec(\dbZ\fracwithdelims[]12)$ which is not locally healthy. We do expect that (at least under some mild conditions) any healthy regular scheme over $\Spec(\dbZ\fracwithdelims[]12)$ is locally healthy. It is a nice problem to check that all healthy regular schemes to be introduced in 3.2.2 5) are locally healthy.
\finishproclaim

\Proclaim{3.2.2. Remarks.} \rm
1) According to [Fa4], if  $e<p-1$, then any  regular formally smooth scheme over $O$ is a healthy regular scheme. As a direct consequence of this and its proof we get that any very healthy regular scheme over a Dedekind ring $D$ (flat over  $\dbZ\fracwithdelims[]12$) is a healthy regular scheme, and, if $D$ is a $\dbZ_{(p)}$-algebra, then it is also a p-healthy regular scheme (see  3.2.17 for a proof of these  statements).

2) Any healthy regular scheme is an almost healthy regular scheme. But we do not know if (or when) an almost healthy normal (regular) scheme  is healthy. However an integral almost healthy regular scheme whose first fundamental group is trivial is a healthy scheme (cf. the classical purity theorem).

3) The role of the Dedekind ring $D$ in the definition of a very healthy regular scheme (over $D$) is essentially just to fix the notations. We can define an abstract very healthy regular scheme to be a flat scheme $Y$ over $\Spec(\dbZ)$ with the property that for any local ring $O_y$ of a point $y$ of $Y$ of positive characteristic $p$, there is a faithfully flat $O_y$-algebra $R_y$, with $R_y$ of the same form as the one in 3.2.1 8). As in 1), any abstract very healthy regular scheme is a healthy regular scheme, and any abstract very healthy regular $\dbZ_{(p)}$-scheme is p-healthy (cf. 3.2.17).

The class of abstract very healthy regular schemes is stable under localization, completion, passage to regular formally smooth schemes, and under pull backs through morphisms defined by ring homomorphisms  between discrete valuation rings of mixed characteristic having the same index of ramification.

4) Let $q\colon Y_1\to Y$ be an \'etale morphism of flat $\dbZ$-schemes. We assume that there is an extensible pair $(Y,U)$ such that $(Y_1,q^{-1}(U))$ is an extensible pair and $q^{-1}(U)$ is an \'etale cover of $U$ (this is equivalent to the fact that $q$ defines an \'etale cover over $Y_{\dbQ}$ and over local rings of $Y$ which are discrete valuations rings of mixed characteristic). We have: 

\smallskip
{\bf A)} If $Y$ is a healthy normal scheme, then $Y_1$ is a healthy normal scheme.

\smallskip
To see this let $(Y_1,U_1)$ be an extensible pair, and let $A_{U_1}$ be an abelian scheme over $U_1$. We can assume that $Y_1$ and  $Y$ are integral.  We can also assume that there is an open subscheme $\tilde U$ of $U$ such that $(Y,\tilde U)$ is an extensible pair and $q^{-1}(\tilde U)=U_1$.

 We consider the abelian scheme over $\tilde U$ obtained from the abelian scheme $A_{U_1}$ through the Weil restriction of scalars (the morphism $U_1\to\tilde U$ is \'etale and finite). It extends to an abelian scheme over $Y$ (as $Y$ is a healthy normal scheme). From this by standard arguments we deduce that $A_{U_1}$ extends to an abelian scheme over $Y_1$.

As a consequence we get: 

\smallskip
{\bf B)} If $q$ is an \'etale cover, then $Y_1$ is a healthy normal scheme iff $Y$ is so. 

\smallskip
This remains true if we replace healthy schemes by any other type of healthy schemes defined in 2), 3), 5), 6) and 8) to 12) of 3.2.1, but
we do not know if (or when) this remains true if the word healthy is replaced by almost healthy.

Even better:

\smallskip
{\bf C)} If $q$ is a pro-\'etale cover, then $Y_1$ is a healthy normal scheme iff $Y$ is so.

\smallskip
To see this we can assume that $Y_1$ and $Y$ are both integral schemes. Let now $A_{U_1}$ be an abelian scheme over an open subscheme $U_1$ of $Y_1$ with the property that $(Y_1,U_1)$ is an extensible pair. There is an \'etale cover $q_2\colon Y_2\to Y$, with $Y_2$ an integral scheme, such that $q$ factors through $q_2$ and the abelian variety over the generic point of $Y_1$ obtained from $A_{U_1}$, is defined over the generic point $\nu$ of $Y_2$. Now the theory of descent implies
that this abelian variety over $\nu$ extends to an abelian scheme over an open subscheme $U_2$ of $Y_2$ with the property that $(Y_2,U_2)$ is an extensible pair. Moreover its extension to $U_1$ is $A_{U_1}$ (we can assume that $U_1$ factors through $U_2$). Now everything results from B).

A similar C) can be stated for the type of healthy schemes introduced in 6), 8) and 10) of 3.2.1. 

5) There are plenty of healthy regular schemes which are not very healthy. Section 3.2.17 is the source of inspiration for such examples. For instance, if $l$ and $p$ are two primes such that $l>p>3$, then the local schemes of whose completion is of the form 
$$
Y=\Spec(W(k)[[x,y,z]]/(x^l+y^2+z^2+p)),
$$ 
with $k$ a perfect field of characteristic $p$, is a healthy regular scheme. This can be easily seen by using Steps A to D of 3.2.17. (Hint: Using Step A we can assume that our local scheme is $Y$ itself. Then we can assume that $\overline{k}=k$ and so that $Y$ is a strictly henselian local scheme. Next we check that the closed subscheme $\Spec(W(k)[[y,z]]/(y^2+z^2+p))$ of $Y$ is a healthy regular scheme.) 
But obviously $Y$ is not a very healthy regular scheme over $W(k)$. It can be checked that $Y$ is also not an abstract healthy regular scheme. 

6) The following definition is not mathematically acceptable, and so it is not  used outside this remark; however we do expect the possibility of defining the class of regular  $O$-schemes it introduces, in terms of different  indices of ramifications of different regular closed subschemes of it. So we do see the possibility of a mathematically acceptable definition of this class, which would lead to a deep understanding of the healthy regular $O$-schemes. 

\Proclaim{\it Definition (tentative).}\rm
 We call a regular $O$-scheme S-healthy (the letter S stands for the word strongly) if the completion of the henselization of an arbitrary  local ring of it of mixed characteristic can be proved to be healthy by using Steps A to D of 3.2.17 (as in the hint of 5)). 
\finishproclaim

Any very healthy regular $O$-scheme is S-healthy, any S-healthy regular $O$-scheme is locally healthy. We do not know what is the relation between $R$-healthy regular schemes (to be defined in 3.2.2.3.1) and S-healthy regular schemes over $O$. In our opinion the most important subclasses of healthy regular schemes over $O$ are: of S-healthy, of locally healthy, and of quasi-compact healthy schemes over $O$ (to which we have to add, in the case when $e<p-1$, the subclasses of abstract very healthy schemes over $O$, of $R$-healthy schemes over $O$, and of regular formally smooth schemes over some DVR $O_1$ which is a faithfully flat $O$-algebra having $e$ as its index of ramification).  
\finishproclaim

\Proclaim{3.2.2.1. Proposition.} 
We assume that $O$ is a henselian DVR. Let $Y$ be a regular scheme over $O$ and let $O\hookrightarrow O_1$ be a formally \'etale inclusion, with $O_1$ a DVR. Then $Y_{O_1}$ is a healthy regular scheme iff $Y$ is a healthy regular scheme.
\finishproclaim

\proof
Obviously $Y_{O_1}$ is a regular scheme. 
If $Y_{\widehat{O_1}}$ is a healthy regular scheme, then from the theory of descent, we deduce that $Y$ and $Y_{O_1}$ are healthy regular schemes. So we can assume that $O_1$ is complete. Using B) of 3.2.2 4) we can assume that $O_1=\widehat{O}$. Let $Y_1:=Y_{O_1}$. 

We assume now that $Y$ is a healthy regular scheme. 
We can assume that $Y$ is an integral $O$-scheme, with a non-empty special fibre. Let $(Y_1,U_1)$ be an extensible pair, and let $A_{U_1}$ be an abelian scheme over $U_1$. There is an extensible pair $(Y,U)$ such that $U_1=U\times_Y Y_1$ (as the special fibers of $Y$ and $Y_1$ are the same).

We treat first the case when $Y$ is an affine (integral) scheme. Then $A_{U_1}$ is defined over a scheme $U_{O_1^{\prime}}$, with $O_1^{\prime}$ a finitely generated $O$-subalgebra of $\widehat{O}$. As $O$ is an excellent ring (as its field of fractions has characteristic zero), we deduce from  [BLR, Th. 12, p. 83] the existence of an $O_1^{\prime}$-algebra $O_2$, smooth over $O$, and such that we have a homomorphism $O_2\to\widehat{O}$ of $O_1^{\prime}$-algebras. Let $O_3$ be the localization of $O_2$ with respect to its prime ideal dominated by the maximal ideal of $\widehat{O}$. 

Remark 3.2.2 4) gives  us the right to assume that the first fundamental group of $Y$ is trivial (and so that $Y$ is an almost healthy regular scheme), and that $O$ is a strictly henselian DVR. From the smoothness of  $O_2$ (over $O$) we deduce the existence of an epimorphism $s_3\colon O_3\twoheadrightarrow O$ of $O$-algebras.
  Now it is easy to check that the resulting (abelian scheme) model of $A_{U_1}$ over $U_{O_3}$ extends to an abelian scheme over $Y_{O_3}$: using the fact that $Y_{O_3}$ is a regular scheme (being the localization of a smooth $Y$-scheme), we can follow entirely the independent section 3.2.17 (the role of $V$ being replaced by $Y$; the existence of $s_3$ guarantees that everything is the same). We deduce that $A_{U_1}$ does extend to an abelian scheme over $Y_1$.

The same argument works for the case when $Y$ is quasi-compact (i.e. a local $O$-algebra $O_3$ as above does exist in this case also). The general case is treated similarly: we can assume that $A_{U_1}$ (cf. 3.2.1.1 16)) is principally polarized; as the moduli stack of principally polarized abelian schemes of a given dimension over $O$-schemes is algebraic over $O$ (and so it is quasi-compact and quasi-separated) we deduce (see also below) the existence of a local $O$-algebra $O_3$ having the same properties as above. The rest of the argument is the same.
  
In fact we can avoid using stacks as follows. Let $V$ be the normalization of $O$ in the field of fractions of $Y$. It is a strictly henselian DVR (as $O$ is so and as the special fibre of the regular scheme $Y$ is non-empty). The generic fibre of $Y$ is geometrically connected over the field of fractions $K(V)$ of $V$. As $V_1:=V\otimes_O O_1$ is a DVR, this implies that $Y_1$ is an integral scheme. Moreover its generic fibre is geometrically connected over the field of factions of $V_1$. Now it is an easy exercise to check (starting from the fundamental exact sequence of [SGA1, p. 253], applied to the generic fibre of $Y$, viewed as a $K(V)$-scheme) that the first fundamental group of $Y_1$ is trivial (we view $Y$ as a $V$-scheme and $Y_1$ as a $V_1$-scheme); as $Y$ is regular we can refer to loc. cit. even if $Y$ is not quasi-compact. From this and the classical purity theorem we deduce that $A_{U_1}$ has level-$N$ structures, for any $N\in\dbN$ satisfying $(N,p)=1$. So we can replace the referred stack, by a Mumford scheme $\scrA_{d(A_{U_1}),1,N}$ (we view it as a quasi-projective smooth scheme over $\dbZ\fracwithdelims[]1N$) (cf. [Mu]). Here $d(A_{U_1})$ is the relative dimension of $A_{U_1}$, while $N>2$ is an integer satisfying $(N,p)=1$. Let $m:U_1\to\scrA_{d(A_{U_1}),1,N}$ be such that $A_{U_1}$ is the pull back through $m$ of the universal abelian scheme over $\scrA_{d(A_{U_1}),1,N}$ (in what follows we do not need to mention polarizations). We can assume $V=O$ and $V_1=O_1$. As any morphism between spectra of fields of characteristic 0 is regular, using descent and [FC, 2.7 of p. 9] we deduce the existence of a subfield $L$ of $O_1[{1\over p}]$ which is finitely generated over $K(V)$ and such that the generic fibre of $m$ factors through $Y\times_O L$. 

Let $C_1$, ..., $C_l$ be a finite cover of $\scrA_{d(A_{U_1}),1,N}$ by affine, open subschemes whose complements are smooth and have pure codimension 1, $l\in\dbN$. For $i\in\{1,...,l\}$, let $Z_i$ be an arbitrary affine, open subscheme of $U$ such that its special fibre is mapped by $m$ into $C_i$. Let $W_i$ be the maximal open subscheme of $Z_i\times_O O_1$ such that the restriction $m_i$ to it of the natural morphism $Z_i\times_O O_1\to\scrA_{d(A_{U_1}),1,N}$ factors through $C_i$; it is affine and contains the special fibre of $Z_i\times_O O_1$. We deduce the existence of a finitely generated $O$-subalgebra $O^i$ of $L\cap O_1$ such that, provided we replace $U$ by the open complement in it of a closed subscheme of its special fibre of codimension in this fibre at least 1, $m_i$ factors through an open, affine subscheme of $Z_i\times_O O^i$ containing the special fibre. As the notation suggests, we can assume that $O^i$ depends only on $i\in\{1,...,l\}$ and not on $Z_i$ (cf. simple arguments on global sections based on [Ha, ex. 2.4, p. 79] and on dilatations; for this part we can assume that the generic fibre of $Y$ is affine). We deduce the existence of a finitely generated $O$-subalgebra $O_1^\prime$ of $O_1$ containing all $O^i$'s, having $L$ as its field of fractions and such that $A_{U_1}$ is the pull back of an abelian scheme over $Z$, with $Z$ an open subscheme of $U_{O_1^\prime}$ containing the generic points of its special fibre, the special fibre of $U\times_O O_1$, $Y\times_O L$ and $\tilde U_{O_1^\prime}$, with $\tilde U$ an affine, dense, open subscheme of $U$. So $Y_{O_1^\prime}\setminus Z$ has codimension at least $2$. 

We consider a projective compactification $\bar\scrA_{d(A_{U_1}),1,N}$ of $\scrA_{d(A_{U_1}),1,N}$ as in [FC, \S 5 of Ch. V]; so the complement of $\scrA_{d(A_{U_1}),1,N}$ in $\bar\scrA_{d(A_{U_1}),1,N}$ is of pure codimension 1. We embed $\bar\scrA_{d(A_{U_1}),1,N}$ in a projective space $\dbP^r$ over $\dbZ\fracwithdelims[]1N$. Let $\scrF$ be the pull back to $Z$ of the canonical line bundle $\scrO(1)$ on $\dbP^r$. Now the existence of $O_2$ and $O_3$ as above and such that the morphism $Z\times_{O_1^\prime} O_3\to\scrA_{d(A_{U_1}),1,N}$ extends to a morphism $U_{O_3}\to\scrA_{d(A_{U_1}),1,N}$, can be easily checked starting from [Ha, 7.1, p. 150] and standard arguments with Weil divisors (cf. also [FC, iv) of 10.1, p. 88], the fact that the pull back of $\scrF$ to $Z\times_{O_1^\prime} O_2$ extends to a line bundle on $Y\times_O O_2$ and [Ma, Th. 38]). This ends the proof of the Proposition.

\Proclaim{3.2.2.2. Lemma.} Let $\Spec(R)$ be a local henselian healthy regular scheme over $O$ of dimension at least 2. Let $Z$ be a  normal $R$-scheme which is a projective limit of smooth schemes of finite type over $\Spec(R)$ such that:

\smallskip
-- each member of the projective system has the property that the only open subscheme of its special fibre (defined by $\pi_O=0$) containing its fibre over the maximal point of $\Spec(R)$ is the special fibre itself;

-- the transition morphisms are dominant modulo $\pi_O$; 

-- either the transition morphisms are \'etale or each connected component of $Z_{O^{\text{sh}}}$ is such that its special fibre has a finite number of points of codimension (in this special fibre) zero.

\smallskip
Then $Z$ is a  healthy normal scheme.  
\finishproclaim

\proof
Remark 3.2.2 4) gives us the right to assume that $\Spec(R)$ is a strictly henselian local scheme, and so that it is an almost healthy scheme, cf. 3.2.2 2). So we can assume that $O=O^{\text{sh}}$ and that $Z$ is connected. It is enough (cf. 3.2.1.1 15)) to prove this Lemma  for the case of a smooth scheme $Y$ over $\Spec(R)$ satisfying the condition that the only open subscheme of its special fibre containing its fibre over the maximal point of $\Spec(R)$ is the special fibre itself. This condition implies the existence of a Zariski dense set of good sections $\Spec(R)\to Y$; here by good we mean that, fixing an open subscheme $U$ of $Y$ such that $(Y,U)$ is an extensible pair, we take only sections $s\colon\Spec(R)\to Y$ such that the pair $(\Spec(R),s^{-1}(U))$ is also an extensible pair.
Now everything is entirely similar to Steps C and D of 3.2.17 (cf. also the proof of 3.2.2.1). This ends the proof.

\smallskip
Let now $O\hookrightarrow O_1$ be an inclusion between two discrete valuation rings which are faithfully flat over $\dbZ_{(p)}$. We assume that it is of index of ramification 1 and that $O$ is a henselian DVR. We recall that a faithfully flat inclusion $O_3\hookrightarrow O_2$ between two discrete valuation rings  is said to be of index of ramification 1, if a uniformizer of $O_3$ is a uniformizer of $O_2$, and if at the level of residue fields we get a separable field extension.

\Proclaim{3.2.2.3. Corollary.}
A) Let $Y$ be a healthy regular $O$-scheme such that one of the following two conditions is satisfied:

\smallskip
a) any maximal point of $Y$ of positive characteristic has a local ring whose henselization is a healthy regular scheme of dimension at least two;

b) any smooth scheme over a DVR of mixed characteristic which is a local ring of $Y$, is a healthy regular scheme.  

\smallskip
We have:

\smallskip
1) For any projective limit $Z$ of smooth schemes of finite type over $O$ having \'etale transition $O$-morphisms whose fibres are dominant morphisms, $Y_Z$ is a healthy normal scheme.

2) $Y_{O_1}$ is a healthy regular scheme.

B) If the completion $\Yhat$ $Y$ of a local henselian healthy regular scheme $Y$ is a projective limit of smooth affine schemes over $Y$, then this completion is a healthy regular scheme.
\finishproclaim

\proof The scheme $Y_Z$ is normal: it is a projective limit of normal schemes with dominant transition morphisms. To check 1) let $(Y_Z,U)$ be an extensible pair, and let $A$ be an abelian scheme over $U$. The conditions a) and b) imply that we can assume that there is an extensible pair $(Y,U(Y))$ such that $U=U(Y)_Z$ (in case a) cf. 3.2.2.2). We can assume that $Z$ is local, and even more using descent and C) of 3.2.2 4), we can assume that $Z$ is strictly henselian. So, as in the last paragraph of the proof of 3.2.2.1 we come back to an essentially finite type picture, i.e. we can assume that $Z$ is a localization of a smooth scheme of finite type over $O$.
So the part of the proof of 3.2.2.1 involving passage to $O^{\text{sh}}$ and taking sections applies: we get that $A$ extends to an abelian scheme over $Y_Z$.

To see 2), we can assume (cf. 3.2.2.1) that
both $O$ and $O_1$ are complete DVR's. Now everything results by using in this order  part 1), 3.2.2.1 and 3.2.2 4), once we remark that $\Spec(O_1)$ is a pro-\'etale cover of the spectrum of a DVR $O_2$, which is the completion of a henselian DVR  of whose spectrum is a projective limit of smooth affine
$O$-schemes whose transition $O$-morphisms are \'etale and have fibres which are dominant morphisms between integral schemes (as the inclusion $O\hookrightarrow O_1$ has index of ramification 1).

Part B) results from 3.2.2.2 if $Y$ is of dimension at least 2 (the case when $Y$ is of dimension 1, i.e. when $Y$ is the spectrum of a DVR, is trivial). The only extra thing we need to add: as $\Yhat$ has a finite number of points of its special fibre of codimension 1 in it, any abelian scheme over $U$, with $(\Yhat,U)$ an extensible pair, is defined over an open subscheme $U_Z$ of a smooth scheme $Z$ of finite type over $Y$, with $(Z,U_Z)$ an extensible pair, and with the natural morphisms $\Yhat\to Y$ and $U\to Y$ factoring through $Z$ and respectively through $U_Z$. This ends the proof of the Corollary.  

\Proclaim{3.2.2.3.1. Definition.} \rm
A healthy regular scheme over $\dbZ_{(p)}$ is called $R$-healthy ($R$ stands to honor the theorem of Raynaud mentioned in 3.2.1.1 2)) if the local rings of the generic points of its special fibre have indices of ramification smaller than $p-1$.
\finishproclaim

\Proclaim{3.2.2.4. Remarks.} \rm
a) Part 2) of 3.2.2.3 A) is in essence the maximum it can be said in full generality for the case of a ring homomorphism $O\to O_1$ of index of ramification 1, as the spectrum of any DVR of mixed characteristic is a healthy regular scheme, and as the condition b) of 3.2.2.3 A) is a natural one (in this context). Of course there are variants of 3.2.2.3 A) when we intermingle the conditions a) and b) of 3.2.2.3.

b) Using def. 3.2.2.3.1, from 3.2.2.3 A) we get (cf. 3.2.2 1)): a regular scheme $Y$ over $O$ is $R$-healthy iff $Y_{O_1}$ is an $R$-healthy regular scheme. 

c) There are $R$-healthy regular schemes which are not abstract healthy regular schemes (see 3.2.2 5)).
\finishproclaim

\smallskip
We start now by clarifying and restating the definitions introduced in
[Mi4, Ch. 2], and commented in the footnote of [Mi3, p. 513]. So the conjecture [Mi4, 2.7] also gets restated (see 3.2.5). 

Let $(G,X)$ define a Shimura variety and let $v$ be a prime
of $E(G,X)$ dividing the rational prime $p$.
Let $H$ be a compact open subgroup of $G(\dbQ_p)$.
We assume now that  $O$ is a faithfully flat $O_{(v)}$-algebra. Let $L$ be its field of fractions. We have $L\supset E(G,X)$. Let $f\colon (G,X)\to (G_1,X_1)$ be a map and let $H_1$ be a compact open subgroup of $G_1(\dbQ_p)$ such that $f$ takes $H$ into $H_1$. Let $v_1$ be the prime of $E(G_1,X_1)$ divided by $v$. Let $O_1$ be a DVR which is a faithfully flat $O_{(v_1)}$-subalgebra of $O$. Let $L_1$ be its field of fractions.

\Proclaim{3.2.3. Definitions.} \rm
1) An integral model of $\Sh_H(G,X)$ over $O$ is a faithfully flat scheme $\scrM$ over $O$
together with a $G(\dbA_f^p)$-continuous action and a
$G(\dbA_f^p)$-equivariant isomorphism
$$
\scrM_L\arrowsim \Sh_H(G,X)_L.
$$
When the $G(\dbA_f^p)$-action on $\scrM$ is obvious, by abuse of language, we say that $\scrM$ is an integral model.

$1^\prime$) By a (map or) morphism from an integral model $\scrM$ of $\Sh_H(G,X)$ over $O$ to an integral model $\scrM_1$ of $\Sh_{H_1}(G_1,X_1)$ over $O_1$ we mean a $G(\dbA_f^p)$-equivariant $O_1$-morphism 
$$
\scrM\to\scrM_1,
$$
whose restriction to generic fibres is the natural $L_1$-morphism $\Sh_H(G,X)_L\to \Sh_{H_1}(G_1,X_1)_{L_1}$ defined by $f$ (to be compared with 2.9).

In particular if $f$ is the identity map of $(G,X)$ we get the definition of morphisms between two integral models of $\Sh_H(G,X)$ over $O$.

2) The integral model $\scrM$ is said to be smooth (resp. normal) if there is a compact open
subgroup $H_0$ of $G(\dbA_f^p)$ such that for any
inclusion $H_2\subset H_1$ of compact open subgroups
of $H_0$, the natural morphism $\scrM/H_2\to \scrM/H_1$ is a finite \'etale morphism
between smooth schemes (resp. between normal schemes) of finite type over $O$. In other words, there is a compact open subgroup $H_0$ of $G(\dbA_f^p)$ such that $\scrM$ is a pro-\'etale cover of the smooth scheme (resp. of the normal scheme) $\scrM/H_0$ of finite type over $O_{(v)}$.

$2^{\prime}$) The integral model $\scrM$ is said to be quasi-projective, projective or proper if for any (it is enough just for one) compact open subgroup
$H_0$ of $G(\dbA_f^p)$ the scheme $\scrM/H_0$ is respectively quasi-projective, projective or proper. 

3) A scheme $T$ over $O$ is said to have the extension property, abbreviated EP (resp. the extended extension property, abbreviated EEP), if, for
any  healthy
regular scheme (resp. for any almost healthy normal scheme) $Y$ over $O$, 
every $L$-morphism $Y_L\to T_L$ extends uniquely to an
$O$-morphism $Y\to T$. Similarly, using $R$-healthy regular schemes instead of healthy regular schemes, we speak about a scheme having the $R$-extension property (abbreviated REP). 

4) A scheme $T$ over $O$ is said to have the weak extension property, abbreviated WEP (resp. the smooth extension property, abbreviated SEP), if, for any abstract very healthy regular scheme $Y$ over $O$ (resp. for any regular formally smooth scheme $Y$ over a DVR which is a faithfully flat $O$-algebra and has the same ramification index as $O$), every $L$-morphism $Y_L\to T_L$ extends uniquely to an $O$-morphism $Y\to T$.

5) A scheme $T$ over $O$ is said to have the quasi extension property, abbreviated QEP (resp. the local extension property, abbreviated LEP), if, for any quasi-compact healthy regular scheme (resp. for any locally healthy regular scheme) $Y$ over $O$, every $L$-morphism $Y_L\to T_L$ extends uniquely to an $O$-morphism $Y\to T$. Similarly we define the quasi extended extension property (abbreviated QEEP).

6) A smooth integral model of $\Sh_H(G,X)$ over $O_{(v)}$ (resp. over its completion $\widehat{O_{(v)}}$) having the EP is called an integral canonical model (resp. a local integral canonical model) of our Shimura variety $\Sh(G,X)$ with respect to $v$ and $H$ (or simply with respect to $H$ as the prime $v$ is determined by the integral model).
\finishproclaim

\Proclaim{3.2.3.0. Remark.} \rm
Other extension properties can be defined starting from quasi healthy schemes, or from locally healthy schemes. Not to be to long, this is not going to be done here.
\finishproclaim 

\Proclaim{3.2.3.1. Remarks.} \rm
0) Allowing $(G,X)$ and $v$ to vary we get that def. 3) to 5) of 3.2.3 make sense for any DVR which is a faithfully flat $\dbZ_{(p)}$-algebra. Moreover 1) and $1^{\prime})$ of 3.2.3 make sense for any compact subgroup $H$ of $G(\dbA_f^p)$ not necessarily open, but for 2) and $2^{\prime})$ of 3.2.3 we do need to assume that $H$ is also open.
 
1) Any scheme over $O$ having the EEP, has the EP (cf. 3.2.2 2)), and any scheme over $O$ having the $EP$, has the WEP (cf. 3.2.2 3)). If $e<p-1$ then any scheme over $O$ having the $WEP$ has the $SEP$. We do not know when the converses of these statements are true.  

2) Any quotient $\scrM/H_0$ (with $H_0$ a compact open subgroup of $G(\dbA_f^p)$) of a normal integral model $\scrM$ of $\Sh_H(G,X)$ over $O$ having the EP, is separated.

To see this we first remark that any DVR of mixed characteristic is a healthy regular scheme. We use the valuative criterion of separatedness. We need to check it just for a DVR of mixed characteristic, as $\scrM$ is a pro-\'etale cover of the normal scheme $\scrM/H_0$ of finite type over $O$ and has a separated generic fibre. Now everything results from the EP.

3) A scheme $Y$ over $O$ has any of the extension properties we defined iff the reduced scheme $Y_{\text{red}}$ attached to it has it. A reduced scheme $Y$ over $O$ has any of the extension properties we defined iff any connected component of its normalization in its ring of fractions has it. This reduces the study of schemes over $O$ having an extension property to the case of integral normal schemes over $O$. All these results from the fact that we defined the different extension properties in terms of normal schemes.

4) Any scheme over $O$ having the EP (resp. EEP) has the LEP and the QEP (resp. has the QEEP). We do not know if (or when) the converse is true.  

5) If $Y$ is a scheme over $O$ having any type of extension property, and if $Y_{1L}$ is a closed reduced subscheme of $Y_L$, then the Zariski closure  $Y_1$ of $Y_{1L}$ in $Y$ also has the same type of extension property. Moreover: the normalization of $Y_1$ in any pro-\'etale scheme over the spectrum  of the ring of fractions of $Y_1$ has the same type of extension property. We will use this trivial fact without any further comment.  

$5^{\prime})$ If $Y$ is an $O$-scheme having the EP, and if $q\colon Y\to Y_1$ is a morphism which is an isomorphism on generic fibres, then $Y_1$ has the EP. This remains true for any of the extension type properties we defined above.  

6) If $Y_1\to Y$ is a pro-\'etale  cover of $O$-schemes, then $Y_1$ has the EP (or QEP, or WEP, or SEP) iff $Y$ has it (for the EP
and QEP this is a consequence of C) of 3.2.2 4); for the WEP and SEP cf. def. 4) of 3.2.3).

7) A regular formally smooth scheme over $O$ having the SEP is uniquely determined by its generic fibre.
\finishproclaim

\Proclaim{3.2.3.2.} \rm
Let $\scrM$ be a smooth integral model of $\Sh_H(G,X)$ over $O$. Let $H_0$ be a compact open subgroups of $G(\dbA_f^p)$ such  that the quotient morphism  $\scrM\to\scrM/H_0$ is a pro-\'etale cover. 

\Proclaim{Proposition.}
a) If $\scrM$ has the SEP (resp. WEP or EP) then $\scrM/H_0$  has the following extension type property:
If $(Y,U)$ is an extensible pair with $Y$ a regular formally smooth scheme over a DVR $O_1$ which is a faithfully flat $O$-algebra having the same index of ramification as $O$ (resp. with $Y$ an abstract very healthy regular scheme, or resp. with $Y$ a healthy regular scheme), then any morphism $U\to\scrM/H_0$ extends uniquely to a morphism $Y\to\scrM/H_0$; 

b) We assume that $\scrM$ has the SEP and satisfies the valuative criterion of properness  with respect to discrete valuation rings of mixed characteristic (for instance these hold if $\scrM$ has the EP and $e<p-1$). We have:

i) Let $\scrM^0/H_0$ be an open closed subscheme of $\scrM/H_0$, and let $q^0\colon\scrM^0/H_0\to Z$ be a proper morphism, with $Z$ a faithfully flat scheme of finite type over $O$, which is an isomorphism on generic fibres. Then the natural map $\scrM^0/H_0(O^{\text{sh}})\to Z(O^{\text{sh}})$ is a bijection;

ii) The morphism $q^0$ of i) is in fact finite.
\finishproclaim

\proof
The proof of a) is a consequence of the classical purity theorem and of the fact that the class of schemes $Y$ mentioned in a) are stable under pro-\'etale covers (cf. rm. 4) of 3.2.2). We now prove b).

i) We can assume that $O=O^{\text{sh}}$, that $Z$ is normal and connected and that $\scrC:=\scrM^0/H_0$ is connected. We first remark that $\scrC$ is a separated scheme (the argument is the same as in 3.2.3.1 2)) of finite type over $O$ (cf. def. 2) of 3.2.3) and has a smooth quasi-projective generic fibre (the generic fibre is a model of the quotient of a Hermitian symmetric domain by an arithmetic subgroup).

$i_A)$ We consider a proper morphism $q\colon\scrC\to Z$, with $Z$ a faithfully flat scheme of finite type over $O$, having the properties mentioned in i) of b). From the smoothening process (cf. [BLR, Th. 3, p. 61]) we deduce the existence of a $Z$-scheme $Z_1$, smooth over $O$, quasi-projective over $Z$, and having the property that the induced map $Z_1(O)\to Z(O)$ is a bijection. Moreover the generic fibre of $Z_1$ is the same as the generic fibre of $Z$ (or of $\scrC$).

$i_B)$ Let $\scrC^1$ be a connected component of $\scrM$ which is a pro-\'etale cover of $\scrC$. So $\scrC^1\times_{\scrC} \scrC_1$ is a pro-\'etale cover of the normal $O$-scheme $\scrC_1$ of finite type. As $q$ is proper, the normalization of any local ring of $Z_1$ which is a DVR in the field of fractions of $\scrC^1$ is a regular ring of dimension 1. This implies the existence of a morphism from the spectrum of any such normalization into $\scrC^1$ (as $\scrM$ satisfies the valuative criterion of properness with respect to discrete valuation rings of mixed characteristic). So there is a rational map $q_1$ from $Z_1$ to $\scrC$ defined on points of codimension 1, and inducing an isomorphism on generic fibres. From the mentioned extension type property of $\scrM/H_0$ (cf. a)), which is also enjoyed by its connected component $\scrC$, we deduce that $q_1$ is in fact a morphism. Moreover the induced maps $Z_1(O)\to\scrC(O)\to Z(O)$ are bijections. This implies that $q_1$ is a surjective morphism (as $\scrC$ is a smooth scheme and as $O=O^{\text{sh}}$). So i) of b) holds.

ii) As $Z$ is normal, we need to show that $q$ is an isomorphism. We just need to show that $q$ is an isomorphism in codimension 1. If this is not so, then there is a connected component $\scrC_p$ of the special fibre of $\scrC$ dominating a reduced closed subscheme $Z_p$ of the special fibre of $Z$ of dimension $d<\dim(\scrC_p)$. So $Z_p$ is a closed subscheme of the non-smooth locus of $Z$. Let $\tilde\scrC$ be the open subscheme of $\scrC$ defined by $\scrC_p$ and the generic fibre of $\scrC$. 

From [BLR, p. 72] we deduce that the morphism $\tilde\scrC\to Z$ lifts to a morphism $\tilde q_p\colon\tilde\scrC\to\tilde Z$, where $\tilde Z$ is obtained from $Z$ through the first blowing up needed to get $Z_1$: we always blow up a reduced connected component of the maximal reduced closed subscheme $S_Z$  of the special fibre of $Z$ having the property that it is included in the non-smooth locus of $Z$ and the points of it with values in the residue field of $O^{\text{sh}}$ which admit lifts (in $Z$) to $O^{\text{sh}}$-valued points, are Zariski dense in it. 
 As $\tilde\scrC$ is smooth and its fibres over $Z$ are proper schemes (over residue fields of points of $Z$), we deduce that $\tilde q_p$ dominates a closed subscheme $\tilde Z_p$ of the special fibre of $\tilde Z$ of the same dimension $d$: the morphism $\tilde\scrC\to\tilde Z$ factors through an open subscheme of $\tilde Z$ which is affine over $Z$, cf. the properties of dilatations [BLR, p. 62]. So $\tilde Z_p$ is included in the non-smooth locus of $\tilde Z$. We can apply induction to get that $q_1$ has a section above $\tilde\scrC$  such that $\scrC_p$ dominates a closed subscheme of the special fibre of $Z_1$ of dimension $d$. Contradiction. We conclude that $q$ is an isomorphism in codimension 1, and so an isomorphism. This ends the proof of the Proposition.

\Proclaim{Expectations.} \rm
Under the hypotheses of b) above we expect that the following statements can be proved without assuming that $\scrM$ is a quasi-projective integral model:

\smallskip
iii) If $\Spec(O)\to\Spec(O_1)$ is a finite Galois cover, with $O_1$ an $O_{(v)}$-subalgebra of $O$, then $\scrM$ is the extension to $O$ of a smooth integral model $\scrM_1$ of $\Sh_H(G,X)$ over $O_1$. Also $\scrM_1$ inherits the properties of $\scrM$ mentioned in b);

iv) The quotient of $\scrM/H_0$ through a finite free action exists as a scheme (not only as an algebraic space).
\finishproclaim

We present the reasons for these expectations.

iii) To prove iii) we can assume that both $O$ and $O_1$ are complete (for instance cf. Raynaud's result mentioned in [BLR, p. 166])). Let $C:=\Gal(O/O_1)=\Gal(k/k_1)$, with $k$ and $k_1$ the residue fields of $O$ and respectively of $O_1$. We view $C$ as a finite \'etale group scheme over $O$. Due to the fact that $\scrM$ has the SEP and that its generic fibre is definable over the field of fractions of $O_1$ (being definable over $E(G,X)$) we deduce the existence of a natural action of $C$ on $\scrM$, compatible with the action of $G(\dbA_f^p)$ on $\scrM$. It provides us with a Galois-descent datum (see [BLR, 6.2]). To prove iii) we just have to show that it is effective. It is enough to work with $\scrM/H_0$ instead of $\scrM$.

$iii_A)$ From [BLR, Lemma 4, p. 155] and [Mu1, p. 112] we deduce the existence of a quasi-projective smooth scheme $U^1$ of finite type over $O_1$ such that $U^1_O$ is an open subscheme of $\scrM/H_0$ containing the generic fibre and all the points of codimension 1. Let $Z$ be a faithfully flat projective scheme over $O_1$ having $U^1$ as an open subscheme. We can assume that its generic fibre is smooth (cf. the resolution of singularities in characteristic 0). We can assume that the generic fibre of $U^1$ is dense in the generic fibre of $Z$. We get a rational map from $\scrM/H_0$ to $Z_O$ defined on the generic fibre and in points of codimension 1. 

$iii_B)$ 
Let $y\colon \Spec(k_2)\hookrightarrow\scrM/H_0$ be an arbitrary maximal point of positive characteristic. Here $k_2$ is a finite field extension of $k_1$. Let $\Spec(O_2)$ be the \'etale cover of $\Spec(O)$ having $\Spec(k_2)$ as its special fibre. Let $z\colon \Spec(O_2)\hookrightarrow\scrM/H_0$ be an arbitrary lift of $y$.
Let $Z_1$ be obtained from $Z$ as above, using the smoothening process. So $Z_1$ is the smooth locus of a scheme $Z_1^{\prime}$ obtained from $Z$ through a sequence of blowings up centered on special fibres. We get a natural bijection $Z_1(O^{\text{sh}})\to Z(O^{\text{sh}})$. As $Z$ is a projective $O$-scheme, we can view $z$ as an $O_2$-valued point $z_2$ of $Z_1$. Let $y_2\colon \Spec(k_1^{\prime})\hookrightarrow Z_1$ be the maximal point of the special fibre of $Z_1$ through which the $k_2$-valued point of $Z_1$ defined by $z_2$, factors. Let $\Spec(O_1^{\prime})$ be the \'etale cover of $O_1$ having $k_1^{\prime}$ as its residue field. Let $\Spec(O_1^{\prime})\hookrightarrow Z_1$ be a lift of $y_2$. Let $\Spec(O^{\prime})$ be the Galois cover of $\Spec(O_1^{\prime})$ generated by $O$. Let $W_1$ be the closed subscheme of $Z_1$ which is the Zariski closure of the closed subscheme of its generic fibre defined by the complement of the generic fibre of $U^1$. Let $Z_2$ be the open subscheme of $Z_1$ defined by the complement of $W_1$. As in $i_B)$  we get a morphism $q\colon Z_{2O}\to\scrM/H_0$, which at the level of generic fibres is an isomorphism.

$iii_C)$ We can assume that $z_2$ factors through $Z_2$. To see this we have to use blowings up centered on special fibres. First we blow up $y_2$ on $Z_1$. We get similarly a point $y_2^{\prime}$ on the resulting scheme $Z_2^{\prime}$. Now we blow up $y_2^{\prime}$ on $Z_2^{\prime}$. After a finite number of operation we achieve the separation of the point $z_2$ from $W_1$. This is possible due to the fact that in characteristic zero we do have such a separation: let $O_{y_2}$ be the local ring of $y_2$ in $Z_2$, and let $n\in\dbN$ be the valuation (with respect to the normalized valuation of $O_2$) of the image in $O_2$ (through the epimorphism $O_{y_2}\twoheadrightarrow O_2$ defined by $z_2$) of an element of $O_{y_2}$ defining $W_1$ in $\Spec(O_{y_2})$; after at most $n$ blowings up we achieve the desired separation.

So we get a  $C$-equivariant morphism $\Spec(O_{y_2}\otimes_{O_1} O)\to\scrM/H_0$. Its image contains the $C$-orbit of $y$ in $\scrM/H_0$. The same is true for any other maximal point of $Z_2$ whose inverse image to $Z_{2O}$ dominates the $C$-orbit of $y$. So this orbit should be contained in an affine open scheme of $\scrM/H_0$. If $O_1^{\prime}=O_1$, this is obvious. The general case should be handleable by standard arguments on local rings: we just need to show that the intersection of the local rings of the points of the $C$-orbit of $y$ is a semi-local ring whose localizations  with respect to maximal ideals are local rings of the points of the $C$-orbit of $y$; this should be provable using the fact that $q$ is an isomorphism above points of $\scrM/H_0$ of codimension 1, starting from [Ma, Th. 38].

We assume now that we were able to get that the $C$-orbit of $y$ is contained in an affine open subscheme of $\scrM/H_0$.
As $y$ was an arbitrary maximal point of the special fibre of $\scrM/H_0$, from [Mu1, p. 112] we deduce that the quotient of $\scrM/H_0$ through the action of $C$ exists as a scheme. This scheme is $\scrM_1/H_0$. Taking its normalization in the ring of fractions of the extension of $\Sh_H(G,X)$ to the field of fractions of $O_1$, we get the desired  integral model $\scrM_1$ of $\Sh_H(G,X)$ over $O_1$ (obviously $\scrM_{1O}=\scrM$). The last part of iii) involving the inheritance property is trivial.   

iv) The above ideas of iii) can be entirely adapted for the case of quotients. The easy details are left as an exercise. We just need to replace the operation of extension of scalars (from $O_1$ to $O$) used above, by the operation of taking the normalization (of a reduced scheme whose ring of fractions is the subring of the ring of fractions $\scrF$ of $\scrM/H_0$ fixed by the action) in $\scrF$.
\finishproclaim

\Proclaim{3.2.3.2.1. Remarks.} \rm
1) We call the part of i) of 3.2.3.2 b) involving $O^{\text{sh}}$-valued points as the maximality property.

2) We think it is possible to prove that $\scrM/H_0$ is a quasi-projective scheme over $O$ by just refining 3.2.3.2. In the case when $(G,X)$ is of preabelian type and $(v,6)=1$ we will prove this in 6.4.1 using the extra fact that different schemes related to $\scrM$ are moduli schemes of abelian varieties (subject to some conditions).  

3) In [Va6] we will develop the general theory of integral canonical models of smooth  schemes of finite type over the field of fractions of a Dedekind domain (of mixed characteristic), starting from 3.2.3.2 and rm. 1) of 6.4.6.

$3^{\prime})$ The ideas and results of 3.2.3.2 can be used in a much larger context (not involving Shimura varieties). For instance for a) we just used the fact that $\scrM/H_0$ has a pro-\'etale cover having some extension type property, while for expectation iii) (resp. iv))  we used (besides the mentioned fact) the fact that the descent (resp. the quotient) we are dealing with is known to be effective at the level of generic fibres.

4) Expectation iii) is not true in the larger context if the finite morphism $\Spec(O)\to\Spec(O_1)$, with $O_1$ a DVR, is not an \'etale cover, as it can be easily seen through examples involving N\'eron models of abelian varieties. 
\finishproclaim

\Proclaim{3.2.3.3. Proposition.}
Let $i_O\colon O\hookrightarrow O_1$ be a faithfully flat inclusion of discrete valuation rings, with $O_1$ having also $e$ as its index of ramification. We have:

1) A scheme $Y$ over $O$ has the WEP or the SEP iff $Y_{O_1}$ has it.

2) If moreover $O$ is a henselian local ring and if $i_O$ is formally \'etale, then a scheme $Y$ over $O$ has the EP (or QEP) iff $Y_{O_1}$ has it.

3) If $i_O$ has index of ramification 1, then a scheme $Y$ over $O$ has the REP iff $Y_{O_1}$ has the REP.
\finishproclaim

\proof
We just need to check that the class of schemes involved in the definition of these extension properties is stable under pull backs via $i_O$ and that any $O_1$-scheme belonging to such given class, as an $O$-scheme also belongs to the given  class. This last part is trivial, while the first part is a direct consequence of def. 4) of 3.2.3 for 1), of  3.2.2.1 for 2), and of 3.2.2.4 b) for 3). 

\Proclaim{3.2.3.4. Remark.} \rm
We do expect that the  condition on $O$ of being a henselian DVR used in 3.2.3.3 2) is not needed. For this we need to prove that for any \'etale morphism $\Spec(O_1)\to\Spec(O)$, with $O_1$ a DVR, a scheme $Y$ over $O$ is a healthy regular scheme iff $Y_{O_1}$ is so.   
\finishproclaim

\Proclaim{3.2.4. Remark.} \rm
We assume that $G$ is unramified over $\dbQ_p$ and that $H$ is a hyperspecial subgroup of $G(\dbQ_p)$. Then, if $p>2$, any (local) integral canonical model $\scrN$ of $\Sh_H(G,X)$ is uniquely determined up to a unique isomorphism (cf. 3.2.3.1 7); i.e $\scrN$ has the SEP as it has the EP: this results from 3.2.2 1) and from [Mi3,
4.7] which shows that $v$ is unramified over $p$). If $p=2$ then we know the uniqueness of an integral canonical model of $\Sh_H(G,X)$ only when $G$ is a torus (cf. 3.2.8).  
\finishproclaim

\Proclaim{3.2.5. Milne's conjecture [Mi4].}
If $G$ is unramified over $\dbQ_p$ and if $H$ is a hyperspecial
subgroup of $G(\dbQ_p)$, then $\Sh_H(G,X)$ has an integral
canonical model with respect to  $v$ and $H$.
\finishproclaim

\Proclaim{3.2.6. Notations and Definitions.} \rm
By $(G,X,H,v)$ we always denote a quadruple where: $(G,X)$ defines a
Shimura variety, $v$ is a prime of $E(G,X)$ dividing a rational prime $p$ such that $G$ is unramified over $\dbQ_p$,
and $H$ is a hyperspecial subgroup of $G(\dbQ_p)$.
The maps from a quadruple $(G,X,H,v)$ into another quadruple $(G_1,X_1,H_1,v_1)$ are defined by maps $f\colon (G,X)\to (G_1,X_1)$ taking $H$ into $H_1$ and inducing
an inclusion $E(G,X)\supset E(G_1,X_1)$ with $v$ dividing
$v_1$. We denote it by $f\colon (G,X,H,v)\to (G_1,X_1,H_1,v_1)$. The map $f$ is called injective, or finite, or a cover if as a map $f\colon (G,X)\to (G_1,X_1)$ of Shimura pairs it is so. If $(G,X,H,v)$ is a quadruple then $(G^{\ad},X^{\ad},H^{\ad},v^{\ad})$ (with $H^{\ad}$ as in the part b) of 3.2.7 2) and with $v^{\ad}$ the prime of $E(G^{\ad},X^{\ad})$ divided by $v$) is called its adjoint quadruple and $(G^{\ab},X^{\ab},H^{\ab},v^{\ab})$ (with $H^{\ab}$ the only hyperspecial subgroup of $G^{\ab}(\dbQ_p)$ and with $v^{\ab}$ the prime of $E(G^{\ab},X^{\ab})$ divided by $v$) is called its toric part quadruple.
We have maps from $(G,X,H,v)$ into its adjoint and toric part quadruples.

By $(G,X,H)$ we always denote a triple which can be extended to a quadruple $(G,X,H,v)$. The definitions of maps between quadruples extend to triples. We also speak about the adjoint and the toric part triple of a triple $(G,X,H)$.

By an integral canonical model of a quadruple $(G,X,H,v)$ we
mean an integral canonical model of $\Sh_H(G,X)$ over $O_{(v)}$. We denote it by $\Sh_v(G,X,H)$.  It is clear what we mean by $\Sh_v(G,X,H)$ having the EEP. Similarly, we speak about integral smooth (or normal) models of $(G,X,H,v)$ over $O$, or about a local integral canonical model of $(G,X,H,v)$. We say that $(G,X,H,v)$ or $(G,X,H)$ is of abelian (preabelian, etc.) type if $(G,X)$ is so.

If all the quadruples $(G,X,H,v)$ extending a triple $(G,X,H)$ have uniquely determined integral canonical models, then we denote by $\Sh_p(G,X,H)$ the model of $\Sh_H(G,X)$ over the normalization of $\dbZ_{(p)}$ in $E(G,X)$, obtained by gluing along their generic fibres the integral canonical models of all these quadruples. We call it the integral canonical model of the triple $(G,X,H)$. Similarly we define a (smooth or normal) integral model over $O$ of $(G,X,H)$. The rm. 2) of 3.2.7 shows that if $\Sh_p(G,X,H)$ exists, then for any other hyperspecial subgroups $H_1$ of $G(\dbQ_p)$, $\Sh_p(G,X,H_1)$ exists and as a scheme it is isomorphic to $\Sh_p(G,X,H)$. This means that it is irrelevant with which hyperspecial subgroup $H$ of $G(\dbQ_p)$  we work and so we sometimes write $\Sh_p(G,X)$ instead of $\Sh_p(G,X,H)$ and $\Sh_v(G,X)$ instead of $\Sh_v(G,X,H)$. We say that $\Sh_p(G,X)$ exists if for a (any) hyperspecial subgroup $H$ of $G(\dbQ_p)$, $\Sh_p(G,X,H)$ exists. We call $\Sh_p(G,X)$ the $\dbZ_{(p)}$-model or the $\dbZ_{(p)}$-canonical model of our Shimura variety $\Sh(G,X)$. We say that $\Sh_p(G,X,H)$ has the EP (or the EEP) if as a $\dbZ_{(p)}$-scheme it has it. 
\finishproclaim

\Proclaim{3.2.7. Remarks.} \rm
1) Milne's conjecture can be reformulated: any quadruple $(G,X,H,v)$ has an integral canonical model.

2) If a quadruple $(G,X,H,v)$ has an integral
canonical model $\scrM$, then any other quadruple of the form
$(G,X,H_1,v)$ has also an integral canonical model,
which is isomorphic to $\scrM$ as an $O_{(v)}$-scheme.
This results from the following fact:

\Proclaim{3.2.7.1.} 
Under the canonical action of $\Aut(Sh(G,X))$ on $G$ (cf. 2.4.3) and so on $G(\dbQ_p)$, the hyperspecial subgroups of $G(\dbQ_p)$ are permuted transitively. 
\finishproclaim

So actually $(G,X,H,v)\arrowsim (G,X,H_1,v)$. To see this we first remark that:

\smallskip
\item{a)}
Any two hyperspecial subgroups of $G(\dbQ_p)$ are conjugate
by an element of $G^{\ad}(\dbQ_p)$ [Ti, p. 47].

\smallskip
\item{b)}
There is a hyperspecial subgroup $H^{\ad}$ of
$G^{\ad}(\dbQ_p)$ normalizing $H$ $(H^{\ad}$ is the group of
$\dbZ_p$-valued points of the quotient 
$G_{\dbZ_p}^{\der}/Z$, where $G_{\dbZ_p}^{\der}$
is the derived subgroup of the reductive group $G_{\dbZ_p}$ over $\dbZ_p$
 having $G_{\dbQ_p}$ as its generic fibre and having $H$ as
its group of $\dbZ_p$-valued points, and where $Z$ is the center
of $G_{\dbZ_p}^{\der}$).

\smallskip
\item{c)}
$G^{\ad}(\dbQ_p)=G^{\ad}(\dbQ)H^{\ad}$ [Mi3, 4.9].

\smallskip
\item{d)}
If $g\in G^{\ad}(\dbQ)$ takes $X$ onto $X$, then
$(G,X,H_1,v)$ has an integral canonical model if and only if
$(G,X,gH_1g^{-1},v)$ has an integral canonical model.

\smallskip
\item{e)}
$G^{\ad}(\dbZ_{(p)}):=G^{\ad}(\dbQ)\cap{H^{\ad}}$ permutes transitively the connected components of $X^{\ad}$ (cf. 3.3.3).

\smallskip
\item{f)}
If an element of $G^{\ad}(\dbR)$ leaves invariant a connected component of $X$, it leaves invariant $X$.

\smallskip
So a), b) and c) imply that there is $g\in G^{\ad}(\dbQ)$ such that $H_1=gHg^{-1}$. From e) we get that we can replace $g$ with $gh$, with  $h\in G^{\ad}(\dbZ_{(p)})$, in such a way that $gh$ takes a (fixed) connected component $X^0$ of $X$ into itself. So f) implies that $gh\in G^{\ad}(\dbQ)$ produces by inner conjugation (of $G$) an isomorphism $(G,X,H,v)\arrowsim (G,X,H_1,v)$.

\smallskip
The integral canonical model $\scrM$ of our Shimura variety $\Sh(G,X)$ with respect to $v$ and $H$, will be often referred to as an integral canonical model of $\Sh(G,X)$, as the prime $v$ is determined by it and as  it is irrelevant with which hyperspecial subgroup we work. Similarly, we will often speak about a local integral canonical model of a Shimura variety, without mentioning the hyperspecial subgroup and the prime with respect to which it is defined.

3) The category $qf-Sh$ ($tr-Sh$) whose objects are quadruples (respectively triples) as in 3.2.6 and whose morphisms are finite maps between them has quasi fibre products (as in 2.4.0). If $f_i\colon (G_i,X_i,H_i,v_i)\to (G_0,X_0,H_0,v_0)$, $i=\overline{1,2}$, are finite maps such that the intersection $X_1\cap X_2$ is not empty (see 2.4.0), then a quasi fibre product of $f_1$ and $f_2$ is described by maps $p_i^j\colon (G_3,X_3^j,H_3,v_3)\to (G_i,X_i,H_i,v_i)$, $i=\overline{1,2}$, where $(G_3,X_3^j)$ is as in 2.4.0, $H_{3}:=(H_{1}\times H_{2})\cap{G_{3}(\dbQ_p)}$, and $v_3$ is uniquely determined as $E(G_3,X_3^j)$ is the composite field of $E(G_1,X_1)$ and $E(G_2,X_2)$. 

If $f_1$ or $f_2$ is a cover then the set $I$ introduced in 2.4.0 has precisely one element; so we speak about the fibre product of $f_1$ and $f_2$.

\smallskip
This allows us to define the standard quadruple  situation of Shimura varieties of preabelian  type (abbreviated SQSPT). For a given quadruple $(G,X,H,v)$ of preabelian type, this is a commutative diagram

$$
\CD
(G_4,X_4,H_4,v_4) @>{p_4}>> (G_3,X_3,H_3,v_3) @>{p_1}>>
(G_1,X_1,H_1,v_1) \\
@VV{p_2}V @VV{p_3}V @VV{f_1}V \\
(G_2,X_2,H_2,v_2) @>{f_2}>> \ (G,X,H,v)
@>{f_0}>> (G^{\ad},X^{\ad},H^{\ad},v^{\ad})
\endCD
$$

such that:

\smallskip
a) all its maps are finite;

\smallskip
b) the two squares are quasi fibre products;

\smallskip
c) $f_2$ is a cover with $E(G,X)=E(G_2,X_2)$ (see 10) below);

\smallskip
d) $G_2^{\der}$ is either a simply connected semisimple group, or is isomorphic to $G_1^{\der}$ (as we  need); in both situations we have $G_4^{\der}=G_2^{\der}$;

\smallskip
e) there is an injective map $f\colon (G_1,X_1,H_1,v_1)\hookrightarrow (GSp(W,\psi),S,K_p,p)$.

\smallskip
To show its existence once we assume the existence of $f$ and $f_1$ (cf. 6.4.2), we just need to modify the map $f_1$ in such a way that the intersection of $X_2$ and $X_1$ (inside $X^{\ad}$) is non-empty. As $G^{\ad}(\dbZ_{(p)}):=G^{\ad}(\dbQ)\cap H^{\ad}$ permutes transitively the connected components of $X^{\ad}$ (cf. 3.3.3), by composing an arbitrary map $f_1$ with an automorphism (cf. 9) below) of $(G^{\ad},X^{\ad},H^{\ad})$, we can always achieve a non-empty intersection $X_1\cap X_2$.

\smallskip
When $G_1^{\der}=G_2^{\der}$, all the quadruples of the above diagram are of abelian type, and then we refer to it as the standard quadruple situation of Shimura varieties of abelian type (abbreviated SQSAT).

4) Let ICM-Sh (ICM-tr-Sh) be the category whose objects are quadruples $(G,X,H,v)$ (resp. triples $(G,X,H)$) having an integral canonical model and satisfying $(v,2)=1$ (resp. satisfying $(p,2)=1$, where $p$ is the prime such that $H\subset G(\dbQ_p)$), and whose morphisms are the  maps between quadruples (resp. triples). Any such integral canonical model is formally smooth over the localization of $\dbZ$ with respect to some prime $p>2$ and has the SEP (cf. 3.2.4). So we have a functor $\scrF$ from ICM-Sh (ICM-tr-Sh) to the category of schemes: it associates to a quadruple $(G,X,H,v)$ (resp. to a triple $(G,X,H)$) its integral canonical model $\Sh_v(G,X,H)$ (resp. $\Sh_p(G,X,H)$, with $p$ as before), and to a map $(G,X,H,v)\to (G_1,X_1,H_1,v_1)$ (resp. $(G,X,H)\to (G_1,X_1,H_1)$) the morphism 
$$
\Sh_v(G,X,H)\to \Sh_{v_1}(G_1,X_1,H_1)
$$ 
(resp. $\Sh_p(G,X,H)\to \Sh_p(G_1,X_1,H_1)$) whose generic fibre is the natural morphism $\Sh_H(G,X)\to \Sh_{H_1}(G_1,X_1)$.

$4^{\prime}$) With the notations and definitions of 1) and $1^\prime$) of 3.2.3, we get the category $SIM(\Sh_H(G,X),O)$ of smooth integral models of $\Sh_H(G,X)$ over $O$. If there is such an integral model having the SEP, then as an object of this category, it is a final object.
  
5) The definition of a healthy or of an almost healthy normal  scheme appeals to abelian schemes, while the definition of an abstract very healthy regular scheme is intrinsic. We could have defined the notion of an integral canonical model of a Shimura variety using the WEP (or SEP) instead of the EP. Defining it using  the WEP instead of EP or even instead of SEP would have been definitely more convenient (and then we would have been speaking about integral canonical models having the EP).  We preferred to work out def. 6) of 3.2.3 using the EP due to the following reasons:

\smallskip
-- it is closer to the spirit of Milne's original (though inadequate, cf. footnote of [Mi3, p. 513]) definition in [Mi4, Ch. 2];  

-- the philosophy of 6) below;

-- it makes sense and works also for $p=2$: the WEP is enjoyed by any scheme over $\dbZ_{(2)}$, and we just hope that the SEP works for $p=2$ (cf. 3.2.1.4 5) and 3.2.9);

-- all integral canonical models of Shimura varieties (of preabelian type) whose existence we are able to prove in this paper (or in [Va2-3] and [Va5]) have the EP (and so they have the WEP and the SEP);

-- the worries that 3.2.3 6) might not work for Shimura varieties which are not  of preabelian type are not so justified (cf. 8) below);

-- the greatest advantage of using the EP instead of the SEP (and even instead of the WEP) consists in the fact that in this way we can get (the simplest way is by extension of scalars; but there are other ways like dealing with cases of bad reduction or like taking quotients of extended integral canonical models to be introduced in 3.5.1) (very often uniquely determined) (smooth or normal) integral models having the EP, of some quotients of Shimura varieties (of preabelian type) over discrete valuation rings which do not have the index of ramification 1 (or some $e\in\dbN$, $e<p-1$) (cf. also rm. 3) of 3.2.3.2.1);

-- it it easy to see, using N\'eron models and the fact that any DVR of mixed characteristic defines a  healthy scheme, that the EP is a stronger property than the WEP or than the SEP (cf. also 3.2.3.1 1)). 

\smallskip
6) In our philosophy (cf. [Va6]), the healthy regular schemes over $\Spec(\dbZ)$ are forming  the largest class $\scrR$  of regular schemes over $\Spec(\dbZ)$ which contains all the smooth schemes over $\Spec(\dbZ\fracwithdelims[]12)$ and it is such that  for any extensible pair $(Y,U)$, with $Y$ a regular scheme (belonging to $\scrR$)  over a Dedekind ring $D$ faithfully flat over a localization of $\dbZ$, every morphism from $U$ to a  familiar smooth moduli scheme over $D$ (such as moduli of semistable curves, of semistable vector bundles of a projective smooth curve, of polarized abelian schemes satisfying some extra conditions, etc.) extends uniquely to a morphism from $Y$ into that moduli scheme over $D$.

7) In  3.2.3 1) we could have defined an integral model $\scrM$ (of $\Sh_H(G,X)$) without requiring that $\scrM$ is faithfully flat over $O$. But we can not see any use of such integral models $\scrM$  which are not faithfully flat: the Zariski closure $\scrM_1$ of $\scrM_L$ in $\scrM$ is ``the only part of $\scrM$ influenced (controlled) by $\scrM_L$". So it makes no sense to say that $\scrM$ is an integral  model of $\scrM_L=\Sh_H(G,X)$.

8) It is well known (cf. \S4) that Shimura varieties of Hodge type are moduli schemes of principally polarized abelian schemes of a given dimension, endowed with a family of Hodge cycles and some level structures, and satisfying some additional conditions. So it looks reasonable to define an integral canonical model of a Shimura variety of preabelian type (cf. Definitions 3 of 2.5) in the way we did. As in this paper we are dealing only with Shimura varieties of preabelian type, we would like to indicate briefly why the def. 6) of 3.2.3 of an integral canonical model of a Shimura variety should work also for Shimura varieties which are not of preabelian type. We have four reasons for this:

\smallskip
\item{a)} We expect the possibility of interpreting a large class of quotients of Shimura varieties of special type over the completion of their reflex fields in finite primes, as moduli schemes of $p$-divisible groups (or of something similar) endowed with tensors (a notion with which we will be dealing extensively in [Va2]; here, for a glimpse of what we have in mind see 5.6.5). Remarks 1) and 3) of 3.2.2, together with the expectations of 3.2.1.4 6), of 3.2.3.4 and of 3.2.1.2, do motivate why we dared to work with the EP instead of the WEP (for a scheme which is a moduli of $p$-divisible groups).

\item{b)} There are generalized Shimura filtered $\sigma$-crystals of special type (cf. [Va2] for the meaning of this). Here we just give an idea: for instance, there are quadruples $(M,\varphi,(u_{\alpha})_{\alpha\in\scrJ}, G_{W(k)})$ as in 5.6.5 satisfying d), and a variant of f) and g) of 5.6.5, and such that $G_{W(k)}^{\ad}$ is a simple adjoint group of $E_7$ Lie type, etc. The local deformation theory of 5.4 remains true for generalized Shimura filtered $\sigma$-crystals (cf. [Va2]).

\item{c)} The philosophy of 6) above.

\item{d)} The philosophy of [Mi1, paragraph 9, p. 343-345]. 

\smallskip
Moreover once we know the existence of local integral canonical models of Shimura varieties of special type, we should be able to get, using the above four reasons (and 6.4.1), the existence of integral canonical models of Shimura varieties of special type.

9) The group $\Aut((G,X,H))$ of automorphisms of a triple $(G,X,H)$ (or of a quadruple $(G,X,H,v)$) is the subgroup of $Aut(G_{\dbZ_{(p)}})(\dbZ_{(p)})$ (it is of finite index if $G$ is an adjoint group) leaving $X$ invariant (cf. 3.1.3.2; here $G_{\dbZ_{(p)}}$ is the reductive group over $\dbZ_{(p)}$ having $G$ as its generic fibre, and such that $G_{\dbZ_{(p)}}(\dbZ_p)=H$, cf. 3.1.3). If $G$ is adjoint and all simple factors of $(G,X)$ are such that [De2, 1.2.8 (ii)] applies, then we have $\Aut((G,X,H))=Aut(G_{\dbZ_{(p)}})(\dbZ_{(p)})$.

10) For any quadruple $(G,X,H,v)$ and for any isogeny (of connected groups) $G_1\to G^{\der}$, there is a cover $(G_0,X_0,H_0,v_0)\to (G,X,H,v)$ with $G_0^{\der}=G_1$ (and if needed also with $E(G^0,X^0)=E(G,X)$). This is a direct consequence of the proof of [MS, 3.4] (i.e. we can take $G_0$ unramified over $\dbQ_p$, if $G$ is unramified over $\dbQ_p$).

11) For any quadruple $(G,X,H,v)$ there are finite maps $f\colon (G_1,X_1,H_1,v_1)\to (G,X,H,v)$ and $f_1\colon (G_1,X_1,H_1,v_1)\to (G_2,X_2,H_2,v_2)$ such that: 

\smallskip
-- $(G_2,X_2,H_2,v_2)$ is a product of quadruples $(G_i,X_i,H_i,v_i)$, $i$ running through the elements of a finite set, with $G_i^{\ad}$ a simple adjoint $\dbQ$--group;

-- they define a quasi fibre product of the natural maps $f_0\colon  (G,X,H,v)\to (G^{\ad},X^{\ad},H^{\ad},v^{\ad})$ and $f_2\colon (G_2,X_2,H_2,v_2)\to (G_2^{\ad},X_2^{\ad},H_2^{\ad},v_2^{\ad})=(G^{\ad}, X^{\ad},H^{\ad},v^{\ad})$;

-- there are injective maps $(G_i,X_i,H_i,v_i)\hookrightarrow (G,X,H,v)$, $i\in I$, producing an isogeny $\prod_{i\in I} G_i^{\der}\to G^{\der}$.

\smallskip
This results from  2.12 1) using an argument similar to the one used in 3.1.4.  

12) The advantage of working with triples instead of quadruples consists in the fact that if $(G,X,H)\to (G_1,X_1,H_1)$ is a finite map between two triples having integral canonical models, with $H\subset G(\dbQ_p)$ for a prime $p>2$, then the natural morphism (cf. 4)) $\Sh_p(G,X,H)\to\Sh_p(G_1,X_1,H_1)$ is (at least) in the majority of cases the composite of a pro-\'etale cover with an open closed embedding (cf. 6.4.5). But the natural morphism $\Sh_v(G,X,H)\to \Sh_{v_1}(G_1,X_1,H_1)$, with $v$ a prime of $E(G,X)$ dividing $p$ and the prime $v_1$ of $E(G_1,X_1)$, is not so if there are other primes (besides $v$) of $E(G,X)$ dividing $v_1$. This together with C) of 3.2.2 4) makes such triples more suitable for passing the EP enjoyed by an integral canonical model of a triple to a smooth integral model of another triple having the same adjoint triple (for instance cf. 6.2.3). 

\Proclaim{3.2.8. Example.} \rm
We consider a Shimura pair $(T,\{h\})$ with $T$ a torus. Let $p$ be a rational prime. Then $T$ is unramified over $\dbQ_p$ iff $T$ splits over an unramified extension of $\dbQ_p$. If this is so then $T(\dbQ_p)$ has a unique hyperspecial subgroup $H_T$. For any compact open subgroup $H_T^p$ of $T(\dbA_f^p)$, $\Sh_{H_T\times H_T^p}(T,\{h\})$ is the scheme associated to a finite product of finite field extensions of $E(T,\{h\})$ which are unramified over $p$ (this results from the reciprocity map 2.6 and from the fact that $T(\dbQ)H_T=T(\dbQ_p)$  [Mi4, 4.11]). So, for every prime $v_T$ of $E(T,\{h\})$ dividing $p$, $(T,\{h\},H_T,v_T)$ has an integral canonical model, obtained by taking the normalization of $O_{(v_T)}$ in $\Sh_{H_T}(T,\{h\})$. This integral canonical model is uniquely determined even for $p=2$.

\Proclaim{3.2.9. Example.} \rm
We consider a Siegel modular variety $\Sh(G\Sp(W,\psi),S)$. Let $g\in\dbN$ be defined by $\dim_{\dbQ}(W)=2g$. Then any quadruple of it $(G\Sp(W,\psi),S,K_p,p)$ has an integral canonical model $\scrM$ over $\dbZ_{(p)}$: as a scheme it parameterizes isomorphism classes of principally polarized abelian schemes of dimension $g$ (over $\dbZ_{(p)}$-schemes) having (compatibly) level-$N$ symplectic similitude structure for any $N\in\dbN$ relatively prime to $p$; we have a natural continuous action of $GSp(W,\psi)(\dbA_f^p)$ on this scheme. 

This can be seen as follows: [De1, 4.21] takes care of the generic fibre of $\scrM$. The results of [Mu] implies the existence and the smoothness of the integral model $\scrM$. The fact that it has the EP is explained in [Mi4, p. 170-171].\finishproclaim
 
The definition of an integral canonical model of a quadruple $(G,X,H,v)$ was inspired by the desire that this example works.

\Proclaim{3.2.10. Definition.} \rm
We call an injective map $(T,\{h\},H_T,v_T)\hookrightarrow (G,X,H,v)$, with $T$ a maximal torus of $G$, a special quadruple of $(G,X,H,v)$.
\finishproclaim

\Proclaim{3.2.11. Lemma.}
Every quadruple has  special quadruples. 
\finishproclaim

\proof
This results easily from an argument similar to the one in 3.1.4. Let $G_{\dbZ_{(p)}}$ be a reductive group having  $G$ as its generic fibre. For any maximal torus ${T_1}_{\dbZ_p}\hookrightarrow G_{\dbZ_p}$, there is a special quadruple $(T,\{h\},H_T,v_T)$ of $(G,X,H,v)$ such that the Zariski closure $T_{\dbZ_p}$ of $T_{\dbQ_p}$ in $G_{\dbZ_p}$ is $G_{\dbZ_p}(\dbZ_p)$-conjugate to ${T_1}_{\dbZ_p}$. 

\smallskip
Similarly, we can impose different conditions on the $G(\dbQ_l)$-conjugacy class of $T_{\dbQ_l}$, for $l$ belonging to a finite set of rational primes different from $p$ (cf. the argument in 3.1.4). We express this property by: every quadruple has plenty of special quadruples.

\Proclaim{3.2.12. The relation between different types of models.} \rm
Let $(G,X,H,v)$ be an arbitrary quadruple. It can have more than one smooth integral model over $O_{(v)}$ (or $\widehat{O_{(v)}}$). Starting with such a smooth integral model, we can cook from it new smooth integral models of it  by using blowings up (dilatations) and by removing a $G(\dbA_f^p)$-invariant closed subscheme of its special fibre, which is not the whole special fibre. If $\dim(X)\Ge 1$ it should be always possible to construct a smooth integral model of our quadruple whose special fibre does have a $G(\dbA_f^p)$-invariant closed subscheme, strictly included in the special fibre of it (cf. [Va2], where this is proved for the case when $(G,X)$ is of preabelian type with $v$ not dividing 2).

\Proclaim{Fact.}
We assume that $(G,X,H,v)$ has an integral canonical model $\scrM$ and that $v$ does not divide 2. If $e<p-1$ then any normal integral model  $\scrM_1$ of it over $O$ having the SEP is isomorphic to $\scrM_O$. 
\finishproclaim

\proof
Let $H_0$ be a compact open subgroup of $G(\dbA_f^p)$ such that for any inclusion $H_2\subset H_1$ of open subgroups of $H_0$, the morphisms $\scrM/H_2\to \scrM/H_1$ and $\scrM_1/H_2\to \scrM_1/H_1$ are \'etale covers. We have a natural $G(\dbA_f^p)$-equivariant morphism $\scrM_O\to \scrM_1$, as $\scrM_1$ has the SEP. It is enough to show that the induced morphism $q\colon \scrM_O/H_0\to\scrM_1/H_0$ is an isomorphism. Due to the EP of $\scrM$, $q$ satisfies the valuative criterion of properness with respect to discrete valuation rings of mixed characteristic. From this and  Nagata's embedding theorem ([Na], [Vo]) we deduce that $q$ is proper. As $e<p-1$ part ii) of 3.2.3.2 b) applies. So $q$ is an isomorphism. This ends the proof of the Fact.
\finishproclaim

\Proclaim{3.2.12.1. Remark.} \rm
If $e\Ge p-1$ and $\dim(X)>0$ we do not know if (or when) $\scrM_O$ has the SEP. 
\finishproclaim

\Proclaim{3.2.13. Fact.}
Let $(G,X)$ be an arbitrary Shimura pair and let $v$ be an arbitrary prime of $E(G,X)$ dividing $p$. Any integral model $\scrM_1$ of $\Sh_{\tilde H}(G,X)$ over $W(\overline{k(v)})$ (with $\tilde H$ a compact open subgroup of $G(\dbQ_p)$) which as a scheme is normal and has a quotient $\scrM_1/{\tilde H}_0$ (with ${\tilde H}_0$ a compact open subgroup of $G(\dbA_f^p)$) of finite type over $W(\overline{k(v)})$, descends to an integral model over an \'etale DVR extension $O_{(v^{\prime})}$ of $O_{(v)}$.  
\finishproclaim

\proof
Claim 3.1.3.1 allows us to descend $\scrM/\tilde H_0$ to a scheme $\scrM^{\text{sh}}/\tilde H_0$ of finite type over $O_{(v)}^{\text{sh}}$. So $\scrM^{\text{sh}}/\tilde H_0$ descends to a scheme $\scrM^{v^{\prime}}/\tilde H_0$ over an \'etale DVR extension $O_{(v^{\prime})}$ of $O_{(v)}$. Now the normalization of  $\scrM^{v^{\prime}}/\tilde H_0$ in the ring of fractions of the extension of $\Sh_{\tilde H}(G,X)$ to the field of fractions $L^{\prime}$ of $O_{(v^{\prime})}$ (there is a natural $G(\dbA_f^p)$-continuous action on this normalization) is an integral model of $\Sh_{\tilde H}(G,X)$ over $O_{(v^{\prime})}$. Obviously its extension to $W(\overline{k(v)})$ is $\scrM_1$. This ends the proof of the Fact.

\Proclaim{3.2.13.1.} \rm
There are variants of descent when we work with an arbitrary DVR $O$ faithfully flat over $O_{(v)}$, instead of $O_{(v)}$. The expectation of 3.2.3.2 iii), if true, implies  that in many cases we can assume that $k(v^\prime)=k(v)$. But we do not know (cf. 3.1.3.1) when we can take $O_{(v^\prime)}=O_{(v)}$. This motivates why we also introduced the notion of local integral canonical models: if a quadruple $(G,X,H,v)$ has an integral canonical model then it has a local integral model, but we do not know (even if $(v,2)=1$) if the converse is true.  
\finishproclaim

\Proclaim{3.2.14. Remark.} \rm
Let $f\colon\Sh(G,X)\hookrightarrow\Sh(G_1,X_1)$ be an injective map and let $p$ be a rational prime such that $G$ and $G_1$ are unramified over $\dbQ_p$. We assume the existence of a hyperspecial subgroup $H$ of $G(\dbQ_p)$ included in a hyperspecial subgroup $H_1$ of $G_1(\dbQ_p)$. Then for any compact open subgroup $H^p$ of $G(\dbA_f^p)$, the natural morphism 
$$
Sh_{H^p\times H}(G,X)\to \Sh_{H^p\times H_1}(G_1,X_1)\times_{E(G_1,X_1)} E(G,X)
$$
is a closed embedding. 

The proof of this is entirely similar to the proof of [De1, 1.15] (being just the $\dbZ_{(p)}$-version of it), starting from 3.3.1. In particular $\Sh_H(G,X)$ is a closed subscheme of $\Sh_{H_1}(G_1,X_1)\times_{E(G_1,X_1)} E(G,X)$.

\Proclaim{3.2.15. Remark.} \rm
Let $f\colon (G,X,H,v)\hookrightarrow (G_1,X_1,H_1,v_1)$ be an injective map between two quadruples having integral canonical models $\scrM$ and respectively $\scrM_1$. We assume that $v$ does not divide 2. Then $\scrM$ is the normalization of the Zariski closure of $\Sh_H(G,X)$ in $\scrM_{1O_{(v)}}$ (due to 3.2.14 this makes sense). 

This results by putting together 3.2.12 and 3.4.1. If we also have $G^{\der}=G_1^{\der}$, then $\scrM$ is an open closed subscheme of $\scrM_1$ and for every compact open subgroup $H_0$ of $G(\dbA_f^p)$, $\scrM/H_0$ is an open closed subscheme of $\scrM_1/H_0$ (we have $E(G,X)=E(G_1,X_1)$, cf. [De1, 3.8], and so $O_{(v)}=O_{(v_1)}$). In this case we do not need to refer to 3.2.12 or 3.4.1: 3.2.14 is sufficient.
\finishproclaim

\Proclaim{3.2.16. Remark.} \rm
Let $(G,X)=(G_1\times G_2,X_1\times X_2)$ define a Shimura variety which is a product of two Shimura varieties defined by $(G_i,X_i)$, $i=\overline{1,2}$. Let $p$ be a prime such that $G$ is unramified over $\dbQ_p$ and let $H=H_1\times H_2$ (cf. 3.1.5) be a hyperspecial subgroup of $G(\dbQ_p)$. Here $H_i\subset G_i(\dbQ_p)$, $i=\overline{1,2}$. Let $v$ be a prime of $E(G,X)$ dividing $p$ and let $v_i$ be the prime of $E(G_i,X_i)$ divided by $v$. If $(G_i,X_i,H_i,v_i)$ has an integral canonical model $\scrM_i$,
$i=\overline{1,2}$, then $(G,X,H,v)$ has an integral canonical model $\scrM$ defined by the product over $O_{(v)}$  of the extensions to $O_{(v)}$ of the two integral canonical models $\scrM_1$ and $\scrM_2$.

\Proclaim{3.2.17. The proof of 3.2.2 1) and 3).} \rm
Let $D$ be a Dedekind ring flat over $\dbZ\fracwithdelims[]12$. Let $(Y,U)$ be an extensible pair, with $Y$ a very healthy regular scheme over $D$. Let $A_U$ be an abelian scheme over $U$. We have to prove that $A_U$ extends to an abelian scheme over $Y$. For this we can assume that $D$ is a DVR faithfully flat over $\dbZ_{(p)}$ (for some prime $p\Ge 3$), that  $Y=\Spec(R)$ is a local regular scheme of dimension $d+1$ (with $d\in\dbN$), that $U=\Spec(R)\setminus \Spec(R/I)$ with $I$ an ideal of $R$ of height at least 2, and that the residue field of $R$ is an algebraic extension of the residue field  of $D$. 

\Proclaim{Step A.} \rm
It is enough to show that $B_U:=(A_U\times A_U^t)^{4}$
extends to an abelian scheme over $Y$ (we can apply [FC, 2.7] to the projectors of $B_U$ onto its factors). The abelian scheme $A_U$ is a projective scheme over $U$ (cf. [FC, 1.10 a)]) and so it is polarizable. The Zarhin's trick [Za] implies that $B_U$ has a principal polarization $p_U$. Let $N\Ge 4$ be an integer relatively prime to $p$. Let $U_0:=B_U[N]$. It is an \'etale cover of $U$. Let $Y_1$ be the normalization of $Y$ in the ring of fractions of $U_1$. From the classical purity theorem we get that $Y_1$ is an \'etale cover of $Y$. Using descent (based on [FC, 2.7]), it is enough to show that $B_{U_1}:=B_U\times_U U_1$ extends to an abelian scheme over $Y_1$. So we can assume that $U_1=U$; so  the principally polarized abelian scheme $(B_U,p_U)$ has a level-N structure. Let $\scrA_{d(B_U),1,N}$ be the moduli scheme over $\dbZ_{(p)}$ parameterizing principally polarized abelian schemes (over $\dbZ_{(p)}$-schemes) (of dimension $d(B_U)$ equal to the relative dimension of $B_U$) endowed with a level-N structure. We get a morphism $q_U\colon U\to\scrA_{d(B_U),1,N}$ corresponding to $(B_U,p_U)$ and its level-N structure. We need to show that $q_U$ extends to a morphism $q_Y\colon Y\to\scrA_{d(B_U),1,N}$. Let $\dbD_w$ have the same meaning as in 3.2.1 8). 

We can replace $R$ by $R_1:=R\otimes_D \dbD_w$ and then we can replace $R_1$ by the completion $R_0$ of a localization of $R_1$ in a point of it  having $\overline{k(w)}$ as its residue field. This admits an argument at the level of extensions of morphisms: to show that $q_U$ extends, it is enough to show that for any $R_0$ as above, the morphism $q_{U_0}\colon U_0\to\scrA_{d(B_U),1,N}$, with $Y_0:=\Spec(R_0)$ and $U_0:=Y_0\setminus \Spec(R_0/IR_0)$, extends to a morphism $q_{Y_0}\colon Y_0\to\scrA_{d(B_U),1,N}$. From the very definition of a very healthy regular scheme, we get that $R_0=V[[x_1,...,x_d]]$, with $V$ a finite flat DVR extension  of $W(\overline{k(w)})$ of degree $e< p-1$. We get an abelian scheme $B_{U_0}$ over $U_0$.
\finishproclaim

In the case of an abstract very healthy regular scheme, the same argument at the level of extensions of morphisms, allows us to reduce the proof of 3.2.2 3) involving healthy schemes to the case of an abelian scheme $B_{U_0}$ over a scheme $U_0$ as above.

Now we forget how $B_{U_0}$ has been obtained and we just use the fact that it is an abelian scheme over $U_0$. The fact that it has a polarization implies that below we can deal with abelian schemes over complete, local, affine schemes  and not over affine, formal schemes. From now on we  follow [Fa4]. Let $K:=V\fracwithdelims[]1p$.

\Proclaim{Step B.} \rm
We assume first that $d=1$. Let $n,m\in\dbN$. Then $B_{U_0}[p^n]$ extends to a finite flat group scheme $G_n=\Spec(O_n)$ (with $O_n$ the ring of global sections of the ring sheaf of the ringed space $B_{U_0}[p^n]$) over $Y_0$ (cf. 3.2.1.1 9)). 

The natural homomorphisms $G_n\to G_{n+m}$ are closed embeddings. To see this let ${G_n}_K$ be the generic fibre of the restriction ${G_n}_V$ of $G_n$ to $R_0/x_1R_0=V$. The finite group ${G_n}_K$ extends uniquely to a finite flat group scheme ${G_n}_V$ over $V$, and so ${G_n}_V$ is the Zariski closure of ${G_n}_K$ in ${G_{n+m}}_V$, cf. [Ra, 3.3.6]; hence the corresponding ring homomorphisms $O_{n+m}\to O_n$ become surjective by tensoring with $V$, and thereby, cf. Nakayama's Lemma, they are epimorphisms. 

Due to the uniqueness of an extension of a flat finite group scheme over $U_0$ (to a flat finite group scheme over $Y_0$) (cf. 3.2.1.1 9)) we get that $G_{n+m}/G_n\arrowsim G_m$. So the $p$-divisible group of $B_{U_0}$ extends to a $p$-divisible group $G_{Y_0}$ over $Y_0$.${}^1$ $\vfootnote{1}{The rest of this Step B does not follow [Fa4].}$ 

Let $\bar\scrA_{d(B_U),1,N}$ be a projective, toroidal compactification of $\scrA_{d(B_U),1,N}$, such that the complement of $\scrA_{d(B_U),1,N}$ in $\bar\scrA_{d(B_U),1,N}$ has pure codimension 1 and there is a semiabelian scheme over it extending the universal abelian scheme over $\scrA_{d(B_U),1,N}$ (cf. [FC]). Let $\tilde Y_0$ be the normalization of the Zariski closure of $U_0$ in $Y_0\times_V \bar\scrA_{{d(B_U),1,N}V}$. It is a projective, normal, integral $Y_0$-scheme having $U_0$ as an open subscheme. Let $\bar B_{\tilde Y_0}$ be the semiabelian scheme over $\tilde Y_0$ extending $B_{U_0}$. As $Y_0$ is strictly henselian and due to the classical purity theorem, $B_{U_0}$ has level-$\bar N$ structure for any $\bar N\in\dbN$ prime to $p$. So from N\'eron--Ogg--Shafarevich criterion we get that $\bar B_{\tilde Y_0}$ is an abelian scheme in codimension at most 1. As the complement of $\scrA_{d(B_U),1,N}$ in $\bar\scrA_{d(B_U),1,N}$ has pure codimension 1, we get that $\bar B_{\tilde Y_0}$ is an abelian scheme. A theorem of Tate implies that the $p$-divisible group of $\bar B_{\tilde Y_0}$ coincides with the pull back $G_{\tilde Y_0}$ of $G_{Y_0}$ to $\tilde Y_0$ in codimension at most 1. So as $\tilde Y_0$ is normal, the $p$-divisible group of $\tilde B_{Y_0}$ is $G_{\tilde Y_0}$. The complement $C_0$ of $U_0$ in $\tilde Y_0$ when endowed with the reduced structure, is a connected, projective scheme over the field of fractions of $Y_0$ (cf. [Hart, 11.3 of p. 279]). From the last two sentences we get that the morphism $q_{\tilde Y_0}:\tilde Y_0\to\scrA_{d(B_U),1,N}$ is constant on $C_0$. We easily get that $q_{\tilde Y_0}$ factors through a morphism $q_{Y_0}:Y_0\to\scrA_{d(B_U),1,N}$. So $B_{U_0}$ extends to an abelian scheme $B_{Y_0}$ over $Y_0$. 
\finishproclaim

\Proclaim{Step C.} \rm
We now treat the general case by induction on $d\in\dbN$. Let now $d\Ge 2$. First we apply the inductive assumption to $R_y:=R_0\fracwithdelims []1y$ (with $y$ an arbitrary regular parameter of $R_0$): $R_y$ is a regular scheme of dimension $d$ (the local rings of the maximal points of $R_y$ are very healthy regular schemes over different DVR's, so the inductive assumption can be applied). So we can assume that $U_0=\Spec(R_0)\setminus  \Spec(R_0/m_0)$ with $m_0$ the maximal ideal of $R_0$. 
\finishproclaim

\Proclaim{Step D.} \rm
Let $x:=x_1$ and $U_x:=\Spec(R_0/xR_0)\setminus \Spec(R_0/m_0)$. By induction $B_{U_0}\times U_x$ extends to an abelian scheme $B_1$ over $\Spec(R_0/xR_0)$. Let $\scrT_{B_1}$ (resp. $\scrT_{B_1^t}$) be the tangent space of $B_1$ (resp. of $B_1^t$). Both are free module over $R_x:=R_0/xR_0$ of dimension $d(B_U)$. The liftings of an abelian scheme over $R_0/x^nR_0$ which is a lift of $B_1$, to an abelian scheme over $R_0/x^{n+1}R_0$, are parameterized by sections of a principal homogeneous space of $\scrT_{B_1}\otimes\scrT_{B_1^t}$. But this free $R_x$-module
has the same sections over $\Spec(R_x)$ as over $U_x$. So there is a unique way of lifting (compatibly) $B_1$ to an abelian scheme $B_{Y_0}$ over $Y_0$ which over $U_0$ is $B_{U_0}$. This completes the induction, and ends the proof of the part of 3.2.2 1) and 3) involving healthy regular schemes.

\smallskip
The above Steps B to D can be easily adapted to get the part of 3.2.2 3) pertaining to $p$-healthy regular schemes. This ends the proof of 3.2.2 1) and 3).
\finishproclaim

\smallskip
\Proclaim{3.3. The complex points of an integral canonical model.} \rm
Let $p$ be a rational prime and let $(G,X,H)$ be an arbitrary triple, with $H$ a hyperspecial subgroup of $G(\dbQ_p)$.

\Proclaim{3.3.1.} \rm
We have 
$$\Sh_H(G,X)(\dbC)=G(\dbZ_{(p)})\setminus {(X\times G(\dbA_f^p))}/{Z(G)^p},$$ 
where $G(\dbZ_{(p)}):=G(\dbQ)\cap H$ and $Z(G)^p$ is the topological closure of $Z(G)(\dbQ)\cap H$ in $G(\dbA_f^p)$ [Mi3, 4.11].

\Proclaim{3.3.2. Lemma.} 
$G(\dbA_f^p)$ permutes transitively the connected components of $\Sh_{H}(G,X)_{\dbC}$.
\finishproclaim

\proof
If $G^{\der}$ is simply connected, this results from 3.3.1 and from [De1, 2.5] (by passage to limit). For an arbitrary $G$, we have to use the well known trick [MS, 3.4] (cf. 3.2.7 9)) for reducing the problem to the case when $G^{\der}$ is simply connected (as described in [Mi4, 4.19]). This ends the proof.

\Proclaim{3.3.3. Corollary.} 
$G(\dbZ_{(p)})$ permutes transitively the connected components of $X$.
\finishproclaim

\smallskip
\Proclaim{3.4. Methods of constructing integral models.} \rm
Let $\Sh(G,X)$ be an arbitrary Shimura variety.
In essence there are four methods of constructing good integral models of quotients of $\Sh(G,X)$:

\smallskip
\item{1)} By proving first that a suitable quotient of $\Sh(G,X)$ is the moduli scheme parameterizing some objects which make sense over $O_{(v)}$-schemes (with $v$ a prime of $E(G,X)$), and that in fact we have a moduli scheme over $O_{(v)}$. Such a moduli scheme over $O_{(v)}$, in a suitable context, is (expected to be) an integral canonical model of $\Sh(G,X)$ (cf. 3.2.8 and 3.2.9).

\item{2)} By taking the normalization of the Zariski closure of a quotient of $\Sh(G,X)$ into a good integral model of a quotient of another Shimura variety $\Sh(G_1,X_1)$  (here we need an injective map $(G,X)\hookrightarrow (G_1,X_1)$) (cf. what follows below).

\item{3)} By taking the normalization of a good integral model of a quotient of $\Sh(G,X)$ into the ring of fractions of a quotient of another Shimura variety $\Sh(G_1,X_1)$ (here we need a finite map $(G_1,X_1)\to (G,X)$) (cf. 6.1.2).

\item{4)} By taking the quotient through a (torsion) group action on a connected component of a good integral model of a quotient of $\Sh(G,X)$ (here the group action is related to a finite map $(G,X)\to (G_1,X_1)$) (cf. 6.2.2).

\smallskip
These methods are supported by well known ideas pertaining to Shimura varieties (like 3.2.14 and 3.2.7 9)).
Variants for 1) are obtained by working over a DVR faithfully flat over $O_{(v)}$ (instead of $O_{(v)}$).
The method 2) is used as well for constructing integral canonical models of a Shimura variety $\Sh(G,X)$ of abelian type for which there is a Shimura variety $\Sh(G_1,X_1)$ of Hodge type with $G^{\der}=G_1^{\der}$ and $(G^{\ad},X^{\ad})=(G_1^{\ad},X^{\ad}_1)$ (cf. 3.2.15, 5.1 and 6.2.3). The method 4) is used for passing from the existence of integral canonical models of these Shimura varieties to the existence of integral canonical models of all Shimura varieties of abelian type (cf. [Mi4, 4.11 and 4.13]; see also 3.4.5 and 6.2.2). The method 3) is used for the passage from the abelian type case to the preabelian type case (cf. 6.1). 

\smallskip
We start with an injective map $f\colon (G,X,H,v)\hookrightarrow (G_1,X_1,H_1,v_1)$. We assume that $(G_1,X_1,H_1,v_1)$ has a normal integral model $\scrM_1$ over $O_{(v)}$. Let $\scrM$ be the normalization of the Zariski closure of $\Sh(G,X)/H$ in $\scrM_1$ (cf. 3.2.14). It has an obvious $G(\dbA_f^p)$-continuous action ($p$ being the rational prime divided by $v$). Let $E:=E(G,X)$.  

\Proclaim{3.4.1. Proposition.}
The integral model $\scrM$ is a normal integral model of $(G,X,H,v)$. It has the EP (or EEP, or WEP, or SEP) if $\scrM_1$ has it.
\finishproclaim

\proof
Obviously $\scrM$ has the EP (or EEP, etc.) if $\scrM_1$ has it.
Let $H_0$ be a compact open subgroup of $G(\dbA_f^p)$
 such that:

\smallskip
\item{i)}
the subgroup $H_0\times H$ of $G(\dbA_f)$ is smooth for $(G,X)$;

\smallskip
\item{ii)}
there is a compact open subgroup $K_0$ of $G_1(\dbA_f^p)$
including $H_0$ and such that for any compact open
subgroup $K_1$ of $K_0$, $\scrM_1/K_1$ is a normal scheme of finite type over $O_{(v)}$ and 
\'etale over $\scrM_1/K_0$.

The existence of such a subgroup $H_0$ is implied by the fact that $\scrM_1$ is a normal integral model and by 2.11.

Let $H_1\subset H_2$ be two open subgroups of $H_0$.
Let $\scrP_i$ be the normalization of the Zariski closure of the generic fibre of $\scrM/H_i$
in $\scrM_1
/H_i$, for $i=\overline{1,2}$
( $\scrM_E/H_i$ is a closed subscheme of
$\scrM_{1E}/H_i$, cf. 3.2.14).
We get the following diagram:
$$
\CD
\scrM @= \scrM @>>> \scrM_1\\
@VVV @VVV @VVV\\
\scrM/H_1 @>{g_1}>> \scrP_1 @>>> \scrM_1/H_1\\
@VVV @VVV @VVV\\
\scrM/H_2 @>{g_2}>> \scrP_2 @>>> \scrM_1/H_2.
\endCD
$$
The conditions \ i) and \ ii) and the fact that $\scrM$ is
$H_i$-invariant imply that the two right squares are
Cartesian. So $\scrM$ is a pro-\'etale cover of $\scrP_1$
and $\scrP_2$.
The generic fibre of $\scrM/H_i$ is a scheme of finite type
over $E$. The scheme $\scrM_1/H_i$ is a projective limit of schemes of the form $\scrM_1/T$ with $T$ an open subgroup of $K_0$
including $H_i$.
So there is an open subgroup $K_i$ of $K_0$,
with $H_i\subset K_i$, such that the morphism
$\scrM_E/H_i\to\scrM_{1E}/K_i$ is a 
closed embedding.
As the morphism $\scrM/H_i\to\scrM_1/K_i$ is integral, we deduce
that $\scrP_i$ is integral over the Zariski closure $\scrS_i$ of $\scrP_{iE}$ in $\scrM_{1E}/K_i$ and has the same generic
fibre as $\scrS_i$. 
As $\scrS_i$ is an excellent scheme (it is of
finite type over $O_{(v)}$), we get that $\scrP_i$ is finite over $\scrS_i$,
and so of finite type over $O_{(v)}$.
Both $\scrP_1$ and
$\scrP_2$ are faithfully flat over $O_{(v)}$.
Moreover $g_1$ and $g_2$ are integral morphisms between flat schemes over $O_{(v)}$ having the same
generic fibre.
The normality of $\scrP_1$ and $\scrP_2$ implies that $g_1$ and $g_2$ are isomorphisms; so $\scrM/H_1\to\scrM/H_2$ is  an \'etale morphism between schemes of finite type
over $O_{(v)}$ (as the morphism $\scrP_1\to\scrP_2$ is so).
We conclude that $\scrM$ is a normal integral model. This ends the proof of the Proposition.

\Proclaim{3.4.1.1. Remark.} \rm
The above Proposition as well as 3.4.2 and 3.4.3 below remain true if $H$ and $H_1$ are just compact open subgroups of $G(\dbQ_p)$ and respectively of $G_1(\dbQ_p)$ satisfying $f(H)=f(G(\dbQ_p))\cap H_1$, or if $O_{(v)}$ is replaced by an arbitrary DVR $O$ faithfully flat over $O_{(v)}$.
\finishproclaim

\Proclaim{3.4.2. Remark.} \rm
The above proof shows that $\scrM$ is a
pro-\'etale cover of a normal scheme $\scrP$ of finite type over $O_{(v)}$.
As $O_{(v)}$ is a universally catenary ring, all the maximal points of $\scrM$ have dimension $d+1$, where $d=\dim\,X$ (as the dimension
formula holds between $O_{(v)}$ and any connected component
of $\scrP$ [Ma, p. 85]).
\finishproclaim

\Proclaim{3.4.3. Remark.} \rm
For any compact open
subgroup $H_0$ of $G(\dbA_f^p)$ small enough, $\scrM/H_0$ is the normalization of 
a closed subscheme of $\scrM_1/H_0$. If $\scrM$ is a subscheme of $\scrM_1$, then we do not need to take any normalization.
\finishproclaim

\Proclaim{3.4.4. Corollary.}
We assume that $\scrM_1$ has the EP. Then $\scrM$ is an integral canonical model iff $\scrM$ (as a scheme) is formally smooth over $O_{(v)}$.
\finishproclaim

\Proclaim{3.4.5. Expectation.} \rm
Let $\scrM$ be a smooth integral model of a quadruple $(G,X,H,v)$ over a DVR $O$. Let $p$ be the rational prime divided by $v$. Let  $H_0$ be a  subgroup of $G(\dbA_f^p)$ such that the subgroup $H_0\times H$ of $G(\dbA_f)$ is smooth for $(G,X)$. We do expect that under some mild conditions (like $H_0\times H$ is $p$-smooth for $(G,X)$ and the index of ramification $e$ of $O$ is less than $p-1$) $\scrM$ is a pro-\'etale cover of $\scrM/H_0$.

This expectation is based on two facts. First we can prove it (under the restriction $e<p-1$ and under some mild assumptions on $H_0$) for the case of a quadruple of preabelian type (for $p\Ge 5$ cf. 6.4.2.1; for $p=3$ cf. [Va2]). Second we have the following considerations.
 
Let $\Htil_0$ be an open subgroup of $H_0$ such that $\scrM$ is a pro-\'etale cover of $\scrM/\Htil_0$ (cf. the definition of a smooth integral model). We can assume that $\Htil_0$ is a normal subgroup of $H_0$. Let $C_0$ be the quotient of $H_0/\Htil_0$ by its subgroup acting trivially on $\scrM/\Htil_0$. It is a finite group. Then $\scrM/H_0$ is the quotient of $\scrM/\Htil_0$ by $C_0$ (cf. the definition of a continuous action). The action of $C_0$ on the generic fibre of $\scrM/\Htil_0$ is free (as $\Sh(G,X)/H\times\Htil_0$ is an \'etale cover of $\Sh(G,X)/H\times H_0$). But then it is expected (cf. 3.4.5.1 below) that the action of $C_0$ on $\scrM/\Htil_0$ is free. If this is so then $\scrM/\Htil_0$ is an \'etale cover of $\scrM/H_0$ (and so $\scrM$ is a pro-\'etale cover of $\scrM/H_0$). 
\finishproclaim

\Proclaim{3.4.5.1. Proposition.}
Let $p$ be a rational prime. Let $V$ be a complete DVR which is a faithfully flat $\dbZ_{(p)}$-algebra, and has an index of ramification $e<p-1$. Let $C$ be a finite (abstract) group acting on a regular formally smooth $V$-algebra $R$ in such a way that it acts freely on $R\fracwithdelims[]1p$. Let $V_1$ be the DVR obtained by adjoining to $V$ a primitive $p$-th root of unity. We assume that either the order of $C$ is relatively prime to $p$, or it is $p$ and the subring $R^C$ of $R$ formed by elements fixed by $C$ is such that the affine scheme $\Spec(R^C\otimes_V V_1)$ is locally factorial. Then $C$ acts freely on $R$.
\finishproclaim

\proof
We assume that we do have a situation with a non-free action. We can assume that $C$ is a finite cyclic group of prime order $l$. Let $\pi_V$ be a uniformizer of $V$ and let $k_V$ be its residue field. We can also assume that $R$ is a local ring.

If $l$ is different from $p$ this is well known. We can assume further on that $V$ is a complete DVR of index of ramification $e<p-1$, that $k_V$ is an algebraically closed field, and that $R=V[[x_1,...x_d]]$ is the ring of formal power series in $d$ variables with coefficients in $V$. We can write $R=\oplus_{\gamma\in \widehat{C}} R^{\gamma}$, with $\widehat{C}$ the dual group of $C$ (i.e. the group of characters of $C$), and with $C$ acting on $R^{\gamma}$ through the character $\gamma\in\widehat{C}$. Now it is trivial to see that if for a non-trivial character $\gamma$ of $C$, $R^{\gamma}$ is different from zero, then the action of $C$ on $R\fracwithdelims[]1p$ is not free (i.e. there is an element $y$ of the maximal ideal $m_R$ of $R$, whose image in $m_R/m_R^2$ is non-zero and is different from the image of $\pi_V$ in $m_R/m_R^2$, and which belongs to an  $R^{\gamma}$, for a non-trivial character $\gamma$; this disturbs the free action of $C$ on $R\fracwithdelims[]1p$ as we can see by induction on $d$). Contradiction. For this part we do not need that $e<p-1$.

Let now $l=p$. We abbreviate the notion of unique factorization domain by UFD. From the theory of tamely totally ramified extensions of $W(\overline{k_V})$, and from the fact that $e$ is smaller than $p-1$, we deduce that the index of ramification of $V_1$ is $e_1$, a multiple of $p-1$ relatively prime to $p$. In fact $e_1=l.c.m.(p-1,e)$. 

Let us recall a well known fact:

\Proclaim{3.4.5.2.}
 Let $M$ be a torsion free $V$-module separated with respect to the $\pi_V$-topology, and let $1_M$ be its identity automorphism. Then any $V$-automorphism $a_M$ of $M$ such that $a_M^p=1_M$ and $a_M$ modulo $\pi_V$ is the identity, is the identity automorphism.
\finishproclaim

\proof  Writing $a_M=1_M+\pi_Vb_M$ with $b_M\in\text{End}(M)$, by induction on $n\in\dbN$, we can check that $b_M$ is of the form $\pi_V^{n-1}c_M$ with $c_M\in\text{End}(M)$. As $M$ is separated with respect to the $\pi_V$-topology, we deduce that $\text{End}(M)$ is separated with respect to this topology. So $b_M=0$. This is the only place where we need that $e<p-1$ (3.4.5.2 is not true if $e\Ge p-1$). This proves 3.4.5.2.

\smallskip
So the case $l=p$ results once we show that the action of $C$ on $R/\pi_VR$ is trivial. 
As $V$ is complete we deduce that the completion of $R$ is of the form $V^{\prime}[[x_1,...,x_d]]$, with $V^{\prime}$ a finite \'etale DVR extension of $V$. 
The $V$-algebra $V^\prime$ is a subring of $R$ ($R$ is normal). The group $C$ acts on $V_1$ trivially (we assumed that the action is non-free). So $V^\prime\subset R^C$. This allows us to replace $V$ by $V^\prime$ and $V_1$ by $V_1^\prime$, where $\Spec(V_1^\prime)$ is a connected open-closed subscheme of $\Spec(V_1\otimes_V V^\prime)$. Not to complicate the notations, we assume that $V=V^\prime$. So $R_1:=R\otimes_V V_1$ is an integral domain.   

Let $\pi_1$ be a uniformizer of $V_1$. Let $\scrO_1$ be the local ring of the generic point of the special fibre of $\Spec(R_1)$. The group $C$ acts on it. Let $R_1^C$ and $\scrO_2:=\scrO^C$ be the subrings of $R_1$ and respectively of $\scrO_1$ formed by elements fixed by $C$. We have $R_1^C=R^C\otimes_V V_1$.
The sets $\pi_1R_1$ and $\pi_1R^C_1$ are prime ideals of $R_1$ and respectively of $R_1^C$. 

Both $\scrO_1$ and $\scrO_2$ are discrete valuation rings  having the same index of ramification equal to $e_1$ (both being $V_1$-algebras). Let $\scrK_i$ be the  field of fractions of $\scrO_i$, $i=\overline{1,2}$.

As $V_1$ contains the $p$-th roots of unity, and as the action of $C$ on $R$ is non-trivial, there is $y\in\scrK_1$ such that $C$ acts on it through a non-trivial character $\gamma$ of $C$. So $y^p\in\scrK_2$, but $y\notin\scrK_2$. By reasons of dimension, we deduce that $\scrK_1$ is a Kummer extension of $\scrK_2$. We get the situation: 

\smallskip\noindent
{\bf (3.4.5.3)} $\scrK_1$ is a Galois extension of $\scrK_2$ of degree $p$, obtained by adjoining a $p$-th root of an element of $\scrK_2$. 

\smallskip
In all that follows $y$ denotes an element of $\scrK_1\setminus\scrK_2$ such that $y^p\in\scrK_2$. We repeatedly replace it by $y_1=yy^C$, with $y^C\in\scrK_2$. 

If the action of $C$ on the residue field of $\scrO_1$ is non-trivial (i.e. if the action of $C$ on $R/\pi_VR$ is non-trivial), then we deduce easily that the residue field $k_1$ of $\scrO_1$ is a Galois extension of the residue field $k_2$ of $\scrO_2$ ($k_2\subset k_1^C$, where $k_1^C$ is the subfield of $k_1$ formed by elements fixed by $C$; but $k_1$ is a Galois extension of $k_1^C$ of degree $p$, and so by reasons of dimension we must have $k_1^C=k_2$). We deduce that $\Spec(\scrO_1)$ is a Galois cover of $\Spec(\scrO_2)$.

The morphism $\Spec(R_1)\to\Spec(R_1^C)$ is \'etale above points of $\Spec(R_1^C)$ of characteristic zero or of codimension 1. So $\Spec(R_1^C)$ is regular in all these points (and so is regular in codimension at most 1).

\Proclaim{Step a).} \rm
From the fact that $R_1^C$ is a local UFD, we deduce that the Picard group of $\Spec(R_1^C)$ is trivial and isomorphic to its divisor class group. This implies that we can assume that $y$ is an invertible element of $R_1$. In other words we can replace $y$ by $y_1:=yy^C$, with $y^C\in\scrK_2$, such that, in any point of $\Spec(R_1)$ of codimension 1, $y_1$ is invertible; so $y_1^p$ is an invertible element of $R_1^C$ (this can be deduced  from  [Ma, Th. 38], as $R_1^C$ is a normal ring). In detail: For any prime $\grp$ of $R_1^C$ of codimension 1, as the ring extension $R_1^C\to R_1$ is \'etale above it, we deduce the existence of an element $y_{\grp}\in\scrK_2$ such that $y^py_{\grp}^p$ is an invertible element of the localization of $R_1^C$ with respect to the prime $\grp$. The elements $y_{\grp}$, with $\grp$ running through all the primes of $R_1^C$ of codimension 1 are defining a Weil divisor. As this Weil divisor is linearly equivalent to the zero divisor, we deduce the existence of an element $y^C\in\scrK_2$ producing this Weil divisor. We can take now $y_1=yy^C$.
  
As a conclusion: the extension $\scrK_1$ of $\scrK_2$ is obtained  by adjoining a $p$-th root (still denoted by  $y$) of an invertible element of $R_1^C$.
\finishproclaim

By dropping the assumption that $R_1^C$ is a UFD, we can assume that $k_V$ is separable closed and that $R_1$ is a complete local ring. So the first fundamental group of $\Spec(R_1)$ is trivial. Moreover $R_1^C$ is a complete local ring (as $R_1$ is so, and as the inclusion $R_1^C\hookrightarrow R_1$ is finite). We deduce that the first fundamental group $\pi_1(R_1^C)$ of $\Spec(R_1^C)$ is trivial $(\pi_1(R_1^C)$ is a subgroup of $C$; but it is not $C$ as the inclusion $R_1^C\hookrightarrow R_1$ of complete local rings having the same residue field, is not \'etale).

We assume now that $e$ divides $p-1$. So $e_1=p-1$.

\Proclaim{Step b).} \rm
 If the image of $y$ in $k_1$ is not in $k_2$, then we deduce by reasons of dimension that $k_1$ is obtained from $k_2$ by adjoining a $p$-th root of an element of $k_2$. We get a contradiction with the fact that $k_1$ is a Galois extension of $k_2$. So the image of $y$ in $k_1$ is in $k_2$. Replacing $y$ with $y_1=yy^C$, with $y^C\in R_1^C$, we can assume that $y$ is congruent to 1 modulo the  ideal of $R_1^C$ generated by $\pi_1$. We can assume that $\pi_1^{p-1}=p$.
\finishproclaim

\Proclaim{Step c).} \rm
Let now $y=1+\pi_1y_0$, with $y_0\in R_1$. So $y^p$ is congruent to $1+p\pi_1(y_0+y_0^p)$ modulo $\pi_1^{p+1}R_1^C$ (or modulo $\pi_1^{p+1}R_1$ as $\pi_1^{p+1}R_1\cap R_1^C=\pi_1^{p+1}R_1^C)$. Let $z_0\in R_1^C$ which modulo $\pi_1R_1^C$ is $y_0+y_0^p$. The equation $x^p+x=z_0$ defines an \'etale $R_1^C$-algebra. As $\pi_1(R_1^C)$ is trivial, we deduce that there is $y^C\in R_1^C$ such that $y^C$ is congruent to $y_0$ modulo $\pi_1R_1$. Replacing $y$ with $y_1=y(1-\pi_1y^C)$, we can assume that $y$ is congruent to 1 modulo $\pi_1^2R_1$ (we have $p>2$ as $e<p-1$).
\finishproclaim

\Proclaim{Step d).} \rm
Now by trivial induction on $n\in\dbN$, we can assume that $y$ is congruent to 1 modulo $\pi_1^{n+1}R_1$ (if $y=1+\pi_1^ny_0$, with $n\in\dbN$ greater than 1 and with $y_0\in R_1$, then $y^p$ is congruent to $1+\pi_1^{n+p-1}y_0$ modulo $\pi_1^{n+p}R_1^C$ or modulo $\pi_1^{n+p}R_1$).
\finishproclaim

\Proclaim{Step e).} \rm
This implies, as $R_1^C$ and $R_1$ are complete with respect to the $\pi_1$-topology, that we can assume that $y^p=1$. As $V_1$ contains the $p$-th roots of unity, this contradicts the fact that $\scrK_1$ is a field. 
\finishproclaim

The case when $e_1$ is not $p-1$ is entirely similar. The only difference is that the above Steps c) and d) have to be applied intermingled. The trivial details are left to the reader.

The contradiction of the Step e) ends the proof of 3.4.5.1.

\Proclaim{3.4.5.4. Remarks.} \rm
1) It is an easy exercise now to see that once we assume in 3.4.5.1 that $\Spec(R^C\otimes_V V_1)$ is a locally factorial scheme, the condition on the order of $C$ (of being $p$) can be weaken: it is enough to assume that $C$ is a $p$-elementary finite group. From the fact that $\Spec(R^C\otimes_V V_1)$ is a locally factorial scheme  we deduce easily that $\Spec(R^C)$ is locally factorial, but we do not know if (or when) the converse to this is true. 

2) Proposition 3.4.5.1 can be formulated for regular formally smooth schemes instead of affine such schemes as the condition of having a free action is local. We have inserted 3.4.5.1 for the case $l=p$ mainly to give an idea how bad the singularities can be for a non-free action (cf. 1)). We hope to use it later on to study the singularities of different quotients of different extended integral canonical models (to be defined in 3.5.1) (cf. 3.5.3). 

3) For the order of $C$ equal to $p$, the Lemma 3.4.5.1 is not true if we do not assume that $\Spec(R^C\otimes_V V_1)$ is a locally factorial scheme, as it can be seen through examples involving smooth schemes $X$ over a DVR $O$ faithfully flat over $\dbZ_{(p)}$ and of index of ramification 1, whose relative dimension is greater than $p-2$. But if the relative dimension of $X$ over such a DVR $O$ is less than $p-1$, then any finite group  acting on it in such a way that it acts freely on its generic fibre, acts freely on $X$. This can be checked starting from 3.4.5.2 and the fact that any representation of a cyclic group of order $p$ of degree less than $p-1$ over such a DVR $O$ is trivial.
\finishproclaim

\Proclaim{3.4.6. Remark.} \rm
We come back to 3.4.1 to 3.4.3. In practice $p$ is different from 2 and then we can take $H_0$ to be a product of its $q$-components $H^q$ ($q$ being an arbitrary prime different from $p$), with $H^2$ a compact open subgroup of $G(\dbQ_2)$ small enough, and with any other component $H^q$ of it a maximal compact subgroup of $G(\dbQ_q)$ (which can be chosen to be a hyperspecial subgroup of $G(\dbQ_q)$ if $G$ is unramified over $\dbQ_q$).  
\finishproclaim

\Proclaim{3.4.7. Corollary.}
We assume that $\scrM_1$ has the EP. If $\Sh(G_1,X_1)$ is a Siegel modular variety and if $p$ is big enough (without an effectively computable lower bound) then $\scrM$ is a closed subscheme of $\scrM_1$.
\finishproclaim

\proof
From 3.2.12 we deduce that $\scrM_1$ is the extension to $O_{(v)}$ of the integral canonical model of $(G_1,X_1,H_1,v_1)$ (see 3.2.9).
Let $\Htil$ be a compact open subgroup of $G(\dbA_f)$ which is a product of its $q$-components $\Htil^q$ (so $\Htil^q$ is a hyperspecial subgroup of $G(\dbQ_q)$, for any big enough prime $q$). We assume that it is smooth for $(G,X)$ and that $\Sh_{\Htil}(G,X)$ is a closed subscheme of the extension to $E$ of $\Sh_{\Ktil}(G_1,X_1)$, with $\Ktil$ a compact open subgroup of $G_1(\dbA_f)$ which is a product of its $q$-components, contains $\tilde H$, and is small enough (cf. 3.2.9 and 4.1) so that $\Sh_{\Ktil}(G_1,X_1)_E$ extends to a  smooth moduli scheme $\scrM_1(\Ktil)$ over $O_E\fracwithdelims[]1{N!}$ (with $N\in\dbN$ big enough and with $O_E$ the ring of integers of $E$).

 Taking $N$ big enough we can assume that the Zariski closure  $\scrM(\Htil)$ of $\Sh_{\Htil}(G,X)$ in $\scrM_1(\Ktil)$ is a smooth scheme over  $O_E\fracwithdelims[]1{N!}$, that $\Htil^q$ is a hyperspecial subgroup of $G(\dbQ_q)$ for any prime $q\Ge N$, and that (cf. the proof of 3.4.1) for any such prime $q$, the normalization of $\scrM(\Htil)_{\dbZ_{(q)}}$ in the ring of fractions of $\Sh_{\Htil^q}(G,X)$ is the integral canonical model of the triple $(G,X,\Htil^q)$. We can take now $p\Ge N$. This ends the proof of the Corollary.

\Proclaim{3.4.8. Definition.} \rm
With the notations of 3.2.3 1) and 2), a smooth (resp. normal) integral model $\scrM$ (of $\Sh(G,X)/H$ over $O$) is said to be strongly smooth (resp. strongly normal) if for any  compact open subgroup $H_0$ of $G(\dbA_f^p)$ such that the subgroup $H_0\times H$ of $G(\dbA_f)$ is smooth for $(G,X)$, $\scrM$ is a pro-\'etale cover of $\scrM/H_0$. Similarly, $\scrM$ is said to be quasi-strongly smooth (resp. quasi-strongly normal) if for any compact, open subgroup $H_0$ of $G(\dbA_f^p)$ such the subgroup $H_0\times H$ of $G(\dbA_f)$ is $p$-smooth for $(G,X)$, $\scrM$ is a pro-\'etale cover of $\scrM/H_0$. 
\finishproclaim

\Proclaim{3.4.8.1. Remark.} \rm
If $\scrM$ is a (quasi-) strongly normal integral model of $\Sh_H(G,X)$ over $O$ having the SEP, then any smooth integral model of $\Sh_H(G,X)$ over $O$ is (quasi-) strongly smooth (cf. rm. $4^{\prime})$ of 3.2.7). In particular, if there is a strongly normal integral model of $\Sh_H(G,X)$ over $O_{(v)}$ having the EP and if $e<p-1$, then any smooth integral model of $\Sh_H(G,X)$ over $O$ is strongly smooth. We do not know if (or when) the condition $e<p-1$ is truly needed.
\finishproclaim

\smallskip
\Proclaim{3.5. Extended integral canonical models.} \rm
Let $(G,X,H,v)$ be an arbitrary quadruple and let $p$ be the rational prime divided by $v$.

\Proclaim{3.5.1. Definition.} \rm
A normal scheme $\scrMtil$ over $O_{(v)}$ together 
with a $G(\dbA_f^p)\times H$-continuous action is called an extended integral canonical model of $(G,X,H,v)$ if:

\smallskip
\item{a)}
There is a $G(\dbA_f^p)\times H$-equivariant isomorphism $\scrMtil_{E(G,X)}\arrowsim \Sh(G,X)$;

\smallskip
\item{b)}
$\scrMtil/H$ is an integral canonical model of $(G,X,H,v)$.

\smallskip
Similarly, we speak about an extended local integral canonical model of $(G,X,H,v)$ or about the extended integral canonical model of $(G,X,H)$.

\Proclaim{3.5.2. Remark.} \rm
The $O_{(v)}$-scheme $\scrMtil$ is determined by the integral canonical model $\scrMtil/H$, being the normalization of $\scrMtil/H$ in the ring of fractions of $\Sh(G,X)$. So $\scrMtil$ exists iff $(G,X,H,v)$ has an integral canonical model. If $v$ is relatively prime to 2, then any extended integral integral model of $(G,X,H,v)$ is uniquely determined up to unique isomorphism, cf. 3.2.7 2).  
\finishproclaim

\Proclaim{3.5.3. Problem.} \rm
For $\tilde H$ a compact open subgroup of $G(\dbA_f^p)\times H$ determine the type of singularities of $\tilde M/\tilde H$.
\finishproclaim

\bigskip
\noindent
{\boldsectionfont \S4. Shimura varieties of Hodge type and special families of tensors}

\bigskip
Let $(G,X)$ be a Shimura pair defining a  Shimura variety of Hodge type.
Let $f\colon(G,X)\hookrightarrow (GSp(W,\psi),S)$ be an embedding of it into a Shimura pair defining a Siegel modular variety. We fix a  family 
$(s_\alpha)_{\alpha\in\scrJ}$ of tensors in spaces of the form
$W^{\otimes m}\otimes W^{*\otimes n}$, $m,n\in\dbN$, such that $G$ is the
subgroup of $GSp(W,\psi)$ fixing its tensors. As $G$ is reductive, the existence of finite such families is implied by [De3, 3.1]. We allow the above family of tensors to be infinite. Let $L$ be a $\dbZ$-lattice of $W$ such that we have a perfect form $\psi\colon L\otimes L\to\dbZ$.

We start by reviewing the interpretation of the complex Shimura variety $\Sh(G,X)_{\dbC}$ as a moduli space with respect to the $\dbZ$-lattice $L$ of $W$  and the above family of tensors.  Then in 4.2 and 4.3 we treat the problem: for a  rational prime $p$ for which $G$ is unramified over $\dbQ_p$, find a $\dbZ$-lattice $L$ and a family of tensors $(s_\alpha)_{\alpha\in\scrJ}$ (subject to the above conditions) which are $\dbZ_{(p)}$-well adapted for using successfully the integral version of Fontaine's comparison theory [Fa3], and so  for proving (cf. \S5) the existence of $\Sh_p(G,X)$.  

\smallskip
\Proclaim{4.1. Shimura varieties of Hodge type as moduli schemes.} \rm
As $G$ contains the group of multiplications by scalars (cf. Definition 1 of 2.5), our tensors are in spaces of the form $(W\otimes
W^*)^{\otimes m}\overset\text{can}\to= W^{\otimes m}\otimes W^{*\otimes m}$, $m\in\dbN$.
If $s_\alpha\in(W\otimes W^*)^{\otimes m(\alpha)}$ then $\deg(s_\alpha)=2m(\alpha)$.
The form $2{\pi}i\psi$ is a bilinear map $W\otimes
W\to\dbQ(1):=2{\pi}i\dbQ$, inducing
an isomorphism $W\simover{\longrightarrow} W^*(1)$. Any $x\in X$ defines a Hodge $\dbQ$--structure on $W$ and on $W^*$, and the above isomorphism $W\simover{\longrightarrow}W^*(1)$ is an isomorphism of Hodge $\dbQ$--structures.
This gives us the right to think of the tensors $s_\alpha$ as
being in spaces of the form $W^{*\otimes 2m}(m)$. 
Let $L^*\subset W^*$ be the dual $\dbZ$-lattice of $L$. What follows is very close to [MS, Ch. 2] except that we do not work in a rational context: we work with principally polarized abelian varieties and not with their isogeny classes.

We consider quadruples of the form
$[A,p_A,(v_\alpha)_{\alpha\in\scrJ},k]$ where:

\smallskip
\item{a)}
$(A,p_A)$ is a principally polarized abelian variety over
$\dbC$;

\smallskip
\item{b)}
$(v_\alpha)_{\alpha\in\scrJ}$ is a family of Hodge cycles of
$A$;

\smallskip
\item{c)}
$k$ is an isomorphism $H_1(A,\dbZ)\otimes\dbA_f
\overset\text{can}\to=
V_f(A)\twiceover{k}{\longrightarrow}
W\otimes\dbA_f$ taking the Betti realization $w_\alpha$ of
$v_\alpha$ into $s_\alpha$,$\forall\,\alpha\in\scrJ$,
mapping $H_1(A,\dbZ)\otimes\dbZhat$ onto $L\otimes\dbZhat$
and inducing a symplectic similitude between
$(H_1(A,\dbZ)\otimes\dbZhat,p_A)$ and
$(L\otimes\dbZhat,\psi)$.

\smallskip
We define $\scrA(G,X,W,\psi)$ to be the set of isomorphism
classes of quadruples of the above form satisfying the
following conditions:

\smallskip
\item{(i)}
there exists a similitude  isomorphism $(H_1(A,\dbQ),p_A)\arrowsim (W,\psi)$ taking the Betti realization
$w_\alpha$ of $v_\alpha$ into $s_\alpha$, $\forall\,\alpha\in\scrJ$;

\smallskip
\item{(ii)}
composing the homomorphism $h_A\colon\dbS\to GSp(H_1(A,\dbR),p_A)$, defined by the
Hodge structure on $H_1(A,\dbR)$, with an isomorphism
$GSp(H_1(A,\dbR),p_A)\arrowsim GSp(W\otimes\dbR,\psi)$, induced by an
isomorphism as in (i), we get an element of $X$.

\smallskip 
We have a right action of $G(\dbA_f)$ on $\scrA(G,X,W,\psi)$
defined by:

$$[A,p_A,(v_\alpha)_{\alpha\in\scrJ},k]\cdot
g=[A',p_{A'},(v_\alpha)_{\alpha\in\scrJ},g^{-1}k].$$

\noindent
Here $A'$ is the abelian variety, from the same isogeny class as $A$, defined
by the $\dbZ$-lattice of $H_1(A,\dbQ)$ induced from
$L\otimes\dbZhat$ through the isomorphism $g^{-1}\circ k:H_1(A,\dbQ)\otimes\dbA_f\tilde\to W\otimes\dbA_f$, while $p_{A'}$ is the
only rational multiple of $p_A$ which produces a principal polarization of
$A'$ (see [De1, 4.7] for the theorem of Riemann used here). Here as well as in e) below we identify a polarization with its Betti realization. 

There is a $G(\dbA_f)$-equivariant bijection 
$$
f_{(G,X,W,\psi)}\colon\Sh(G,X)(\dbC)\arrowsim\scrA(G,X,W,\psi)
$$ 
defined as follows.
To $[h,g]\in\Sh(G,X)(\dbC)=G(\dbQ)\setminus X\times G(\dbA_f)$ we
associate the quadruple
$[A,p_A,(v_\alpha)_{\alpha\in\scrJ},k]$ where:

\smallskip
d) $A$
is associated to the Hodge structure $(W,h)$ and the $\dbZ$-lattice $H_1(A,\dbZ)$ of $W$
induced from the $\dbZ$-lattice $L$ of $W$ through $k\colon V_f(A)=W\otimes\dbA_f\twiceover{g^{-1}}{\longrightarrow} W\otimes\dbA_f$ (i.e. $k(H_1(A,\dbZ)\otimes\dbZhat)=L\otimes\dbZhat$);

e) $p_A$ is the only (rational) multiple of
$\psi$ which gives birth to a principal polarization of $A$;

f)  $\forall\,\alpha\in\scrJ$, the Betti realization
of $v_\alpha$ is $s_\alpha$.

\smallskip
The inverse $g_{(G,X,W,\psi)}$ of $f_{(G,X,W,\psi)}$ is defined as follows. Let $[A,p_A,(v_\alpha)_{\alpha\in\scrJ},k]\in
\scrA(G,X,W,\psi)$.
We choose a similitude isomorphism $i_A\colon (H_1(A,\dbQ),p_A)
\arrowsim (W,\psi)$ as in (i). It produces an isomorphism $\tilde i_A\colon GSp(H_1(A,\dbQ),p_A)\arrowsim GSp(W,\psi)$. We
define $h\in X$ to be $\tilde i_{A\dbR}\circ h_A$ ($h_A$ being the homomorphism $\dbS\to GSp(H_1(A,\dbR),p_A)$ defining the Hodge structure of $A$) and $g\in G(\dbA_f)$ to be the composite map $W\otimes\dbA_f\twiceover{k^{-1}}{\longrightarrow} V_f(A)=H_1(A,\dbQ)\otimes\dbA_f\twiceover{i_A\otimes 1}{\longrightarrow} W\otimes\dbA_f$.
Then 
$$g_{(G,X,W,\psi)}([A,p_A,(v_\alpha)_{\alpha\in\scrJ},k])=[h,g].$$

\smallskip
Taking $(G,X)\!=\!(GSp(W,\psi),S)$ and $\scrJ\!=\!\phi$, we get a
bijection between the set 
$\Sh(GSp(W,\psi),S)$$(\dbC)$ and the set of isomorphism
classes of principally polarized abelian varieties over
$\dbC$ of dimension $g_W$ (with $2g_W=\dim_{\dbQ}(W)$) having (compatibly) level-$N$ symplectic similitude
structure for any $N\in\dbN$. So to give a $\dbC$-valued point of $\Sh(GSp(W,\psi),S)$ is the same  as to  give a  triple $[A,p_A,(l_N)_{N\in\dbN}]$, where $(A,p_A)$ is a principally polarized abelian variety over $\dbC$ of dimension $g_W$, for which we have a compatible system of similitude isomorphisms $l_N\colon (L/NL,\psi)\arrowsim (H_1(A,\dbZ/{N\dbZ}),p_A)$ ($N\in\dbN$). The compatibility means that if $N,M\in\dbN$ are such that $N|M$, then $l_N$ is obtained from $l_M$ by tensoring with $\dbZ/N\dbZ$. Below we identify as well $L/NL$ (resp. $L\otimes\dbZhat$) with a finite \'etale (resp. with a pro-\'etale) scheme over any base scheme.
\finishproclaim

\Proclaim{4.1.0.} \rm
For $N\in\dbN$ let $K(N):=\{g\in GSp(W,\psi)(L\otimes\dbZhat)\mid g\, \text{mod $N$ is the identity}
\}$.
Then the set $\Sh_{K(N)}(GSp(W,\psi),S)(\dbC)$ is in one to one
correspondence with the set of isomorphism classes of principally 
polarized abelian varieties over $\dbC$ having a level-$N$
symplectic similitude structure.
This implies (cf. [De1, 4.21]) that $\Sh(GSp(W,\psi),S)$ is the $\dbQ$--scheme representing the functor that sends a
$\dbQ$--scheme $T$ to the set of isomorphism classes of
principally polarized abelian schemes (of dimension $g_W$) over $T$, having
(compatibly) level-$N$ symplectic similitude structure for
any $N\in\dbN$ (see [Mu] for why this functor is
representable).
So $\Sh(G,X)$ is the closed subscheme (cf. [De1, 1.15 and 5.9])
of $\Sh(GSp(W,\psi),S)_{E(G,X)}$ whose complex points are those
triples $[A,p_A,(l_N)_{N\in\dbN}]$ for which:

\parindent=32pt
\smallskip
\Item{(4.1.1)}
the isomorphism $k^{-1}\colon L\otimes\dbZhat\arrowsim
H_1(A,\dbZ)\otimes\dbZhat$, defined by the fact that mod $N$ it is $l_N$, $\forall N\in\dbN$, when tensored with $\dbQ$, 
takes $s_\alpha$ to the Betti
realization $w_\alpha$ of a Hodge cycle $v_\alpha$ of $A$
($\forall\,\alpha\in\scrJ$);

\smallskip
\Item{(4.1.2)}
$H_1(A,\dbQ)$ together with $p_A$ and the family of tensors $(w_\alpha)_{\alpha\in\scr J}$ satisfies the
above two conditions (i) and (ii).
\finishproclaim

\parindent=25pt
\Proclaim{4.1.3. Lemma.}
Let $Z=\hbox{\rm \Spec}(R)$ be an integral affine scheme
over $\dbC$ and let $(A,p_A)$ be a principally polarized
abelian scheme over $Z$ having (compatibly) level-$N$
symplectic similitude structure (defined by an isomorphism $l_N\colon L/NL_Z\arrowsim A[N]$) for any $N\in\dbN$.
Let $g_Z\colon Z\to\Sh(G\Sp(W,\psi),S)$ be the morphism induced by the above data.
For every $\alpha\in\scrJ$, we assume the existence of a
cycle $t_\alpha\in F^0((H_{dR}^1(A/Z)\otimes
H_{dR}^1(A/Z)^*)^{\otimes m(\alpha)})$ (we recall that 
$2m(\alpha)=\deg(s_\alpha)$), annihilated by  the Gauss--Manin
connection $\nabla$ (of $A$).
Let $f_1,f_2\colon\hbox{\rm \Spec}(\dbC)\to Z$ be two
complex points.
If the quadruple $[A,p_A,(t_\alpha)_{\alpha\in\scrJ},k]$ (with $k\colon\text{proj.lim.}_{N\in\dbN} A[N]\arrowsim L\otimes\dbZhat$ such that its inverse mod $N$ is  $l_N$; here we identify a Hodge cycle with its de Rham component) becomes a
quadruple of $\scrA(G,X,W,\psi)$ in the point $f_1$, then it
becomes a quadruple of $\scrA(G,X,W,\psi)$ in the point $f_2$
also ($i.e.$ the morphism $g_Z\circ f_2\colon\hbox{\rm \Spec}(\dbC)\to
\Sh(GSp(W,\psi),S)$ factors through
$\Sh(G,X)$).
\finishproclaim

\proof
There is an integral affine scheme $Y=\Spec(T)$ of finite
type over $\dbC$, with $T$ a subring of $R$, such that 
$(A,p_A)$ and its cycles $(t_\alpha)_{\alpha\in\scrJ}$ descend to $(B,q)$ and
cycles $(u_\alpha)_{\alpha\in\scrJ}$.

We have $\nabla u_\alpha=0$, $\forall\alpha\in\scrJ$.
Let $T\hookrightarrow T_1$ be an injective ring
homomorphism, with $T_1$ a smooth integral $\dbC$-algebra, such that $\Spec(T_1)(\dbC)\to Y(\dbC)$ is
surjective (cf. the resolution of singularities).
Let $h_1,h_2\colon \Spec(\dbC)\to Y_1=\Spec(T_1)$ be two
points such that the diagram
$$
\spreadmatrixlines{2\jot}
\matrix
\Spec(\dbC) &\doublemaprights{30}{30}{h_2}{h_1}
  &Y_1\\
\lrdoublemapdown{f_1}{f_2} &&\rmapdown{j}\\
Z &\longrightarrow &Y
\endmatrix
$$
is commutative (the morphisms $j$ and $Z\to Y$ are associated to the inclusions  $T\hookrightarrow T_1$ and respectively $T\hookrightarrow R$).

We denote by $(B_1,q_1)$ and $(u_\alpha^1)_{\alpha\in\scrJ}$ the pull back through $j$ of
$(B,q)$ and $(u_\alpha)_{\alpha\in\scrJ}$.
Let $h\colon B_1\to Y_1$ be the morphism defining the abelian scheme $B_1$.
We get that $\nabla(u_\alpha^1)=0$, and so $u_\alpha^1\in (R^1h_*(\dbC)\otimes (R^1h_*(\dbC))^*)^{\otimes m(\alpha)}$, $\forall\alpha\in\scrJ$. As $u_\alpha^1$ is rational in
$h_1$, we deduce that 
$u_\alpha^1\in(R^1h_*(\dbQ)\otimes R^1h_*(\dbQ)^*)^{\otimes
m(\alpha)}$, $\forall\alpha\in\scrJ$. So $u_\alpha^1$ is rational in $h_2$, $\forall\alpha\in\scrJ$.
From [De3, p. 36] we deduce that the tensors  $(h_2^*(u_\alpha^1))_{\alpha\in\scrJ}$ are 
de Rham components of Hodge cycles $(v_\alpha)_{\alpha\in\scrJ}$
of $A_2:=A\times_{Z}{}_{f_2}\Spec(\dbC)$
(their \'etale components are uniquely determined).
As $Y_1(\dbC)$ is connected, we easily deduce that $A_2$
together with $(v_\alpha)_{\alpha\in\scrJ}$ satisfy the condition (4.1.2).
The isomorphisms $(l_N)_{N\in\dbN}$ are producing an isomorphism
$k_2\colon H_1(A_2,\dbZ)\otimes\dbZhat\arrowsim
L\otimes\dbZhat$.
The fact that $k_2^{-1}$ carries $s_\alpha$ to the Betti
realization of $v_\alpha$ (condition (4.1.1)) can be seen working
mod $N$ (for any $N\in\dbN$).
Multiplying by a natural number big enough all $v_\alpha$
and $s_\alpha$, we can work with families $(\vtil_\alpha)_{\alpha\in\scrJ}$ and
$(\stil_\alpha)_{\alpha\in\scrJ}$ assumed to be integral with respect to
$H_1(A,\dbZ)\otimes\dbZhat$ and $V(\dbZ)\otimes\dbZhat$.
The fact that $k_2(\vtil_\alpha)=\stil_\alpha$, $\forall\alpha\in\scrJ$, results from
the analogue property of the isomorphism $k_1\colon H_1(A_1,\dbZ)\otimes\dbZhat\arrowsim L\otimes\dbZhat$ (with $A_1:=A\times_Z {}_{f_1}\Spec(\dbC)$) and from the fact that a level-$N$
symplectic similitude structure on $Z$ can be descended
to an integral affine $Y$-scheme $Y_N$ of finite type over
$\dbC$ (i.e. for any given $N\in\dbN$ we can assume that the isomorphism $l_N$ is defined over $Y$, and so over $Y_1$).
From the characterization of $\Sh(G,X)(\dbC)$ (cf. 4.1.0), we deduce
that the morphism $g_Z\circ f_2\colon \Spec(\dbC)\to\Sh(G\Sp(W,\psi),S)$
factors through $\Sh(G,X)$ (i.e.
$[A_2,p_{A_2},(v_\alpha)_{\alpha\in\scrJ},k_2]\in\scrA(G,X,W,\psi)$).

\Proclaim{4.1.4. Remark.} \rm
A similar result can be proved if, instead of
$\Sh(G\Sp(W,\psi),S)$ and
$\Sh(G,X)$, we work with
$\scrM:=\Sh_{K_p}(G\Sp(W,\psi),S)$ and
$\scrN:=\Sh_{H_p}(G,X)$, where
$K_p:=\{g\in G\Sp(W,\psi)(\dbQ_p)\mid
g(L\otimes\dbZ_p)=L\otimes\dbZ_p\}$ and $H_p:=K_p\cap
G(\dbQ_p)$ ($p$ being a fixed rational prime).
This follows from the fact that a situation of the form
$\Spec(\dbC)\doublemaprights{30}{30}{f_2}{f_1} Z\to\scrM$, with
$f_1$ factoring through $\scrN$, can be lifted to a situation 
$$
\Spec(\dbC)
\doublemaprights{30}{30}{f_{20}}{f_{10}}  Z_1\to
\Sh(G\Sp(W,\psi),S),
$$ 
with $f_{10}$ factoring through
$\Sh(G,X)$ and with $Z_1$ an integral affine $Z$-scheme.
\finishproclaim

\Proclaim{4.1.5. Remark.} \rm
Later on we need a formal version of 4.1.3 and 4.1.4. We work under the hypotheses of 4.1.3 with $R=\dbC[[z_1,...,z_n]]$ a ring of formal power series over $\dbC$ and with $f_1$ the complex point of $Z$ associated to the surjective ring homomorphism $R\twoheadrightarrow\dbC$ taking all $z_i$ to zero. But instead of assuming that $t_{\alpha}$'s are parallel with respect to $\nabla$, we assume just that $t_{\alpha}$'s are annihilated by all $\frac{\delta}{\delta z_i}$. Then the generic point $w$ (this replaces the point $f_2$ of 4.1.3) of $Z$ is mapped through $g_Z$  into $\Sh(G,X)$, i.e. the cycles $t_{\alpha}$ become (in $w$) de Rham components of Hodge cycles of $A_w$ (the fibre of $A$ over $w$) and the \'etale components of these Hodge cycles are related to $v_{\alpha}$ (through the family of isomorphisms $(l_N)_{N\in\dbN}$) as expected. 

It is enough to see the first part, i.e. that $t_{\alpha}$ becomes in $w$ the de Rham component of a Hodge cycle of $A_w$, $\forall\alpha\in\scrJ$ (the second part involving the expected relation is entirely the same as in the above proof of 4.1.3).
This is a result of Faltings. The proof of this is entirely analogous to the proof of its integral version [Fa3, rm iii) after Th. 10]. The only difference is that now we have to use the strictness property of  maps between Hodge structures, instead of the strictness property of maps between objects of $\scrM\scrF(V_0)$ (cf. [Fa1] for the definition of $\scrM\scrF(V_0)$; here $V_0$ is a Witt ring over a perfect field). 
\finishproclaim

\Proclaim{4.1.6. Remark.} \rm
Sometimes it is more convenient to work with families $(s_{\alpha})_{\alpha\in\scrJ}$ such that $G$ is the subgroup of $GL(W)$ (and not of $GSp(W,\psi)$) fixing its tensors. This has the advantage that we can be loose about  mentioning alternating forms (like $\psi$) or different Tate-twists (to be compared with 5.2.9). In particular, in such a situation, the form $\psi$ is uniquely determined by an isomorphism as in (i) of 4.1, up to scalar multiplication with a rational number; so it is more natural to denote the set $\scrA(G,X,W,\psi)$ just by $\scrA(G,X,W)$.  
\finishproclaim

\smallskip
\Proclaim{4.2. Digression on reductive Lie algebras.}\rm
Till the end of \S4 the notations to be introduced are independent of the ones in 4.1. Let $W$ be a finite vector space over an arbitrary field of
characteristic zero. All the reductive Lie subalgebras of $\grg\grl(W)$ considered in 4.2 are assumed to satisfy the following condition: their centers are generated by semisimple endomorphisms of $W$ and the trace forms on them are perfect.

Let $\grg\subset\grg\grl(W)$ be (such) a reductive Lie subalgebra. It is known that the above assumption implies that the restriction to $\grg$ of the trace
form $\Tr$ on $\grg\grl(W)$ is perfect (for $a,b\in\grg\grl(W)$, $\Tr(a,b)$ is the trace of the endomorphism $ab$ of $W$): one just needs to point out that the restriction of $\Tr$ to the semisimple part $[\grg,\grg]$ of $\grg$ is automatically non-degenerate (argument: the kernel of this restriction is an ideal of $\grg$ and so a semisimple Lie algebra; based on Cartan's solvability criterion in characteristic 0, it is as well solvable and so trivial). For any vector subspace $\grm$ of $\grg\grl(W)$ let 
$$\grm^\bot:=\{x\in\grg\grl(W)\mid\Tr(xy)=0,\,\forall\,y\in\grm\}.$$
In particular we get a direct sum decomposition $\grg\grl(W)=\grg\oplus \grg^\bot$.

\Proclaim{4.2.1. Convention.} \rm                                              Any time we have a situation as above,
we  denote by $\pi(\grg)$ (or by $\pi_W(\grg)$) the projector of $\grg\grl(W)$ defined
by $\pi(\grg)(x)=x$ if $x\in\grg$ and $\pi(\grg)(x)=0$ if $x\in\grg^\bot$.
\finishproclaim

The Lie subalgebra of $\grg\grl(W)$ centralizing $\pi(\grg)$ under the adjoint representation is of the
form $\grg\oplus \gru$, where
$$
\gru:=\{y\in \grg^\bot\mid[\grg,y]\subset\grg,\,
[\grg^\bot,y]\subset \grg^\bot\}=
\{y\in\grg^\bot\mid [\grg,y]=\{0\}\}.
$$
The last equality is due to the fact that $[\grg,\grg^\bot]\subset
\grg^\bot$ and $\Tr([a,b],c)=\Tr(a,[b,c])$, $\forall\, a,b,c\in\grg\grl(W)$.

\Proclaim{4.2.2. Proposition.}
Let $\grg\subset\grh\subset\grg\grl(W)$ be inclusions of reductive Lie algebras.
We consider reductive Lie algebras $\grg_1$ satisfying : \ a) $\grg\subset\grg_1\subset\grh$ and
\ b) $[\grg,\grg]=[\grg_1,\grg_1]$.
They form a set $\scrS$.
Then an element $\grg_1$ of $\scrS$ is maximal under the
relation of inclusion if and only if
$\grg_1=\grh\cap\{\text{the Lie subalgebra of }\grg\grl(W)$
centralizing $\pi(\grg_1)\}$.
\finishproclaim

\proof
If $\grg_1=\grh\cap\{\text{the Lie subalgebra of }\grg\grl(W)$
centralizing $\pi(\grg_1)\}$ then
$$
\grh\cap\grg_1^\bot\cap\{\text{centralizer of }\grg_1\text{
in }\grg\grl(W)\}=\{0\}.
$$
This implies that there is no reductive Lie subalgebra of
$\grh$ strictly containing $\grg_1$ and having the same semisimple
part as $\grg_1$.
So $\grg_1$ is a maximal element of $\scrS$.

Let now $\grg_1$ be a maximal element of $\scrS$.
We deduce that the centralizer $\grc$ of
$\grg_1$ in $\grh$ has no reductive Lie subalgebra included in
$\grg_1^\bot$. But $\grc$ is a reductive Lie subalgebra of $\grg\grl(W)$. Argument: the centralizer of $\grg_1$ in $\grg\grl(W)$ is the Lie algebra of a reductive group (this can be seen moving to an algebraically closed field and using irreducible representations) and it is of the form $\grc\oplus {\bar{\grc}}$ with ${\bar{\grc}}$ a subspace of $\grh^\bot$; so the trace form on $\grc$ is perfect. This implies  that $\grc\cap\grg_1^\bot$ is zero. So $\grg_1$ is the subalgebra of $\grh$
centralizing $\pi(\grg_1)$. This ends the proof of the Proposition.

\Proclaim{4.2.3. Remark.} \rm
Let $f\colon (G,X)\hookrightarrow (G\Sp(W,\psi),S)$ be an injective map. If in 4.2.2 we take $\grg=\Lie(G)$ and $\grh=\gsp(W,\psi)$, then for any maximal element $\grg_1$ of $\scrS$ there is a uniquely determined (up to isomorphism) Shimura variety $\Sh(G_1,X_1)$ for which there are injective maps $f_0\colon(G,X)\hookrightarrow(G_1,X_1)$ and $f_1\colon(G_1,X_1)\hookrightarrow (G\Sp(W,\psi),S)$ such that $f=f_1\circ f_0$ and  $df_1(\Lie(G_1))=\grg_1$.
\finishproclaim

\smallskip
\Proclaim{4.3. Special families of tensors.}\rm
\finishproclaim

\Proclaim{4.3.1. Definitions.} \rm
Let $(G,X)$ define an arbitrary Shimura variety.
A pair $(G_1,X_1)$ is called an enlargement of $(G,X)$ if
there is an injective map $f\colon(G,X)\hookrightarrow
(G_1,X_1)$ such that $f(G^{\der})=G_1^{\der}$ and
$f(G)\not=G_1$.
If $i\colon(G,X)\hookrightarrow (G_2,X_2)$ is an injective
map, by an enlargement of $(G,X)$ in $(G_2,X_2)$ we mean a
pair $(G_1,X_1)$, with $G\subsetneqq G_1\subset G_2$, $G^{\der}=G_1^{\der}$ and
$X\subset X_1\subset X_2$.
We say $(G,X)$ is saturated in $(G_2,X_2)$ if it
has no enlargement in $(G_2,X_2)$.
\finishproclaim

\Proclaim{4.3.1.1.} \rm
Let now $(G,X)$ be of Hodge type and let
$f\colon(G,X)\hookrightarrow (G\Sp(W,\psi),S)$ be an injective map.
From 4.2.2 and 4.2.3 we deduce that
either $(G,X)$ is saturated in $(G\Sp(W,\psi),S)$ or
there is an enlargement of $(G,X)$ in $(GSp(W,\psi),S)$ which is saturated in
$(G\Sp(W,\psi),S)$.

The advantage of having injective maps
$(G,X)\hookrightarrow(G\Sp(W,\psi),S)$ with $(G,X)$ saturated in
$(G\Sp(W,\psi),S)$ is: $\Lie(G)$ is the Lie
subalgebra of $\gsp(W,\psi)$  centralizing (just one tensor of degree $4$ which is a projector of $\grg\grl(W)$) $\pi_{W}(\grg)$.
\finishproclaim

\Proclaim{4.3.2.} \rm
We consider now the following situation.
Let $(W,\psi)$ be a symplectic space over a field of
characteristic zero. Let $G_0$ be a semisimple subgroup of $G\Sp(W,\psi)$ and let $\grg_0:=\Lie(G_0)$. Let $G$ be a reductive subgroup of $G\Sp(W,\psi)$ having $G_0$ as its derived subgroup and such that its Lie
algebra $\grg$ is the Lie subalgebra of $\gsp(W,\psi)$
centralizing $\pi(\grg)$ (cf. 4.2.2). 
We now list some useful tensors fixed by the group
$G$.

We have $\gsp(W,\psi)=\grg\oplus\grh$, with
$\grh:=\gsp(W,\psi)\cap\grg^\bot$.
Let $\grh=\operatornamewithlimits{\oplus}\limits_{i\in I}\grh_i$ be a direct sum 
decomposition of $\grh$ in (non-zero) irreducible $\grg$-modules.
Let $\grm_i$ be the kernel of the representation
$\grg\to\grg\grl(\grh_i)$.
We deduce the existence of a reductive Lie
subalgebra $\grg_i$ of $\grg$ such that $\grg$ is the direct
sum of Lie algebras $\grg=\grg_i\oplus\grm_i$ [Bou1, p. 57].
We have faithful irreducible representations
$\grg_i\hookrightarrow\grg\grl(\grh_i)$.
Associated to the direct sum decomposition
$\grg\grl(W)=\grg\oplus\operatornamewithlimits{\oplus}\limits_{i\in
I}\grh_i\oplus\gsp(W,\psi)^\bot$ we consider the
projectors $p_i\colon\grg\grl(W)\to\grg\grl(W)$, the image of $p_i$
being $\grh_i$, $\forall\,i\in I$.
For every $i\in I$, let $r_i$ be the projection of $\grg\grl(W)$ on
$\grg_i$ associated to the direct sum decomposition
$\grg\grl(W)=\grg_i\oplus\grm_i\oplus\grh\oplus \gsp(W,\psi)^\bot$.

For $i\in I$, let $k_i$ be the Casimir element of the
representation $\grg_i\hookrightarrow\grg\grl(\grh_i)$ (we have
$\grg_i\not=0$, as $\grg$ is the Lie subalgebra of
$\gsp(W,\psi)$ centralizing $\pi(\grg)$).

The element $k_i$ induces a linear map  $q_i\colon\grg\grl(W)\to\grg\grl(W)$ such
that $q_i|\grh_i\colon \grh_i\to \grh_i$ is an isomorphism.
We choose a  linear combination of $(q_i)_{i\in I}$ with coefficients
in $\dbZ$ such that the resulting  linear map $q\colon\grg\grl(W)\to\grg\grl(W)$ 
has the property that $q|\grh\colon \grh\to \grh$ is an isomorphism (using induction,
it is enough to handle the case when $I$ has two elements; but this case is obvious, as $\dbZ$ is infinite).
Let $\qbar\colon\grg\grl(W)\to\grg\grl(W)$ be the linear map  such that
$\qbar$ is zero on $\grg\oplus\gsp(W,\psi)^\bot$ and
$\qbar|\grh\colon \grh\to \grh$ is $(q\vert \grh)^{-1}$.

For  $i\in I$, let $t_i\colon\grg\grl(W)\to\grg\grl(W)^*$ be the
linear map such that $t_i$ is zero on
$\grm_i\oplus \grh\oplus\gsp(W,\psi)^\bot$ and
$t_i|\grg_i\colon\grg_i\to\grg_i^*$ is the isomorphism
induced by the restriction to $\grg_i$ of the trace form $\Tr_{\grh_i}$ on $\grg\grl(\grh_i)$.
Explicitly: if $x\in\grg_i$, then
$t_i(x)(y)=\Tr_{\grh_i}(x,y)$.
For $i\in I$, let $s_i\colon\grg\grl(W)^*\to\grg\grl(W)$ be the
linear map which is zero
on $(\grm_i\oplus \grh\oplus\gsp(W,\psi)^\bot)^*$ and
$s_i\vert\grg_i^*\colon\grg_i^*\to\grg_i$ is
$(t_i\vert\grg_i)^{-1}$.

Let $t\colon\grg\grl(W)\to\grg\grl(W)^*$ and
$s\colon\grg\grl(W)^*\to\grg\grl(W)$ be linear maps defined in the
same manner as $t_i$ and $s_i$ but for the representation
$\grg\hookrightarrow\grg\grl(W)$.

Let 
$$
B\colon\grg\grl(W)\to\grg\grl(W)^*
$$ 
be the linear map which is zero on $\grg_0^\bot$ and $B|\grg_0\colon\grg_0\to\grg_0^*$ is the isomorphism induced by the Killing form on $\grg_0$. Let 
$$
B^*\colon\grg\grl(W)^*\to\grg\grl(W)
$$ 
be the linear map obtained from $B$ in the same manner as the tensors $s_i$ were obtained from $t_i$.

The tensors $\pi(\grg_0)$, $\pi(\grg)$, $B$, $B^*$, $\qbar$, $t$ and $s$, as well as the tensors 
$p_i$, $r_i$, $s_i$ and $t_i$,  $i\in I$, are centralized by $\grg$, and so
fixed by the group $G$.

\Proclaim{4.3.3. Notation.} \rm
Let $\tilde W$ be a finite vector space over a field $k$ of characteristic zero and let $\tilde\grg$ be the Lie algebra of a  semisimple subgroup $\tilde G$ of $GL(\tilde W)$. 
We call an $\grs\grl_2$ Lie subalgebra of $\tilde\grg\otimes {\kbar}$ standard if with respect to a Weyl direct sum decomposition ${\tilde\grg\otimes {\kbar}}={\grt\opoplus_{\alpha\in\Phi}\grg_{\alpha}}$, with $\Phi$ the system of roots associated to a maximal torus $T$ of $\tilde G_{\kbar}$ (so $\grt=\Lie(T)$), is generated by $\grg_{\alpha}$ and  $\grg_{-{\alpha}}$, for some $\alpha\in\Phi$.
We denote by $s(\tilde\grg,\tilde W)$ the maximum dimension which appears among the irreducible subrepresentations
of $\tilde  W\otimes {\kbar}$ of any standard $\grs\grl_2$ Lie subalgebra of  $\tilde\grg\otimes {\kbar}$. 
\finishproclaim

\Proclaim{4.3.4. Definitions.} \rm
Let $O$ be a discrete valuation ring, let $\pi$ be a
uniformizer of it, and let $K$ be its field of fractions.
Let $(W,\psi)$ be a symplectic space over $K$.
Let $(s_\alpha)_{\alpha\in\scrJ}$ be a family of tensors in
spaces of the form $W^{\otimes m}\otimes W^{*\otimes n}$. The family of tensors $(s_{\alpha})_{\alpha\in\scrJ}$ is called essentially finite, if the $O$-submodule of the tensor algebra of $W\oplus W^*$ generated by its tensors, is a free $O$-module of finite rank. Let $R$ be a faithfully flat integral ring over $O$. A free $R$-module $M$ satisfying $M\fracwithdelims[]1\pi=W\otimes_K R\fracwithdelims[]1\pi$, is said to envelop the above family of tensors with respect to $\psi$, if $\psi$  induces a perfect form $\psi\colon M\otimes M\to R$ and if all the tensors of the family $(s_\alpha)_{\alpha\in\scrJ}$ are in spaces of the form $M^{\otimes m}\otimes
M^{*\otimes n}$.
Let $H$ be a  reductive subgroup of $G\Sp(W,\psi)$ fixing the tensors of the above family.
The family of tensors $(s_\alpha)_{\alpha\in\scrJ}$ is said to be
$O$-well positioned with respect to $\psi$ for the group $H$
if the following condition is satisfied:

\smallskip
(4.3.5) For any  faithfully flat integral ring $R$ over
$O$ and for any free $R$-module $M$, satisfying $M\fracwithdelims[]1\pi=W\otimes_K R\fracwithdelims[]1\pi$ and enveloping the family of tensors
$(s_\alpha)_{\alpha\in\scrJ}$ with respect to $\psi$, the
Zariski closure of $H_{R\fracwithdelims[]1\pi}$ in $G\Sp(M,\psi)$ is a reductive group scheme $H_R$ over $R$. 

\smallskip
In addition, if there is an $O$-lattice $M_O$ of $W$ enveloping the family of tensors ${(s_\alpha)}_{\alpha\in\scrJ}$ with respect to $\psi$, then we say that our family of tensors is $O$-very well positioned with respect to $\psi$ for the group $H$.

We have variants, depending on the class of $O$-algebras we use in 4.3.5. If we use the class of normal integral faithfully flat $O$-algebras (resp. of reduced faithfully flat $O$-algebras) we obtain the notion of weakly (resp. strongly) $O$-well (or $O$-very well) positioned families of tensors with respect to $\psi$ for the group $H$.   
\finishproclaim

\Proclaim{4.3.6. Remarks.} \rm
0) Warning: if $H$ extends to a reductive group over $O$, we do not require that the extension of it to $R$ is isomorphic to $H_R$.

1) If the family of tensors $(s_\alpha)_{\alpha\in\scrJ}$ is $O$-well
positioned (resp. $O$-very well positioned) with respect to $\psi$ for the group $H$, then
the family $(s_\alpha)_{\alpha\in\scrJ}$ is $O$-well
positioned (resp. $O$-very well positioned) with respect to $\psitil$ for the group $\Htil$,
where $\psitil\colon W^*\otimes W^*\to K$ is the
perfect alternating form on $W^*$ obtained from
$\psi$ through the isomorphism $\tilde {f}\colon W\arrowsim W^*$
canonically induced by $\psi$, ($\tilde {f}(x)(y)=\psi(x,y)$) and
where $\Htil$ is the subgroup of $G\Sp(W^*,\psitil)$ corresponding
to $H$ under the canonical identification of $GL(W)$ with
$GL(W^*)$ produced by $\tilde f$.

2) The family of tensors $(s_{\alpha})_{\alpha\in\scrJ}$ is $O$-well positioned with respect to $\psi$ for $H$ iff it is well positioned with respect to $\psi$ for $H^{\der}$ and for the toric part of $Z(H)$ (cf. 3.1.6.1). The same remains true in a weakly (this is obvious) or strongly (as the proof of 3.1.6 applies) context. 

3) Let $R$ be a noetherian, reduced, faithfully flat, local $O$-algebra. Let $M$ be a free $R$-module of finite rank. Let $H^\prime$ be a reductive subgroup of $GL(M)_{R\fracwithdelims[]1\pi}$. We assume there is a noetherian, integral, local $O$-subalgebra $R_0$ of $R$ such that $R$ is the strict henselization of $R_0$, $M$ is obtained by extension of scalars from a free $R_0$-module $M_0$ and $H^\prime$ is the pull back of a reductive subgroup $H^\prime_0$ of $GL(M_0)_{R_0[{1\over \pi}]}$. We also assume the existence of a projector $\Pi$ of $\End(M)$ having the properties:

\smallskip
-- its image $\Lie$ is such that $\Lie[{1\over \pi}]=\Lie(H^\prime)$;

-- as a projector of $\End(M[{1\over {\pi}}])$ it is fixed by $H^\prime$.

\smallskip
 So we have a direct sum decomposition
$\End(M)=\text{Lie}\oplus\text{Ker}(\Pi)$ preserved after inverting $\pi$ by the action of $H^\prime$ on $\End(M[{1\over {\pi}}])$. Let $l$ (resp. $l_R$) be the residue field of $O$ (resp. of $R$). Let $\tilde k$ be the algebraic closure of the field of fractions $\tilde k_0$ of $R_0$.

Let $\{I_i|i\in\{1,...,r\}\}$, with $r\in\dbN$, be a set of ideals of $R$ which are intersections of prime ideals of $R$ of codimension 0. We assume that $\cap_{i=1}^{i=r} I_i=0$ and that the Zariski closure of $H^\prime_{R/I_i\fracwithdelims[]1\pi}$ in $GL(M\otimes R/I_i)$ is a reductive group $H_i^\prime$  over $R/I_i$, $\forall i\in\{1,...,r\}$. 

\Proclaim{Definition.} \rm
A faithful representation $\rho_T:T\hookrightarrow GL(V)$ of a split torus over an arbitrary field $\tilde k$ is called $p$ Lie recoverable, if $T$ is generated by a finite family $(T_j)_{j\in J}$ of subtori of it such that for any $j\in J$, the representation of $T_j$ on $V$ does not involve two distinct characters of $T_j$ whose difference inside the group of characters of $T_j$ (viewed additively) is divisible by $p$. If $J$ has just $1$ element, then we also say $\rho_T$ is elementary $p$ Lie recoverable.
\finishproclaim

With this definition, we can state:

\Proclaim{Proposition.} 
The Zariski closure of $H_0^\prime$ in $GL(M_0)$ is a reductive group if moreover any one of the following three conditions hold:
\medskip
i) $l$ is of characteristic $p$, $H^\prime$ is a torus and the representation $\rho_{\tilde k}$ of ${H^\prime_0}_{\tilde k}$ on $M_0\otimes_{R_0} \tilde k$ is $p$ Lie recoverable;
\smallskip
ii) $l$ is of characteristic $0$;
\smallskip
iii) $H^\prime_0$ is non-trivial and semisimple, the characteristic of $l$ is an odd prime $p\ge s(\Lie({H_0^\prime}_{\tilde k}),M_0\otimes_{R_0} \tilde k)$ ($s(\Lie({H_0^\prime}_{\tilde k}),M_0\otimes_{R_0} \tilde k)\in\dbN$ is defined as in 4.3.3 even if the characteristic of $O$ is not $0$), and moreover the semisimplification $AUT^{\text{ss}}$ of the connected component $AUT^0$ of the origin of $Aut(\Lie\otimes_R l_R)_{\text{red}}$ is smooth of dimension equal to the relative dimension of $H^\prime_0$.
\finishproclaim

\proof
We start with a warning: part of the proof below uses Part 2 of 4.3.10; so it should be read only after reading the mentioned Part 2. This Proposition is used in what follows just to get some complements in 4.3.7 4) and 5) and in 4.3.10.2; for the Part 2 of 4.3.10, from 4.3.7 4) and 5) only the first paragraph of 4.3.7 4) is needed. So no vicious circle is created.

It is enough to show $H_R^\prime$ is a reductive subgroup of $GL(M)$. We start showing that the pull backs of $H_{i_1}^\prime$ and $H_{i_2}^\prime$ to $\Spec(R/I_{i_1}\cap I_{i_2})$ coincide, $\forall i_1, i_2\in\{1,...,r\}$; we refer to this property as $IND$. We can assume $(i_1,i_2)=(1,2)$.

The first step is to check that the pull backs of ${H^\prime_1}_{R/I_1+I_2}$ and ${H^\prime_2}_{R/I_1+I_2}$ to ${\Spec(R/I_1+I_2)}_{\text{red}}$ coincide. For ii) this is a consequence of [Bo, 7.1]. For i) this can be checked as follows. The generic fibre $GF$ of $H^\prime_0$ splits over a Galois field extension $\tilde k_1$ of $\tilde k_0$ contained in the residue field of any point of $\Spec(R/I_i)$ of codimension $0$, $\forall i\in I$. So (by replacing $R_0$ with $R\cap \tilde k_1$) we can assume $GF$ is a split torus. So, based on similar arguments to the ones of 3.1.6 pertaining to $T_R^2$, we can assume $\rho_{\tilde k}$ is elementary $p$ Lie recoverable; but this case follows from very definitions. In other words: we can recover the maximal direct sum decomposition of $M\otimes_R R/I_i$ normalized by $H^\prime_i$ from the representation of $\Lie\otimes_R R/I_i$ on $M\otimes_R R_i$; moreover, the character of the representation of $H^\prime_i$ on any fixed member of it can be read out from the representation of $GF$ on $M_0\otimes_{R_0} \tilde k_0$ ($i=\overline{1,2}$).

For the time being we assume that the pull backs of ${H^\prime_1}_{R/I_1+I_2}$ and ${H^\prime_2}_{R/I_1+I_2}$ to ${\Spec(R/I_1+I_2)}_{\text{red}}$ coincide even for iii); so let ${H^\prime_{12}}_{\text{red}}$ be the reductive subgroup of $GL(M\otimes_R {(R/I_1+I_2)}_{\text{red}})$ with which these two pull backs coincide. Let $R_{12}:=R/I_1+I_2$ and let $N$ be its ideal of nilpotent elements. Let $n\in\dbN$ be such that $N^n=\{0\}$. We consider a maximal torus $T_i$ of ${H^\prime_i}_{R_{12}}$, $i=\overline{1,2}$. We can assume that $T_1$ and $T_2$ lift the same maximal torus $T_{12}$ of ${H^\prime_{12}}_{\text{red}}$, cf. [SGA3, Vol. II, 3.6 of p. 48]. Loc. cit. implies as well the existence of $g\in GL(M\otimes_R R_{12})$ which modulo $N$ is the identity and such that $gT_1g^{-1}=T_2$ as tori of $GL(M\otimes_R R_{12})$.

The second step is to check that $\Lie(T_1)=\Lie(T_2)$ implies  $T_1=T_2$. For i) this is a consequence of Definition: as above, we can recover the maximal direct sum decomposition of $M\otimes_R R_{12}$ normalized by $T_1$ as well as the characters through which $T_1$ acts on its members from $\Lie(T_1)$ (and $GF$). The same applies to ii). 

For iii) we have to proceed as in Part 2 of 4.3.10. We have two subcases: 

\smallskip
-- if $p>2s(\Lie({H_0^\prime}_{\tilde k}),M_0\otimes_{R_0} \tilde k)$, then we can apply the part of the previous paragraph referring to i) (argument: we need to consider as in the mentioned place ``standard" $1$ dimensional split subtori of $T_1$ in order to get that we are in the context of elementary $p$ Lie recoverable representations of their fibres over $\tilde k$; so the representation of $T_{1\tilde k}$ on $M_0\otimes_{R_0} \tilde k$ is $p$ Lie recoverable); 

-- if $s(\Lie({H_0^\prime}_{\tilde k}),M_0\otimes_{R_0} \tilde k)\le p\le 2s(\Lie({H_0^\prime}_{\tilde k}),M_0\otimes_{R_0} \tilde k)$, then we have to use as in the mentioned Part 2 standard ${\grs}{\grl}_2$ Lie subalgebras of $\Lie\otimes_R R_{12}$ in order to be able to recover the maximal direct sum decomposition of $M\otimes_R R_{12}$ normalized by $T_1$, from the representations of $\Lie(T_1)$ and of $\Lie\otimes_R R_{12}$ on it.

\smallskip
The third step is to show, by induction on $j\in\{1,...,n\}$, that we can assume that the $R_{12}$-submodules $\Lie(T_1)$ and $\Lie(T_2)$ of $\End(M\otimes_R R_{12})$, coincide modulo $N^j$. The passage from $j$ to $j+1$ goes as follows. We can assume $g$ modulo $N^j$ is the identity (this is just a variant of the second step, with $R_{12}$ replaced by $R_{12}/N^j$). The existence of $g$ implies the existence of $A\in N^j\End(M\otimes_R R_{12})$ such that under the automorphism $1_{\End(M\otimes_R R_{12})}+{\ad}(A)$ of $\End(M\otimes_R R_{12})$, $\Lie(T_1)$ modulo $N^{j+1}$ is mapped onto $\Lie(T_2)$ modulo $N^{j+1}$. But due to the existence of $\Pi$, modulo $N^{j+1}$ we can assume $A\in N^j\Lie\otimes_R R_{12}$. As the Lie algebra of $H_1^\prime$ is $\Lie\otimes_R R/I_1$, we can assume $g$ modulo $N^{j+1}$ is an $R_{12}/N^{j+1}$-valued point of $H_1^\prime$. So, by replacing $T_1$ by an $H_1^\prime(R_{12})$-conjugate of it, we can assume that the restrictions modulo $N^{j+1}$ of $\Lie(T_1)$ and $\Lie(T_2)$ coincide. So, by induction, we can assume $\Lie(T_1)=\Lie(T_2)$.

What follows next for ii) and iii) is very much the same as Part 2 of 4.3.10: using exponential maps (for iii) this is allowed as $p\ge s({H^\prime}_{0\tilde k},M_0\otimes_{R_0} \tilde k)$), we can recover from $T_1$ and $\Lie\otimes_R R_{12}$ the $\dbG_a$ subgroups of ${H^\prime_i}_{R_{12}}$ normalized by $T_1$ and so, based on [SGA3, Vol. III, 4.1.2 of p. 172] we get that an open subscheme of ${H^\prime_1}_{R_{12}}$ coincides with an open subscheme of ${H^\prime_2}_{R_{12}}$ and so ${H^\prime_1}_{R_{12}}={H^\prime_2}_{R_{12}}$. This paragraph forms the fourth step. 

So to end the proof of the fact that ${H^\prime_1}_{R_{12}}={H^\prime_2}_{R_{12}}$, we need to argue the first step for iii). The natural homomorphism from ${H^\prime_i}_{l_R}$ to $AUT^{\text{ss}}$ is an isogeny (it is nothing else but the central isogeny ${H^\prime_i}_{l_R}\to {H^{\prime{\ad}}_i}_{l_R}$; as $p$ is odd this can be checked easily using standard ${\grs}{\grl}_2$ Lie subalgebras of $\Lie\otimes_R l_r$). We easily get: $\Lie(T_1)\otimes_{R_{12}} l_R$ is the Lie algebra of a maximal torus of ${H^\prime_2}_{l_R}$ and so we can assume it is the Lie algebra of $\Lie(T_2)\otimes_{R_{12}} {l_R}$. As in the above part referring to two subcases, we get ${T_1}_{l_R}={T_2}_{l_R}$. As in the previous paragraph we get that ${H^\prime_1}_{l_R}$ and ${H^\prime_2}_{l_R}$ coincide. But the role of $l_R$ is just to fix some notations: the same applies to any point $\Spec(\tilde l_R)\to \Spec(R_{12})$, with $\tilde l_R$ an algebraically closed field containing $l_R$, lifting a geometric point of the special fibre of $\Spec(O)$; we just need to mention that $\Lie\otimes_R \tilde l_R$ and $\Lie\otimes_R l_R\otimes_{l_R} \tilde l_R$  are isomorphic (cf. the uniqueness theorem of [SGA3, Vol. III, p. 313-314] and the fact that $H_0^\prime$ is over $R_0[{1\over {\pi}}]$).

From now we forget about $R_0$ and $H_0^\prime$ and just use $IND$ to get that $H_R^\prime$ is reductive. By induction we can assume $r=2$. As ${H^\prime_1}_{R_{12}}={H^\prime_2}_{R_{12}}$ and as the fibres of $H^\prime_i$, $i=\overline{1,2}$, are connected, we get: the fibres of $H^\prime_R$ are connected. The ring $R/I_i$ is also strictly henselian and so $H^\prime_i$ is a split group. This implies that $H^\prime$ itself is split: as ${H^\prime_1}_{R_{12}}={H^\prime_2}_{R_{12}}$, we get this by first applying (we are in a local, strictly henselian context) [SGA3, Vol. III, 1.5 of p. 329] over $R/I_1+I_2$ and then by making $\pi$ invertible. Let $H^{''}_R$ be a split reductive group over $R$ having $H^\prime$ as its generic fibre. As $H^{''}_R$ is smooth and as ${H^\prime_1}_{R_{12}}={H^\prime_2}_{R_{12}}$, the amalgamated sum of $H^\prime_1$ and $H^\prime_2$ along $H^\prime_{12}$ is a reductive group over $R$ which can be identified (cf. the uniqueness theorem of [SGA3, Vol. III, p. 313-4]) with $H^{''}_R$. We get a natural homomorphism $q\colon H^{''}_R\to GL(M)$ factoring through $H^\prime_R$. As $q$ is a closed embedding over $R/I_i$ we deduce that $q$ itself is a closed embedding. So $H_R$ can be identified with $H^\prime_R$. This ends the proof.
\finishproclaim

\Proclaim{4.3.7. Remarks.} \rm
1) We could have worked out 4.3.4 without the relative context, i.e. with respect to $\psi$. The relative context is all we need for applications to Shimura varieties of Hodge type. When the role of $\psi$ is irrelevant (for instance in 4.3.10 b)) we do not mention with respect to $\psi$. 

2) The definition of $O$-well positioned families of tensors presented here is different from the one in [Va1, 3.7.4], where we also asked that the subgroup of $GSp(M,\psi)$ fixing $v_{\alpha}$, $\forall\alpha\in\scrJ$, is a group scheme whose connected components of the origins of its fibres are (reductive groups defined by) the fibres of $H_R$.

3) Let $R_0$ be an integral ring and let $M_{R_0}$ be a free $R_0$-module of finite rank. Let $\scrK_0$ be the field of fractions of $R_0$ and let $G_{\scrK_0}$ be a subgroup of $GL(M\otimes\scrK_0)$. It is not always true (cf. [BT, 3.2.15]) that the Zariski closure $G_{R_0}$ of $G_{\scrK_0}$ in $GL(M_{R_0})$ is a group subscheme of $GL(M_{R_0})$. However, $G_{R_0}$ is a group subscheme of $GL(M_{R_0})$ if it is a flat scheme over $R_0$.

So, in $4.3.5$, the fact that the Zariski closure $H_R$ of $H_{R\fracwithdelims[]1\pi}$ in $G\Sp(M,\psi)$ is a group subscheme of $GSp(M,\psi)$ is part of the requirements on a family of tensors $(s_{\alpha})_{\alpha\in\scrJ}$ in order to be $O$-well positioned with respect to $\psi$ for the group $H$. To show that $H_R$ is a reductive group scheme over $R$, we need to check two things:

\smallskip
a) that $H_R$ is flat over $R$ (and so a group subscheme of $GSp(M,\psi)$);

\smallskip
b) that the fibres of $H_R$ over $\Spec(R)$ are reductive groups (over fields).
 
\smallskip
Simple arguments at the level of tangent spaces and of reduction to the case $R$ noetherian, show that a) follows from b). 

\smallskip
4) If the family of tensors $(s_{\alpha})_{\alpha\in\scrJ}$ is essentially finite, for proving that it is $O$-well positioned with respect to $\psi$ for the group $H$, it is enough to check $4.3.5$ only for integral rings $R$ which are faithfully flat and of finite type over $O$ (and so noetherian). To see this, let $R$ and $M$ be as in $4.3.5$.  We choose  a basis $\scrB$ of $M$. It naturally produces  a basis of the tensor algebra of $M\oplus M^*$. Let $R_1$ be a finitely generated $O$-subalgebra of $R$ such that $\scrB$ is included in $W\otimes_K R_1\fracwithdelims[]1\pi$. Let $R_2$ be the $O$-subalgebra of $R$ generated by $R_1$ and by the coefficients of all $s_{\alpha}$ with respect to the above basis of the tensor algebra of $M\oplus M^*$. The $O$-algebra $R_2$ is finitely generated as the family of tensors $(s_{\alpha})_{\alpha\in\scrJ}$ is essentially finite. Let $M_2$ be the free $R_2$-submodule of $M$ generated by the elements of $\scrB$. We have $M_2\fracwithdelims[]1\pi=W\otimes_K R_2\fracwithdelims[]1\pi$. Moreover $M_2$ envelopes the family of tensors $(s_{\alpha})_{\alpha\in\scrJ}$. So, if the Zariski closure of $H_{R_2\fracwithdelims[]1\pi}$ in $GSp(M_2,\psi)$ is a reductive group scheme over $R_2$, then, by pull back, the Zariski closure of  $H_{R\fracwithdelims[]1\pi}$ in $G\Sp(M,\psi)$ is a reductive group scheme over $R$.

We assume now that the extra assumptions of 4.3.6 3) pertaining to iii) of 4.3.6 3) hold; so in particular, there is a projector of $\text{End}(W)$ on $\Lie(H)$ fixed by $H$ and which is part of our family of tensors. Localizing $R$, replacing it by a quotient $R_2$ of $R^{\text{sh}}$ dominating $R$, or by $R_1$, where $\Spec(R_1)$ is an integral finite flat scheme over $\Spec(R)$  (the operation of taking the Zariski closure of $H_{R\fracwithdelims[]1{\pi}}$ in $GSp(M,\psi)$ is well behaved --in connection to checking 4.3.5-- with respect to these operations, cf. a) and b) of 3) and iii) of 4.3.6 3)) we can assume, for checking 4.3.5, that:

\smallskip
c) $R$ is a noetherian strictly henselian integral local ring with  an algebraically closed residue field, and $H_{R^0}:=H_R\times R_0 $ is a reductive group scheme over $R^0$, where $R^0$ is the open subscheme of $\Spec(R)$ defined by the complement of the maximal ideal of $R$.

\smallskip
d) This allows us to pass from $O$ to its strict henselization $O^{\text{sh}}$ and so we can assume that $O$ is a strictly henselian DVR.

\smallskip
e) If moreover $K$ is of characteristic zero (so $O$ is an excellent ring), we have to deal only with excellent rings (as the set of excellent rings is stable under the operations performed in this remark).

\smallskip
f) If $K$ is of positive characteristic and if $O$ is a Nagata ring, we have to deal only with Nagata noetherian rings (as the set of such rings is stable under the operations performed in this remark, cf. [Ma, Ch. 12]).

\smallskip
5) If the family of tensors $(s_{\alpha})_{\alpha\in\scrJ}$ is essentially finite, if the extra conditions needed to get c) to f) above (i.e. the assumptions of 4.3.6 3) pertaining to iii) of 4.3.6 3)) are satisfied, and if $K$ is of characteristic zero (so $O$ is an excellent ring), then, for checking $4.3.5$, we can assume that $R$ is an integral noetherian complete local ring having an algebraically closed residue field. In other words we can replace $R$ (with $R$ the localization of an  integral finitely generated $O$-algebra with respect to a prime dominating the maximal ideal of $O$) by its completion $\widehat{R}$. Argument: $\widehat{R}$ is a reduced ring (as $R$ is an excellent
ring) and so the statement follows from iii) of 4.3.6 3).
So we can replace $O$ with the completion of $O^{\text{sh}}$ (cf. also to d) of 4)), i.e. we can assume that $O$ is a strictly henselian complete DVR.

$5^\prime)$ If in 4) and 5) we work in a weakly (resp. strongly) $O$-well positioned context, we do not have to make any assumption on the existence of a good projector of $\text{End}(W)$ as part of the family of tensors. We get c) to f) and 5) above, but always assuming that we have normal integral domains (resp. reduced  rings) instead of integral rings.

6) All concrete families of tensors used  in this paper are essentially finite and fit in the strongly context. Any essentially finite family of tensors in spaces of the form $W^{\oplus m}\oplus W^{*\oplus n}$, with $m,n\in\dbN$, is of bounded degree, but the converse to this is not true.   

7) Any time we can replace $O$ with another DVR $O^1$ (faithfully flat over $O$), we can replace the family of tensors $(s_{\alpha})_{\alpha\in\scrJ}$ with the family of tensors $(s_{\alpha^1})_{\alpha\in\scrJ^1}$ (of the tensor algebra of $(W\oplus W^*)\otimes_O O_1$) formed by linear combinations (with coefficients in $O^1$) of the tensors of $(s_{\alpha})_{\alpha\in\scrJ}$. If the family of tensors $(s_{\alpha})_{\alpha\in\scrJ}$ is essentially finite, then the family of tensors $(s_{\alpha^1})_{\alpha\in\scrJ^1}$ is also essentially finite.

8) To check 4.3.5 for a noetherian ring $R$, we can assume (cf. 3)) that it is local, and that $H_{R^0}$ is a reductive group scheme over $R^0$, where $R^0$ is the open subscheme of $\Spec(R)$ defined by the complement of the maximal ideal of $R$, even if the family of tensors $(s_{\alpha})_{\alpha\in\scrJ}$ is not essentially finite. In 4.3.4 and 4.35 we could have worked with $M$ a projective (instead of free) $R$-module but this would have made no significant difference.    

9) The role of $O$ is mostly just to fix up notations. For the greatest part of 4.3.4-17 it can be replaced by any other integral noetherian scheme $Z$, and then the role of $R$ is replaced by an arbitrary  integral flat $Z$-scheme. We will not stop to state the results in this generality, as they can be immediately deduced from the ones stated. 
\finishproclaim

\Proclaim{4.3.8. Remark.} \rm
The tensors which give a lot of information about the modules enveloping them, are  projections and isomorphisms.
\finishproclaim

\Proclaim{4.3.9. Remark.} \rm
If $H_1$ is a reductive subgroup of $H$ with $H_1^{\der}=H^{\der}$, then any weakly $O$-well positioned family of tensors (with respect to $\psi$) for $H$ is also a weakly $O$-well positioned family of tensors (with respect to $\psi$) for $H_1$. This results easily from 3.1.6 and from the fact that the Zariski closure in a torus $T_R$ (over a normal ring $R$ as in 4.3.5) of a subtorus of the generic fibre of $T_R$, is a torus over $R$: this is a local statement for the \'etale topology of $\Spec(R)$; so we can assume that $T_R$  is split and then we can use characters of $T_R$. The same thing remains true for weakly $O$-very well positioned families of tensors.
\finishproclaim

\Proclaim{4.3.10. Proposition.}
With the notations of 4.3.2, if $W$ is a vector space over
$\dbQ$, then:

\smallskip
a) there is $N\in\dbN$, such that for
any prime $p$ not dividing $N$, the family of tensors formed by $\pi(\grg)$, $\qbar$, and by $p_i$, $r_i$, $s_i$ and $t_i$, $i\in I$, is strongly $\dbZ_{(p)}$-very well
positioned with respect to $\psi$ for the group $G$;

b) for any odd prime $p\Ge s(\grg_0,W)$, the family of three tensors formed by $\pi(\grg_0)$, $B$ and $B^*$ is strongly $\dbZ_{(p)}$-well positioned for the group $G_0$.
\finishproclaim

\proof
Let $L$ be a $\dbZ$-lattice in $W$ such that $\psi$ induces a perfect form $\psi\colon L\otimes L\to\dbZ$.
As the family of tensors of a) is finite, we deduce the existence
of a number $N\in\dbN$, such that for any prime $p$ not dividing $N$, $L\otimes
\dbZ_{(p)}$ envelopes the family of tensors of a) with respect to $\psi$. So a) follows once we show the strongly $\dbZ_{(p)}$-well positioned part.
We fix a prime  $p$ not dividing $N$, for the case a) respectively an odd prime $p\Ge s(\grg_0,W)$, for the case b). Let $R$ be a reduced faithfully flat $\dbZ_{(p)}$-algebra and let $S:=R\fracwithdelims[]1p$. Let $M$ be a free $R$-module, with $M\otimes_R S=W\otimes_{\dbQ} S$, enveloping the family of tensors of a) with respect to $\psi$, respectively enveloping the family of tensors of  b). We have to show that the Zariski closure $G(M)$ of $G_S$ in $GSp(M,\psi)$ in case a), and respectively that the Zariski closure $G_0(M)$ of ${G_0}_S$ in $GL(M)$ in case b), are reductive groups over $R$. We can assume that $R$ is a local reduced noetherian ring (cf. 4.3.7 5')).
Let $m$ be its maximal ideal.

\Proclaim{Case a).} \rm
Let $\grg(M):=(\grg\otimes S)\cap\grg\grl(M)$.
We have $\grg(M)=\pi(\grg)(\grg\grl(M))$ and so $\grg(M)$ is direct summand of $\gsp(M,\psi):=\Lie(GSp(M,\psi))$.
Let
$$
A:=((\gsp(M,\psi)/\grg(M))\otimes R/m)^{\grg(M)}.
$$
(the upper right index refers to the operation of taking the elements annihilated by $\grg(M)$).
\finishproclaim
 
\Proclaim{Claim 1.}
We have $A=\{0\}$.
\finishproclaim

This results from the following facts.

\smallskip
\item{(a)}
The fact that the family of tensors $(p_i,t_i,s_i,r_i)_{i\in I}$ is enveloped by $M$
implies that the trace form  on $\grg_i(M):=r_i(\grg\grl(M))$ associated
to its representation on $\grh_i(M):=p_i(\grg\grl(M))$ is
perfect.
So the Casimir element $k_i$ of this representation induces a linear map
$\grg\grl(M)\to\grg\grl(M)$.

\smallskip
\item{(b)}
The fact that $\qbar$ is enveloped by $M$ implies that the
linear combination of $k_i$ used in the formation of
$\qbar$, induces an endomorphism $q\colon\grg\grl(M)\to\grg\grl(M)$ such
that its restriction to $\grh(M):=\opoplus_{i\in I}\grh_i(M)$
is an isomorphism $\grh(M)\arrowsim \grh(M)$.

\smallskip
\item{(c)}
Any element of $A$ is annihilated by $\qbar$ (as $\qbar$ is the
endomorphism induced by a sum of Casimir elements).

\smallskip
As $A=\{0\}$, the Lie subalgebra of $\gsp(M/mM)$
centralizing the reduction of $\pi(\grg)$ modulo $m$, is
$\grg(M)/m\grg(M)$.
This implies that the scheme $G(M)$ has smooth fibres.

Moreover it is smooth in the $R$-valued point defining its origin. To check this, let $R_0(G(M))$ be the ring of the completion of $G(M)$ in the origin, and let $R[[\grg(M)]]$ be the ring of formal power series defined by the free $R$-module $\grg(M)$. We get naturally an epimorphism $i_0(R)\colon R[[\grg(M)]]\twoheadrightarrow R_0(G(M))$. If $R$ is integral, by reasons of dimension, we get that $i_0(R)$ is an isomorphism. As $R$ is reduced, this implies that $i_0(R)$ is an isomorphism: the  kernel of $i_0(R)$ is included in $\scrP[[\grg(M)]]$, for any prime ideal $\scrP$ of $R$ of codimension 0.

The fact that $\pi(\grg)$ is enveloped by $M$ implies that
the trace form on $\grg(M)/m\grg(M)$ is perfect and so the
Lie algebra of the nilpotent radical of the connected
component of the origin of $G(M)\times_R \Spec(R/m)$ is zero
(cf. [Bou1, p. 41]).
From this we deduce easily (cf. [SGA3, Vol. 3, p. 12] and [Ti, 3.8.1]) that the connected component of the origin of any fibre of $G(M)$ is a reductive group scheme. From 3.1.2.1 c) and [Hart, ex. 4.11 p. 107] we deduce that all fibres of $G(M)$ are connected. From this and the fact that $G(M)$ is smooth in the origin  we deduce that $G(M)$ is a smooth scheme over $R$. 

We conclude that $G(M)$ is a reductive group scheme over $R$
and so condition $4.3.5$ (for reduced rings) is satisfied.
This proves a).

\Proclaim{Case b).} \rm
 We can assume, cf. 4.3.7 $5^\prime)$, that $R$ is a  
noetherian excellent strictly henselian  local ring, that $R/m$ is an algebraically closed field, and that $G_0(M)_{R^0}$ is a semisimple group over the open subscheme  $R^0$ of $\Spec(R)$ defined by the complement of the maximal point $\Spec(R/m)$ of $\Spec(R)$. From the properties  implied by the excellence property we need just that $R$ is an $N$-1 ring (cf. def. of [Ma, 31.A]), i.e. that the normalization $R_n$ of $R$ in its ring of fractions is a finite $R$-module and so a noetherian ring.   
\finishproclaim

\Proclaim{Part 1.} \rm 
The integrality of $\pi(\grg_0)$ gives us a direct sum decomposition $\grg\grl(M)=\grg_0(M)\oplus\grg_0(M)^\bot$ and the integrality of $B$ and $B^*$ implies that the Killing form $b(M)$ on $\grg_0(M)$ is perfect. Let $Aut(\grg_0(M))$ be the group scheme (of finite type) over $R$ defined by the Lie algebra automorphisms of $\grg_0(M)$,  and let $G_0(M)^{\ad}$ be the connected component of the origin of $Aut(\grg_0(M))$, defined as the Zariski closure in $Aut(\grg_0(M))$ of the connected components of the origins of the  fibres of $Aut(\grg_0(M))$ over points of $\Spec(R)$ of codimension 0. 
\finishproclaim

\Proclaim{Claim 2.} 
The subscheme $G_0(M)^{\ad}$ is a subgroup of $Aut(\grg_0(M))$. It is an adjoint group over $R$ having $\grg_0(M)$ as its Lie algebra. 
\finishproclaim

\proof
We first remark that for any algebraically closed field $\kbar$ which is an $R$-algebra, $\Lie(Aut(\grg_0(M))_{\kbar})$ is the Lie algebra of differentiations of $\grg_0(M)\otimes\kbar$; the same argument --based on the fact  that  the Killing form of $\grg_0(M)\otimes\kbar$ is perfect-- as in the characteristic zero case, gives us $\Lie(Aut(\grg_0(M)_{\kbar}))=\grg_0\otimes\kbar$. So, by
reasons of dimension, the tangent space in the origin of $G_0(M)^{\ad}$ is also $\grg_0\otimes\kbar$. This implies that $Aut(\grg_0(M))$ is smooth (over $R$) in the origin (the argument for this is the same as the one used in Claim 1, in a similar situation) and that every fibre of it is a smooth group, which is the extension of an adjoint group by a finite \'etale group. The finite \'etale group corresponds to outer automorphisms of the Lie algebra of the semisimple part of this extension.

As $R$ is a strictly henselian ring, we deduce from the smoothness of $Aut(\grg_0(M))$ in the origin, by using translations, that $Aut(\grg_0(M))$ is smooth over $R$ in any point of the connected component of the origin of a fibre of it. All these points belong to $G_0(M)^{\ad}$, and by reasons of dimension, they are smooth points of $G_0(M)^{\ad}$.

But $G_0(M)^{\ad}$ has all its fibres connected: an inner automorphism of a semisimple Lie algebra can not specialize to an outer automorphism. To see this, we first remark that $\grg_0(M)$ is defined over a subring of $R$ which is finitely generated over $\dbZ$. So everything comes down to checking this in the case of a complete DVR, having an algebraically closed residue field. If $R$ is such a ring, then the open subscheme of $G_0(M)^{\ad}$ defined by putting together the connected components of the origins of its fibres is a semisimple group, and so everything results from 3.1.2.1 c).   

So $G_0(M)^{\ad}$ is a smooth subgroup of $Aut(\grg_0(M))$ and has connected fibres. So $G_0(M)^{\ad}$ is a semisimple group over $R$ (cf. the above statement on the fibres of $Aut(\grg_0(M))$). We have $\Lie(G_0(M)^{\ad})=\Lie(Aut(\grg_0(M))=\grg_0(M)$. The group $G_0(M)^{\ad}$ is adjoint as its fibres over points of $\Spec(R)$ of codimension 0 are so. This ends the proof of Claim 2.

\Proclaim{Part 2.} \rm
Let $\grg_0(M)=\grt\opoplus_{\alpha\in\Phi}\grg_{\alpha}$ be a Weyl direct sum decomposition of $\grg_0(M)$ with respect to a system of roots $\Phi$ associated to the Lie algebra $\grt$ of a maximal split torus $T^{\ad}$ of $G_0(M)^{\ad}$ ($T^{\ad}$ exists as $R$ is a strictly henselian local ring; warning below, despite notations, $T^{\ad}$ is a torus and not the adjoint group of the group $T$ to be introduced later). For any $\alpha\in\Phi$ let $\dbG_{a,\alpha}$ be the $\dbG_a$ subgroup of $G_0(M)^{\ad}$ having $\grg_{\alpha}$ as its Lie algebra. The inequality $p\Ge s(\grg_0,W)$ implies that for any $\alpha\in\Phi$, every $x\in\grg_{\alpha}$, as an endomorphism of $M$, satisfies $x^p=0$. Let  $\alpha$ be an arbitrary element of $\Phi$. Let $V(\grg_{\alpha})$ be the affine scheme over $R$ defined by the $R$-module $\grg_{\alpha}$ (for an $R$-algebra $R_1$, $V(\grg_{\alpha})(R_1)=\grg_{\alpha}\otimes R_1$). There is a natural identification $V(\grg_{\alpha})=\dbG_{a,\alpha}$.

The homomorphism $exp\colon V(\grg_{\alpha})\to GL(M)$, defined on an $R$-valued point $x\in\grg_{\alpha}$ by $exp(x)=\sum_{i=0}^{p-1}\frac {x^i}{i!}$ (the above sum is an isomorphism of $M$ as $x$ is a nilpotent endomorphism of $M$), is an isomorphism onto its image $\dbG_{a,\alpha}$: at the Lie algebra level we get an isomorphism $\Lie(V(\grg_{\alpha}))\arrowsim\grg_{\alpha}$. We deduce that $\dbG_{a,\alpha}(R)\subset GL(M)(R)$
and so the groups $\dbG_{a,\alpha}$ can be considered as subgroups of $GL(M)$.

We treat first the special case when $R$ is a complete $DVR$ with an algebraically closed  residue field. Let ${G_0^{\sc}}_S$ be the simply connected semisimple group cover of ${G_0(M)}_S^{\ad}$. Using [Ti, 3.1.1] we get that the  subgroup of ${G_0^{\sc}}_S$$(S)$ generated by the subgroups $\dbG_{a,\alpha}(R)$ is hyperspecial. It is mapped under the composite homomorphism ${G_0^{\sc}}_S$$\to {G_0(M)}_S\to GL(M)_S$ into $GL(M)(R)$.
From 3.1.2.1 a) and c) we deduce that $G_0(M)$ is a reductive (and so semisimple) subgroup of $GL(M)$.

\smallskip
We come back to the general case. As in the special case we ``recovered" $G_0(M)$ from $\grg_0(M)$, we get directly that:

\smallskip
d) The reduced subscheme of the connected component of the origin of any  fibre of $G_0(M)$ is a semisimple  group scheme.

\Proclaim{Claim 3.} There is a subtorus $T$ of $GL(M)$ having $\grt$ as its Lie algebra.
\finishproclaim

\proof
First we remark that $T_S$ is well defined (it is the inverse image of $T^{\ad}_S$ under the natural homomorphism $G_0(M)_S\to G_0(M)^{\ad}_S$). So $T_S$ is a split torus. Let $C$ be the set of characters of $T_S$ through which it acts on $M\otimes S$. We consider the direct sum decomposition $M\otimes S=\oplus_{\gamma\in C} M_S^{\gamma}$ associated to the faithful representation $T_S\hookrightarrow GL(M\otimes S)$. So $T_S$ acts on $M_S^{\gamma}$ through the character $\gamma$. We need to show that the above direct sum decomposition of $M\otimes S$  extends to a direct sum decomposition of $M$, i.e. that the natural $R$-linear map 
$$
i_T\colon\oplus_{\gamma\in C} M^{\gamma}\to M, 
$$
with $M^{\gamma}:=M\cap M_S^{\gamma}$, is an isomorphism. 

To see this, let $\scrB(\Phi)$ be a basis of roots of $\Phi$. Let $\alpha\in\scrB(\Phi)$. Let $\grs\grl_2(\alpha)$ be the Lie subalgebra of $\grg_0(M)$ generated by $\grg_{\alpha}$ and $\grg_{-\alpha}$. As $\grg_0(M)$ is the Lie algebra of the adjoint group $G_0(M)^{\ad}$, and as $p>2$, we deduce that it is an $\grs\grl_2$ Lie algebra over $R$; so the notation is justified. As an $R$-module, it is isomorphic to $R^3$. We choose a standard basis $\{h_{\alpha},x_{\alpha},y_{\alpha}\}$ of it. So $x_{\alpha}\in\grg_{\alpha}$, $y_{\alpha}\in\grg_{-\alpha}$, $h_{\alpha}\in [\grg_{\alpha},\grg_{-\alpha}]$, and the formulas $h_{\alpha}=[x_{\alpha},y_{\alpha}]$, $[h_{\alpha},x_{\alpha}]=2x_{\alpha}$ and $[h_{\alpha},y_{\alpha}]=-2y_{\alpha}$ hold. The element $h_{\alpha}$ is a semisimple  element of $\grt$. Over $S$ it generates the Lie algebra of a subtorus $T_{S\alpha}$ of $GL(M\otimes S)$. It is a split torus, as it is a subtorus of the split torus $T_S$.

The key fact is: as $p\Ge s(\grg_0,W)$, we deduce that the eigenvalues of $h_{\alpha}$, as a semisimple endomorphism of $M$, are integers in the set $A(\alpha):=\{-p+1,-p+2,...,p-1\}$. For any $i\in A(\alpha)$ let $M(i)$ be the $R$-submodule of $M$ formed by elements on which $h_{\alpha}$ acts as the multiplication with $i$. So if any two such integers are not congruent mod $p$ (and so they are not congruent modulo $m$) (this is the case if $p>2s(\grg_0,W)$) then $M=\oplus_{i\in A(\alpha)} M(i)$. To see that this remains true even when two distinct eigenvalues are congruent mod $p$ we have to use $x_{\alpha}$ and $y_{\alpha}$. 

We need to show that for any $i\in\{1,...,p-1\}$, if $v(p-i)\in M(p-i)$ and $v(-i)\in M(-i)$ are such that 
$$v(p-i)+v(-i)\in mM,\leqno (1)$$
then $v(p-i)\in mM(p-i)$ and $v(-i)\in mM(-i)$. We can assume that $p-i\Ge i$ and we will prove the statement by induction on $i\in\{1,...,\frac{p-1}2\}$. 

So let us first treat the case when $i=1$. Applying first $x_{\alpha}$ to the relation $(1)$ a couple of times $v(-1)$ gets annihilated. Applying then $y_{\alpha}$ to the result the same number of times (to bring the things back)  we get something which is a multiple of $v(p-1)$ by an integer which is non-zero mod $p$. But what we get is in $mM$. In fact it is in $mM(p-1)$: we get this by  applying first $p-1$ times $x_{\alpha}$ to $(1)$  and then applying (backwards) $p-1$ times  $y_{\alpha}$ to $(1)$. So $v(p-1)\in mM(p-1)$. Similarly we get that $v(-1)\in mM(-1)$. We deduce that $M(p-1)$ and $M(-1)$ are direct summands of $M$, and so they are free ($R$ being a local ring). For $j\in\{1,...,p-1\}$ let $M_j(p-1):=x_{\alpha}^j(M(p-1))$. It is a submodule of $M(p-1-2j)$. Let $M_0(p-1):=M(p-1)$ and let 
$$\Mtil(p-1):=\oplus_{j=0}^{p-1}   M_j(p-1).$$
Using the fact that $p$ is greater then all eigenvalues of the endomorphism $h_{\alpha}$ of $M$ we deduce that $y_{\alpha}(M_i(p-1))=M_{i-1}(p-1)$, $\forall i\in\{1,...,p-1\}$. This implies that $\Mtil(p-1)$ is a direct summand of $M$, and so a free $R$-module. 

To proceed further on we just have to repeat everything for $i=2$ and for the quotient $\grs\grl_2(\alpha)$-module $M/\Mtil(p-1)$. Then we repeat everything for $i=3$ and the new  $\grs\grl_2(\alpha)$-module which is the  quotient of $M/\Mtil(p-1)$ (by a similarly constructed $\Mtil(p-2)$ submodule), etc. The induction becomes obvious.

We conclude that $M$ is a direct sum of submodules on which $h_{\alpha}$ acts diagonally. 
 This implies that $T_{S\alpha}$ extends to a subtorus $T_{\alpha}$ of $GL(M)$. 

Let $\tilde T:=\prod_{\alpha\in\scrB(\Phi)} T_{\alpha}$. As the subtori $T_{\alpha}$ of $GL(M)$, $\alpha\in\scrB(\Phi)$, commute, we get a group homomorphism $i_{\tilde T}\colon\tilde T\to GL(M)$, obtained by taking the product of homomorphisms $T_{\alpha}\hookrightarrow GL(M)$. Over $S$, $i_{\tilde T}$ factors through $T_S$. Let $T$ be the quotient of $\tilde T$ by the finite flat group subscheme (over $R$) of $\tilde T$, which over $S$ is the kernel of the factorization $\tilde T_S\to T_S$; this finite flat group scheme is the kernel of $i_{\tilde T}$. The notation is justified, i.e. the fibre of $T$ over $S$ is indeed the torus $T_S$ we previously considered. We get a homomorphism $T\to GL(M)$. The torus $T$ over $R$ is split as $R=R^{\text{sh}}$. 

The group of characters of $T$ is the same as the group of characters of $T_S$. So $T$ acts on $M$ through the characters $\gamma\in C$, achieving a direct sum decomposition of $M$ on submodules on which it acts diagonally through the characters of $C$. This proves that $i_T$ is an isomorphism and that $T$ is a subtorus of $GL(M)$. This ends the proof of Claim 3.
   
\smallskip
Let now 
$$
U(M):=T\times \prod_{\alpha\in\Phi} \dbG_{a,\alpha}.
$$ 
Let $u_M\colon U(M)\to GL(M)$ be the morphism defined by taking the product (as in [SGA3, Vol. 3, p. 172]) of the inclusions of the factors of $U(M)$ in $GL(M)$. It factors through $G_0(M)$. We have:

\smallskip
e) $u_M$ is injective on points with values in fields (i.e. it is radicial); 

\smallskip
f) In any $R/m$-valued point of the group scheme $U(M)$, $u_M$ induces an injection at the level of tangent spaces, producing a surjection at the level of cotangent spaces;

\smallskip
g)  At the level of completions of local rings  (defined by an $R/m$-valued point of $U(M)$), $u_M$ induces an epimorphism.

\smallskip
Part e) is a direct consequence of d) (cf. [Bo, 14.14] and the particular case). As over $R/m$ the morphism $u_M$ is a locally closed immersion (cf. the special case), using translates ($U(M)$ being smooth over $R$), it is enough to check part f) in the origin of $U(M)$. But in this case it results from the fact that   the tangent space of $U(M)$ in the origin $\Spec(R)\hookrightarrow U(M)$ is $\grg_0(M)$ (cf. the definition of the factors of $U(M)$ and of the expression of their Lie algebras), and from the fact that $\grg_0(M)$ is a direct summand of $\grg\grl(M)$ (as $\pi(\grg_0$) is enveloped by $M$). Part g) is a direct consequence of f) and of the fact that $R/m$ is an algebraically closed field.

We consider the simply connected group cover $G_0(M)^{\sc}$ of $G_0(M)^{\ad}$. $T$ fixes $\pi(\grg_0)$, as $T_S$ does. So $T$ acts under the adjoint representation on $\grg_0(M)$. We get a homomorphism $m_T\colon T\to Aut(\grg_0(M))$. As $Aut(\grg_0(M))$ is a subgroup of $GL(\grg_0(M))$, the kernel of $m_T$ is the same as the kernel of the representation of $T$ on $\grg_0(M)$. But any linear representation of a split torus (over $R$) is a direct sum of irreducible one-dimensional representations (associated to characters). So $\ker(m_T)$ is a finite flat group scheme over $R$. The quotient of $T$ by it is a subtorus of $GL(\grg_0(M))$, and so a subtorus  of $G_0(M)^{\ad}$: any torus over a field is a geometrically connected variety. This subtorus of $G_0(M)^{\ad}$ is nothing else but the subtorus $T^{\ad}$ we considered in the first paragraph of Part 2. 

The inverse image of $T^{\ad}$ under the natural isogeny $G_0(M)^{\sc}\to G_0(M)^{\ad}$ is a maximal torus $T^{\sc}$ of $G_0(M)^{\sc}$. We get an isogeny $i^{\sc}\colon T^{\sc}\to T^{\ad}$ of split tori over $R$. Its kernel is the center of $G_0(M)^{\sc}$. Moreover $i^{\sc}$  factors through $T$, as this happens over $S$. We get another isogeny of split tori $T^{\sc}\to T$. Let $C_T$ be its kernel. It is a finite flat group scheme over $R$, contained in the center of $G_0(M)^{\sc}$. 
Let $\tilde G_0(M)$ be the semisimple group over $R$ which is the quotient of $G_0(M)^{\sc}$ by $C_T$. From the very construction of $C_T$ we get that $\tilde G_0(M)_S$ is $G_0(M)_S$.

We want to show that $G_0(M)$ is $\tilde G_0(M)$. We have a morphism $l_0\colon\tilde G_0(M)_S\to G_0(M)$. We view it as a rational map from $\tilde G_0(M)$ to $G_0(M)$. We also view it, keeping the same notation,  as a rational map from $\tilde G_0(M)$ to $GL(M)$.

We have a canonical homomorphism $\tilde G_0(M)\to G_0(M)^{\ad}$.
The $R$-scheme $U(M)$ is an open subscheme of $\tilde G_0(M)$: each factor of $U(M)$ (i.e. $T$ and each  $\dbG_{a,\alpha}$, $\alpha\in\Phi$) are subgroups of $\tilde G_0(M)$. This is obvious for $\dbG_{a,\alpha}$, i.e. the subgroup $\dbG_a$ of $\tilde G_0(M)$, corresponding to an element $\alpha\in\Phi$, is mapped isomorphically into the subgroup $\dbG_{a,\alpha}$ of $G_0(M)^{\ad}$ (we are dealing only with central isogenies). For $T$ this is obvious from its construction. So we can apply [SGA3, Vol. 3, p. 172]; we get: 

\smallskip
h) The rational map $l_0$ is defined in codimension 1. 

\smallskip
We  first assume that $R$ is a normal ring, i.e. that $R=R_n$. From h) and from [BLR, Th. 1 of 4.4] we deduce that $l_0$ can be extended to a morphism $l_1\colon\tilde G_0(M)\to GL(M)$. The morphism $l_1$ is a group homomorphism, as $\tilde G_0(M)$ is a smooth scheme over $R$ and as the fibre of $l_1$ over $S$ is a group homomorphism. From the special case we deduce that all fibres of $l_1$ are closed immersions. But $l_1$ is proper (as its fibre over $S$ is proper, this results from the valuative criterion of properness, cf. 3.1.2.1 c)) and so it is a finite morphism. From Nakayama's Lemma we deduce  that $l_1$ is a closed immersion and so $\tilde G_0(M)=G_0(M)$. This ends the proof in the case $R=R_n$. 

We would like to point out  that if $R$ is as in the special case (i.e. it is a complete DVR with an algebraically closed field), from 3.1.2.1 c) we get directly that $l_1$ is a closed embedding. This represents a second proof of the special case without a reference to [Ti, 3.1.1] but based on the elementary result [BLR, Th. 1 of 4.4]: the facts e) to g) above, obtained based on d) above, were not used to get h). 

We now come back to the general case (i.e. we do not assume anymore that $R=R_n$). From the fact that the result is known for $R_n$ and from the fact that $R_n$ is a finite $R$-module, we deduce the existence of a finite morphism $\tilde G_0(M)_{R_n}\to GL(M)$. It factors through $G_0(M)$, producing a finite dominant morphism $\tilde G_0(M)_{R_n}\to G_0(M)$. We deduce that:

\smallskip
i)  The reduced scheme defined by $G_0(M)_{R/m}$ is a semisimple group having as its Lie algebra $\grg_0(M)\otimes R/m$. 

\smallskip
This implies that the localization of $u_M$ in the $R$-valued point defining the origin of $GL(M)$, is a finite morphism. From g) we deduce that it is a closed embedding. This implies that around this origin, $u_M$ is a closed embedding. We deduce that $G_0(M)$ is smooth in the origin. As any $R/m$-valued point of $U(M)$ has a lift to $R$ (as $R$ is a strictly henselian ring), using translations with $R$-valued points of $G_0(M)$, we deduce from i) above, that $G_0(M)$ is smooth in all its $R/m$-valued points. As $G_0(M)_{R^0}$ is smooth over $R^0$, we deduce that $G_0(M)$ is a smooth scheme over $R$, and so it is a subgroup of $GL(M)$. From the fact that  $G_0(M)_{R^0}$ is a semisimple group over $R^0$ and from i) above, we get that $G_0(M)$ is a semisimple group over $R$. This ends the proof of the case b) and so of the Proposition.

\Proclaim{4.3.10.1. Remarks.} \rm
1) Proposition 4.3.10 b) remains true if instead of $\dbQ$ (and $\dbZ_{(p)}$) we work with an arbitrary field $K$ of characteristic zero, which is the field of fractions  of a DVR $O$ of mixed characteristic (and with $O$), and if, instead of $\pi(\grg_0)$, we work with any other projector $\pi_0$ of $\grg\grl(W)$ on $\grg_0$
centralized by $\grg_0$ (the role of $\pi(\grg_0)$ was just to produce a direct sum decomposition $\grg\grl(W)=\grg_0\oplus\grg_0^{\bot}$). 

Moreover, the condition $p>2$ is not needed: If $p=2=s(\grg_0,W)$ then 4.3.10 b) remains true as it can be easily checked. Of course in the majority of cases for $p=2$ we get a non-perfect Killing form on $\grg_0(M)$.${}^1$ $\vfootnote{1}{Unfortunately, one can check that for $p=2$ we never have a perfect Killing form.}$  However: 

$1^{\prime})$ Part 2 of the proof of 4.3.10 b) is a result independent of Part 1 (we just needed that there is an adjoint group over $R$ whose Lie algebra is $\grg_0(M)$). It is a result on representations of the Lie algebra of an adjoint group and so it remains true even if the Killing form (or the trace form) on $\grg_0(M)$ is not perfect.

Part 1 of the proof of 4.3.10 a) is a result on the existence of adjoint groups having a prescribed Lie algebra which is subject to the condition that its Killing form is perfect.

2) Proposition 4.3.10 a) remains true if instead of $\dbZ$ we work with any other Dedekind domain $D$ of characteristic zero having an infinite number of maximal ideals (the number $N$ being replaced by a non-zero ideal of $D$).

3) Proposition 4.3.10 admits versions in positive characteristic. Of course, some precautions have to be taken. For instance the restriction of the trace form on $\grg\grl(W)$ to $\grg_0$ (or $\grg$) might not be perfect. Concentrating just on 4.3.10 b) we can state:
\finishproclaim

\Proclaim{4.3.10.2.}
 With the notations of 4.3.2, we assume that there is a projection $\pi_0$ of $\grg\grl(W)$ on $\grg_0$, annihilated by $\grg_0$, and that the Killing form on $\grg_0$ is perfect. If $s(\grg_0,W)$ is not greater than the characteristic $p$ of the residue field of $O$, if this residue field is perfect, and if $p>2$, then the family of tensors formed by
$\pi_0$, $B$ and $B^*$ (as the Killing form on $\grg_0$ is perfect, we can define $B$ and $B^*$ as in 4.3.2) is strongly $O$-well positioned for the group $G_0$.
\finishproclaim

The proof of this is entirely analogous to the proof of 4.3.10 b). We just have to check --it is easy-- that the condition $s(\grg_0,W)\Le p$ can be used in the same manner as in the proof of 4.3.10 b) (instead of e) of 4.3.7 4) we have to use d) and f) of 4.3.7 4)). Also, it can be easily checked that the condition on the residue field being perfect is not needed.
\finishproclaim
 
\Proclaim{4.3.10.3.} \rm
The family of tensors of 4.3.10 a) is not so suited for explicit computations, while the one of 4.3.10 b) is. The advantage offered by the family of tensors of 4.3.10 a) is: it cuts out of $\gsp(M,\psi)$ the Lie algebra of the group $G(M)$ (cf. Claim 1 of 4.3.10) by using only one tensor $\pi(\grg)$. However we do not use it in the rest of the paper. There are variants of 4.3.10 a) when $GSp(W,\psi)$ is replaced by another reductive subgroup of $GL(W)$.
\finishproclaim  

\Proclaim{4.3.11. Example.} \rm
We consider the case of Shimura varieties of PEL type, to emphasize that the (incipient) idea of using $\dbZ_{(p)}$-very well positioned families of tensors goes back to [LR]. We use the situation and notations used in [Ko, Ch. 5]. For simplicity we denote the nondegenerate $\dbQ$--valued alternating form on $V$ by $\psi$. If $p=2$, we assume that $4$ does not divide $\dim_{\dbQ}(V)$; this can be weaken, but the fast argument we present for the following Claim does require this assumption.

\Proclaim{Claim.} 
The elements of $\scrO_B$ form a family of tensors which is strongly $\dbZ_{(p)}$-very well positioned with respect to $\psi$ for the group $G$. 
\finishproclaim

\proof
The conditions imposed on B imply that the group scheme over $\dbZ_{(p)}$ defined by invertible elements of $\scrO_B$ is reductive.

We get that the group scheme $C$ over $\dbZ_{(p)}$ defined as the centralizer of $\scrO_B$ in $GL(L)$ is reductive. This is a property of linear representations of semisimple algebras over discrete valuation rings of mixed characteristic. In our case, passing from $\dbZ_{(p)}$ to $W(\dbF)$, $\scrO_B\otimes W(\dbF)$ is a finite product of algebras of the form $\text{End}(N)$, with $N$ a finite free $W(\dbF)$-module. So, inside $V\otimes_{\dbQ} W(\dbF)\fracwithdelims[]1p$ we can find a $W(\dbF)$-lattice $M$ such that $M$ is a direct sum of irreducible representations of $\scrO_B\otimes W(\dbF)$ (and so $M/pM$ is a direct sum of irreducible representations of $\scrO_B\otimes\dbF$). Using the fact that the determinants (as defined in loc. cit.) of $\scrO_B$ with respect to $L\otimes W(\dbF)$ and with respect to $M$ are the same, we deduce that the two representations of $\scrO_B\otimes\dbF$ on $M/pM$ and on $L\otimes\dbF$ are isomorphic. So $C$ is indeed a reductive group scheme over $\dbZ_{(p)}$. It is defined by invertible elements of a $\dbZ_{(p)}$-order of a semisimple $\dbQ$--algebra. 

Moreover there is $n\in\dbN$ such that $\frac1n$ times the bilinear form $\grb$ on $\Lie(C)$ induced by the trace form $Tr$ on $\grg\grl(L)$ is perfect. This can be read out from the end of [Sh, 2.1]. For instance, with the terminology and notations of the loc. cit., we can take $n=\frac m2$ if $L$ is of type $I$, $II$ or $III$, etc. Here we use that $B$ is a simple $\dbQ$--algebra. 

We assume first that $2p$ does not divide $\dim_{\dbQ}(V)$. This is equivalent to: the standard trace form on $\Lie(GSp(L,\psi))$ is perfect.

The fact that $\scrO_B$ is self dual with respect to $\psi$ implies that $\Lie(C)=\grc\oplus\grc^\bot$, with $\grc:=\Lie(C)\cap\Lie(GSp(L,\psi))$ and with $\grc^\bot:=\Lie(C)\cap\Lie(GSp(L,\psi))^\bot$ (here $\Lie(GSp(L,\psi))^\bot$ refers to perpendicularity with respect to the trace form). So the Zariski closure $G_{\dbZ_{(p)}}$  of the connected component $G$ of the origin  of the intersection of $C_{\dbQ}$ with $GSp(V,\psi)$ in $GSp(L,\psi)$ is a reductive group scheme over $\dbZ_{(p)}$. 

To see this  let $\scrG$ be the connected component of the origin of the special fibre of $G_{\dbZ_{(p)}}$. The above direct sum decomposition of $\Lie(C)$ implies that the dimension of $\Lie(\scrG)$ is equal to the dimension of $\scrG$, and so $\scrG$ is a smooth group over $\dbF_p$. So $G_{\dbZ_{(p)}}$ is smooth over $\dbZ_{(p)}$ in the points of $\scrG$. From the fact that $\frac1n\grb$  is a perfect form on $\Lie(C)$, we deduce that $\scrG$  is a  reductive group over $\dbF_p$. This results from the fact that the Lie algebra $\grn$ of the nilpotent radical of $\scrG$ is zero as it is included in the null space of the restriction mod $p$ of the symmetric bilinear form $\frac1n\grb$: [Bou1, p. 41] implies that $\grn$ is perpendicular to $\grc$, while from the definition of $\grb$ and $\grc^\bot$ we get that $\grn$ is perpendicular to $\grc^\bot$; here perpendicularity is with respect to $\frac1n\grb$. From 3.1.2.1 c) we deduce easily that $G_{\dbZ_{(p)}}$ is a reductive group over $\dbZ_{(p)}$ (i.e. its special fibre is $\scrG$). 

Using 3.1.2.1 a) and c) and the determinant condition of [Ko, Ch. 5], the same things remain true if we work with an arbitrary reduced ring $R$ which is faithfully flat over $\dbZ_{(p)}$, with a free $R$-module $M$ which satisfies $M\fracwithdelims[]1p=W\otimes R$ and envelopes the elements of $\scrO_B$ with respect to $\psi$ (i.e. with the same arguments we get that the group scheme over $R$ defined by the invertible elements of $\scrB\otimes_{\dbZ_{(p)}} R$ is reductive, that its centralizer in $GL(M)$ is a reductive group scheme $C_R$ over $R$, and that $G_{R\fracwithdelims[]1p}$ extends to a reductive subgroup $G_R$ of $GL(M)$). 

If $2p$ divides $\dim_{\dbQ}(V)$, then we have to work all the above with $Sp(L,\psi)$ instead of $GSp(L,\psi)$; this is supported by 3.1.6. This ends the proof of the Claim.
\finishproclaim

\Proclaim{4.3.12. Remark.} \rm
We start with an injective map $(G,X)\hookrightarrow (G\Sp(W,\psi),S)$. Let $B$ we the subalgebra of $\text{End}(W)$ formed by elements fixed by $G$;  it is a semisimple $\dbQ$--algebra. The connected component $G_1$ of the origin of the subgroup of $GSp(W,\psi)$ fixing B contains $G$ and we get an injective map $(G,X)\hookrightarrow (G_1,X_1)$
(with $X_1$ determined naturally by $X$). The pair $(G_1,X_1)$ defines a Shimura variety (it is easy to see that the axiom SV3 of 2.3 is satisfied) of PEL type (cf. their def.; see [Mi4, p. 161]). We call it the PEL-envelope of $(G,X)$ with respect to the injective map $f$. The tensors of degree 2 does not allow us to distinguish $(G,X)$ from $(G_1,X_1)$. So we were forced in 4.3.10 to use tensors of degree 4, to be able to conclude that the Zariski closure of $G$ in a $GL(L_{(p)})$ (for  some particular $\dbZ_{(p)}$-lattices $L_{(p)}$ of $W$) is a reductive group scheme over $\dbZ_{(p)}$.
\finishproclaim

\Proclaim{4.3.13. The case of a torus.} \rm
We consider a situation of the form $T\hookrightarrow GL(M)$ with $T$ the maximal torus of the center of a reductive subgroup $G$ of $GL(M)$, and with $M$ a free module of finite rank over a discrete valuation ring $O$. Let $B$ be the subalgebra of $\text{End}(M)$, formed by endomorphisms fixed by $G$. Then $B$ forms a family of tensors which is strongly $O$-very well positioned for $T$. To see this we can assume (cf. 4.3.7 $5^\prime$)) that $O$ is a strictly henselian DVR. Then $T$ is a split torus and $M=\oplus_{\alpha\in\scrJ} M_{\alpha}$, with $\scrJ$ a set of characters of $T$, and with $T$ acting on $M_{\alpha}$ through the character $\alpha$ ($\forall\alpha\in\scrJ$). Now the subfamily of $B$ (cf. 4.3.7 7)) formed by the projections of $M$ on  $M_{\alpha}$ ($\alpha\in \scrJ$) (they are fixed by $G$) associated to the above direct sum decomposition, is obviously strongly $O$-very well positioned for $T$.  
\finishproclaim

\Proclaim{4.3.14. Remark.} \rm
Let $G_0\hookrightarrow G\hookrightarrow G\Sp(W,\psi)$ be monomorphisms between reductive groups over $\dbQ$. Let $p$ be a rational prime. We assume the existence of a family of tensors $(v_{\alpha})_{\alpha\in\scrJ_0}$  in spaces of the form $W^{\otimes m}\otimes W^{*\otimes n}$ which is $\dbZ_{(p)}$-very well positioned with respect to $\psi$ for the group $G$. We assume also the existence of a $\dbZ_{(p)}$-lattice $L$ of $W$ enveloping the above family of tensors with respect to $\psi$ and such that there is a torus $T$ of the Zariski closure $G_{\dbZ_{(p)}}$ of $G$ in $G\Sp(L,\psi)$ having as its centralizer in $G_{\dbZ_{(p)}}$, a reductive group having the Zariski closure of $G_0$ in $G\Sp(L,\psi)$ as a reductive subgroup containing its derived subgroup. 

\Proclaim{Fact.}
The family of tensors $(v_{\alpha})_{\alpha\in\scrJ_0}$ can be enlarged (by adding only tensors of degree 2) to a family of tensors $(v_{\alpha})_{\alpha\in\scrJ}$, with $\scrJ\supset\scrJ_0$, which is $\dbZ_{(p)}$-very well positioned with respect to $\psi$ for the group $G_0$. 
\finishproclaim

This is a direct consequence of 4.3.13 (cf. also 4.3.9). A similar result can be stated for strongly or weakly $\dbZ_{(p)}$-very well positioned families of tensors. 
\finishproclaim

\Proclaim{4.3.15. Remark.} \rm
Let $O$ be a DVR and let $(M,\psi)$ be a symplectic space over its field of fractions $K$. Let $O_1$ be a DVR which is an \'etale cover of $O$ and let $K_1$ be its field of fractions. Let $G$ be a reductive subgroup of $G\Sp(M,\psi)$. If there is a family of tensors $(s_{\alpha})_{\alpha\in\scrJ_1}$ in spaces of the form $M^{\otimes m}\otimes M^{*\otimes n}\otimes O_1$ which is strongly (resp. weakly) $O_1$-very well positioned for the group $G_{K_1}$, and if there is an $O$-lattice $L$ of $M$ such that $L\otimes O_1$ envelopes the above family of tensors with respect to $\psi$, then there is a family of tensors $(w_{\beta})_{\beta\in\scrJ}$ of degree not bigger than the maximal degree of the tensors of $(s_{\alpha})_{\alpha\in\scrJ_1}$, situated in spaces of the form $M^{\otimes m}\otimes M^{*\otimes n}$, which is enveloped by $L$ and strongly (resp. weakly) $O$-very well positioned  with respect to $\psi$ for $G$. This is so due to the fact that the tensors of $L^{\otimes m}\otimes L^{*\otimes n}\otimes O_1$ fixed by the reductive group $G_{O_1}$, the Zariski closure of $G_{K_1}$ in $G\Sp(L\otimes O_1,\psi)$, are linear combinations with coefficients in $O_1$ of tensors of $L^{\otimes m}\otimes L^{*\otimes n}$ fixed by the reductive group $G_O$, the Zariski closure of $G$ in $G\Sp(L,\psi)$. We can now take the family of tensors $(w_{\beta})_{\beta\in\scrJ}$ showing up in such linear combinations of tensors of $(s_{\alpha})_{\alpha\in\scrJ_1}$. The same thing remains true when we do not work in the relative context (i.e. when we replace $GSp(M,\psi)$ by $GL(M)$ and there is no alternating form $\psi$ on $M$).  
\finishproclaim

\Proclaim{4.3.15.1.} \rm
Remark 4.3.15 remains true if instead of $O_1$ we work with the completion of $O$ (the argument is the same). 
\finishproclaim

\Proclaim{4.3.16. The relative PEL situation.} \rm
Let $O$ be a discrete valuation ring of mixed characteristic and let $K$ be its field of fractions. Let $M$ be a free module of finite rank over $O$. Let $G$ be a reductive subgroup of $GL(M)$ and let $L\subset\text{End}(M)$ be a semisimple algebra over $O$. So $L\otimes O^{\text{sh}}$ is a product of algebras of the form $\text{End}(P)$ with $P$ a free module over $O^{\text{sh}}$. We assume that the subgroup $C(L)$ of $GL(M)$ fixing $L$ is a reductive group over $O$ and that the connected component $G_0$ of the origin of $C(L)\cap G$ (defined as the Zariski closure in $G$ of the connected component the origin of the generic fibre of $C(L)\cap G$) is a reductive group over $O$, containing the maximal torus of the center of $G$. We assume that the bilinear  form on $\grg:=\Lie(G^{\der})$ induced by the trace form $Tr$ on $\text{End}(M)$ is perfect, and that $\pi(\grg)$ leaves invariant $\Lie(C(L))$. We also assume that one of the following two conditions is satisfied:

\smallskip
1) There is a torus $T$ of $G$ such that $G_0$ is contained in the centralizer $G^0$ of $T$ in $G$, $G_0^{\ab}=G^{0\ab}$ and the monomorphism $G_0^{\der}\hookrightarrow G^{0\der}$ becomes over $O^{\text{sh}}$ the diagonal embedding of $G_0^{\der}$ in a product of a finite number of copies of $G_0^{\der}$, which are permuted transitively (under conjugation) by the invertible elements of $L\otimes O^{\text{sh}}$;

2) A rational multiple of $Tr$ restricts to a perfect form on $\Lie(C(L))$. 

\smallskip
Let $(s_{\alpha})_{\alpha\in\scrJ}$ be a family of tensors of the tensor algebra of $M\oplus  M^*$ fixed by $G$, which is enveloped by $M$ and is $O$-very well positioned for $G$. Then the family of tensors formed by $(s_{\alpha})_{\alpha\in\scrJ}$, $\pi(\grg)$, and all the tensors of degree 2 fixed by $G_0$ and enveloped by $M$ (the elements of $L$ are examples of such tensors), is $O$-very well positioned for $G_0$. The proof of this presents no difficulty, being just an extended version of 4.3.11 and 4.3.14. The same remains true in a strongly or weakly context. 

We refer to the above situation as the relative PEL situation defined by the triple $(G,L,T)$ (resp. by the pair $(G,L)$) if condition 1) (resp. condition 2)) above is satisfied.
When 2) above is satisfied we get the relative PEL situation generalizing 4.3.11. We would like to remark that in 4.3.11 the tensor $\pi(\grg)$ is still present in disguise, cf. the connection between $\psi$ and $\scrO_B$ (see [Ko, Ch. 5]).
\finishproclaim

\Proclaim{4.3.17. Remark.} \rm
If in 4.3.4 we have $W=W_1\oplus W_2$ and $\psi=\psi_1\oplus\psi_2$ (with $(W_i,\psi_i)$ a symplectic space over $K$), if $H_i$ is a reductive subgroup of $GSp(W_i,\psi_i)$, and if $(s_{\alpha})_{\alpha\in\scrJ_i}$ is a family of tensors of the tensor algebra of $W_i\oplus W_i^*$ which is $O$-well positioned with respect to $\psi_i$ for $H_i$, $i=\overline{1,2}$, then the family of tensors $(s_{\alpha})_{\alpha\in\scrJ_1\cup\scrJ_2\cup\{1\}}$ (of the tensor algebra of $W\oplus W^*$; here $s_1$ is the projection of $W$ on $W_1$ having $W_2$ as its kernel) is $O$-well positioned with respect to $\psi$ for $H:=H_1\times H_2$. The same thing remains true for $O$-very well positioned families of tensors, or in a context without $\psi$, or in a strongly or weakly context, or if $H$ is replaced by a reductive subgroup $\tilde H$, obtained from $H$ in the same manner as we got $G_3$ from $G_1\times G_2$ in Example 3 of 2.5.

\bigskip
\noindent
{\boldsectionfont \S5. The basic result}

\bigskip
We present our procedure for proving the existence of integral canonical models of Shimura varieties of Hodge type.

\smallskip
\Proclaim{5.1. Theorem.}
Let $(G,X)$ define a Shimura variety of Hodge type and let $p>2$ be a prime such that $G$ is unramified over $\dbQ_p$.
We assume that the pair $(G,X)$ satisfies the following condition with respect to the prime $p$:

\smallskip
\parindent=35pt
\Item{\rm (*)}
There is an injective map $f\colon(G,X)\hookrightarrow(G\Sp(W,\psi),S)$ such that there is a  family of $G$-invariant tensors $(v_\alpha)_{\alpha\in\scrJ_0}$ in spaces of the form $(W\otimes W^*)^{\otimes m}$ (with $m\in\dbN$) and of degree not bigger than $2(p-2)$, which is
$\dbZ_{(p)}$-very well positioned with respect
to $\psi$ for the group $G$. 

\parindent=25pt
\smallskip\noindent
Then $\Sh_p(G,X)$ exists and has the EEP.
\finishproclaim

\proof
For the sake of clarity we divide the proof into steps.

\Proclaim{5.1.1. Step 0. Preliminaries.} \rm
Let $f\colon(G,X)\hookrightarrow(G\Sp(W,\psi),S)$ be an
injective map for which there is a family of $G$-invariant tensors in spaces of the form $(W\otimes W^*)^{\otimes m}$ and of degree not bigger than $2(p-2)$, which is $\dbZ_{(p)}$-very well positioned with respect to $\psi$ for the group $G$.
We fix such a family $(v_\alpha)_{\alpha\in\scrJ_0}$ and a prime $v$ of $E:=E(G,X)$ dividing $p$. 

Let $L$ be a $\dbZ$-lattice of $W$ such that $L_p:=L\otimes
\dbZ_{(p)}$ envelopes the family $(v_\alpha)_{\alpha\in\scrJ_0}$
and we have a perfect form $\psi\colon L\otimes L\to\dbZ$. This implies (cf. def. 4.3.4) that the Zariski closure $G_{\dbZ_{(p)}}$ of $G$ in $GSp(L_p,\psi)$ is a reductive group scheme over $\dbZ_{(p)}$. So the group $H:=\{g\in G(\dbQ_p)\mid g(L_p\otimes\dbZ_p)=L_p\otimes\dbZ_p\}$ is a hyperspecial subgroup of $G(\dbQ_p)$. Due to 3.2.7.1 it is enough to work with $(G,X,H,v)$.
Let $K_p:=\{g\in G\Sp(W,\psi)(\dbQ_p)\mid
g(L\otimes\dbZ_p)=L\otimes\dbZ_p\}$.
It is a hyperspecial subgroup of $G\Sp(W,\psi)(\dbQ_p)$.

The fact that $G$ is unramified over $\dbQ_p$ implies that $v$ is
unramified over $p$ [Mi3, 4.7].
Let $\dbF:=\overline{k(v)}$. Let $\scrM$ be the extension to $O_{(v)}$ of
the integral canonical model $\Sh_p(G\Sp(W,\psi),S)$ of $(G\Sp(W,\psi),S,K_p,p)$ (cf. 3.2.9).
Let $\scrN$ be the normalization of the Zariski closure of $\Sh_H(G,X)$ in
$\scrM$.
Let $V_0:=W(\dbF)$ and let $K_0$ be its field of fractions.
Let $\scrNbar:=\scrN_{V_0}$ and $\scrMbar:=\scrM_{V_0}$.

We claim that $\scrNbar$ is formally smooth over $V_0$.
For this it is enough to show that the completion of
the local ring of $\scrNbar$ in a point
$\Spec(\dbF)\to\scrNbar$ is $V_0[[X_1,\dotsc,X_d]]$,
with $d:=\dim\,X$. This is achieved at the end of Step 5 (of 5.5 below). 
\finishproclaim

\Proclaim{5.1.2. Step 1. The moduli setting.} \rm
We start with an arbitrary point 
$$
y\colon\Spec(\dbF)\hookrightarrow\scrNbar.
$$
From the definition of $\scrNbar$ we deduce (cf. 3.4.2) the
existence of a morphism 
$$
m_V\colon\Spec(V)\to\scrNbar
$$
lifting $y$, with
$V$ the normalization of $V_0$ in a finite field extension $K$ of
$K_0$.

Using the interpretation
of $\Sh_p(G\Sp(W,\psi),S)$ as a moduli scheme (working with the lattice
$L$) (see 3.2.9 and 4.1),  we get a universal principally polarized abelian scheme $(\scrA,\scrP_{\scrA})$ over $\Sh_p(G\Sp(W,\psi),S)$ (of relative dimension equal to  half the dimension of $W$ over $\dbQ$), having (compatibly) level-N symplectic similitude structure  for any $N\in\dbN$ satisfying $(p,N)=1$. Let $(\scrA_{\scrM},\scrP_{\scrM})$ and  $(\scrA_{\scrN},\scrP_{\scrN})$ be its pull backs to $\scrM$ and respectively to $\scrN$. 

The morphism $m_V$ gives birth to a
principally polarized abelian scheme $(A,p_A)$ over $V$, having (compatibly)
level-$N$ symplectic similitude structure for any $N\in\dbN$ satisfying $(p,N)=1$.
We fix an embedding $j\colon\Kbar\hookrightarrow\dbC$. We still denote by $j$ its restriction to $K$, $V$, $K_0$ or $V_0$.
The morphism $\Spec(\Kbar)\to\Sh_H(G,X)_{K_0}=\scrN_{K_0}$ induced
by $m_V$ can be lifted to a morphism
$u\colon\Spec(\Kbar)\to\Sh(G,X)_{K_0}$ such that the
point $x\in\Sh(G,X)(\dbC)$ induced from $u$ through the
inclusion $j$, is of the form $[h,a]$ with the $p$-component
of $a$ equal to one (i.e. $a\in G(\dbA_f^p)$ as 
we have $\dbA_f=\dbA_f^p\times\dbQ_p$).
This results from the fact that $G(\dbQ_p)=G(\dbQ)H$
[Mi3, 4.9].

The subgroup of $G\Sp(W,\psi)$ fixing
$(v_\alpha)_{\alpha\in\scrJ_0}$ might not be $G$.
Let $(v_\alpha)_{\alpha\in\scrJ}$, with
$\scrJ\supset\scrJ_0$, be an enlarged family of tensors
such that $G$ is the subgroup of $G\Sp(W,\psi)$ fixing them. If $\scrJ\setminus\scrJ_0$ is finite (we can assume this, cf. [De3, 3.1], but it is irrelevant for what follows), then the family of tensors $(v_\alpha)_{\alpha\in\scrJ}$ is essentially finite.

We think of $\Sh(G,X)_\dbC$ as the moduli scheme associated
to the injective map $f$, the lattice $L$ and the family of
tensors $(v_\alpha)_{\alpha\in\scrJ}$ (cf. 4.1).
Using d), e) and f) of 4.1 for the point $x$, we deduce that:

\smallskip
\item{a)}
the isogeny class of $A_\dbC$ is given by the pair $(W,h)$;

\smallskip
\item{b)}
$A_\dbC$ has a family $(t_\alpha)_{\alpha\in\scrJ}$ of Hodge
cycles, the Betti realization of $t_\alpha$ being
$v_\alpha$;

\smallskip
\item{c)}
the linear map $V_f(A_\dbC)=W\otimes\dbA_f@>{a^{-1}}>>
W\otimes\dbA_f$ induces a  similitude isomorphism
$(H_1(A_\dbC,\dbZ)\otimes\dbZhat,p_A)\arrowsim (L\otimes\dbZhat,\psi)$ ($p_A$ being the polarization of $A$).

\smallskip
As $a\in G(\dbA_f^p)$ we deduce that
$H_1(A_\dbC,\dbZ)\otimes\dbZ_p=L\otimes\dbZ_p$ (this identification is unique up to an isomorphism of $L\otimes\dbZ_p$ induced by an element of $G(\dbZ_{(p)}):=G(\dbQ)\cap H$) and that (under this identification)
$p_A=\beta_{(p)}\psi$, with $\beta_{(p)}\in\dbG_m(\dbZ_{(p)})$.
 Let
$$
H_{\text{\'et}}^1:=H_{\text{\'et}}^1(A_\dbC,\dbZ_p)=
H_{\text{\'et}}^1(A_{\Kbar},\dbZ_p);
$$ 
it is identified with $L^*\otimes\dbZ_p$.
So there is a  family of tensors
$(v_\alpha)_{\alpha\in\scrJ}$ in spaces of the form
$(H_{\text{\'et}}^1\otimes H_{\text{\'et}}^{1*})^{\otimes m}
\otimes \dbQ_p$ such that:

\smallskip
\parindent=35pt
\Item{(5.1.3)}
$\forall\,\alpha\in\scrJ$, $v_\alpha$ is the $p$-component
of the \'etale component of $t_\alpha$.

\smallskip
\Item{(5.1.4)}
There is a cycle $\psitil\colon H^1_{\text{\'et}}\otimes
H^1_{\text{\'et}}\to\dbZ_p(-1)$, which is a perfect
alternating form (it comes from the polarization of the Hodge $\dbQ$--structure on
$W^*=H^1(A_\dbC,\dbQ)$ induced from the polarization $\psi$ of the Hodge $\dbQ$--structure on $W$ defined by $h$).
The cycle $\psitil$ differs from the
perfect alternating form $\ptil_V\colon
H_{\text{\'et}}^1\otimes H_{\text{\'et}}^1\to\dbZ_p(-1)$ (induced by
the principal polarization $p_A$) just by multiplication with a $\dbZ_p$-unit.
It is fixed by the Galois group $\Gal(\Kbar/K)$.

\smallskip
\item{(5.1.5)}
For any integral ring $R$ which is faithfully
flat over $\dbZ_p$, and for every free $R$-module $M_R$ satisfying  $M_R\fracwithdelims[]1p=H^1_{\text{\'et}}\otimes_{\dbZ_p} R\fracwithdelims[]1p$
and enveloping the family of tensors $(v_\alpha)_{\alpha\in\scrJ_0}$ with
respect to $\psitil$, the Zariski closure of $G_{R\fracwithdelims[]1p}$ in $G\Sp(M_R,\psitil)$ is a reductive group over $\Spec(R)$.

\smallskip
\item{(5.1.6)}
The subgroup of $G\Sp(H^1_{\text{\'et}}\otimes\dbQ_p,\psitil)$
fixing $v_{\alpha}$, $\forall\alpha\in\scrJ$, is exactly
$G_{\dbQ_p}$.

\smallskip
\item{(5.1.7)}
The Galois group $\Gal(\Kbar/K)$ fixes $v_{\alpha}$, $\forall\alpha\in\scrJ$.
\parindent=25pt

\smallskip
In 5.1.5 and 5.1.6 we think of $G$ as a subgroup of $GL(W^*)$. Property 5.1.5 results from 4.3.5 and 4.3.6 1) as the family $(v_\alpha)_{\alpha\in\scrJ_0}$ is
$\dbZ_{(p)}$-well positioned with respect
to $\psi$ for the group $G$. Properties 5.1.3, 5.14 and 5.1.6 are trivial.

\Proclaim{5.1.8.} \rm
Property 5.1.7 results from the fact that the family $(t_\alpha)_{\alpha\in\scrJ}$ of Hodge cycles of $A_\dbC$ is defined over $K$: from  the fact that the abelian variety $A$ over $V$ has level-$N$ structure for any $N\in\dbN$ relatively prime to $p$, we deduce that the $l$-components of the \'etale components of the Hodge cycles of this family are defined over $K$ (here  $l$ is an arbitrary prime different from $p$).
\finishproclaim

\smallskip
\Proclaim{5.2. Step 2. Crystalline machinery.} \rm
\finishproclaim

Sections 5.2.1 to 5.2.10, 5.2.13 and 5.2.14 we follow closely [Fa3]. The new things are 5.2.1.1 and the use of the ring $\Rtil e$. ${}^1$ $\vfootnote{1}{Only in its final version [Fa3] changed to rings as $\tilde Re$ of 5.2; most of the previous versions of loc. cit. were stated in terms of rings as $Re$ of 5.2. So, though all statements of 5.2 are correct as they stand (including the one of 5.2.4 on $gr_F$'s; cf. also [Fa3, ending paragraph]), they do not match entirely with the references to the published [Fa3] in 5.2;  this paper and [Fa3] were refereed in the same period of time (ending 1998), and so we were not able to adjust this paper to the last minute changes in [Fa3].}$

\Proclaim{5.2.1.} \rm 
Let $\pi$ be a uniformizer of $V$.
As $V$ is totally ramified over $V_0$, there is an Eisenstein
polynomial $f_e(T)\in V_0[T]$ of degree $e:=[K:K_0]$ such that $f_e(\pi)=0$
is a minimal equation of $\pi$ over $V_0$.
Denoting $R:=V_0[[T]]$, we get $V=R/f_eR$.

Let $S_e$ be the subring of $K_0[[T]]$ generated by $R$ and
 divided powers $\frac{f_e^n}{n!}$, $n\in\dbN$.
As $\frac{p^n}{n!}\in R$, $\forall n\in\dbN$, and as $f_e$ is an Eisenstein polynomial,
this is the same as the subring of $K_0[[T]]$ generated by
$R$ and divided powers $\frac{(T^e)^n}{n!}$, $n\in\dbN$.
Let $Re$ be the $p$-adic completion of $S_e$ and let $\Rtil e$
be the completion of $S_e$ with respect to the (decreasing)
filtration given by its ideals $I_n:=I^{[n]}$, $n\in\dbN\cup\{0\}$, where
$I:=(p,f_e(T))=(p,T^e)$. So $\Rtil e$ is the projective limit of artinian rings $S_e/I_n$, $n\in\dbN\cup\{0\}$.
We recall that $I^{[n]}$ is the ideal generated by elements of
the form $\frac{\beta_1^{a_1}}{a_1!}\ldots\frac{\beta_m^{a_m}}{a_m!}$, with
$m,a_1,...,a_m$ non-negative integers such that $a_1+\cdots+a_m\ge n$, and with $\beta_1,\dotsc,\beta_m\in I$.

We get that the $V_0$-algebra $\Rtil e$ (resp. $Re$) is contained in $K_0[[T]]$ and consists
of power series $\Sigma_{n\Ge 0} a_nT^n$ such that the sequence $b_n=a_n{\fracwithdelims[]{n}{e}}!$, $n\in\dbN\cap\{0\}$, is integral, i.e. $b_n\in V_0$, $\forall n\in\dbN\cup\{0\}$, (resp. it is integral and convergent to zero). Here we used $p>2$. 

Let $\Phi$ be the Frobenius lift of $S_e$, $Re$ or $\Rtil e$ extending the Frobenius automorphism of $V_0$ and such that
$\Phi(T)=T^p$.
Decreasing filtrations are defined on $Re$ and $\Rtil e$ by the rules: For $m\in\dbN\cup\{0\}$, $F^m(Re)$
is the ideal of $Re$ obtained as the $p$-adic completion of the ideal of $S_e$ generated by divided powers $\frac{f_e^n}{n!}$ with $n\ge m$, while $F^m(\Rtil e)=I_m\Rtil e$.

We have ring epimorphisms $S_e \twoheadrightarrow V_0$, $Re\twoheadrightarrow
V_0$, $\Rtil e\twoheadrightarrow V_0$ defined by the rule $\Sigma_{n\ge 0}
a_nT^n\to a_0$. We have also a $V_0$-epimorphism $\Rtil e\twoheadrightarrow V$,
sending $T$ to $\pi$.
\finishproclaim

\Proclaim{5.2.1.1. Remark.} \rm
The ring $S_e/pS_e$ is a local ring with the property that any element of its maximal ideal is nilpotent. Its residue field is $\dbF$. So any reductive group over $S_e$ is a split group, and so any reductive group over $Re$ is also a split group.
\finishproclaim

\Proclaim{5.2.2.} \rm 
Keeping the notations of 5.1, let $(M,\Phi_M,\nabla)$ be the
evaluation at the thickening attached naturallly to the closed embedding $\Spec(V/pV)\hookrightarrow \Spec(Re)$ of the Frobenius crystal over $V/pV$ defined by taking the dual
of the Lie algebra of the universal vector extension of the
abelian scheme $A$ (or of the $p$-divisible group of $A$) (see [Me] and [Fa3]).

So $M$ is a free $Re$-module of dimension $\dim_{Re}(M)=\dim_\dbQ (W)$ endowed with an $Re$-submodule $F^1(M)$, $\nabla$ is
an integrable connection  (nilpotent mod $p$) on $M$, and $\Phi_M$ is a $\nabla$-parallel $\Phi$-linear endomorphism of $M$.
The
restriction of $\Phi_M$ to $F^1(M)$ is divisible by $p$ and we
have an isomorphism $\left(M+\frac1p
F^1(M)\right)\otimes_{Re} {}_{\Phi}Re\twiceover{\Phi_M}{\longrightarrow}
M$. We have $M/F^1(Re)M=H_{dR}^1(A/V)$.
The submodule $F^1(M)$ of $M$ is the inverse image of the Hodge
filtration of $H_{dR}^1(A/V)$ defined by $A$, under the surjective map 
$M\twoheadrightarrow M/F^1(Re)M=H_{dR}^1(A/V)$. So $F^1(Re)M\subset F^1(M)$.

Using $F^1(M)$ the tensor algebra of $M\oplus M^*$ gets a natural filtration indexed by integers. In particular, $F^0(M^*)$ contains $F^1(Re)M^*$ and the quotient $F^0(M^0)/F^1(Re)M^*$ is naturally identified with $\Hom_V(M/F^1(M),V)$. Moreover, for $n\in\dbN$ we speak about $F^n(M^{\otimes 2n})$.
\finishproclaim

\Proclaim{5.2.2.1.} \rm
Let $(M_0,\varphi_0):=(M,\Phi_M)\otimes_{Re}V_0$.
It is the contravariant Dieudonn\'e-module of 
$A_{\dbF}$. There is an isomorphism 
$$
(M\fracwithdelims[]1p,\Phi_M)\arrowsim (M_0,\varphi_0)\otimes
Re\fracwithdelims[]1p
$$ 
of Frobenius isocrystals [Fa3, Ch. 6].
\finishproclaim

\Proclaim{5.2.3.} \rm 
Let $\Vbar$ be the integral closure of $V$ in $\Kbar$ and let
$\Vbar{^{\wedge}}$ be its $p$-adic completion.
Let $S_0$ be the ring consisting of
sequences $(x_n)_{n\in\dbN\cup\{0\}}$, with $x_n\in\Vbar/p\Vbar$ and
$x_{n-1}=x_n^p$, $\forall n\in\dbN$. $\Gal(\Kbar/K)$ acts naturally on $S_0$.

The $\Gal(\Kbar/K)$-module $\dbQ_p(1)$ can be identified with sequences $(\mu_n)_{n\in\dbN\cup\{0\}}$
 of $p$-power roots of unity (these are elements of $\Vbar$)
such that $\mu_{n-1}=\mu_{n}^p$, $\forall n\in\dbN$.
Taking such sequences modulo $p$, we get a group homomorphism
$\gamma\colon\dbQ_p(1)\to\dbG_m(S_0)$ respecting the Galois actions.
For an element $z\in\Vbar$, we choose a sequence $(z(n))_{n\in\dbN\cup\{0\}}$ of
elements of $\Vbar$ such that $z(0)=z$ and $z(n-1)=z(n)^p$, $\forall n\in\dbN$. Taking this sequence  modulo $p$ we obtain an element 
$\underline{z}\in S_0$, well
defined by $z$ up to multiplication with an element of  $\gamma(\dbZ_p(1))$.

Let $W(S_0)$ be the ring of Witt vectors of $S_0$.
Let $\theta\colon W(S_0)\twoheadrightarrow \Vbar{^{\wedge}}$ be
the ring epimorphism defined by $\theta((x_0,x_1,\ldots))=\sum_{n\ge 0} p^nx_{n,n}$,
where $(x_{n,m})_{m\in\dbN\cup\{0\}}$ is the sequence (of elements of $\Vbar/p\Vbar$) used for  defining
$x_n\in S_0$.
Let $\xi:=f_e((\underline{\pi},0,0,\ldots))$.
It is a generator of the kernel $I_{\theta}$ of $\theta$.

Let $B^+(V)$ be the Fontaine's ring defined as the completion of the divided power hull of the ideal $I_{\theta}$ of $W(S_0)$.
The ring $B^+(V)$ is a $W(S_0)$-algebra and so also a $V_0$-algebra, as $W(S_0)$ is a $V_0$-algebra. It has a (decreasing) filtration $F^n(B^+(V))$ by  divided
powers: $F^n(B^+(V))$ is the completion of $I_{\theta}^{[n]}$, $n\in\dbN\cup\{0\}$.
We have $B^+(V)/F^1(B^+(V))=\Vbar^\wedge$.
The Frobenius automorphism of $W(S_0)$ extends to a Frobenius lift $\Phi$ of
$B^+(V)$ (it makes sense to still denote it by $\Phi$, cf. 5.2.4 below).
The Galois group $\Gal(\Kbar/K)$ acts in an obvious manner on
$B^+(V)$, respecting its filtration.

There is a well defined homomorphism 
$$
\beta\colon\dbZ_p(1)\to F^1(B^+(V)), 
$$
obtained by taking log of the homomorphism obtained by composing the Teichm\"uller map $\dbQ_p(1)\to\dbG_m(W(S_0))$ (obtained from $\gamma$) with the canonical homomorphism $\dbG_m(W(S_0))\to\dbG_m(B^+(V))$.
We have $\Phi\circ\beta=p\beta$.
We also denote by $\beta$ the image of a fixed generator of
$\dbZ_p(1)$ through this log map $\beta$.
\finishproclaim

\Proclaim{5.2.3.1.} \rm
Let $B_{dR}^+(V)$ be the completion of $B^+(V)\otimes \dbQ_p$ in the
filtration topology.
We have $B_{dR}^+(V)/F^1(B_{dR}^+(V))=\Kbar{^{\wedge}}:=\Vbar{^{\wedge}}\fracwithdelims[]1p$.
Let $B_{dR}(V):=B_{dR}^+(V)\fracwithdelims[]1\beta$. It has a decreasing filtration $(F^i(B_{dR}(V))_{i\in\dbZ}$ obtained from the filtration of $B^+_{dR}(V)$ by declaring $\frac1{\beta}\in F^{-1}(B_{dR}(V))$.
As $K$ is separable over $K_0$ and so formally smooth over it,
we can lift the inclusion $K\hookrightarrow \Kbar^\wedge=B_{dR}^+(V)/
F^1(B_{dR}^+(V))$ to a $K_0$-monomorphism $K\hookrightarrow B_{dR}^+(V)$.
\finishproclaim

\Proclaim{5.2.4.} \rm 
There is an injective homomorphism of filtered rings 
$$i_V\colon Re\hookrightarrow B^+(V)
$$ 
defined by:
$T\to\underline{\pi}$. It respects Frobenius lifts.
We have $F^n(B^+(V))\cap Re=F^n(Re)$, $\forall\,n\in\dbN\cup\{0\}$,
and an isomorphism of graded $\Vbar{^{\wedge}}$-algebras  $gr_F(Re)\otimes_V
\Vbar{^{\wedge}}\arrowsim gr_F(B^+(V))$ induced by $i_V$, cf. [Fa3, Ch. 4].
As $(M,\nabla)$ defines a crystal over $V/pV$, the tensor product
$M\otimes_{Re}B^+(V)$ acquires a canonical $\Gal(\Kbar/K)$-action.
\finishproclaim

\Proclaim{5.2.5.} \rm
Let us return to the situation of 5.1. The integral version of Fontaine's comparison theorem [Fa3, Th.
7] provides us with an injective linear map of filtered $B^+(V)$-modules 
$$
\rho\colon M
\otimes_{Re}B^+(V)\hookrightarrow H_{\text{\'et}}^1\otimes_{\dbZ_p}
B^+(V).
$$
The filtration on $M\otimes_{Re}B^+(V)$ is the tensor product one, while the filtration on $H_{\text{\'et}}^1\otimes_{\dbZ_p} B^+(V)$ is the one induced by the filtration of $B^+(V)$.
We list the properties of $\rho$ we need.

\smallskip\noindent
\parindent=40pt
\Item{(5.2.6)}
The map $\rho$ respects Frobenius lifts and the Galois actions.

\Item{(5.2.7)}
Inverting $p\beta$, we obtain an isomorphism to be denoted by $\rho_1$; in what follows we still denote by $\rho_1$ the resulting isomorphism at the level of tensor algebras. 

\Item{(5.2.8)}
A tensor $v_\alpha\in (H_{\text{\'et}}^1\otimes H_{\text{\'et}}^{1*})^{\otimes r(\alpha)}
\otimes \dbQ_p$, $\alpha\in\scrJ_0$ (resp.
$\alpha\in\scrJ\setminus \scrJ_0)$ corresponds through $\rho_1$ 
to an element $w_\alpha\in F^0\bigl((M\otimes M^*)^{r(\alpha)}\bigr)$
(resp. to an element $w_\alpha\in F^0\bigl((M\otimes M^*)^{\otimes r(\alpha)}
\fracwithdelims[]1p\bigr)$), with $r(\alpha):=\frac12\deg(v_\alpha)$.

\Item{(5.2.9)}
We have $\Phi_M(w_\alpha)=w_\alpha$ and $\nabla w_\alpha=0$,
$\forall\,\alpha\in\scrJ$.

\Item{(5.2.10)}
Under the identification $M/F^1(Re)M\fracwithdelims[]1p=
H_{dR}^1(A_K/K)$, the tensor $w_\alpha$ is mapped into the de Rham
component of the Hodge cycle $t_\alpha$,
$\forall\,\alpha\in\scrJ$.

\item{(5.2.11)}
The bilinear maps
$\ptil_V$, $\psitil\colon H_{\text{\'et}}^1\otimes
H_{\text{\'et}}^1\to
\dbZ_p(-1)$ are inducing bilinear maps 
$M\otimes M\to Re(1):=\beta Re$
which become perfect alternating forms $\ptil_M$,
$\psitil_M\colon M\otimes M\to Re$.

\Item{(5.2.12)}
The subscheme $\Gtil_{Re}$ of $G\Sp(M,\psitil_M)$ obtained by taking the Zariski closure of the subgroup $G_{Re\fracwithdelims[]1p}$ of $GSp(M\fracwithdelims[]1p,\psitil_M)$ fixing  $w_{\alpha}$, $\forall\alpha\in\scrJ$, is a reductive group scheme over
$\Spec(Re)$, isomorphic to $G_{Re}$.

\parindent=25pt
\smallskip
Properties 5.2.6 and 5.2.7 are part of Fontaine's
comparison theory.
The existence of $\ptil_M$, $\psitil_M$ (of 5.2.11) results either from 5.2.13
below or from the functorial aspect of Fontaine's comparison
theory (for the category of $p$-divisible groups over $V$).
The fact that $\ptil_M$, $\psitil_M$ are perfect can be
seen looking at their restriction modulo $F^1(Re)$.

The bilinear map $\psitil$ induces an isomorphism
$(H_{\text{\'et}}^1)^*\arrowsim H_{\text{\'et}}^1(1)$ (of Galois
modules) and $\psitil_M$ induces an isomorphism $M\arrowsim M^*(1)$
(of filtered Frobenius crystals). So we get isomorphisms
$(H_{\text{\'et}}^1\otimes H_{\text{\'et}}^{1*})^{\otimes m}\arrowsim
H_{\text{\'et}}^{1\otimes 2m}(m)$ and $(M\otimes M^*)^{\otimes
m}\arrowsim M^{\otimes 2m}(-m)$, $\forall m\in\dbN$.

Using these isomorphisms, 5.2.9 and the part of 5.2.8
involving the family of tensors
$(v_\alpha)_{\alpha\in\scrJ\setminus\scrJ_0}$, result from
Fontaine's comparison theory.
If $n\in\dbN$, then:
\finishproclaim

\Proclaim{5.2.13.}
The $\dbQ_p$-vector spaces spanned by a Galois-invariant class
$w_{\text{\'et}}\in H_{\text{\'et}}^1
\fracwithdelims[]1p^{\otimes 2n}(n)$ (also called an
\'etale Tate-cycle) are in one to one correspondence, through the
isomorphism $\rho_1$, with the $\dbQ_p$-vector spaces spanned by a class
$w\in F^n\left(M^{\otimes 2n}\fracwithdelims[]1p\right)$ which
is annihilated by $\nabla$ and fixed by $\Phi_M^n:=\Phi_M/p^n$.
The correspondence is achieved through the formula $\rho_1(w\otimes
1)=w_{\text{\'et}}\otimes\beta^n$.
\finishproclaim

\noindent
We stated 5.2.13 in terms of $\dbQ_p$-vector spaces as there is no natural choice for $\beta$ (cf. 5.2.3); we could have also stated it in terms of free $\dbZ_p$-modules of rank one.
The part of 5.2.8 concerning the family of tensors
$(v_\alpha)_{\alpha\in\scrJ_0}$, results from the fact that
$(v_\alpha)_{\alpha\in\scrJ_0}$ are integral (i.e. they are tensors of the tensor algebra of $H_{\text{\'et}}^1\oplus {H_{\text{\'et}}^{1*}}$) with
 $\deg(v_\alpha)\le
2(p-2)$, $\forall\,\alpha\in\scrJ$, and from the following key
supplement of 5.2.13 [Fa3, Cor. 9]:

\Proclaim{5.2.14.}
If $n\le p-2$, then in the correspondence of 5.2.13,
$w_{\text{\'et}}$ is integral (i.e an element of
$H_{\text{\'et}}^{1\otimes 2n}(n)$) if and only if $w$ is
integral (i.e. iff $w\in F^n(M^{\otimes 2n}))$.
\finishproclaim

\Proclaim{5.2.15.} \rm
We now prove 5.2.10.
This property results from the following observations.

\smallskip
{\bf 1)}\enspace
Tensoring the isomorphism $\rho_1$ with $B_{dR}(V)$ (using
the canonical inclusions $B^+(V)\hookrightarrow B_{dR}^+(V)\hookrightarrow B_{dR}(V)$), we
get an isomorphism, which can be written in the form
$$
\rho_2\colon H_{dR}^1(A/V)\otimes_V B_{dR}(V)\arrowsim
H_{\text{\'et}}^1\otimes_{\dbZ_p} B_{dR}(V)
$$
(the inclusion $V\hookrightarrow B_{dR}(V)$ is induced by the inclusion $K\hookrightarrow B_{dR}^+(V)$ of 5.2.3.1; we used the canonical identification $M\otimes_{Re}V
=H_{dR}^1(A/V)$ and  5.2.2.1).

\smallskip
{\bf 2)}\enspace
The isomorphism $\rho_2$ is nothing else but the isomorphism
which comes up in the de Rham conjecture, proved in [Fa1] (with slight correction in the unpublished [Fa2]) (i.e. the comparison map for the $p$-divisible group of an abelian variety over $V$ is the same as the comparison map for the abelian variety itself, cf. [Fa3, introd. to Ch. 6]).

\smallskip
{\bf 3)}\enspace
The Hodge cycles $(t_\alpha)_{\alpha\in\scrJ}$ are de Rham
cycles.
This means that $\rho_2$ takes the de Rham component ${t_{\alpha}}_{dR}$ of
$t_\alpha$ into the $p$-component of
the \'etale component $v_{\alpha}$ of $t_\alpha$.

\smallskip
Part 1) is obvious.
For proving \ 2) it is enough to show that the isomorphism
$$
\rho_3\colon M_0\otimes_{V_0}B(V)\arrowsim
H_{\text{\'et}}^1\otimes_{\dbZ_p} B(V)
$$
(with $B(V):=B^+(V)\fracwithdelims[]1{p\beta}$) obtained
from $\rho_1$ through the isomorphism $M\fracwithdelims[]1p
\arrowsim M_0\otimes Re\fracwithdelims[]1p$ of 5.2.2.1), is exactly the isomorphism 
$$\rho_4\colon M_0\otimes_{V_0} B(V)\arrowsim H^1_{\text{\'et}}\otimes_{\dbZ_p} B(V)$$ 
of the crystalline cohomology, defined for the abelian
variety $A$ over $V$ (see [Fa1, 5.6]).

\smallskip
To see this it is enough to show that the isomorphism
$\rho_3\circ\rho_4^{-1}$ of $H_{\text{\'et}}^1\otimes
B(V)$:

\smallskip
i) preserves the $F^0$-filtrations (i.e. preserves $H_{\text{\'et}}^1\otimes
F^0(B(V))$),

\smallskip
ii) and becomes identity on $H_{\text{\'et}}^1\otimes
F^0(B(V))/F^1(B(V))$.

\smallskip
This is so due to the fact that there is no element of
$\text{End}(H_{\text{\'et}}^1)\otimes F^1(B(V))$ fixed by
 the Frobenius endomorphism of $B(V)$ induced naturally by $\Phi$.

\smallskip
 Due to the way in which $H_{\text{\'et}}^1$ can be recovered from $H_{dR}^1(A_K/K)$ through the comparison map, we get that $\rho_4^{-1}$ takes  $H_{\text{\'et}}^1$
into $F^0(M_0\otimes_{V_0} B(V))$. This implies i).

For ii) it is enough to check that the Hodge--Tate structures defined on  $H_{\text{\'et}}^1\otimes
F^0(B(V))/F^1(B(V))$ by the two isomorphisms $\rho_3$ and $\rho_4$ are the same. This is done (by direct computation) in [Fa5, the proof of Th. 4] for abelian varieties over $\Spec(V[x])$ (with $x$ an independent variable), and so, due to  functoriality (under the morphism $\Spec(V[x])\to\Spec(V)$), for abelian varieties over $V$.

The proof of \ 3) is almost contained in [Bl].
The extra ingredient needed is an improvement in Principle B
of [Bl, 3.1], as our abelian variety $A_K$ might not be definable over
$\bar{\dbQ}$, a condition required in [Bl, 3.1] (of course as $\scrN$ is defined over $O_{(v)}$, we can select the
lift $m_V$: $\Spec(V)\to\scrNbar^0$ of $y$ in such a way that
$A_K$ is definable over
$\bar{\dbQ}$).
This improvement in Principle B is achieved by the trick of Lieberman. 

We can think of the de Rham component $w_\alpha$ of
$t_\alpha$ as a tensor of $H_{dR}^1(A/V)^{\otimes 2r(\alpha)}$, and so as a tensor of $H_{dR}^{2r(\alpha)}$ ($A^{r(\alpha)}/V$),
where $A^{r(\alpha)}$ is the product of $A$ over $V$ taking $r(\alpha)$ times
(for instance $A^2=A\times_V A$).
If $r(\alpha)=1$, there is nothing to be proved ($\rho_2$ respects algebraic cycles of degree 2). For $r(\alpha)\ge 2$ we get 3) above for $t_{\alpha}$ as a consequence of the following general principle. 
\finishproclaim

\Proclaim{5.2.16. Principle B.}
Let $\scrL$ be the field of fractions of a complete DVR of mixed characteristic having a perfect residue field of characteristic $p$. 
Let $Y$ be a geometrically connected smooth variety over
$\scrL$, and let $\Pi\colon B\to Y$ be an abelian scheme over
$Y$.
Let $n\Ge 2$.
Then if a pair $v=(v_{\text{\'et}},v_{dR})$, with
$v_{\text{\'et}}\in H^0(Y,R^1\Pi_*(\dbQ_p)^{\otimes 2n}(n))$
and with $v_{dR}\in H^0(Y,R_{dR}^1(B/Y)^{\otimes 2n})$, is a de
Rham cycle in a point $z_1\in Y(\overline{\scrL})$, then it is a de
Rham cycle in any other point $z_2\in Y(\overline{\scrL})$.
\finishproclaim

\smallskip
\proof
Let $B^n$ be the $n$-times product of $B$ over $Y$.
All the spectral sequences connecting the cohomology of
$B^n$ with the cohomology of $Y$ degenerate (this is called
the trick of Lieberman).
This results from the fact that $B^n$ has many endomorphisms
over $Y$ --which have to respect the spectral sequences--
defined by multiplying with integers the different factors
$B$ of $B^n$.
For every pair $(r,q)$ of non-negative integers we obtain commutative diagrams (which are part of these spectral sequences) 
$$
\spreadmatrixlines{1\jot}
\CD
E_2^{r,q} @>{d}>> E_2^{r+2,q-1}\\
@V{\tilde a}VV @VV{\tilde b}V\\
E_2^{r,q} @>{d}>> E_2^{r+2,q-1},
\endCD
$$
where $\tilde a$ and $\tilde b$ are multiplications with some integers $n_1$
and $n_2$.
For a suitable choice of multiplications on the factors $B$
of $B^n$, we can achieve $n_1\not=n_2$.

This implies that the $\scrL$-linear map 
$$
H_{dR}^{2n}(B^n/\scrL)\to H^0(Y,R_{dR}^1(B/Y)^{\otimes 2n})
$$ 
is surjective.
The rest of the proof is exactly as in [Bl, 3.2].

\Proclaim{5.2.17.} \rm
We are left with the proof of 5.2.12. We first remark that once we know that $\Gtil_{Re}$ is a reductive group scheme over $\Spec(Re)$, the fact that it is isomorphic to $G_{Re}$ ($Re$ is a $\dbZ_{(p)}$-algebra) is a direct consequence of the fact that  $\Gtil_{Re}$ and $G_{Re}$ are both split reductive groups (cf. 5.2.1.1) and of properties 5.2.7 and 5.2.8 (which guarantee that they are isomorphic over $\Spec(B^+(V)\fracwithdelims[]1p$) (cf. the uniqueness of a split reductive group associated to a given root datum; see [SGA3, Vol. 3, p. 305]).
\finishproclaim

\Proclaim{5.2.17.1.} \rm
To prove that $\Gtil_{Re}$ is indeed a reductive group over $Re$ we can move over the faithfully flat $Re$-algebra $Re^1:=Re\otimes_{V_0} V$. The ring $Re^1$ is integral as it is a subring of $K[[t]]$. It is also a $V$-algebra.
We have an isomorphism
$$\rho_0\colon M\otimes_{Re} Re^1\fracwithdelims[]1p\arrowsim H_{dR}^1(A_K/K)\otimes_K Re^1\fracwithdelims[]1p$$
taking $w_{\alpha}$ into ${t_{\alpha}}_{dR}$, $\forall\alpha\in\scrJ$, and taking $\ptil_M$ into the perfect form $p_A\colon H_{dR}^1(A/V)\otimes H_{dR}^1(A/V)\to V$. It is defined  starting from the isomorphism $M\fracwithdelims[]1p\arrowsim M_0\otimes_{V_0} Re\fracwithdelims[]1p$ of 5.2.2.1 and from the isomorphism $M/F^1(Re)M\arrowsim H_{dR}^1(A/V)$.

The fact that $\rho_0$  takes $w_{\alpha}$ into ${t_{\alpha}}_{dR}$, $\forall\alpha\in\scrJ$, results from  5.2.15 3), as the extension of $\rho_0^{-1}$ to $B_{dR}(V)$ (we have a natural inclusion $Re^1\hookrightarrow B_{dR}(V)$, cf. 5.2.3) when composed with the extension of $\rho_1$ to $B_{dR}(V)$ is nothing else but the isomorphism $\rho_2$ of 5.2.15 (cf. the way $\rho_0$ and $\rho_2$ are defined). Obviously $\rho_0$ takes $\ptil_M$ into $p_A$ (cf. the def. of $\tilde p_M$ and the functoriality of 5.2.2.1).
\finishproclaim

\Proclaim{5.2.17.2.} \rm
We have an isomorphism
$$H_{\text{\'et}}^1\otimes K\arrowsim H_{dR}^1(A_K/K)$$
taking $v_{\alpha}$ into ${t_{\alpha}}_{dR}$, $\forall\alpha\in\scrJ$, and taking $\psitil$ into $p_A\colon H_{dR}^1(A_K/K)\otimes H_{dR}^1(A_K/K)\to K$.

To see this, we first remark that such an isomorphism exists over the field $L_{dR}(u)$ obtained from the field of fractions of $B_{dR}(V)$ by adjoining a square root of $\beta$. This results from 5.2.13 1) and 2): $\rho_2$ takes $p_A$ into $\beta\psitil$ (cf. 5.2.11), and so over $L_{dR}(u)$, by changing the extension  of $\rho_2$ (to $L_{dR}(u)$) by a scalar factor $u$, we get rid of the scalar $\beta$. Now everything results from the well known fact: $H_{ff}^1(K,G_K)=\{0\}$, as $V$ is a complete DVR with an algebraically closed residue field. Here the right lower index $ff$ refers to the faithfully flat topology.
\finishproclaim

\Proclaim{5.2.17.3.} \rm
From 5.2.17.1-2 we get an isomorphism
$$\rho_5\colon H_{\text{\'et}}^1\otimes Re^1\fracwithdelims[]1p\arrowsim M\otimes_{Re} Re^1\fracwithdelims[]1p$$
taking $v_{\alpha}$ into $w_{\alpha}$, $\forall\alpha\in\scrJ$, and taking $\psitil$ into $\ptil_M$. From 5.1.5 we deduce (cf. 5.2.8) that the Zariski closure of the subgroup of $GSp(M\otimes_{Re} Re^1\fracwithdelims[]1p,\ptil_M)$ fixing $w_{\alpha}$, $\forall\alpha\in\scrJ$, in $GSp(M\otimes Re^1,\ptil_M)$ is a reductive group scheme over $\Spec(Re^1)$. So $\Gtil_{Re}$ is a reductive group scheme over $\Spec(Re)$. This ends the proof of 5.2.12.
\finishproclaim

\Proclaim{5.2.18. Remark.} \rm
It is an easy exercise to see that under the identifications
$$H_{dR}^1(A_K/K)=M\otimes_{Re} K\arrowsim M_0\otimes_{V_0} Re\otimes K=M_0\otimes K$$
(defined by inverting $p$ in the isomorphisms of 5.2.2 and 5.2.2.1), $\forall\alpha\in\scrJ$, ${t_{\alpha}}_{dR}$ is an element of the tensor algebra of $(M_0\oplus M_0^*)\fracwithdelims[]1p$; so we could have avoided the replacement of $Re$ by $Re^1$ in 5.2.17.1 to 5.2.17.3.
\finishproclaim

\Proclaim{5.2.19. Remark.} \rm
In applications of 5.1 to the proof of the main results of 6.4, we  use only families  of $G$-invariant tensors $(v_\alpha)_{\alpha\in\scrJ_0}$ in spaces of the form $(W\otimes W^*)^{\otimes m}$ (with $m\in\dbN$) and  of degree not bigger than $2(p-2)$, which are
$\dbZ_{(p)}$-very well positioned for the group $G$ (i.e. we do not work in the relative situation with respect to $\psi$). Moreover we can choose the family $(v_{\alpha})_{\alpha\in\scrJ}$ such that $G$ is the subgroup of $GL(W)$ fixing it. This simplifies the things, in the sense that we do not really have to keep track of all bilinear forms ($\psi$, $p_A$, $\ptil_M$, etc.) we came across. 
\finishproclaim

\smallskip
\Proclaim{5.3. Step 3. The existence of a 
good morphism $\hbox{\rm \Spec}(V_0)\to\overline{\scrN}$ lifting $y$.} \rm
\finishproclaim

\smallskip\noindent
{\bf 5.3.1.}\enspace 
Let $(\Mtil,\tilde {\nabla}):=(M,\nabla)\otimes_{Re}\Rtil e$ be obtained by extension of scalars.
It is a natural evaluation of the $F$-crystal defined by the dual of the Lie algebra of the universal vector
extension of the $p$-divisible group over $V(p):=V/pV$ associated to
the abelian scheme $A_1:=A_{V(p)}$.
At the level of filtrations we have $F^1(M)\otimes_{Re}\Rtil
e\subset F^1(\Mtil)$, where $F^1(\Mtil)$ is the pull back of the Hodge filtration of
$H_{dR}^1(A_1/V(p))$ defined by $A_1$, through the surjective map
$\Mtil\twoheadrightarrow\Mtil\otimes_{\Rtil e}V(p)=M\otimes_{Re}V(p)=H_{dR}^1(A_1/V(p))$.

Let $G_{\Rtil e}$ and $G_V$ be the reductive groups obtained from $G_{Re}$ through the
canonical $V_0$-homomorphisms $Re\hookrightarrow \Rtil e\twoheadrightarrow V$. 
Let $M_V:=\Mtil\otimes_{\Rtil e}V=M\otimes_{Re}V=H_{dR}^1(A/V)$
and let $F^1(M_V)$ be its Hodge filtration defined by $A$.
It has the property that ${F^1(M_V)\otimes_{V}}_j\dbC$ is the $F^{1,0}$ summand of the Hodge direct sum decomposition 
$H_{dR}^1(A/V)\otimes_V {}_j\dbC=F^{1,0}\oplus F^{0,1}$ (here $j$ is the inclusion of 5.1.2).
Let $\mu_\dbC\colon \dbG_m\to GL(M_V\otimes_j\dbC)$ be the
cocharacter such that $\gamma\in\dbG_m(\dbC)$ acts as identity on
$F^{0,1}$ and as the multiplication with $\gamma^{-1}$ on $F^{1,0}$.
The cocharacter $\mu_\dbC$ factors through $G_V\times_V {}_j\dbC$.
Let $\mu_1\colon\dbG_m\to G_V$ be a cocharacter which over $\dbC$
becomes conjugate to $\mu_\dbC$.
Let $M_V=F^1_V\oplus F^0_V$ be the direct sum decomposition obtained from
$\mu_1\colon\gamma_V\in\dbG_m(V)$ acts as the multiplication with
$\gamma_V^{-1}$ on $F^1_V$ and as identity on $F^0_V$.

Let $P^1$ be the parabolic subgroup of $G_V$ which leaves invariant
$F^1(M_V)$ and let $P^2$ be the parabolic subgroup of $G_V$
which leaves invariant $F^1_V$.
As $\mu_1$ and $\mu_\dbC$ are conjugate over $\dbC$, we deduce
that $P_K^1$ and $P_K^2$ become conjugate over $\dbC$ and so they
are conjugate over $\Kbar$.
As $P_K^1$ and $P_K^2$ are defined over $K$, we deduce from [Bo, 20.9] 
that they are conjugate over $K$, i.e. there is an element
$g\in G_V(K)$ such that  $gP_{K}^2g^{-1}=P_K^1$.
From the Iwasawa decomposition [Ti, 3.3.2] we
deduce that $G_V(K)=G_V(V)P^2(K)$.
This implies the existence of an element $g_0\in G_V(V)$ such
that $g_0P^2g_0^{-1}=P^1$.
We get a direct sum decomposition $M_V=\Fbar^1\oplus\Fbar^0$, with
$\Fbar^1:=g_0(F^1)$, associated to the
cocharacter $\mu:=g_0\mu_1g_0^{-1}\colon\dbG_m\to G_V$.
The parabolic subgroup of $G_V$ which leaves invariant $\Fbar^1$ is
$P^1$.
This implies $F^1(M_V)=\Fbar^1$. To check  this it suffices to show that $F^1(M_V)\otimes\dbC=\Fbar^1\otimes\dbC$.
There is an element $g_1\in G_V(\dbC)$ such that
$g_1(F^1(M_V)\otimes\dbC)=\Fbar^1\otimes\dbC$ and so
$g_1P_\dbC^1 g_1^{-1}=P_\dbC^1$.
We deduce that $g_1\in P^1(\dbC)$ (cf. [Bo, 11.16]) and so
$F^1(M_V)\otimes\dbC=\Fbar^1\otimes\dbC$.

\Proclaim{5.3.2. Lemma.}
The cocharacter $\mu\colon\dbG_m\to G_V$ lifts to a
cocharacter $\mutil\colon\dbG_m\to G_{\Rtil e}$.
\finishproclaim

\proof
Let $\grb_0$ be the ideal of $\Rtil e$ generated by the divided
powers of $f_e$.
Let $\grb_n:=\grb_0+I_n\Rtil e$, $n\in\dbN$.
Let $S_n:=\Spec(\Rtil e/I_n\Rtil e)$ and $T_n:=\Spec(\Rtil e/\grb_n)$.
So $T_n$ is a closed subscheme of $S_n$ and $T_{n+1}$ and $S_n$ is a closed subscheme of $S_{n+1}$, $\forall n\in\dbN$.
We have $T_{n+1}\cap S_n=T_n$.

We first remark that 
$$
\Rtil e=\text{proj.lim.}\Rtil e/I_n\Rtil e
$$ 
and 
$$V=\text{proj.lim. }\Rtil e/\grb_n
$$ 
(as $p>2$ and as $\Rtil e/I_n\Rtil e$ is $p$-adically complete); the projective systems are indexed by $n\in\dbN$. Second we show: if $\mu_n\colon\dbG_m\to G_{S_n}$ (with $n\in\dbN$) is a cocharacter such
that $\mu_n\vert T_n=\mu\vert T_n$, then there is a cocharacter
$\mu_{n+1}\colon\dbG_m\to G_{S_{n+1}}$ such that $\mu_{n+1}\vert
T_{n+1}=\mu\vert T_{n+1}$ and $\mu_{n+1}\vert S_n=\mu_n$ (here if $i_0\colon Y_0\hookrightarrow Y$ is a closed embedding and if $\nu$ is a morphism between two $Y$-schemes, we denote by $\nu\vert Y_0:=i_0^*(\nu)$).

To prove this, let $\mutil_{n+1}\colon\dbG_m\to G_{S_{n+1}}$ be any (group homomorphism)
lift of $\mu_n$ (cf. [SGA3, Vol. 2, p. 48]).
Now $\mutil_{n+1}\vert
T_{n+1}$ and $\mu \vert T_{n+1}$ are two lifts of $\mu\vert
T_n$.
From loc. cit.  we deduce the existence of an element
$h_n\in\ker(G_{\Rtil e}(T_{n+1})\to G_{\Rtil e}(T_n))$ such that
$h_n(\mutil_{n+1}\vert T_{n+1})h_n^{-1}=\mu \vert T_{n+1}$.
From the smoothness of $G_{\Rtil e}$ we deduce the existence of
an element 
$h_0\in\ker(G_{\Rtil e}(S_{n+1})\to G_{\Rtil e}(S_n))$
such that under the homomorphism $G_{\Rtil e}(S_{n+1})\to G_{\Rtil
e}(T_{n+1})$  it goes to $h_n$.
Then $\mu_{n+1}=h_0\mutil_{n+1}h_0^{-1}$ satisfies the required
conditions.

We start with a cocharacter $\mu_1\colon\dbG_m\to G_{S_1}$ lifting $\mu|T_1$. We build up inductively $\mu_n\colon\dbG_m\to G_{S_n}$ as above. Conclusion: we can choose $\mutil$ in such a way that $\mutil\vert S_n=\mu_n$, $n\in\dbN$.
Obviously $\mutil|V=\mu$. This ends the proof of the Lemma.

\smallskip\noindent
{\bf 5.3.3.}\enspace
Let now $\mutil\colon\dbG_m\to G_{\Rtil e}$ be a cocharacter
lifting $\mu$.
It achieves a direct sum decomposition $\Mtil=\Ftil^1\oplus\Ftil^0$
with the property that $\Ftil^1\otimes_{\Rtil e}V=F^1(M_V)$.

As $\Spec(\Rtil e)$ is a projective limit of nilpotent thickenings
of $V(p)$, from the deformation theory of principally
polarized abelian schemes (cf. [Me]; see also [FC, p. 14]) we deduce the existence of  a
principally polarized abelian scheme $(\Atil,p_{\Atil})$ over $\Rtil e$ 
associated to the filtered $F$-crystal $(\Mtil,\Ftil^1,\Phi_M,\tilde{\nabla})$ (we still denote by $\Phi_M$ the Frobenius endomorphism of $\Mtil$ induced from the one of $M$ by extension of scalars, as the ring homomorphism $Re\to\Rtil e$ respects the Frobenius lifts) and the symplectic form $\ptil_{\Mtil}$ on $\Mtil$ (obtained from $\ptil_M$ by extension of scalars; it guarantees that we get things over $\Spec(\Rtil e)$ and not only over $Spf(\Rtil e)$), such that under
the epimorphism $\Rtil e\twoheadrightarrow V$, it becomes $(A,p_A)$ (cf. the fact that $F^1(M_V)=\Ftil^1\otimes_{\Rtil e}V$ and that
$\ptil_V$ is obtained from $\ptil_{\Mtil}$ by tensorization).

\Proclaim{5.3.3.1. Lemma.}
The morphism $\mtil\colon\Spec(\Rtil
e)\to\bar{\scrM}$ associated to $(\Atil,p_{\Atil})$ and its
level-$N$ symplectic similitude structures (lifting those of $A_1$), factors through the Zariski closure of $\scrN_{K_0}$ in
$\bar{\scrM}$. 
\finishproclaim

\proof 
We can move from $\Rtil e$ to $R_{\dbC}:=\dbC[[T]]$ under the composition 
$$
\gtil\colon\Rtil e\hookrightarrow\Rtil e\otimes_{V_0} V\operatornamewithlimits{\hookrightarrow}\limits^{\gtil_0}\Rtil e\otimes_{V_0} V\operatornamewithlimits{\hookrightarrow}\limits^{\gtil_1}\dbC[[T]]
$$ 
(the first inclusion being the natural one). Here $\gtil_0$ is the affine transformation taking $T$ into $\pi^{e-1}T+\pi$, with $\pi$ the uniformizer of $V$ used in 5.2.1. This is well defined as the series $\sum_{n=0}^{\infty}\frac{\pi^{en}}{n!}$ is convergent in $V$ (as $p>2$). The homomorphism $\gtil_1$ is the inclusion defined  by the inclusion $j\colon V\hookrightarrow\dbC$ (of 5.1.2) and by the fact that it takes $T$ into $\frac{T}{j(\pi)^{e-1}}$.  Under the canonical surjective map $\Omega_{\Rtil e/V_0}\twoheadrightarrow \overline{\Omega}_{\Rtil e/V_0}$, with $\overline{\Omega}_{\Rtil e/V_0}$ the free module over $\Rtil e$ generated by $dT$, the Gauss--Manin connection of $\Mtil$ (defined by $\Atil$), becomes $\tilde{\nabla}$. This implies that $\frac{\delta}{\delta T}$ annihilates $w_{\alpha}$, $\forall\alpha\in\scrJ$. The principally polarized abelian variety over $\dbC$, obtained from the principally polarized abelian variety over $\Spec(R_{\dbC})$ induced from $(\Atil,p_{\Atil})$ through $\gtil$, by taking (in $R_{\dbC}$) $T=0$, is the extension of $(A,p_A)$ to $\dbC$ via $j$. We have $\frac {\delta T}{\delta (\pi^{e-1}T+\pi)}=\pi^{1-e}$ and $\frac {\delta T}{\delta {\pi}^{e-1}T}=\frac 1{\pi^{e-1}}$. Now everything results from 4.1.5. This ends the proof of the Lemma.

\Proclaim{5.3.4.} \rm
Let $\Rtil e^n$ be the normalization of $\Rtil e$ in its field of fractions. The natural surjection $\Rtil e\twoheadrightarrow V_0$  factors through $\Rtil e^n$ (due to the graded structure of $\Rtil e^n$ inherited as a subring of $K_0[[T]]$, cf. 5.2.1) producing a natural surjection $\Rtil e^n\twoheadrightarrow V_0$. From 5.3.3.1 and the definition of $\scrN$ we get a morphism 
$$\Spec(\Rtil e^n)\to\scrNbar.
$$ 
So we get a morphism 
$$m_0\colon\Spec(V_0)\to\scrNbar
$$ 
lifting $y$. 
It gives birth to:

\smallskip\noindent
(5.3.5)\enspace
a principally polarized abelian scheme $(A_0,p_0)$ over
$\Spec(V_0)$ (it is obtained from $(\scrA_{\scrN},\scrP_{\scrN})$ by pull back) having (compatibly) level-$N$ symplectic similitude
structure for any $N\in\dbN$ satisfying $(N,p)=1$ (defined by a similitude  isomorphism   $k_N\colon (L\otimes \dbZ/N\dbZ,\psi)\arrowsim (A[N],p_0)$ of principally quasi-polarized finite flat group schemes over $V_0$);

\smallskip\noindent
(5.3.6)\enspace
a family $(t_\alpha^0)_{\alpha\in\scrJ}$ of Hodge cycles of
$A_{0K_0}$ (we recall that $K_0=V_0\fracwithdelims[]1p$).

\smallskip
We have:

\smallskip\noindent
(5.3.7)\enspace
The quadruple $[A_{0\dbC},
p_{0\dbC},(t_\alpha^0)_{\alpha\in\scrJ},k]$ is a class
of $\scrA(G,X,W,\psi)$ (see 4.1) (here $k$ is induced as in 4.1.1 from the the isomorphisms $k_N$, $N\in\dbN$, while the embedding $V_0\hookrightarrow\dbC$ is the inclusion $j$ of 5.1.2; $k_N$, with $N$ a power of $p$, is obtained via the identification of 5.1.2).

\smallskip\noindent
(5.3.8)\enspace
Under the identifications 
$$H_{dR}^1(A_0/V_0)=M_0=H_{\text{crys}}^1(A_{0\dbF}/V_0)=\Mtil\otimes V_0$$
the de Rham component $u_\alpha$ of $t_\alpha^0$ is obtained from $w_{\alpha}$ through the epimorphism $\Rtil e\twoheadrightarrow V_0$, $\forall\alpha\in\scrJ$,  and is
a tensor of  $(M_0\otimes M_0^*)^{\otimes
r(\alpha)}\fracwithdelims[]1p$ if $\alpha\in\scrJ\setminus\scrJ_0$ and
a tensor of  $(M_0\otimes M_0^*)^{\otimes r(\alpha)}$ if
$\alpha\in\scrJ_0$.

\smallskip\noindent
(5.3.9)\enspace
If $\varphi_0$ is the Frobenius endomorphism of $M_0$, we have
$\varphi_0(u_\alpha)=u_\alpha$, $\forall\,\alpha\in\scrJ$.

\smallskip\noindent
(5.3.10)\enspace
The polarization $p_0$ induces a perfect
alternating form $\psi_0\colon$\break
$M_0\otimes M_0\to V_0(1)$
(i.e. $\psi_0(\varphi_0(t),\varphi_0(z))=
p\sigma(\psi_0(t,z))$, $\sigma$ being the Frobenius automorphism of $V_0$).

\smallskip\noindent
(5.3.11)\enspace
There is a direct sum decomposition $M_0=F^1\oplus F^0$, with
$F^1$ as the Hodge filtration of
$M_0=H_{dR}^1(A_0/V_0)$ defined by $A_0$, such that $u_\alpha$ belongs to the $F^0$-filtration of $(M_0\otimes M_0^*)^{\otimes r(\alpha)}\fracwithdelims[]1p$ defined naturally by $F^1$, $\forall\,\alpha\in\scrJ$.

\smallskip\noindent
(5.3.12)\enspace
The subgroup  of $G\Sp(M_0,\psi_0)$ obtained by taking the Zariski closure of $G_{K_0}$ (the subgroup of $GSp(M_0\fracwithdelims[]1p,\psi_0)$ fixing $u_{\alpha}$, $\forall\alpha\in\scrJ$) is the reductive group
scheme $G_{V_0}$, and the  decomposition $M_0=F^1\oplus F^0$ is associated to a  cocharacter $\mu_0\colon\dbG_m\to G_{V_0}$, with $\beta_0\in\dbG_m(V_0)$ acting through $\mu_0$ as the multiplication with $\beta_0^{-i}$ on $F^i$, $i=\overline{0,1}$. 
\finishproclaim

All these things result from the analogue properties (see 5.2.8 to 5.2.12) of the family of 
tensors $(w_\alpha)_{\alpha\in\scrJ}$ (situated in spaces of
the form $(M\otimes M^*)^{\otimes m}\fracwithdelims[]1p$)
(5.3.9 results from 5.2.9 and the isomorphism $M\fracwithdelims[]1p\arrowsim M_0\otimes
Re\fracwithdelims[]1p$ of 5.2.2.1); in connection to $\mu_0$, cf. 5.3.2 and 5.3.3.

\smallskip
\Proclaim{5.4. Step 4. Local deformation.} \rm
\finishproclaim

\Proclaim{5.4.1.} \rm
Let $\Rbar:=V_0[[z_1,\dotsc,z_e]]$ be a ring of formal  
power series with coefficients in $V_0$, and let $\Phi_{\Rbar}$ denote the
Frobenius lift on $\Rbar$ which extends the Frobenius automorphism $\sigma$ of
$V_0$ and sends $z_i\to z_i^p$.
Let $\Abar$ be an abelian scheme over $\Spec(\Rbar)$.
Let $M(\Abar):=H_{dR}^1(\Abar/\Rbar)$. It is a free $\Rbar$-module of rank twice the relative dimension $d(\Abar)$ of $\Abar$. Let $F^1(M(\Abar))$ be its Hodge filtration. 
We have:

\smallskip
(a)
The $\Rbar$-module $F^1(M(\Abar))$ is a direct summand in $M(\Abar)$ and free of rank $d(\Abar)$.

\smallskip
(b)
There is a $\Phi_{\Rbar}$-linear endomorphism $\Phi_A\colon
M(\Abar)\to M(\Abar)$ whose restriction to $F^1(M(\Abar))$ is divisible
by $p$ and such that it induces a $\nabla(\Abar)$-parallel isomorphism
$$
\Phi_A\colon\bigl(M(\Abar)+\frac1p F^1(M(\Abar))\bigr)
\otimes_{\Rbar} {}_{\Phi_{\Rbar}}\Rbar\arrowsim
M(\Abar).
$$ 

Here the connection $\nabla(\Abar)$ on $M(\Abar)$ is induced from the Gauss--Manin
connection $\nabla_{\Abar}$ (of $\Abar$) on $M(\Abar)$, through the canonical surjective map $\Omega_{\Rbar/V_0}\twoheadrightarrow\overline{\Omega}_{\Rbar/V_0}$, with $\overline{\Omega}_{\Rbar/V_0}$ the free $\Rbar$-module having as a basis $dz_1,...,dz_e$. The connection on the domain of $\Phi_A$ is induced naturally by $\nabla(\Abar)$. We refer to the quadruple 
$$
(M(\Abar),F^1(M(\Abar)),\Phi_A,\nabla(\Abar))
$$ 
as the $p$-divisible object of the Fontaine's category $\scrM\scrF_{[0,1]}(\Rbar)$ defined by $\Abar$ (this category is defined in the same manner as for smooth $V_0$-algebras; see [Fa1]). 

The above facts are just a variant of Grothendieck--Messing's theory, cf. [Me].
\finishproclaim

\Proclaim{5.4.2.} \rm
Let now $\Spec(\Rbar_0)$ be the completion of
$\Sp(M_0,\psi_0)$ in the origin.
We have an isomorphism $\Rbar_0\arrowsim
V_0[[z_1,\dotsc,z_{\ebar}]]$, with $\ebar:=2\ell^2+\ell$ for $\ell:=\frac12\dim_{\dbQ}(W)$.
Let $\Spec(R_0)$ be the completion  of the derived subgroup
$G_{V_0}^{\der}$ of $G_{V_0}$ in the origin.
We have $R_0\arrowsim V_0[[z_1,\dotsc,z_{e_1}]]$, with
$e_1:=\dim\,G^{\der}$.
The inclusion $G_{V_0}^{\der}
\hookrightarrow\Sp(M_0,\psi_0)$ produces a
surjection $r_0\colon\Rbar_0\twoheadrightarrow R_0$.
We choose identifications $\Rbar_0=V_0[[z_1,\dotsc,z_{\ebar}]]$ and
$R_0=V_0[[z_1,\dotsc,z_{e_1}]]$ such that the epimorphism
$r_0$ of $V_0$-algebras is defined by: $z_i\to z_i$ if $i\le e_1$, and $z_i\to 0$
if $i>e_1$.
Let now $\Phi_{\Rbar_0}$ and $\Phi_{R_0}$ be the
Frobenius lifts of $\Rbar_0$ and respectively $R_0$ such that they take
$z_i\to z_i^p$ and are compatible with $\sigma$.
\finishproclaim

\Proclaim{5.4.3.} \rm
Let $\scrO_y$ be the local ring of $y$ in $\bar{\scrM}$, let
$\widehat{\scrO_y}$ be its completion and let $(A_y,p_{A_y})$ be
the principally polarized abelian scheme over
$\Spec(\widehat{\scrO_y})$ obtained from  $(\scrA_{\scrM}, \scrP_{\scrM})$ through the composite morphism
$\Spec(\widehat{\scrO_y})\to\bar{\scrM}\to\scrM$.
We fix an isomorphism $\widehat{\scrO_y}\arrowsim
V_0[[z_1,\dotsc,z_{e_2}]]$, with $e_2:=\dim\,\Sh(G\Sp(W,\psi),S)$,
such that the epimorphism $\widehat{\scrO_y}\twoheadrightarrow V_0$, 
associated to
the morphism $\Spec(V_0)\to\bar{\scrM}$ defined by $m_0$, is identity on
$V_0$ and sends $z_i$ to zero.
Let $\Phi_y$ be the Frobenius lift on $\widehat{\scrO_y}$,
such that it extends the Frobenius automorphism of $V_0$ and sends
$z_i$ to $z_i^p$.
Let $(M_y,F^1_y,\Phi,\nabla_y)$ be the $p$-divisible object of $\scrM\scrF_{[0,1]}(\widehat{\scrO_y})$ defined by $A_y$.
The principal polarization $p_{A_y}$ induces a perfect
alternating form $\psi_y\colon M_y\otimes M_y\to
\widehat{\scrO_y}$.
\finishproclaim

\Proclaim{5.4.4.} \rm
We consider now the triple
$(M_{\Rbar_0},F^1_{\Rbar_0},\Phi_0)$ defined by
$M_{\Rbar_0}:=M_0\otimes_{V_0}\Rbar_0$,
$F^1_{\Rbar_0}:=F^1\otimes_{V_0}\Rbar_0$ and
$\Phi_0:=g_{Sp}(\varphi_0\otimes\Phi_{\bar R_0})$, with $g_{Sp}$ the universal
element of $\Sp(M_0,\psi_0)(\Rbar_0)$ defined by the natural morphism $\Spec(\Rbar_0)\to Sp(M_0,\psi_0)$.

From [Fa3, Th. 10] we deduce easily the existence of an abelian scheme
$A_{\Rbar_0}$ over $\Spec(\Rbar_0)$, with $A_0=A_{\Rbar_0}\times_{\Rbar_0} V_0$ (the surjection $\Rbar_0\twoheadrightarrow V_0$ is the identity on $V_0$ and sends all $z_i$ to $0$), and
such that the $p$-divisible object of $\scrM\scrF_{[0,1]}(\Rbar_0)$ defined by $A_{\Rbar_0}$ is
exactly $(M_{\Rbar_0},F^1_{\Rbar_0},\Phi_0,\nabla_0)$ (the connection $\nabla_0$ on $M_{\Rbar_0}$ is uniquely determined by the considered triple, cf. loc. cit.).

There is a unique principal polarization $p_{\Rbar_0}$ on
$A_{\Rbar_0}$ (that is why we get an abelian scheme over $\Spec(\Rbar_0)$ and not only over $Spf(\Rbar_0)$) corresponding to $\psi_0$ and lifting the principal polarization $p_0$ of $A_0$ (cf. the theory of
deformations of principally polarized abelian schemes). The principally polarized abelian scheme 
$(A_{\Rbar_0},p_{\Rbar_0})$ endowed  with the level-$N$
(symplectic similitude) structures lifting those of $A_0$, is obtained from $(A_y,p_{A_y})$ (and its level-$N$ symplectic similitude structures obtained from  those of $(\scrA_{\scrM},\scrP_{\scrM})$ by pull back) through  a morphism corresponding to a ring homomorphism $\alpha_y\colon\scrOhat_y\to\Rbar_0$. Here $N\in\dbN$, $(N,p)=1$. Warning: $\alpha_y$ might not respect the two Frobenius lifts
$\Phi_y$ and $\Phi_{\Rbar_0}$.
\finishproclaim

\Proclaim{5.4.4.1.} \rm
If $(A_{R_0},p_{_{R_0}})$ is the principally polarized abelian scheme over
$\Spec(R_0)$ obtained from $(A_y,p_{A_y})$ through  $r_0\circ\alpha_y$, then the $p$-divisible object of $\scrM\scrF_{[0,1]}(R_0)$ defined by  $A_{R_0}$
(together with $p_{_{R_0}}$) can be identified with $(M_{R_0},
F_{R_0},\Phi_1,\nabla_1)$ (together with $\psi_0$), where
$M_{R_0}:=M_0\otimes_{V_0}R_0$, $F^1_{R_0}:=F^1\otimes_{V_0}R_0$,
$\Phi_1:=g_{G^{\der}}(\varphi_0\otimes\Phi_{R_0})$, with $g_{G^{\der}}$ the universal
element of $G_{V_0}^{\der}(R_0)$ (this results from the fact that $r_0$
respects the Frobenius lifts), and with $\nabla_1$ the unique integrable connection on $M_{R_0}$ such that $\Phi_1$ is $\nabla_1$-parallel ([Fa3, Th. 10]).

From the uniqueness of such a connection $\nabla_1$, we deduce (cf. [Fa3, rm. ii) after Th. 10]) that it
respects the $G_{V_0}^{\der}$-action.
This means that $\nabla_1$ is of the form $\delta_0+\gamma_{R_0}$, with $\delta_0$ as the connection annihilating $M_0$ and with
$\gamma_{R_0}\in\text{Lie}(G_{V_0}^{\der})\otimes\overline{\Omega}_{R_0/V_0}$. Here $\overline{\Omega}_{R_0/V_0}$ is the free module over $R_0$ having as a basis $dz_1,...,dz_{e_1}$.
As $G_{V_0}^{\der}$ is a subgroup of $G_{V_0}$, we deduce that $\nabla_1(u_\alpha)=0$,
$\forall\,\alpha\in\scrJ$. As the Gauss--Manin connection on $M_{R_0}$ associated to $A_{R_0}$ becomes under the canonical surjection $\Omega_{R_0/V_0}\twoheadrightarrow\overline{\Omega}_{R_0/V_0}$ the connection $\nabla_1$, we deduce that $\frac {\delta}{\delta z_i}$ annihilates  $u_{\alpha}$, $\forall\alpha\in\scrJ$ ($i=\overline{1,e_1}$).
We have $(A_{R_0},p_{R_0})\otimes_{R_0}V_0=(A_0,p_0)$ (as $\alpha_y$ takes the ideal $(z_1,\dotsc,z_{e_2})$ into the ideal 
$(z_1,\dotsc,z_{\ebar})$).
\finishproclaim

\Proclaim{5.4.5.} \rm
The morphism
$\Spec(R_0)@>{q_0}>>\scrMbar$ associated to
$(A_{R_0},p_{_{R_0}})$ and its level-$N$ symplectic similitude
structures, $N\in\dbN$ such that $(N,p)=1$, induced from those
of $(A_0,p_0)$ ($R_0$ is a strictly Henselian ring) factors through the Zariski closure of $\Sh_H(G,X)_{K_0}$ in $\scrMbar$ (moving from $V_0[[z_1,\dotsc,z_{e_1}]]$ to $\dbC[[z_1,\dotsc,z_{e_1}]]$, this results from 5.4.4.1 and 4.1.5), and so it factors through $\scrNbar$ ($R_0$ being a normal ring).
We denote this factorization by $q_1\colon \Spec(R_0)\to\scrNbar$.
\finishproclaim

\Proclaim{5.4.6.} \rm
The Lie algebra $\grg:=\Lie(G_{V_0})$ is the Lie subalgebra
of $\gsp:=\Lie(GSp(M_0,\psi))$ centralizing $u_{\alpha}$, $\forall\alpha\in\scrJ$. So $\grg\otimes K_0$ is left invariant by $\varphi_0$. Let
$F^0(\grg):=\{x\in\grg\mid x(F^1)\subset F^1\}$ and
$F^{1}(\grg):=\{x\in \grg\mid x(F^1)=0\}$. Similarly we define for $i=\overline{0,1}$, $F^i(\gsp)$. The $V_0$-module $F^i(\grg)$ is the intersection of $\grg$ with  $F^i(\gsp)$, $i=\overline{0,1}$.
This implies that $F^i(\grg)$ are direct summands in $\grg$. We deduce easily that the quadruple
$$(\grg,\varphi,F^0(\grg),F^1(\grg))
$$ 
is a $p$-divisible object of $\scrM\scrF_{[-1,1]}(V_0)$, i.e. we have
$$\varphi(\frac1p F^{1}(\grg)+F^0(\grg)+p\grg)=\grg$$ 
(this Frobenius transform is included in $\grg$ and is a direct summand of $\gsp$, cf. the existence of $\mu_0$ in 5.3.12; so it is $\grg$). We call this quadruple the (Shimura) filtered Lie $\sigma$- crystal attached to the $V_0$-lift $m_0$ of $y$. Forgetting the filtration we get the (Shimura) Lie $\sigma$-crystal  $(\grg,\varphi)$ attached to the point $y$.

Similarly we get that $\grg_0:=\Lie(G_{V_0}^{\ad})$ gets a filtration and that $\grg_0\fracwithdelims[]1p$ gets a Frobenius automorphism (still denoted by $\varphi$), resulting in a $p$-divisible object of $\scrM\scrF_{[-1,1]}(V_0)$. So we similarly speak about the (Shimura) adjoint Lie $\sigma$-crystal attached to $y$, etc. 
\finishproclaim

\Proclaim{5.4.7.} \rm
From 5.4.5 and 5.4.6, we deduce the existence of
a commutative diagram of $V_0$-schemes
$$
\spreadmatrixlines{1\jot}
\matrix
T_1 &\overset{t_0}\to{\rightsubsetarrow{30}} &T_0\\
\lcapmapdown{i_1} &&\rcapmapdown{i_0}\\
\Spec(R_0) &\rightsubsetarrow{30} &\Spec(\Rbar_0)\\
\lmapdown{q_1} &&\rmapdown{\qbar_1}\\
\scrNbar &\longrightarrow &\scrMbar
\endmatrix
$$
and of a morphism $m_1\colon \Spec(V_0)\to T_1$ such that:

a)\enspace
$T_0=\Spec(V_0[[z_1,\dotsc,z_{e_2}]])=\Spec(\scrOhat_y)$
and $T_1=\Spec(V_0[[z_1,\dotsc,z_d]])$ (we recall that $d=\dim\,X=\dim\,\Sh(G,X)$);

\smallskip
b)\enspace
$\qbar_1$ is the morphism associated to $\alpha_y\colon\scrOhat_y\to\Rbar_0$;

\smallskip
c)\enspace
$t_0$, $i_0$ and $i_1$ are closed immersions;

\smallskip
d)\enspace
the tangent space of $T_0$ (in $t_0\circ m_1$) is a direct supplement of
$F^0(\grs\grp(M_0,\psi_0))$ in $\grs\grp(M_0,\psi_0)$;

\smallskip
e)\enspace
the tangent space of $T_1$ (in $m_1$) is a direct supplement of
$F^0(\grg)$ in $\grg$;

\smallskip
f)\enspace
$q_1\circ i_1\circ m_1=m_0$.

\smallskip\noindent
We have $d=\dim_{V_0}(\grg/F^0(\grg))$
and $e_2=\dim_{V_0}(\grs\grp(M_0,\psi_0)/F^0(\grs\grp(M_0,\psi_0)))$ (to justify these formulas it is enough to remark that these dimensions are computing the dimension of the (compact) dual Hermitian symmetric space of a connected component of $X$ and respectively of $S$; this can be seen moving over $\dbC$ and using 5.3.1). Here we identify the tangent space of $\Spec(R_0)$ (resp. of $\Spec(\Rbar_0)$) (in the $V_0$-valued point obtained by taking all $z_i=0$) with the Lie algebra of $G_{V_0}^{\der}$ (resp. of $Sp(M_0,\psi_0)$).
\finishproclaim

\Proclaim{5.4.8. Lemma.}
The ring homomorphism $\widehat{\scrO_y}@>{q_y}>>\widehat{\scrO_y}$
associated to $\qbar_1\circ i_0$ is an isomorphism.
\finishproclaim

\proof
It is enough to show that the tangent map of $q_y$ is an
isomorphism.
If this is not true, we deduce the existence of an epimorphism
$\widehat{\scrO_y}\overset{a_C}\to{\longtwoheadedrightarrow{15}}C:=\dbF[\varepsilon]/(\varepsilon^2)$ such that the composition
$b_C:=a_C\circ q_y$ factors through $\dbF$, i.e. $b_C=i\circ pr$, where $pr\colon\widehat{\scrO_y}\twoheadrightarrow\dbF$
is the homomorphism of $V_0$-algebras taking $z_i$ into $0$, and $i\colon\dbF\hookrightarrow C$ is the natural inclusion. But the Kodaira--Spencer map of the $F$-crystal over $\Spec(C)$ attached to the abelian scheme over $\Spec(C)$ obtained from $A_y$ through $b_C$ is injective (cf. 5.4.7). On the other hand, as $b_C=i\circ pr$, it is zero. We reached a contradiction. 
This proves the Lemma.

\smallskip
This Lemma details the last sentence of [Fa3, rm. iii) after Th. 10].

\smallskip
\Proclaim{5.5. Step 5. End of proof.} \rm
Let $\scrO_y^0$ be the local ring of $y$ in $\scrNbar$. From 5.4.7 and 5.4.8 we deduce that the ring homomorphism
$n\colon\scrO_y^0\to\scrO:=V_0[[z_1\dotsc,z_d]]$,
associated to the morphism $q_1\circ i_1\colon T_1\to
\scrNbar$, induces by completion an epimorphism $r\colon\widehat{\scrO_y^0}\twoheadrightarrow\scrO$. But $\widehat{\scrO_y^0}$ and $\scrO$ are local excellent normal rings of the same dimension. This implies that $r$ is an isomorphism. As $y$ was an arbitrary point of $\scrNbar$, we conclude that $\scrNbar$ is formally smooth over $V_0$ and so $\scrN$ is formally smooth over $O_{(v)}$. From 3.4.4 we deduce that $\scrN$ is an integral canonical model of the quadruple $(G,X,H,v)$ having the EEP. This ends the proof of 5.1.
\finishproclaim

\Proclaim{5.5.1. Remark.} \rm
From 5.5 and 5.4.7 we deduce that we can identify the $V_0$-valued points of $\Spec(\widehat{\scrO_y^0})$ with the $V_0$-valued points of the completion of the quotient $G_{V_0}/P_{V_0}$ in the $V_0$-valued point of it defined by the origin of $G_{V_0}$ (here $P_{V_0}$ is the parabolic subgroup of $G_{V_0}$ having $F^0(\grg)$ as its Lie algebra). 
\finishproclaim

\smallskip
\Proclaim{5.6. Comments.} \rm
\finishproclaim

\Proclaim{5.6.1. Corollary.} 
If $H_0$ is a compact open subgroup of $G(\dbA_f^p)$ small
enough, then $\scrN/H_0$ is the normalization of the Zariski closure of
$\Sh_{H_0\times H}(G,X)$ in $\scrM/H_0$, and is a quasi-projective scheme. The morphism $\scrNbar/H_0\to\scrMbar/H_0$ is a formal immersion in any point of $\scrNbar/H_0(\dbF)$.
\finishproclaim

The quasi-projectiveness part is a consequence of the fact that $\scrM$ is a pro-\'etale cover of a quasi-projective smooth scheme over $O_{(v)}$ (cf. [Mu]).

\Proclaim{5.6.2. Corollary.}
The integral canonical model $\Sh_p(G,X,H)$ of the triple $(G,X,H)$
is obtained by taking the normalization of the Zariski closure of $\Sh_H(G,X)$ in the extension to the normalization of $\dbZ_{(p)}$ in $E(G,X)$ of the integral canonical model of the quadruple $(G\Sp(W,\psi),S,K_p,p)$. It has the EEP.
\finishproclaim

\Proclaim{5.6.3. Example.} \rm
Using 4.3.11 we recover (for primes $p\Ge 3$) the well known results (cf. [Ko]) concerning the existence of integral canonical models of Shimura varieties of PEL type.
\finishproclaim

\Proclaim{5.6.4. Remark.} \rm
Morally $\scrN$ should be a closed subscheme of $\scrM$. To see why this should be so, we can move to $V_0$. We start with two $V_0$-valued points of $\scrNbar$, $x_0$ and  $x_1$, giving birth to the same $\dbF$-valued point $y$ of the special fibre of $\scrMbar$, and which give birth to two different $K_0$-valued points of $\scrNbar$, $z_0$ and $z_1$. Using a prime $l$ different from $p$, and using the
level-$l^N$ structures for any $N\in\dbN$, we get that the two families of tensors of the tensor algebra of $H_{\text{\'et}}^1(A_{\dbF},\dbQ_l)\oplus {H_{\text{\'et}}^1(A_{\dbF},\dbQ_l)}^*$ (here $A_{\dbF}$ is the abelian variety over $\dbF$ obtained from $\scrA_{\bar{\scrM}}$ through the point $y$) defined by  the two families of  $l$-components of  \'etale components of the Hodge cycles with which the  two abelian varieties over $K_0$ (obtained from $\scrA_{\bar{\scrN}}$ through the points $z_0$ and $z_1$) are naturally endowed, are the same. 

This should imply that the two families of tensors of the tensor algebra of $(M_0\oplus M_0^*)\fracwithdelims[]1p$ (with $M_0:=H_{crys}^1(A_{\dbF},V_0)$) defined by the de Rham components of the above two families of Hodge cycles, are the same (this is true if we have only cycles of degree 2, as they come from endomorphisms of $A$). If this is true, then we easily get that actually $x_0$ and  $x_1$ give birth to the same $\dbF$-valued point of $\scrNbar$ (cf. 5.4.7 and 5.4.8; see also 5.5.1). At least in the context of the PEL situation [Ko, Ch. 5], we do regain the well-known fact that $\scrN$ is a closed subscheme of $\scrM$. 

However if $p$ is a rational prime big enough, $\scrN$ is a closed subscheme of $\scrM$ (cf. 3.4.7).
In [Va2] we show how the validity of the Langlands--Rapoport conjecture (mentioned in 1.7) for $\scrN$ implies that $\scrN$ is a closed subscheme  of $\scrM$.
\finishproclaim

\Proclaim{5.6.5. Remark.} \rm
Section 5.3.4 remains true for any $V_0$-valued point of $\scrN$. More generally, for any $W(k)$-valued point of $\scrN$ (with $k$ an algebraically closed  field of characteristic $p$) we get:

\smallskip
\item{a)} 
A principally polarized abelian scheme $(A,p_A)$ over
$W(k)$ (obtained from ($\scrA_{\scrN},\scrP_{\scrN}$) by pull back) having (compatibly) level-$N$ symplectic similitude
structure for any  $N\in\dbN$ satisfying $(N,p)=1$ (defined by a similitude isomorphism  $k_N\colon (L\otimes \dbZ/N\dbZ,\psi)\arrowsim (A[N],p_A)$ of principally quasi-polarized finite flat group schemes over $W(k)$);

\smallskip
\item{b)}
A family $(t_\alpha)_{\alpha\in\scrJ}$ of Hodge cycles of
$A_{B(k)}$ (with $B(k):=W(k)\fracwithdelims[]1p$).

\smallskip
We have:

\smallskip
\item{c)}
Under the identification of
$H_{dR}^1(A/W(k))=M=H_{\text{crys}}^1(A_k/W(k))$
the de Rham component $u_\alpha$ of $t_\alpha$ belongs to  $(M\otimes M^*)^{\otimes r(\alpha)}\fracwithdelims[]1p$ if $\alpha\in\scrJ\setminus\scrJ_0$, and
to $(M\otimes M^*)^{\otimes r(\alpha)}$ if
$\alpha\in\scrJ_0$.

\smallskip
\item{d)} 
$\varphi(u_\alpha)=u_\alpha$, $\forall\,\alpha\in\scrJ$, $\varphi$ being the Frobenius endomorphism of $M$.

\smallskip
\item{e)}
The polarization $p_A$ induces a perfect
alternating form $\psi\colon$\break
$M\otimes M\to W(k)(1)$
(we have $\psi(\varphi(t),\varphi(z))=
p\sigma(k)(\psi(t,z))$, $\sigma(k)$ being the Frobenius automorphism of $W(k)$).

\smallskip
\item{f)}
There is a direct sum decomposition $M=F^1\oplus F^0$, with
$F^1$ as the Hodge filtration of
$H_{dR}^1(A/W(k))=M$ defined by $A$, such that $u_{\alpha}$ belongs to the $F^0$-filtration of $(M\otimes M^*)^{\otimes
r(\alpha)}\fracwithdelims[]1p$ defined naturally by $F^1$, $\forall\,\alpha\in\scrJ$.

\smallskip
\item{g)}
The subgroup of $GSp(M\otimes B(k),\psi)$ fixing $u_\alpha$, $\forall\,\alpha\in\scrJ$, is (reductive and identified with) $G_{B(k)}$.
The subgroup $G_{W(k)}$ of $G\Sp(M,\psi)$, obtained by taking the Zariski closure of $G_{B(k)}$,  is a reductive group
scheme over $W(k)$ and the decomposition $M=F^1\oplus F^0$ is associated to a cocharacter $\mu_{W(k)}\colon \dbG_m\to G_{W(k)}$, with $\beta_0\in\dbG_m(W(k))$ acting through it as the multiplication with $\beta^{-i}_0$ on $F^i$, $i=\overline{0,1}$. 

\smallskip
\item{h**)}
There is an isomorphism 
$$H_{\text{\'et}}^1(A_{\overline {B(k)}},\dbZ_p)\otimes_{\dbZ_p} W(k)\arrowsim H_{dR}^1(A/W(k))$$
taking the $p$-component of the \'etale component of $t_{\alpha}$ into (de Rham component) $u_{\alpha}$ (of $t_{\alpha}$), for any $\alpha\in\scr J$.

\smallskip
Properties a) to g) are just a reformulation of 5.3.4 for a $W(k)$-valued point of $\scrN$.  A proof of h) will be given in [Va2]. Its proof solves positively the following conjecture of Milne (slight restatement):
\finishproclaim

\Proclaim{5.6.6. Conjecture ([Mi5, 0.1]).} 
Let $\tilde A$ be an abelian scheme over the ring $W(k)$ of Witt vectors of an algebraically  closed field $k$ of characteristic $p$ and let $B(k):=W(k)\fracwithdelims[]1p$. Let $(s_{\delta})_{\delta\in I}$ be a family of Hodge cycles of $\tilde A$, including a polarization. We assume that the Zariski closure in $GL(L_p)$, with $L_p:=H_{\text{\'et}}^1$$(\tilde A_{\overline {B(k)}},\dbZ_p)$,  of the subgroup of $GL(L_p\otimes \dbQ_p)$ fixing the $p$-component of the \'etale component of $s_{\delta}$, $\forall\delta\in I$, is reductive. We also assume that $p$ is big enough with respect to the dimension of $\tilde A$. Then, for some (any) faithfully flat $W(k)$-algebra $R(k)$, there is an   isomorphism of $R(k)$-modules
$$L_p\otimes_{\dbZ_p} R(k)\arrowsim H_{dR}^1(\tilde A/W(k))\otimes_{W(k)} R(k)$$
mapping, for any $\delta\in I$, the $p$-component of the \'etale component of $s_{\delta}$ into  de Rham component of $s_{\delta}$.
\finishproclaim

\Proclaim{5.6.7. Remark.} \rm
The well known results for integral canonical models of  Siegel modular varieties (pertaining to  universal principally polarized abelian schemes over them and) concerning the existence of an ordinary isogeny type in positive characteristic and the existence of canonical lifts of ordinary abelian varieties (over perfect fields), remain valid  for our model $\scrN$. We get results pertaining to the principally polarized abelian scheme $(\scrA_{\scrN},\scrP_{\scrN})$ over it (cf. 1.6 and  [Va2]); we call special any such principally polarized abelian scheme over $\scrN$.
\finishproclaim

\Proclaim{5.6.8. Remark.} \rm
In [Va2] we will see that in the majority of cases the whole of 5.6.5 remains true without assuming that the (perfect) field $k$ is algebraically closed.  
\finishproclaim

\Proclaim{5.6.9. Remark.} \rm
We can work out 5.1 with a family of tensors which is $\dbZ_p$-very well positioned instead of a family of tensors which is $\dbZ_{(p)}$-very well positioned. The only thing needed to be changed is: we get $\dbQ_p$-linear combinations of (components of) Hodge cycles instead of (components of) Hodge cycles. Even better: in 5.1 it is enough to assume the existence of a family of tensors (of degrees not bigger than $2(p-2)$) enveloped by $L_p\otimes V_0$ and which is $V_0$-well positioned for $G_{K_0}$. This is a consequence of the proof of 5.1: we needed condition 4.3.5 to be satisfied for rings of the form $Re^1$; but they are $V_0$-algebras. However this often boils down to an enlarged family of tensors (of degrees not bigger than $2(p-2)$) of the tensor algebra of $W\oplus W^*$, which is $\dbZ_{(p)}$-very well positioned with respect to $\psi$ for $G$. For instance, this is so,  if we are dealing with strongly $V_0$-well positioned families of tensors (cf. 4.3.15 and 4.3.15.1): this is the case we will encounter in 6.5 and 6.6 (cf. 4.3.10 and 4.3.13); however we will not bother to mention strongly in 6.5 and 6.6 (as we think it is irrelevant).
\finishproclaim

\Proclaim{5.6.10. Remark.} \rm
We could have worked out the proof of 5.1 working at some finite level, i.e. working with some quotients $\scrN/H_0$ (with $H_0$ as in 5.6.1) and $\scrM/K_0^p$ (with $K_0^p$ a compact open subgroup of $GSp(W,\psi)(\dbA_f^p)$ properly chosen). This would have just slightly complicated the presentation. In [Va2] we refine the things: we work in such a finite level context, with points in perfect fields (here we worked with algebraically closed fields of characteristic $p$).  
\finishproclaim 

\smallskip
\Proclaim{5.7. A practical form of the basic result.}
\finishproclaim

\Proclaim{5.7.1. Theorem.}
Let $(G,X)\hookrightarrow (G\Sp(W,\psi),S)$ be
an injective map 
and let $p\Ge 5$ be a rational prime. We assume the existence of a $\dbZ_{(p)}$-lattice $L$ of $W$  such  that $\psi$ induces a perfect form $\psi\colon\L\otimes L\to \dbZ_{(p)}$ and the Zariski closure of $G$ in $G\Sp(L,\psi)$ is a reductive group $G_{\dbZ_{(p)}}$ over $\dbZ_{(p)}$ (so $G$ is unramified over $\dbQ_p$). If the  Killing  form on  $\Lie(G_{\dbZ_{(p)}}^{\der})$ and the form $\scrT$ on  $\Lie(G_{\dbZ_{(p)}}^{\der})$ induced (by restriction) by the trace form on $\text{End}(L)$ are both perfect, then $\Sh_{p}(G,X)$ exists and has the EEP.
\finishproclaim

\proof
This is a direct consequence of  4.3.10 b), 4.3.13, 3.1.6 and 5.1. We present the details.

Let $G_0:=G^{\der}$ and let $\grg_0:=\Lie(G_0)$. We have  
$$
s(\grg_0,W)=2.
$$ 
This can be easily checked starting from [De2, 1.3.7] (i.e. starting from the fact that all weights given irreducible subrepresentations of $W\otimes\dbC$ of a simple Lie algebra factor of $\grg_0\otimes\dbC$ are minimal weights --poides minuscules--, cf. [Bou2, Ch. VIII, \S7.3]). The fact that the Killing  form and the trace form $\scrT$ on $\grg_0$ are both perfect, can be restated (with the notations of 4.3.2): the tensors (of degree 4) $\pi(\grg_0)$, $B$ and $B^*$ (can be viewed --cf. 4.1-- as  tensors) of the tensor algebra of $W\oplus W^*$ (and) are enveloped by the $\dbZ_{(p)}$-lattice $L$. So the family of tensors formed by $\pi(\grg_0)$, $B$ and $B^*$ is $\dbZ_{(p)}$-well positioned for $G_0$ (cf. 4.3.10 b)). Now 4.3.13 guarantees the existence of a family of endomorphisms $(v_{\alpha})_{\alpha\in\scrJ_1}$ of $L$ fixed by $G$, which is $\dbZ_{(p)}$-well positioned with respect to the maximal torus of $Z(G)$. Let $(v_{\alpha})_{\alpha\in\scrJ_0}$ be the family of tensors formed by putting $\pi(\grg_0)$, $B$, $B^*$ and $(v_{\alpha})_{\alpha\in\scrJ_1}$ together. So $\scrJ_1\subset\scrJ_0$.

The family of $G$-invariant tensors $(v_{\alpha})_{\alpha\in\scrJ_0}$ is enveloped by $L$ and $\dbZ_{(p)}$-well positioned with respect to $G$ (cf. 4.3.6 2)).

For any $\alpha\in\scrJ_0$ we have $\deg(v_{\alpha})\in\{2,4\}$ and so $\deg(v_{\alpha})$ is not bigger than $2(p-2)$ (as $p$ is at least 5). Now everything results from 5.1.  This ends the proof of the Theorem.

\Proclaim{5.7.2. Notations.} \rm
Let $\Gtil_0=\prod_{i\in\scrK} \Gtil_i$ be a product of simple adjoint groups of classical Lie type over a field. Let 
$$B(\Gtil_0):=\prod_{i\in\scrK} B(\Gtil_i)$$ 
where, for any $i\in\scrK$, $B(\Gtil_i)$ is $6(l+1)$ if $\Gtil_i$ is of $A_l$ or $C_l$ Lie type, $6(l-1)$ if $\Gtil_i$ is of $D_l$ Lie type, and $6(2l-1)$  if $\Gtil_i$ is of $B_l$ Lie type with $l\Ge 2$.

Let $(G_0,X_0)$ be an adjoint Shimura variety of abelian type with $G_0$ a simple $\dbQ$--group. Let $f\colon (G,X)\hookrightarrow (G\Sp(W,\psi),S)$ be an injective map with $(G^{\ad},X^{\ad})=(G_0,X_0)$. Let $\grg_0$ be the Lie algebra of $G_0$ (or of $G^{\der}$). Let $\grh_0$ be a non-compact simple factor of ${\grg_0}\otimes {\dbR}$. We denote by $A(G_0,X_0,W)$ the number of elements of the set $I$ defined by an isomorphism $W\otimes\dbR\arrowsim W_0\oplus\oplus_{i\in I} W_i$ of $\grh_0$-modules, with $\grh_0$ acting trivially on $W_0$ and with each $W_i$ as an irreducible non-trivial $\grh_0$-module. It depends only on the representation of $\grg_0$ on $W$, and not on the choice of $G$ or of $\grh_0$ (cf. [De2, 2.3.4]). So the notation $A(G_0,X_0,W)$ is justified. 
\finishproclaim

\Proclaim{5.7.2.1. Lemma.} 
The factor $\delta_0$, that relates the Killing form $\scrK$ on a split simple Lie algebra over $\dbZ\fracwithdelims[]1{B(G_0)}$ of the same Lie type as $G_0$ and the trace form $\scrT$ on it associated to the irreducible representation of it given by a weight $w_i$ corresponding (cf. [De2, 1.3.7]) to the representation $W_i$ of $\grh_0$ (it does not depend on the element $i\in I$!) (so $\scrK=\delta_0\scrT$), is an invertible element of $\dbZ\fracwithdelims[]1{B(G_0)}$. Moreover $\scrK$ and $\scrT$ are perfect forms.
\finishproclaim

This is an easy computation, using the coroots of the classical Lie algebras (they are described in [Bou2, Ch. 8, \S13]) starting from the fact that any two $\grg$-invariant perfect bilinear forms on an absolutely simple Lie algebra $\grg$ over a field of characteristic zero differ one from another just by multiplication with a non-zero element of the field. It should be also compared with the explicit form of the Killing form of the (complex) classical Lie algebras [He, formulas (5), (16) and (22) of Ch. 3 \S8]. The extra thing needed besides these formulas is the fact (implied by the mentioned computation) that over an algebraically closed field of characteristic zero the trace forms on a ${\grs\gro}(n)$ Lie algebra defined by the representations associated to the fundamental weights corresponding to the roots $\alpha_{l-1}$ and $\alpha_l$, with $l=[\frac n2]$, are equal (here $\alpha_{l-1}$ and $\alpha_l$ are having the usual meaning; cf. [Bou2, Ch. 8, \S13] page 193 if $n\in\dbN$ is odd and page 208 if $n$ is even).

\Proclaim{5.7.3. Remark.} \rm
The conditions (in 5.7.1) that $p\Ge 5$ and  the above two bilinear forms on $\Lie(G_{\dbZ_{(p)}}^{\der})$ are perfect, are equivalent to: $p$ does not divide the product 
$$B(G^{\ad})\prod_{i\in\scrK} A(G_i^{\ad},X_i^{\ad},W),$$
where $(G^{\ad},X^{\ad})=\prod_{i\in\scrK} (G_i^{\ad},X_i^{\ad})$, with all $G_i^{\ad}$ as simple $\dbQ$--groups.  
Here the numbers $A(G_i^{\ad},X_i^{\ad},W)$ are computed starting from an injective map $(G_i,X_i)\hookrightarrow (GSp(W,\psi),S)$ factoring through the injective map $(G,X)\hookrightarrow (GSp(W,\psi),S)$, cf. 2.12 1).  
\finishproclaim

\Proclaim{5.7.4. Remark.} \rm
In 5.7.1 we can use instead of the bilinear form on $\Lie(G_{\dbZ_{(p)}}^{\der})$ induced by the trace form on $\grg\grl(L)$, any other bilinear form induced by a bilinear form on $\grg\grl(L)$ which is fixed by $G_{\dbZ_{(p)}}$ (cf.   4.3.10.1 1)). Even better: it is enough that such a bilinear form on $\grg\grl(L)$ is  defined only over $V_0=W(\overline{\dbZ/p\dbZ})$ (cf. 5.6.9). 
\finishproclaim

\Proclaim{5.7.5. Example: Classical Spin modular varieties of odd dimension (and rank two).} \rm
Let $l\Ge 3$ be an integer. Let $G:=SO(2,2l-1)$ be the $\dbQ$--group whose points in a $\dbQ$--algebra $R$ are those matrices $g$ in $SL(2l+1,R)$ which leave invariant the quadratic form $-x_1^2-x_2^2+x_3^2+...+x_{2l+1}^2$, i.e. ${}^{t}gI_{2,2l-1}g=I_{2,2l-1}$, with $I_{2,2l-1}$ the diagonal matrix of order $2l+1$ having $-1$ on the first two lines and $+1$ on the others. 

Let $\Sh(G,X)$ be the adjoint Shimura variety with $X$ a double copy of the Hermitian symmetric domain of $BD$ $I_{(p=2,q=2l-1)}$ type (cf. the classification of symmetric domains [He, p. 518]). The group $G$ is an absolutely simple adjoint group of $B_l$ Lie type which splits over $\dbQ(i)$. We have $\dim\,X=2l-1$ and $E(G,X)=\dbQ$. 

Let $f\colon (G_1,X_1)\hookrightarrow (G\Sp(W,\psi),S)$ (with $(G_1^{\ad},X_1^{\ad})=(G,X)$) be the injective map defined by the Spin representation of the simply connected cover $G_1^{\der}$ of $G$ (this representation is defined over $\dbQ$ as $G_1^{\der}$ splits over $\dbQ(i)$). We have $\dim_{\dbQ}(W)=2^l$ if $l$ mod 4 is 1 or 2, and $\dim_{\dbQ}(W)=2^{l+1}$ if $l$ mod 4 is 0 or 3  (cf. [Sa, p. 458]). The group $G_1^{\der}=\text{Spin}(2,2l-1)$ is a Spin group and $Z(G_1)=\dbG_m$ acts on $W$ by multiplication with scalars (so $G_1^{\ab}=\dbG_m$). For any prime $p$, $G_1$ is unramified over $\dbQ_p$. We have: $A(G,X,W)$ is $1$ or $2$ depending on the fact that $l$ mod $4$ is or is not $1$ or $2$. We call $\Sh(G_1,X_1)$ the classical Spin modular variety of dimension $2l-1$ (and rank two) (cf. [Va4] for the use of the word classical). 

Let $\grh:=\Lie(G_1^{\der})$ and let $\pi_W(\grh)$ be the projection of $\grg\grl(W)$ on $\grh$ associated to the direct sum decomposition $\grg\grl(W)=\grh\oplus\grh^\bot$ (here $\grh^\bot$ is the subspace of $\grg\grl(W)$ perpendicular to $\grh$ with respect to the trace form on $\grg\grl(W)$). Let  $B\colon\grg\grl(W)\to\grg\grl(W)^*$ be the linear map which is zero on $\grh^\bot$ and  $B|\grh\colon\grh\to\grh^*$ is the isomorphism induced by the Killing form on $\grh$, and let $B^*\colon\grg\grl(W)^*\to\grg\grl(W)$ be the linear map which is zero
on ${(\grh^\bot)}^*$ and
$B^*\vert\grh^*\colon\grh^*\to\grh$ is
$(B\vert\grh)^{-1}$. If $l$ mod 4 is 1 or 2, then $\Lie(G_1)$ is the Lie subalgebra of $\grg\grl(W)$ centralizing $\pi_W(\grh)$ due to the fact that the representation ${G_1^{\der}}_{\dbC}\to GL(W_{\dbC})$ is irreducible. So  $(G_1,X_1)$ is saturated in $(G\Sp(W,\psi),S)$. If $l$ mod 4 is 0 or 3 then the maximal connected subgroup $G_2$ of $GL(W)$ fixing $\pi_W(\grh)$ contains $G_1$, $G_2^{\der}$ is isogenous to $G_1^{\der}$ times a form of an $SL_2$-group, and $G_2^{\ab}$ is a torus of dimension 1 (the representation ${G_1^{\der}}_{\dbR}\to GL(W_{\dbC})$ is not irreducible; see [Sa, p. 458]). So $(G_1,X_1)$ is not saturated in $(G\Sp(W,\psi),S)$.

Let now $p$ be a prime not dividing $6(2l-1)$ and let $L$ be a $\dbZ_{(p)}$-lattice of $W$ such  that $\psi$ induces a perfect form $\psi\colon L\otimes L\to\dbZ_{(p)}$ and the Zariski closure of $G_1$ in $G\Sp(L,\psi)$ is a reductive group over $\dbZ_{(p)}$ (the existence of such a $\dbZ_{(p)}$-lattice results from the fact that the Spin representation of $G_1$ has a $\dbZ_{(p)}$-version).

Now the family of tensors formed by $\pi_W(\grh)$, $B$ and $B^*$ is integral with respect to $L$ (i.e. it is enveloped by $L$) (for instance, for $\pi_W(\grh)$ this means that it is a projector of $\grg\grl(L)$) and is $\dbZ_{(p)}$-very well positioned for the group $G_1$ (cf. 5.7.1 to 5.7.3). This implies that the Killing form on the Lie algebra $\grh_L:=\grh\cap\grg\grl(L)$ and the restriction to $\grh_L$ of the trace form on $\grg\grl(L)$ are both perfect. Let $K_p:=\{g\in G\Sp(W,\psi)(\dbQ_p)|g(L\otimes\dbZ_p)=L\otimes\dbZ_p\}$ and let $H_1:=K\cap G_1(\dbQ_p)$. So $K_p$ is a hyperspecial subgroup of $G\Sp(W,\psi)(\dbQ_p)$ and $H_1$ is a hyperspecial subgroup of $G_1(\dbQ_p)$. The normalization of the Zariski closure of $\Sh_{H_1}(G_1,X_1)$ in the integral canonical model $\scrM$ of $(G\Sp(W,\psi),S,K_p,p)$ is an integral canonical model $\scrN$  of $(G_1,X_1,H_1,p)$ (cf. 5.7.1 and 5.6.2). The universal (principally polarized) abelian scheme over $\scrM$ (obtained by choosing a $\dbZ$-lattice $L_{\dbZ}$ such that $L=L_{\dbZ}\otimes\dbZ_{(p)}$ and $\psi\colon L_{\dbZ}\otimes L_{\dbZ}\to\dbZ$ is perfect) gives birth to a principally polarized abelian scheme $(\scrA_{\scrN},\scrP_{\scrN})$ over $\scrN$ of dimension ${\frac 12}\dim_{\dbQ}(W)$. The pro-scheme $\scrN$ admits plenty of smooth toroidal compactifications and the abelian scheme $\scrA_{\scrN}$  extends to semi-abelian schemes over  these smooth toroidal compactifications of $\scrN$ (cf. [Va3]).

If $l=3$ then $\dim_{\dbQ}(W)=16$ and we obtain abelian schemes of dimension 8. If $l=4$ then $\dim_{\dbQ}(W)=32$ and we obtain abelian schemes of dimension 16. If $l=5$ then $\dim_{\dbQ}(W)=32$ and we obtain abelian schemes of dimension 16. If $l=6$ then $\dim_{\dbQ}(W)=64$ and we obtain abelian schemes of dimension 32.  
\finishproclaim

\Proclaim{5.7.6. Remark.} \rm
For $l=10$ we get the Shimura variety $\Sh(G_1,X_1)$ associated to  the moduli space of complex $K3$ surfaces.
\finishproclaim

For more examples, including the case of classical Spin modular varieties of even dimension (and rank 2), see [Va4].

\smallskip
\Proclaim{5.8. Integral good embeddings in a Siegel modular variety.} 

\Proclaim{5.8.1. Definition.} \rm
Let the pair $(G,X)$ define a Shimura variety of Hodge type. Let $p$ (resp. $p\Ge5$) be a rational prime such that $G$ is unramified over $\dbQ_p$. We say that $(G,X)$ (or $\Sh(G,X)$) has a good embedding (resp. a very good embedding) (in a Siegel modular variety) with respect to $p$, if there is an injective map $f\colon (G,X)\hookrightarrow (G\Sp(W,\psi),S)$ such that the hypotheses of 5.1 (resp. of 5.7.1) are satisfied. Similarly, we speak about an injective map $(G,X)\hookrightarrow (G\Sp(W,\psi),S)$ as being a good embedding or a very good embedding with respect to $p$. 
\finishproclaim

\Proclaim{5.8.2. Remark.} \rm
If $(G,X)$ defines a Shimura variety of Hodge type, if $p$ is a rational prime such that $G$ is unramified over $\dbQ_p$, and if $(G,X)$ has a good embedding with respect to $p$, then  $\Sh_p(G,X)$ exists (cf. 5.1) and we can study its points in fields of positive characteristic using the machinery of crystalline cohomology (cf. the proof of 5.1 and [Va1] and [Va2]).

\Proclaim{5.8.3. Definition.} \rm
Let $f\colon (G,X)\hookrightarrow (G\Sp(W,\psi),S)$ be an injective map and let $p$ be a prime such that $G$ is unramified over $\dbQ_p$. A $\dbZ_{(p)}$-lattice $L$ of $W$ is called good with respect to  $f$ if $\psi$ induces a perfect form $\psi\colon L\otimes L\to\dbZ_{(p)}$ and if the Zariski closure of $G$ in $G\Sp(L,\psi)$ is a reductive group over $\dbZ_{(p)}$.

\Proclaim{5.8.4. Proposition.}
Let $f\colon (G,X)\hookrightarrow (G\Sp(W,\psi),S)$ be an injective map with $G^{\ad}$ a simple $\dbQ$--group. Let $l$ be the rank of a simple factor of $G_{\dbC}^{\ad}$ (i.e. $G^{\ad}$ is of $A_l$, $B_l$, $C_l$ or $D_l$ Lie type) and let $N(G^{\ad})$ be the number of non-compact simple factors of $G_{\dbR}^{\ad}$. Let 
$$
p\Ge\text{max}(5,2l,\frac {\dim_{\dbQ}(W)}{2lN(G^{\ad})})
$$ 
be a rational prime. If there is a $\dbZ_{(p)}$-lattice of $W$ good with respect to  $f$, then $f$ is a very good embedding with respect to p.
\finishproclaim

\proof 
This results from 5.7.1 and 5.7.3. We have just to remark that $\dim_{\dbQ}(W)$ is at least $2lN(G^{\ad})A(G^{\ad},X^{\ad},W)$ (with equality only for $G=G\Sp(W,\psi)$) (this is an easy consequence of [De2, 2.3.7 b)]; for $m,n\Ge 2$ positive integers we have $mn\Ge m+n$) and that all the prime factors of $B(G^{\ad})$ are smaller than $\text{max}(5,2l)$ (cf. 5.7.2).

\Proclaim{5.8.5. Remark.} \rm
If in 5.8.4 we concentrate in just one Lie type of rank $l$ we can obtain even better estimates than the estimate of 5.8.4 which works for all Lie types of rank $l$. For instance, if $(G^{\ad},X^{\ad})$ is of $D_l^{\dbR}$ type, with $l\ge 5$, then we need $p\Ge\text{max}(5,l,\frac {\dim_{\dbQ}(W)}{2^{l-1}N(G^{\ad})})$. If $G^{\ad}$ is of $B_l$ Lie type, $l\ge 1$, then we need $p\Ge\text{max}(5,2l,\frac {\dim_{\dbQ}(W)}{2^lN(G^{\ad})})$, etc. In the mentioned cases, these estimates are a consequence of the dimension formula of the Spin representation of a split orthogonal Lie algebra (over $\dbC$) (see [Bou2, Ch. 8, \S13]).

\Proclaim{5.8.6. Corollary.} 
Let $f\colon (G,X)\hookrightarrow (G\Sp(W,\psi),S)$ be an injective map. Let $p$ be a prime greater or equal to $\text{max}(5,2+\dim_{\dbQ}(W)/2)$ (resp. greater or equal to $\text{max}(5,\dim_{\dbQ}(W)/2)$). If there is a $\dbZ_{(p)}$-lattice of $W$ good with respect to  $f$, then $f$ is a very good embedding (resp. is a good embedding) with respect to $p$.
\finishproclaim

\proof
If $p-2\Ge\text{max}(3,\dim_{\dbQ}(W)/2)$ then this is a consequence of 5.8.4 nd 5.8.5. If $p\Ge 5$ and $2p\in\{\dim_{\dbQ}(W),\dim_{\dbQ}(W)+2\}$, and if $f$ is not a very good embedding with respect to $p$, then either $G=G\Sp(W,\psi)$ or $2p=\dim_{\dbQ}(W)$ and $G^{\ad}$ is an absolutely simple $\dbQ$--group of $A_{p-1}$ Lie type. In both these cases we get immediately that we are in the context described in 4.3.11; so 5.6.3 applies.  

\Proclaim{5.8.7. Corollary.}
Let $f\colon (G,X)\hookrightarrow (G\Sp(W,\psi),S)$ be an injective map. Then there is $N(G,X)\in\dbN$ effectively computable such that  $f$ is a (very) good embedding with respect to any prime $p\Ge N(G,X)$ with the property that  $G$ is unramified over $\dbQ_p$.
\finishproclaim

\proof
Let $L$ be a $\dbZ$-lattice of $W$ such that we get a perfect form $\psi\colon L\otimes L\to\dbZ$. There is a number $N(G,L,f)\in\dbN$ such that for any prime $p\Ge N(G,L,f)$ the Zariski closure of $G$ in $GSp(L\otimes\dbZ_{(p)},\psi)$ is a reductive group scheme over $\dbZ_{(p)}$. It is effectively computable (for instance cf. 4.3.10 b)). Now we can take $N(G,X)=\text{max}\bigl(N(G,L,f),\dim(W)/2\bigr)$, cf. 5.8.6.

\Proclaim{5.8.8. Corollary.}
We assume that 5.6.5 h) holds. Then the Milne's conjecture (see 5.6.6) is true if the prime $p$ is bigger than $\text{max}(5,\dim(A))$.  
\finishproclaim

\proof
We use the notations of 5.6.6.
Let $f\colon (G,X)\hookrightarrow (GSp(W,\psi),S)$ be an injective map  defined by $(A,p_A)$ (here $p_A$ is the polarization of $A$ defined by some $s_{\delta(0)}$, $\delta(0)\in I$) and the reductive family $(s_{\delta})_{\delta\in I\setminus\{\delta(0)\}}$ (cf. 2.12 3)) with respect to $p_A$. From the hypotheses of 5.6.6 we deduce that there is a $\dbZ_{(p)}$-lattice $L$ of $W$ good for $f$. If $p_A$ is a principal polarization then the Corollary is a direct consequence of 5.8.6 and of 5.6.5 h) (cf. definitions 5.8.1 and 5.8.3). If $p_A$ is not a principal polarization, then we have to apply the Zarhin's trick [Za]: replacing $A$ by $(A\times A^t)^4$ the numbers $A(G_i,X_i,W)$ defined in 5.7.1 to 5.7.3 for the injective map $f$, are replaced by  numbers which are 8 times bigger. As we are taking $p\Ge 5$, this does not change anything (cf. the proof of 5.8.4), and so we do not have to replace $\dim(W)/2$ by $4\dim(W)$. It is easy to see that the Zarhin's trick does not destroy the $\dbZ_p$-\'etale reductiveness part. This ends the proof of the Corollary.
 
\smallskip
Actually we do not need to assume that $A$ is polarized (as 5.6.6  speaks about) (cf. [Va2]). For better estimates than $\text{max}(5,\dim(A))$ see [Va2].

\Proclaim{5.8.9. Remark.} \rm
If in 5.8.6 to 5.8.8 we concentrate just on one specific type of Shimura varieties, we can obtain much better estimates, cf. 5.8.5.
\finishproclaim

\bigskip
\noindent
{\boldsectionfont \S6. The existence of integral canonical models}

\bigskip
First we complete (cf. 6.1, 6.2 and 6.8) the steps (introduced in 3.4) needed to construct integral canonical models  of Shimura varieties of preabelian type. Then we digress very briefly (cf. 6.3) on conjugates of such models. The main results are gathered in 6.4, while their proofs spread till the very end of 6.8. Besides the tools developed in the previous chapters we rely heavily on [De2]. In particular, as a main new idea, we build up an integral version (6.5.1.1) of [De2, 2.3.10]. Sections 6.4.2, 6.5 and 6.6 are independent of 6.1 and 6.2; so in 6.2.2 E) to G) and 6.2.2.1 we refer to 6.5 and 6.6. Also, the proof of 6.4.5 b) depends only on 6.2.2 B), C) and D) and so we refer to it in 6.2.2.1. As a conclusion, the right order to read \S6 is: first 6.4.2 and its proof in 6.5 and 6.6, then 6.1, 6.2.3, 6.2.3.1, 6.2.1, 6.2.2 and 6.2.4, then 6.4.1, 6.4.1.1, 6.4.2.2 and 6.4.5, then 6.2.2.1 and 6.2.2.2, then 6.4.2.2, 6.4, 6.4.5.1 and 6.4.6 to 6.4.11, and finally 6.7 and 6.8 (6.3 can be read out at any time).

\smallskip
\Proclaim{6.1. The going up between finite maps.}

\Proclaim{6.1.1.} \rm
Let $\Sh(G,X)$ be a Shimura variety of Hodge type and
let $f\colon (G,X)\hookrightarrow (G\Sp(W,\psi),S)$ be an injective map. Let $p\Ge 3$ be a prime. We assume the existence of a $\dbZ_{(p)}$-lattice $L$ of $W$ which is good for $f$. Let $(G,X,H,v)$ be a quadruple of $(G,X)$ having an
integral canonical model $\scrM$, with $v$ dividing $p$. 

\Proclaim{6.1.2.* Theorem.}
We consider a finite map $f\colon (G_1,X_1,H_1,v_1)\to (G,X,H,v)$.
Then $(G_1,X_1,H_1,v_1)$ has an integral canonical model $\scrM_1$ having the EEP and obtained by taking the normalization of $\scrM_{O_{(v_1)}}$ in the ring of fractions $\scrR$ of $\Sh_{H_1}(G_1,X_1)$. The scheme $\scrM_1$ is a pro-\'etale cover of an open closed subscheme of $\scrM_{O_{(v_1)}}$. 

If $\Sh_p(G,X,H)$ exists, then $\Sh_p(G_1,X_1,H_1)$ also exists, has the EEP, and is the normalization of $\Sh_p(G,X,H)$ in $\scrR$.  
\finishproclaim

A complete proof of 6.1.2 will be presented in [Va3]. For a discussion and a proof in many cases, see 6.8. 

\Proclaim{6.1.2.1. Warning.} \rm
The results below (as well as 6.1.2) whose numbers have a right $*$, in the case of Shimura pairs $(G,X)$ of preabelian type which are not of abelian type, are fully proved in this paper only in the generic situation, i.e. working with a prime (or primes in some cases, like in 6.4.4) $p$ which is (or are) big enough, with an upper bound depending only on $(G,X)$ (cf. 6.8.5). See 6.8 for an explanation.  As 6.8.0 to 6.8.2 explain how we  prove (in [Va2] and in [Va3]) 6.1.2 in the remaining cases (see also 6.8.6), we felt it is appropriate to state the main results and remarks in the way we did. The labeled results are fully proved here in the abelian type case.  
\finishproclaim

\smallskip
\Proclaim{6.2. The going down between finite maps.}

\Proclaim{6.2.1.} \rm
Let $f\colon (G,X)\to (G_1,X_1)$ be a cover such that $E(G,X)=E(G_1,X_1)$ (cf. [MS, 3.4]). Let $E:=E(G,X)$. We consider a map $(G,X,H,v)\to (G_1,X_1,H_1,v)$ defined by $f$, with $v$ a prime of $E(G,X)$ dividing a rational prime $p\Ge 2$. Let $V_0:=W(\overline{k(v)})=W(\dbF)$ and let $A$ be the kernel of the homomorphism $G\to G_1$. We recall (cf. 2.4) that $A$ is a torus such that $H^1(\Gal(\overline{k}/k),A(\overline{k}))=0$ for any field $k$ of characteristic zero. Let $B:=G^{\ab}$.
\finishproclaim

\Proclaim{6.2.2. Theorem.}
We assume that $(G,X,H,v)$ has an integral canonical model $\scrM$ and that $\scrM_{V_0}$ has the EEP. We also assume that either

a) $p$ is relatively prime to the order $Q$ of the center of the simply connected semisimple group cover of $G^{\der}_1$ and $\scrM$ is a quasi-projective integral model, or

b) there is an injective map $f_2\colon (G_2,X_2)\hookrightarrow (GSp(W,\psi),S)$ which is a good embedding with respect to the prime $p>2$ and we have $G_2^{\der}=G^{\der}$ and $(G_2^{\ad},X^{\ad}_2)=(G^{\ad},X^{\ad})$.

Then $(G_1,X_1,H_1,v)$ has an integral canonical model $\scrM_1$. Moreover the natural morphism $\scrM\to\scrM_1$ is a pro-\'etale cover.
\finishproclaim

\proof
As the proof is quite long we itemize the steps (ideas).

{\bf A)} [Mi4, 4.11 and 4.13] contains all that is needed to see how to construct an integral model $\scrM_1$ of $(G_1,X_1,H_1,v)$ over $O_{(v)}$, as a quotient of $\scrM$. We just need to remark that such a quotient always exists as a scheme: $\scrM$ is a quasi-projective integral model (in case b)  cf. 5.6.1, 5.8.1 and 5.8.2; see also below), and so we can quote [Mu1, p. 112]. We want to prove that the morphism $\scrM\to \scrM_1$ is a pro-\'etale cover (this implies that $\scrM_1$ is a smooth integral model) and that $\scrM_1$ has the EP.

{\bf B)}  Let $S_1$ be an integral healthy  regular scheme over $O_{(v)}$ and let $q\colon{S_1}_E\to\scrM_1$ be a morphism. Let $S_0$ be the normalization of $S_1$ in the ring of fractions of ${S_1}_E\times_{\scrM_1} \scrM$. For proving that $\scrM_1$ has the EP, we need to show that $q$ extends to a morphism $S_1\to \scrM_1$. For seeing this it is enough to show that $S_0$ is a pro-\'etale cover of $S_1$ (as $\scrM$ has the EP and as a pro-\'etale cover of a healthy regular $O_{(v)}$-scheme is also a healthy regular $O_{(v)}$-scheme, cf. 3.2.2 4)). From the classical purity theorem we get: it is enough to work with $S_1$ the spectrum of a discrete valuation ring $O$ faithfully flat over $\dbZ_{(p)}$. We can assume that $O$ is complete with an algebraically closed residue field, and so that it is a $V_0$-algebra.

{\bf C)} The key fact for checking that $\scrM_1$ has the EP is:

\Proclaim{Fact.}  A connected component of ${\scrM_1}_{V_0}$ is the quotient of a connected component $\scrC^0$ of $\scrM_{V_0}$ by a commutative group $C_p$ which is a $Q^2$-torsion group. 
\finishproclaim

\proof
${\scrM_1}_{V_0}$ is the quotient of $\scrM_{V_0}$ by the group $A(\dbA_f^p)/\overline{A(\dbZ_{(p)})}$, where $\overline{A(\dbZ_{(p)})}$ is the topological closure of $A(\dbZ_{(p)}):=A(\dbQ)\cap H$ in $A(\dbA_f^p)$: this is an easy consequence of   [Mi4, 4.13]. We assume first that $G^{\der}$ is simply connected. So (cf. [De1, 2.4 and 2.5]) the set of connected components of $\scrM_{V_0}$ is in one to one correspondence to the set $B(\dbA_f^p)/\overline{B(\dbZ_{(p)})}$, with $\overline{B(\dbZ_{(p)})}$ having the analogue meaning of $\overline{A(\dbZ_{(p)})}$. If moreover $G_1=G_1^{\ad}$, we just have to add (cf. the Sublemma below) that the canonical homomorphism $A\to B$ has finite kernel of order a divisor of $Q$. 

\Proclaim{Sublemma.}
Let $t\colon T_1\to T_2$ be an isogeny of $\dbQ$--tori. Let $T_0$ be its kernel. Let $p$ be a prime such that $T_2$ is unramified over $\dbQ_p$. Let $H(T_i)$ be the hyperspecial subgroup of $T_i(\dbQ_p)$, $i=\overline{1,2}$. Let $T_i(\dbZ_{(p)}):=H(T_i)\cap  T_i(\dbQ)$; we denote by $\overline{T_i(\dbZ_{(p)})}$ its topological closure in $T_i(\dbA_f^p)$, $i=\overline{1,2}$. Let $Q(t)$ be the least common multiple of the orders of elements of the group $T_0(\dbC)$. Then the kernel of the natural homomorphism $t_p\colon T_1(\dbA_f^p)/\overline{T_1(\dbZ_{(p)})}\to T_2(\dbA_f^p)/\overline{T_2(\dbZ_{(p)})}$ is a $Q(t)^2$-torsion group.
\finishproclaim

\proof
Let $a\in\ker(t_p)$. Let $\tilde{a}\in T_1(\dbA_f^p)$ representing it. There is a sequence $(b_n)_{n\in\dbN}$ of elements of $T_2(\dbZ_{(p)})$ converging to $t(\tilde{a})\in T_2(\dbA_f^p)$: $T_2(\dbA_f^p)$ is a topological group having a countable basis of neighborhoods of its identity element. Let $c_n\in T_1(\dbZ_{(p)})$ be such that its image in $T_2(\dbZ_{(p)})$ is $b_n^{Q(t)}$. As $T_2(\dbA_f^p)$ is a locally compact group and as $T_0(\dbA_f^p)$ is a compact group, we deduce the existence of a subsequence $(c_n)_{n\in\dbN(1)}$, with $\dbN(1)$ an infinite subset of $\dbN$, converging to an element $a_1\in\overline{T_1(\dbZ_{(p)})}$. Obviously $\tilde{a}^{Q(t)^2}a_1^{-Q(t)}\in\overline{T_1(\dbZ_{(p)})}$ is the identity element. So $a^{Q(t)^2}=1$. This proves the Sublemma.

For proving the above Fact in the general case it is enough to remark that:

\smallskip
-- there is a cover $(G_0,X_0,H_0,v_0)\to (G,X,H,v)$ with $G_0^{\der}$ a simply connected semisimple group (cf. rm. 10) of 3.2.7) and so we can apply the previous argument involving only connected components (we do not need to assume that $(G_0,X_0,H_0,v_0)$ has an integral canonical model, as the argument on connected components can be performed over $\dbC$) for the induced cover $(G_0,X_0,H_0,v_0)\to (G_1,X_1,H_1,v)$;

-- the proof of Lemma 6.2.3 allows us to shift the situation to the case when $G_1=G_1^{\ad}$ (even for $p=2$).

\smallskip
{\bf D)} The scheme $S_0$ is a disjoint union of integral schemes. As $C_p$ is a $Q^2$-torsion group, we get that $S_0$ has the property that any abelian scheme $\scrA$ over the generic fibre of a connected component $S_0^0$ of $S_0$, having level-$l^N$ structures for any $N\in\dbN$ (with $l$ a rational prime relatively prime to $p$), extends to an abelian scheme over a finite integral cover of $S_0^0$, and so $S_0$ is an almost healthy normal scheme over $O_{(v)}$. For checking this, we can assume that $\scrA$ is defined over the field of fractions $K_1$ of a finite flat DVR extension $O_1$ of $O$. The Galois-representation on $H_{\text{\'et}}^1$$(\scrA_{K_1},\dbZ_l)$ has an image a $Q^2$-torsion group and so it has a finite image (cf. [Se, 1.3] and  the structure of $l$-adic Lie groups). So the N\'eron--Ogg--Shafarevich criterion applies to get that $\scrA$ extends to an abelian scheme over $O_1$.

Due to the EEP enjoyed by $\scrM_{V_0}$, we get a morphism $S_0\to \scrM$. This implies that $q$ extends to a morphism $S_1\to\scrM_1$, and so, provided the morphism $\scrM\to\scrM_1$ is a pro-\'etale cover, $S_0$ is a pro-\'etale cover of $S_1$ and 
$\scrM_1$ has the EP. 

{\bf E)} In case a), $Q^2$ is relatively prime to $p$. So the smoothness of $\scrM_1$ is a consequence of 3.4.5.1 and of [Mi4, 4.11 and 4.13]. 

In case b) for checking the smoothness of $\scrM_1$ we have to work harder. Let $\scrM_2$ be the integral canonical model of a quadruple $(G_2,X_2,H_2,v_2)$, with $v_2$ a prime of $E(G_2,X_2)$ dividing the same prime of $E(G^{\ad},X^{\ad})=E(G_2^{\ad},X_2^{\ad})$ as $v$ (cf. 5.8.1 and 5.8.2). We choose a $\dbZ_{(p)}$-lattice $L_p$ of $W$ such that there is a family of tensors of degrees not bigger than $2(p-2)$ and situated in $\dbZ_{(p)}$-modules of the form $(L_p\oplus L_p^*)^{\otimes m}$ ($m\in\dbN$), which is $\dbZ_{(p)}$-very well positioned with respect to $\psi$ for $G_2$. We can assume that $H_2=G_2(L\otimes\dbZ_p)$ (cf. 3.2.7.1). 

We can choose the connected component $\scrC^0$ of $\scrM_{V_0}$ such that over an embedding of $V_0$ into $\dbC$, its set of complex points contains those defined by equivalence classes of the form $[x,1]$, with $x$ running through the points of a fixed connected component $X^0$ of $X$ (cf. 3.3). The Lemma 6.2.3 allows us to identify $\scrC^0$ with the connected component $\scrC_2$ of $\scrM_{2V_0}$ of whose complex points (under the same embedding of $V_0$ in $\dbC$) are defined by equivalence classes of the same form $[x_2,1]$, with $x_2$ running through the points of a connected component of $X_2$ which can be (cf. 3.3.3) identified with $X^0$.

{\bf F)} We can assume $G_1$ is adjoint (the action of a subgroup of a group acting freely, is free); so $(G_1,X_1)=(G_2^{\ad},X_2^{\ad})$. Based on 6.2.3 and 6.4.2 we can assume $G_1$ is $\dbQ$--simple. The only cases of b) not covered by a) are those in which $p$ is odd and $G_2^{\ad}$ is of $A_{pm-1}$ Lie type. So based on 6.6.5.2 we can assume there is a $\dbZ_{(p)}$-subalgebra $\scrB_2$ of $\End(L_p)$ such that the Zariski closure $G_{2\dbZ_{(p)}}$ of $G_2$ in $GL(L_p)$ is the subgroup of $GSp(L_p,\psi)$ fixing each element of $\scrB_2$.

We consider an injective map $(G_2,X_2)\hookrightarrow (G_2^\prime,X_2^\prime)$, with $G_2^\prime$ as the subgroup of $GL(W)$ generated by $G_2$ and by the center of the centralizer of $G_2$ in $GL(W)$. The Zariski closure of $G_2^\prime$ in $GL(L_p)$ is a reductive group $G^\prime_{2\dbZ_{(p)}}$ and moreover we get a cover 
$$q_2^\prime:(G_2^\prime,X_2^\prime)\to (G_2^{\ad},X_2^{\ad});$$ 
the first part can be seen immediately inside $GL(L_p\otimes_{\dbZ_{(p)}} V_0)$, starting from the classification provided by [Ko2, top of p. 375 and p. 395], while the second part is a consequence of the fact that the center $Z(G_2^\prime)$ of $G_2^\prime$ is the group scheme defined by invertible elements of an \'etale $\dbQ$--algebra $AL$ of endomorphisms of $W$. The scheme $\scrM_2$ is an open closed subscheme of the integral canonical model $\scrM_2^\prime$ of the Shimura quadruple $(G_2^\prime,X_2^\prime,G^\prime_{2\dbZ_{(p)}}(\dbZ_p),v_2)$ (cf. 6.2.3 and 3.2.15).

So, based on [Mi4, 4.13], we can assume that a connected component of $\scrM_{1V_0}$ is the quotient of $\scrC^0=\scrC_2$ by a group of automorphisms $GA$ of $\scrC_2$ which are defined by right translation by elements of a subgroup of the group $GR$ of $\dbA_f^p$-valued points of the center $Z(G_2^\prime)$ of $G_2^\prime$; based on [De2, 2.1.12], in fact we can replace (cf. also 3.3.1) $GR$ by the familiar group of the class field theory
$$Z(G_2^\prime)(\dbA_f^p)/Z(G_2^\prime)(L_p).$$
We can assume $AL\subset\scrB_2\otimes_{\dbZ_{(p)}}\dbQ$.

We consider an element $h\in Z(G_2^\prime)$ defining an element of $GA$. We assume it fixes a point $y\in\scrC_2(\dbF)$. $y$ gives birth to a quadruple 
$$Q_y=(A_y,p_{A_y},\scrB_2,(k_N)_{N\in\dbN,\,(N,p)=1}),$$ 
where $(A_y,p_{A_y})$ is a principally polarized abelian variety over $\dbF$ of dimension ${1\over 2}\dim_{\dbQ}(W)$, endowed with a family of $\dbZ_{(p)}$-endomorphisms (still denoted by $\scrB_2$) and having (in a compatible way) level-$N$ symplectic similitude structure $k_N$, $N\in\dbN$ with $(N,p)=1$; $Q_y$ satisfies some axioms, cf. the standard interpretation of $\scrM_2$ as a moduli scheme (to be compared also with [Ko2, Ch. 5]; see 4.1 for the rational context). We have a similar modular interpretation for $\dbF$-valued and $V_0$-valued points of $\scrM_2^\prime$, provided we work in a $\dbZ_{(p)}$-context; in such a context we speak about principally $\dbZ_{(p)}$-polarized abelian varieties, $\dbZ_{(p)}$-isomorphisms (i.e. isomorphisms up to $\dbZ_{(p)}$-isogenies), etc. So the translation of $y$ by $h$ gives birth to a similar quadruple 
$Q_y^\prime=(A_y^\prime,p_{A^\prime_y},\scrB_2,(k_N^\prime)_{N\in\dbN,\,(N,p)=1}),$
with $(A_y^\prime,p_{A^\prime_y})$ a principally $\dbZ_{(p)}$-polarized abelian variety over $\dbF$ which is $\dbZ_{(p)}$-isogenous to $(A_y,p_{A_y})$. Due to this $\dbZ_{(p)}$-isogeny, we can identify $H^1_{crys}(A_y^\prime/V_0)$ with $M:=H^1_{crys}(A_y/V_0)$.
\smallskip
The fact that $h$ fixes $y$ means that these quadruples are isomorphic, under a $\dbZ_{(p)}$-isomorphism $a:A_y\tilde\to A_y^\prime$. The automorphism $a_M$ of $M$ we get (via $a$ and the mentioned identification), as an element of ${\End}(M)$, belongs to the $\dbQ$--vector space generated by crystalline realizations of $\dbQ$--endomorphisms of $A_y$ defined naturally by elements of ${\Lie}(Z(G_2^\prime))$: this can be read out from the \'etale context with $\dbQ_l$-coefficients, where $l$ is an arbitrary prime different from $p$. So $a_M$ leaves invariant any Hodge filtration of $M$ defined (as $F^1$ in (5.3.11)) by a $V_0$-valued point $z$ of $\scrC_2$ lifting $y$. Such a lift is determined by the mentioned filtration, cf. the deformation theory (see [Me, Ch. 4 and 5]) of polarized abelian varieties endowed with endomorphisms; see also 5.6.4. So, based on the modular interpretation of $\scrM_2$, we get: $h$ fixes all these lifts. So $h$ acts trivially on $\scrC_2$. We conclude: $\scrC_2$ is a pro-\'etale cover of $\scrC_2/GA$. This ends the proof of b) and so of the Theorem.

\medskip
{\bf G)} The use of $\dbZ_{(p)}$-isogenies in F) can be entirely avoided. This goes as follows. Let $G_3:=G_2^{\ad}\times\dbG_m$; identifying $\dbG_m$ with the quotient of $G_2$ by its subgroup $G_2^0$ fixing $\psi$, we get naturally an epimorphism $q_2:G_2\twoheadrightarrow G_3$. Let $(G_3,X_3)$ be the Shimura pair such that $q_2$ defines a finite map 
$$q_2:(G_2,X_2)\to (G_3,X_3).$$  

\medskip\noindent
{\bf Fact.} {\it Let $A:={\Ker}(q_2)$. For any field $l$ of characteristic $0$, the group $H^1({\Gal}(l),A(\bar l))$ is a $2$-torsion group (which in general --like for $l=\dbR$-- is non-trivial).}

\medskip
\proof: From the structure of $G_2^0$ (for instance, see [Ko2, Ch. 5 and 7]) we get that $A$ is a product of Weil restriction of scalars from some totally real number fields to $\dbQ$ of rank one tori which over $\dbR$ are compact (if $\scrB_2\otimes_{\dbZ_{(p)}} \dbQ$ is a simple $\dbQ$--algebra, then the mentioned product has only one factor). The Fact follows.

\medskip
So $q_2$ is not a cover but from the point of view of free actions (see 3.4.5.1) it is ``close enough". In other words, the image of $\scrC_{2\dbC}$ in the quotient of ${\Sh}_{H_2}(G_2,X_2)_{\dbC}$ by $A(\dbA_f^p)$ is a (potentially infinite) Galois cover of the image of $\scrC_{2\dbC}$ in ${\Sh}_{H_3}(G_3,X_3)_{\dbC}$, whose Galois group is a $2$-torsion group, cf. also 3.3.1; here $H_3:=H_2^{\ad}\times\dbG_m(\dbZ_p)$ and $\scrC_{2\dbC}$ is obtained via extension of scalars through an arbitrary $O_{(v_2)}$-monomorphism $V_0\hookrightarrow\dbC$. So, referring to F), eventually by replacing $h$ by $h^2$, we can assume $h\in G_2(\dbA_f^p)$. So F) can be performed ``in terms" of $q_2$ and not of $q_2^\prime$. However, the context of $q_2^\prime$ is more convenient for generalizations (see below).  

\smallskip
We refer to 6.2.2 D) and E).

\Proclaim{6.2.2.1. Proposition.}
Let $g\in\Aut((G_2,X_2,H_2))$ and let $H_2^p$ be a compact subgroup of $G_2(\dbA_f^p)$ such that $g$ belongs to the normalizer of $H_2^p\times H_2$ in $G_2(\dbA_f)$ and $\scrM_2/H_2^p$ is smooth over $O_{(v_2)}$. We assume that the universal principally polarized abelian scheme over $\scrM_2$ obtained through the map $f_2$ and lattice $L$ (cf. 5.1.2), descends to a principally polarized abelian scheme over $\scrM_2/H_2^p$, having a level-$N$ symplectic similitude structure for some $N\in\dbN$, $N\Ge 3$ and relatively prime to $p$ (i.e. we assume that $H_2^p$ is small enough). We also assume that $g$ is inner (i.e. its image in $Aut(G_2^{\ad})(\dbQ)$ belongs to $G^{\ad}_2(\dbQ)$), that a power of $g$ acts trivially on $\scrM_2/H_2^p$ and that $p$ does not divide the torsion number $t(G_2^{\ad})$ (defined in 2.11.1). If $g$ fixes an $\dbF$-valued point $y$ of $\scrM_{2V_0}/H_2^p$, then it fixes a $V_1$-valued point of $\scrM_{2V_0}/H_2^p$ specializing to $y$, with $V_1$ a DVR finite flat extension of $V_0$. 
\finishproclaim   

\proof
We need to show that $g$ does not act freely on the generic fibre of the local ring of $y$ in $\scrM_{2V_0}/H_2^p$. From 3.4.5.1 we deduce that we can assume that $g^p$ acts trivially on $\scrM_{2V_0}/H_2^p$. 

Let $(M_0,\varphi_0)$ and $(\grg_0,\varphi_0)$ be the (Shimura) $\sigma$-crystal and respectively the Shimura adjoint Lie $\sigma$-crystal attached to $y$ (and  the map $f_2$) (cf. 5.4.6; the assumption that the universal abelian scheme over $\scrM_2$ descends to $\scrM_2/H_2^p$ allows us to define them as in 5.4.6). Here $\grg_0$ is the Lie algebra of an adjoint group $G_{2V_0}^{\ad}$ whose generic fibre is $G_{2K_0}^{\ad}$ (cf. 5.4.6). Writing $G_{2V_0}^{\ad}$ as a product of simple adjoint groups, $\varphi_0$ permutes cyclically the Lie algebras of these factors. This allows us to write $(\grg_0,\varphi_0)$ as a product of whose factors correspond to the cycles of the permutation (of the set of simple factors of $G_{V_0}^{\ad}$) we get. We group together the factors of  this product whose Lie algebras are not included in the $F^0$-filtration defined by an arbitrarily chosen $V_0$-lift $z_0$ of $y$ (cf. 5.6.4). We get what we call the non-trivial part $(\grg_0^{nt},\varphi_0)$ of the (Shimura) adjoint Lie $\sigma$-crystal $(\grg_0,\varphi_0)$ (we still denote by $\varphi_0$ its restriction to $\grg_0^{nt}\fracwithdelims[]1p$). Let $G_{2V_0}^{\adnt}$ be the factor of $G_{2V_0}^{\ad}$ whose Lie algebra is $\grg_0^{nt}$. Let $G_{2V_0}^{\adnc}$ be its factor whose simple factors have the property that their Lie algebras are not included in the $F^0$-filtration of $\grg_0$ defined by $z_0$. Let $P^{\adnt}_{2\dbF}$ be the parabolic subgroup of $G_{2\dbF}^{\adnt}$ whose Lie algebra is the natural $F^0$-filtration of $\grg_0^{nt}/p\grg_0^{nt}$. Let $P_{2V_0}$ (resp. $P_{2V_0}^{\adnt}$) be the parabolic subgroup of $G_{2V_0}$ (resp. of $G_{2V_0}^{\adnt}$) leaving invariant the $F^1$-filtration of $M_0$ defined by the  chosen $V_0$-lift $z_0$ of $y$ (resp. defined as the image of $P_{2V_0}$ in $G_{2V_0}^{\adnt}$ under the canonical quotient homomorphism $G_{2V_0}\to G_{2V_0}^{\adnt}$). 
For a presentation of this in a more general and adequate context cf. [Va2].

\Proclaim{The key fact is:}
$g$ gives birth to an isomorphism $g_0$ of $(\grg_0^{nt},\varphi_0)$, with $g_0^p$ acting trivially.
\finishproclaim

For checking this let $(G_2,X_2,H_2)\to (G_2\times G_2,X_2\times X_2,H_2\times H_2)$ be the map defined by the inclusion of $G_2$ into $G_2\times G_2$ whose composite with the two projections of $G_2\times G_2$ are the identity and respectively the automorphism $g$ of $G_2$. It factors through a Hodge quasi product $(G_3,X_3,H_3)$ of $(G_2,X_2,H_2)$ with itself (to be compared with Example 3 of 2.5, where this is detailed for pairs). Composing this factorization with a Segre embedding we get a map $f_3\colon (G_2,X_2,H_2)\to (GSp(W\oplus W,\psi\oplus\psi),S_2,GSp((L_p\oplus L_p)\otimes\dbZ_p))$ which is still a good embedding with respect to $p$. Using the fact that $f_3$ factors through $(G_3,X_3)$ we deduce that the (Shimura) adjoint Lie $\sigma$-crystal $(\grg_1,\varphi_1)$ attached to $y$ (and the map $f_3$) is a Lie subcrystal of the product of $(\grg_0,\varphi_0)$ with itself. As above we define $(\grg_1^{nt},\varphi_1)$. Moreover the first projection (of $G_3$ on $G_2$) allows us to identify $(\grg_1^{nt},\varphi_1)$ with $(\grg_0^{nt},\varphi_0)$, while the second projection gives us the desired isomorphism $g_0$ of $(\grg_0^{nt},\varphi_0)$.

So $g_0$ can be viewed as an element of $G_{2V_0}^{\adnt}(V_0)$ acting on its Lie algebra by conjugation; to see why $g_0$ it is not an outer automorphism of $\grg_0^{nt}$ we just have to remark that:

\smallskip
-- it leaves invariant the simple factors of $G^{\adnc}_{2V_0}$ (this can be seen moving to $\dbC$: $X^0$ is a product of simple Hermitian symmetric domains, indexed by the simple factors of $G^{\adnc}_{2V_0}$; $g$, as an automorphism of $X^0$, is a product of automorphisms of such factors of $X^0$);

-- it commutes with $\varphi_0$ (and so it leaves invariant $P_{2\dbF}^{\adnt}$). 
\smallskip
Moreover $g_0^p$ belongs to any parabolic subgroup of $G_{2V_0}^{\adnt}$ lifting $P_{2\dbF}^{\adnt}$ (as $g^p$ acts trivially on $\scrM_2/H_2^p$). This implies  that the components of $g_0^p$ corresponding to the non-compact simple factors of $G_{2V_0}^{\adnt}$ (i.e. to simple factors of $G_{2V_0}^{\adnc}$)  are trivial. As $g_0^p$ commutes with $\varphi_0$ we deduce that $g_0^p$ is the identity element of $G_{2V_0}^{\adnt}(V_0)$. But $g_0$ is not identity (as otherwise $g$ fixes the connected component of $\scrM_{2V_0}/H_2^p$ through which $f$ factors). So $p$ divides $t(G_2^{\ad})$. Contradiction. This ends the proof.

\Proclaim{6.2.2.1.1. Corollary.} 
We assume $G^{\ad}_{V_0}(V_0)=G^{\ad}_{2V_0}(V_0)$ has no element of order $p$. If $g$ does not act trivially on $\scrM_2/H_2^p$ but $g^p$ does, then $g$ does not fix any $\dbF$-valued point $y$ of $\scrM_{2}/H_2^p$.
\finishproclaim

\Proclaim{6.2.2.2. Corollary.} 
We assume $G^{\ad}_{V_0}(V_0)=G^{\ad}_{2V_0}(V_0)$ has no element of order $p$. Let $H^{2\ad}$ be a compact open subgroup of $G_2^{\ad}(\dbA_f^p)$ such that $H^{2\ad}\times H_2^{\ad}$ is smooth for $(G_2^{\ad},X_2^{\ad})$. Let $\scrM_2^{\ad}$ be the integral canonical model of $(G_2^{\ad},X_2^{\ad},H_2^{\ad},v_2^{\ad})$. Then $\scrM_2^{\ad}$ is a pro-\'etale cover of $\scrM_2^{\ad}/H^{2\ad}$.
\finishproclaim

\proof: Let $H^2$ be a compact open subgroup of $G_2(\dbA_f^p)$ normalized by $H^{2\ad}$ and such that $\scrM_2$ is a pro-\'etale cover of $\scrM_2/H^2$. Corollary follows once we remark that each connected component of $\scrM_{2V_0}^{\ad}/H^{2\ad}$ is the quotient of a connected component of $\scrM_{2V_0}/H^2$ by a group of automorphisms defined by elements of ${\Aut}((G_2,X_2,H_2))$ (this can be seen over $\dbC$, cf. [De2, 2.1.7]).

\smallskip
From now on we assume for the sake of simplicity that $p>2$.

\Proclaim{6.2.3. Lemma.} 
Let $(G_i,X_i,H_i,v_i)$, $i=\overline{1,2}$, be two quadruples with $G^{\der}=G_1^{\der}$ and such that they have the same adjoint quadruple $(G_0,X_0,H_0,v_0)$. Let $p$ be the rational prime divided by $v_0$. Then $\Sh_p(G_1,X_1,H_1)$ exists and has the EP iff  $\Sh_p(G_2,X_2,H_2)$ exists and has the EP. Assuming the existence of these integral models, the connected components of the extension to $O_{(v_1)}^{\text{sh}}$ of the integral canonical model of $(G_1,X_1,H_1,v_1)$ are isomorphic to the connected components of the extension to $O_{(v_2)}^{\text{sh}}=O_{(v_1)}^{\text{sh}}$ of the integral canonical model of $(G_2,X_2,H_2,v_2)$.
\finishproclaim

\proof
We can assume that we have a finite map $f\colon (G_1,X_1,H_1)\to (G_2,X_2,H_2)$ (cf. rm. 3) of 3.2.7). We first assume that $\Sh_p(G_2,X_2,H_2)$ exists and has the EP. Using the toric part triple of $(G_1,X_1,H_1)$ we can assume (cf. 3.2.8) that $f$ is injective. So $\Sh_{H_1}(G_1,X_1)$ is an open closed subscheme of $\Sh_{H_2}(G_2,X_2)$, cf. 3.2.14 and 3.2.15. As $E(G_1,X_1)=E(G_2,X_2)$, we deduce that the Zariski closure of $\Sh_{H_1}(G_1,X_1)$ in $\Sh_p(G_2,X_2,H_2)$ is the integral canonical model $\Sh_p(G_1,X_1,H_1)$. Obviously $\Sh_p(G_1,X_1,H_1)$ has the EP. 

We assume now that $\Sh_p(G_1,X_1,H_1)$ exists and has the EP. Let $E(G_i,X_i)_{(p)}$ be the normalization of $\dbZ_{(p)}$ in $E(G_i,X_i)$, $i=\overline{1,2}$. From [Mi3, 4.7] we deduce that the affine scheme $\Spec(E(G_i,X_i)_{(p)})$ is an \'etale cover of $\Spec(\dbZ_{(p)})$. Let $\scrC$ be a connected component of the image of the natural morphism $m\colon \Sh_{H_1}(G_1,X_1)\to\Sh_{H_2}(G_2,X_2)$. Let $\scrH$ be the subgroup of $G_2(\dbA_f^p)$ leaving invariant $\scrC$. From 3.3.2 we deduce that it is enough to show that $\scrC$ is the generic fibre of a regular formally smooth $E(G_2,X_2)_{(p)}$-scheme $\scrC_p$ having the EP, and on which $\scrH$ acts continuously so that the resulting $\scrH$-action on $\scrC$ is the natural one, and there is a compact open subgroup $\scrH_0$ of $\scrH$ such that $\scrC_p$ is naturally a pro-\'etale cover of the smooth quasi-compact $E(G_2,X_2)_{(p)}$-scheme $\scrC_p/\scrH_0$. As $\Sh_p(G_1,X_1,H_1)$ exists we deduce the existence of $\scrC^\prime_p$, defined as $\scrC_p$, but working over $E(G_1,X_1)_{(p)}$ instead of over $E(G_2,X_2)_{(p)}$; it is an open closed subscheme of $\Sh_p(G_1,X_1,H_1)$. Let $\Spec(E_{(p)})$ be the Galois extension of $\Spec(E(G_2,X_2)_{(p)})$ generated by $\Spec(E(G_1,X_1)_{(p)})$. Let $C:\Gal(\Spec(E_{(p)})/\Spec(E(G_2,X_2)_{(p)}))$ be the resulting Galois group. Due to the EP enjoyed by the extension $\scrC^{''}_p$ of $\scrC^\prime_p$ to $E_{(p)}$, we have a natural Galois-descent datum: $C$ acts on $\scrC^{''}_p$. The extension of $m$, viewed as an $E(G_2,X_2)$-morphism, to $K_0$ identifies each connected component of $\Sh_{H_1}(G_1,X_1)\times_{E(G_2,X_2)} K_0$ with a connected component of $\Sh_{H_2}(G_2,X_2)_{K_0}$, cf. 3.2.14, 3.2.15 and the fact that each connected component of $\Sh_{H_1}(G_1,X_1)_{E(G_2,X_2)} K_0$ is geometrically connected over $K_0$ (as $\Sh_p(G_1,X_1,H_1)$ exists). This together with [Mu1, p. 112] implies that the Galois-descent datum is effective, and so that $\scrC_p$ exists: it has the EP as $\scrC^\prime_p$ has it and as $\Spec(E(G_1,X_1)_{(p)})$ is an \'etale cover of $\Spec(E(G_2,X_2)_{(p)})$ (so B) of 3.2.2 4) applies).

The last part of the Lemma involving connected components over $O_{(v_2)}^{\text{sh}}$ is trivial. This proves the Lemma.    

\Proclaim{6.2.3.1. Remark.} \rm
From the proofs of 6.2.2 and 6.2.3 we deduce that for any finite map $(G^1,X^1,H^1)\to (G^2,X^2,H^2)$ a connected component of $\Sh_{H^1}(G^1,X^1)_{\dbC}$ is a Galois cover of a connected component of $\Sh_{H^2}(G^2,X^2)_{\dbC}$, with a Galois group which is an $M$-torsion abelian pro-finite group, with $M$ equal to the second power of the least common multiple of orders of elements of the center of the simply connected group cover of $G^{2\der}_{\dbC}$ (we can assume that $G^2$ is an adjoint group and that $G^{1\der}$ is simply connected; now everything results from the Step C) of the proof of 6.2.2).
\finishproclaim   

\Proclaim{6.2.4. Corollary.}
Let $(G,X,H)$ be a triple having an integral canonical model $\scrM$. We assume that it has the EP, and that its extension to $V_0$ has the EEP. We also assume that either 

a) the prime $p$ (such that $H\subset G(\dbQ_p)$) is relatively prime to the order of the center of the simply connected semisimple group cover of $G^{\der}$ and $\scrM$ is a quasi-projective integral model, or

b) there is a pair $(G_2,X_2)$ for which condition b) of 6.2.2 is satisfied.

Then any other triple $(G_1,X_1,H_1)$ such that $(G^{\ad},X^{\ad})=(G_1^{\ad},X_1^{\ad})$ and there is an isogeny $G^{\der}\to G_1^{\der}$, has also an integral canonical model $\scrM_1$ having the EP.
\finishproclaim

\proof
This is a direct consequence of 6.2.2, 6.2.3 and 3.2.7 10). 

\Proclaim{6.2.4.1.* Corollary.}
Under the assumptions of 6.2.2 b), any integral canonical model $\scrM_3$ of a Shimura quadruple $(G_3,X_3,H_3,v_3)$ having the same adjoint quadruple as $(G,X,H,v)$ is a quasi-strongly smooth integral model (cf. def. 3.4.8). Moreover, if $p$ does not divide $t(G^{\ad})$, then $\scrM_3$ is in fact a strongly smooth integral model.
\finishproclaim

\proof
We first deal with the case when $p$ does not divide $t(G^{\ad})$. Let $H_1^p\subset H_2^p$ be two compact open subgroups of $G_3(\dbA_f^p)$ such that the morphism $\scrM_3\to\scrM_3/H_1^p$ is a pro-\'etale cover and the generic fibre of the finite morphism $q\colon\scrM_3/H_1^p\to\scrM_3/H_2^p$ is a Galois cover. We need to show that $q$ itself is a Galois cover. This is just a problem of connected components. We use the notations of 6.2.2. So we can move over $V_0$. We can assume that we are dealing with a connected component $\scrC_3$ of $\scrM_{3V_0}$ which over an embedding of $V_0$ into $\dbC$ corresponds to complex points defined by equivalence classes of the form $[x,1]$, with $x$ running through the points of a connected component of $X_3$ (cf. 3.3.2 and 2.3).

We first treat the case when  there is an isogeny $G_2^{\der}\to G_3^{\der}$. 
Using a cover $(G_4,X_4,H_4,v_4)\to (G_3,X_3,H_3,v_3)$, with $G_4^{\der}=G_2^{\der}$, the arguments of [Mi4, 4.11 and 4.13] allow us (cf. 6.2.3, 5.8.1 and 5.8.2) to assume that $G_3^{\der}=G_2^{\der}$. But this case results from 6.2.2.1 (cf. the part of the proof of 6.2.2.2 referring to automorphisms).

To see the general case, the same argument using a cover allows us to assume that $G_3^{\der}$ is the  simply connected group cover of $G^{\der}_2$ (cf. 6.1.2 and 6.2.3). We consider (cf. 3.2.7 10)) a cover $f_5\colon (G_5,X_5,H_5,v_5)\to (G_2,X_2,H_2,v_2)$ such that $G_3^{\der}=G_5^{\der}$. Let $\scrC_5$ be a connected component of the extension to $V_0$ of the integral canonical model of $(G_5,X_5,H_5,v_5)$ dominating $\scrC_2$ and such that its complex points can be described in a similar manner as the complex points of $\scrC_2$ or of $\scrC_3$. We can assume that $H_1^p$ is as small as we want. This together with 6.2.3 allow us to shift our attention to quotients of $\scrC_5$. We get everything in the following context: 

\smallskip
a) we have a compact subgroup $H_{i0}^p$ of $G_i(\dbA_f^p)$, $i\in\{2,5\}$, acting freely on $\scrC_i$ and producing a quotient $\scrC_i/H_{i0}^p$ of finite type; moreover $f_5(H_{50}^p)\subset H_{20}^p$;

\smallskip
b) the natural morphism $\scrC_5/H_{50}^p\to\scrC_2/H_{20}^p$ is an \'etale cover (cf. also 6.1.2);

\smallskip
c) we have a finite group $C(2)$ which is the quotient of a subgroup of the group
$\Aut((G_2^{\ad},X_2^{\ad},H^{\ad}_2))$ leaving invariant $\scrC_i$ and normalizing $H_{i0}^p$, $i\in\{2,5\}$, through a subgroup of it acting trivially on $\scrC_5/H_{50}^p$.    

\smallskip
We need to prove: if $C(2)$ acts freely on the generic fibre of $\scrC_5/H_{50}^p$ then it acts freely on $\scrC_5/H_{20}^p$. This is easy: We can assume that $C(2)$ is a cyclic group of order $p$ (cf. 3.4.5.1); as $C(2)$ also acts on $\scrC_2/H_{20}^p$ such that the \'etale morphism $\scrC_5/H_{50}^p\to\scrC_2/H_{20}^p$ (cf. c)) is $C(2)$-equivariant, the statement follows from 6.2.2.1.1 and from b) above. So $C(2)$ does act freely on $\scrC_5/H_{50}^p$. This ends the proof of the Corollary for the case when $p$ does not divide $t(G^{\ad})$.      
\smallskip
We now assume that $p| t(G_3^{\ad})$ and $H_{03}\times H_3$ is $p$-smooth for $(G_3,X_3)$. We need to show that $\scrM_3$ is a pro-\'etale cover of $\scrM_3/H_{03}$. As this is a problem of connected components of $\scrM_{3V_0}$, we can assume (cf. def. 2.11), that there is a prime $l$ different from $p$ and such that the image of $H_{03}$ in $G_3^{\ad}(\dbQ_l)$ is contained in a compact, open subgroup $H_{03l}^{\ad}$ having no pro-$p$ subgroups. Let $H_{03}^{l{\ad}}$ be a compact, open subgroup of $G_3^{\ad}(\dbA_f^{p,l})$ containing the image of $H_{03}$ in it; here $\dbA_f^{p,l}$ denotes the ring of finite ad\`eles whose both $p$- and $l$-components are omitted.  

As the natural $O_{(v_3)}$-morphism $\scrM_3\to\scrM_{3O_{(v_3)}}^{\ad}$ is pro-\'etale (see 6.4.5 b)), and as the quotient $O_{(v_3)}$-morphism $\scrM_3\to\scrM_3/H_{03}$ factors through the natural morphism $\scrM_3\to\scrM_{3O_{(v_3)}}^{\ad}/\tilde H_{03}^{\ad}$, with 
$$\tilde H_{03}^{\ad}:=H_{03}^{l{\ad}}\times H_{03l}^{\ad},$$ 
we can assume $G_3$ is adjoint. Based on 3.4.5.1, it is enough to show that the morphism $\scrM_3^{\ad}\to\scrM_3^{\ad}/H_{03}^{l\ad}$ is a pro-\'etale cover; here the role of $H_{03}^{l{\ad}}$ is that of an arbitrary compact subgroup of $G_3^{\ad}(\dbA_f^{p,l})$. Based on this and on 3.2.16, we can assume $G_3$ is $\dbQ$--simple. We can assume (cf. 6.5.1.1 i) and the first part of 6.6.5.1):

\medskip
-- the $f_2$ of 6.2.2 b) is such that centralizer $C_{\dbZ_{(p)}}$ of the Zariski closure of $G_2$ in $GL(L_p)$ is reductive; 

-- we have an injective map $(G_2,X_2,H_2,v_2)\hookrightarrow (\tilde G_2,\tilde X_2,\tilde H_2,v_2)$, with $\tilde G_2$ as the subgroup of $GL(W)$ generated by $G_2$ and by the center of the generic fibre $C$ of $C_{\dbZ_{(p)}}$ and with $\tilde H_2=\tilde G_2(\dbQ_p)\cap GL(L_p)(\dbZ_p)$;

-- we are dealing with a subgroup $\tilde H_{03}^{\ad}$ of $\tilde G_2^{\ad}(\dbA_f^p)$ which has the above shape and properties (so in particular, $\tilde H_{03}^{\ad}\times\tilde H_2^{\ad}$ is $p$-smooth for $(\tilde G_2^{\ad},\tilde X_2^{\ad})$).

\medskip
The natural map $(\tilde G_2,\tilde X_2,\tilde H_2,v_2)\to (\tilde G_2^{\ad},\tilde X_2^{\ad},\tilde H_2^{\ad},v_2^{\ad})=(G_2^{\ad},X_2^{\ad},H_2^{\ad},v_2^{\ad})$ is a cover and so we can assume that a connected component of $\scrM_{3V_0}^{\ad}/H_{03}^{l\ad}$ is a quotient of a connected component $\scrC^0$ of the integral canonical model $\scrMtil_{2V_0}$ of $(\tilde G_2,\tilde X_2,\tilde H_2,v_2)$ through a group of automorphisms $GA$ of $\scrMtil_{3V_0}$ leaving invariant $\scrC^0$ and defined by translations by a subgroup of $\tilde G_2(\dbA_f^p)$ whose image in $\tilde G_1^{\ad}(\dbQ_l)$ is trivial, cf. [Mi1, 4.13] and 3.3.1. If $h\in GA$ fixes $y\in\scrC^0(\dbF)$, then as in 6.2.2 F) we get that $h$ acts trivially on $\scrC^0$ (we need to work precisely with our present $l$). Warning: here we dot need to bother about Hodge cycles which are not defined by endomorphisms, i.e. in connection to the $\dbZ_{(p)}$-automorphism $a$ we get (as in 6.2.2 F)) we are bothered just about $\dbZ_{(p)}$-polarizations, level structures (and if one desires, about $\dbZ_{(p)}$-endomorphisms). So $\scrC^0$ is a pro-\'etale cover of $\scrC^0/GA$. This ends the proof.

\medskip\noindent
{\bf 6.2.4.2*. Variant.} What follows is a natural extension of 6.2.2 G) and so provides a variant of the last paragraph of 6.2.4.1; so we can assume $G^{\ad}$ is $\dbQ$--simple. If $(G^{\ad},X^{\ad})$ is of some $A_n$ Lie type with $p|n+1$, then we can proceed as in 6.2.2 G) to get that we can assume that moreover $h\in G(\dbA_f^p)$. If $(G^{\ad},X^{\ad})$ is not of $A_n$ Lie type with $p|n+1$, then $q_{G^{\ad}}$ (see 2.3.5.2) is relatively prime to $p$ and so for any $\tilde h\in G^{\ad}(\dbA_f^p)$ there is $\tilde q\in\dbN$, with $(\tilde q,p)=1$, such that $\tilde h^{\tilde q}$ belongs to the image of $G^{\der}(\dbA_f^p)$ in $G^{\ad}(\dbA_f^p)$. So, as in 6.2.2 G), we can assume that $h\in G(\dbA_f^p)$; we conclude: 

\medskip
{\it Regardless of how $G^{\ad}$ is, in the last paragraph of 6.2.4.1, the use of ``$\dbZ_{(p)}$-" in front of polarizations (and isogenies) can be entirely avoided.}

\Proclaim{6.2.5. Remark.} \rm
There are examples of almost healthy normal schemes which are not noetherian. Such examples can be constructed by taking the normalization of a DVR in an infinite Galois extension of its field of fractions, having a Galois group of finite exponent. 
\finishproclaim

\Proclaim{6.2.6. Remarks.} \rm
{\bf 1)} There are variants for 6.2.2, 6.2.3 and 6.2.4 (which might be useful in the case of Shimura varieties of special type). For instance:

\smallskip
-- in 6.2.2 if we do not assume that $E(G,X)=E(G_1,X_1)$ then we have to work with triples instead of quadruples (to be compared with 6.2.3);

-- in 6.2.3 or 6.2.4 we can work with quadruples but then we either have to restrict to smooth integral models having a weaker extension property (like the WEP or REP) or we need to find extra arguments to be able to shift the EP.

\smallskip
Also there are variants for 6.2.3 and 6.2.4 for $p=2$. The limitations for $p=2$ come only from the fact that we can not presently prove 6.1.2 for $p=2$ and from the the fact that we do not know the uniqueness of an integral canonical with respect to a prime dividing 2 (cf. 3.2.4). These variants will be stated in [Va5].

{\bf 2)*} The integral canonical models of 6.2.4 are quasi-projective as $\scrM$ is so (cf. its proof; see also the proof of 6.4.1).
\finishproclaim

\Proclaim{6.2.7. Warning.} \rm
Any attempt to try to prove 6.1.2 directly (using arguments similar to the ones in 3.4.5.1 and 6.2.2) is meaningless (cf. the two examples below). So we can not handle 6.1.2 just by using geometrically connected components and using 3.2.11 (which provides with $V_0$-valued points). However see 6.8.
\finishproclaim

\Proclaim{\it Example 1.}\rm
Let $Y:=V_0[x]\fracwithdelims[]1{1-px^2}$, and let $Y_1:=Y[y]/(y^2+2pxy+p)$. So $\Spec(Y_1)$
is a finite cover of $\Spec(Y)$, which becomes an \'etale cover by inverting $p$. Moreover the generic fibre of $\Spec(Y_1)$ is geometrically connected over $K_0$. Obviously $Y_1$ is a regular ring which is not an \'etale $Y$-algebra.
\finishproclaim

\Proclaim{\it Example 2.}\rm
Let $Y:=V_0[x]\fracwithdelims[]1{p^{p-1}(1-x)^{p-1}-x^p(p-1)^{p-1}}$ and let $Y_1:=Y[y]/(y^p+pxy+p(1-x))$. The situation is as above. The extra nice thing is that $\Spec(Y_1)$ has plenty of $V_0$-valued points (which is not the case in the above example), as it can be easily checked. 
\finishproclaim

\smallskip
\Proclaim{6.3. Conjugates of integral canonical models of Shimura varieties.}\rm
We use the notations pertaining to conjugates of Shimura varieties used in [Mi1, p. 335-356]. Let $(G,X,H,v)$ be a quadruple having an integral canonical model $\scrM$ over $O_{(v)}$ and let $p$ be the rational prime divided by $v$. Let $\tau$ be an automorphism of $\dbC$ and let $x$ be a special point of $X$. We denote by ${\tau}v$ the prime of ${\tau}E(G,X)$  such that $O_{({\tau}v)}$ is ${\tau}O_{(v)}$. Let ${}^{\tau,x}H$ be the image of $H$ under the isomorphism $G(\dbQ_p)\to {}^{\tau,x}G(\dbQ_p)$ defined by $sp_{p}(\tau)$. It is a hyperspecial subgroup of ${}^{\tau,x}G(\dbQ_p)$.

\Proclaim{6.3.1. Lemma.} 
The integral model ${\tau}\scrM$ is an integral canonical model of $({}^{\tau,x}G,{}^{\tau,x}X,{}^{\tau,x}H,{\tau}v)$ (having EEP if $\scrM$ does).
\finishproclaim

\proof
Here ${\tau}\scrM$ is defined in the same manner as ${\tau}E(G,X)$.
Obviously ${\tau}\scrM$ has the EP. It has the EEP if $\scrM$ does have it. The scheme ${\tau}\scrM$ has a ${}^{\tau,x}G(\dbA_f^p)$-continuous action due to the fact that $\scrM$ has a $G(\dbA_f^p)$-continuous action and due to [Mi1, Ch. 2, 4.2 b) and 5.5 b)]. Using again the loc. cit. and the smoothness of $\scrM$, we get that ${\tau}\scrM$ is also a smooth model (over $O_{({\tau}v)}$). This ends the proof of the Lemma.
\finishproclaim

\smallskip
\Proclaim{6.4. The main results.}\rm

\Proclaim{6.4.1.* Theorem.} 
Let $\Sh(G,X)$ be a Shimura variety of preabelian type. Let $p\Ge 5$ be a prime  such that $G$ is unramified over $\dbQ_p$. Then $\Sh_{p}(G,X)$ exists and has the EP. As a scheme it  is a pro-\'etale cover  of a quasi-projective smooth scheme over (the normalization in $E(G,X)$ of) $\dbZ_{(p)}$.
\finishproclaim

\proof
Let $(G,X,H,v)$ be a quadruple of preabelian type with $v$ dividing a rational prime $p\Ge 5$. From 6.4.2 below we deduce the existence of an injective map $f\colon (G_1,X_1)\hookrightarrow (GSp(W,\psi),S)$ which is a good embedding with respect to $p$ and such that $(G_1^{\ad},X_1^{\ad})=(G^{\ad},X^{\ad})$. We use the notations of the SQSPT introduced in 3.2.7 3). From  3.2.7 2) and 5.8.2 (cf. def. 5.8.1), we deduce that $(G_1,X_1,H_1)$ has an integral canonical model having the EEP. From [Mu, p. 139] and 5.6.2 we deduce that as a scheme it is a pro-\'etale cover of a quasi-projective smooth scheme over $\dbZ_{(p)}$. The statement of 6.1.2 implies that $(G_4,X_4,H_4)$ has an integral canonical model having the EEP, which as a scheme is a pro-\'etale cover of a quasi-projective smooth scheme over $\dbZ_{(p)}$. From 6.2.3 we deduce that $(G_2,X_2,H_2)$ has an integral canonical model which as a scheme is a pro-\'etale cover  of a quasi-projective smooth scheme over $\dbZ_{(p)}$. It has the EP and its extension to $V_0$ has the EEP. From 6.2.2 b) we deduce that $(G,X,H)$ has an integral canonical model $\scrM$.  As the quotient of a quasi-projective smooth scheme through a free action of a finite group is still a quasi-projective smooth scheme (cf. [Mu, p. 112]), we deduce that $\scrM$ is a pro-\'etale cover of a quasi-projective smooth scheme over $\dbZ_{(p)}$. From  3.2.2 4) we deduce that it also has the EP.
This ends the proof of the Theorem.

\smallskip
If $(G,X,H,v)$ is of abelian type then we can use a SQSAT with $G_1^{\der}$ not depending on $i\in\{1,2,3,4\}$ (cf. 3) and 10) of 3.2.7 and 6.4.2). So we can use 6.2.3 (instead of 6.1.2) for concluding that $(G_4,X_4,H_4)$ has an integral canonical model having the EEP and which as a scheme is a pro-\'etale cover of a quasi-projective smooth scheme over $\dbZ_{(p)}$ (as $(G_1,X_1,H_1)$ has an integral canonical model having these properties).

\Proclaim{6.4.1.1. Remarks. 1)*} \rm
From 6.4.1  we deduce that any integral canonical model of a quadruple $(G,X,H,v)$ of preabelian type, with $(v,6)=1$, is a quasi-projective integral model.

{\bf 2)} We refer to 6.4.1 with  $\Sh(G,X)$ of compact type. It is expected that $\Sh_p(G,X)$ is a pro-\'etale cover of a projective smooth scheme over $\dbZ_{(p)}$.

From the proof of 6.4.1 (see also 6.8) we deduce that for seeing this, we can assume that we have an embedding $(G,X)\hookrightarrow (GSp(W,\psi),S)$ good  with respect to $p$. As different quotients of $\Sh_p(GSp(W,\psi),S)$  have (plenty of smooth projective) toroidal compactifications (cf. [FC]) which are moduli of semi-abelian varieties, we deduce that different quotients of $\Sh_p(G,X)$ admit compactifications (obtained by taking the normalization of some Zariski closures in the mentioned toroidal compactifications), which are projective schemes and moduli of semi-abelian varieties. One needs to show that, in our case, these quotients of $\Sh_p(G,X)$ are in fact identical to their compactifications. This is equivalent to showing that over these compactifications we have in fact abelian schemes (and not just semi-abelian schemes). It is expected that this is an easy consequence of [FC, iv) of 10.1, p. 88].${}^1$ $\vfootnote{1}{Using a slightly different approach, in a manuscript to be made available in August 2003 it is checked that $\Sh_p(G,X)$ is a pro-\'etale cover of a projective smooth scheme over $\dbZ_{(p)}$ if each simple factor $(G_0,X_0)$ of $(G^{\ad},X^{\ad})$ is such that $G_{0\dbR}$ has compact, simple factors.}$ 

{\bf 3)} Theorem 6.4.1 fulfills the expectation of [Mi4, 2.17].
\finishproclaim

\Proclaim{6.4.2. Theorem.}
Let $\Sh(G,X)$ be an adjoint Shimura variety of abelian type. Let $p\Ge 5$ be a prime such that $G$ is unramified over $\dbQ_p$. Then there is a Shimura variety $\Sh(G_1,X_1)$ of Hodge type having $\Sh(G,X)$ as its adjoint variety and having a good embedding in a Siegel modular variety with respect to p, and such that for any other Shimura variety $\Sh(G_2,X_2)$ of abelian type having $\Sh(G,X)$ as its adjoint variety, there is an isogeny $G_1^{\der}\to G_2^{\der}$.
\finishproclaim

The proof of 6.4.2 is presented in 6.5 and 6.6.

\Proclaim{6.4.2.1.* Corollary.}
Any integral canonical model $\scrM$ of a Shimura quadruple $(G,X,H,v)$ of preabelian type, with $(v,6)=1$, is a quasi-strongly smooth integral model. If $p$ does not divide $t(G^{\ad})$, then $\scrM$ is a strongly smooth integral model.
\finishproclaim

This is a direct consequence of 6.4.1, 6.4.2 and 6.2.4.1. We would like to remark that if $(G,X)$ is of abelian type then we do not need to use 6.1.2 (cf. the proofs of 6.2.4.1 and 6.4.2). 

This Corollary implies that many other smooth integral models are strongly smooth, cf. 3.4.8.1.

\Proclaim{6.4.2.2.* Corollary.} If in 6.4.2.1 above there is a quadruple $(G_1,X_1,H_1,v_1)$ having the same adjoint quadruple as $(G,X,H,v)$, admitting an embedding $(G_1,X_1,H_1,v_1)\hookrightarrow (GSp(W,\psi),S),K_p,p)$, and such that there is an isogeny $G^{\der}\to G_1^{\der}$, then $\scrM_{O_{(v)}^{\text{sh}}}$ has the EEP.
\finishproclaim

\proof
This is a consequence of 6.2.2 b), 6.2.3 and 6.1.2 (cf. 6.4.1 and the def. of the EEP). We just need to add that in 6.2.2 D) it was essential just that  over $\scrM_2$ we have a principally polarized abelian scheme which is the pull back of a universal one and it did not matter that it is special in the sense of 5.6.7; 3.2.7 4) and 3.2.9 imply that we have a similar principally polarized abelian scheme over the integral canonical model of $(G_1,X_1,H_1,v_1)$. If the pair $(G,X)$ is of abelian type then we do not need to use 6.1.2.  

\Proclaim{6.4.3.} \rm
Let $(G,X)$ define a Shimura variety of preabelian type. Let $\scrS$ be the set of primes whose elements are 2, the primes $p$ for which $G$ is  ramified over $\dbQ_p$, and 3 if $G$ is unramified over $\dbQ_3$ but $\Sh_3(G,X)$ does not exist (if a quadruple $(G_1,X_1,H_1,v_1)$ with $v_1$ dividing a rational prime $p\Ge 3$, has an integral canonical model, then we expect that $\Sh_p(G_1,X_1)$ does exist; this is motivated by rm. 8) of 3.2.7 and by the proof of 5.1, where was irrelevant  with which prime  of the reflex field dividing $p$ we were working). Let $A_f^{\scrS}$ be the ring of finite ad\`eles with all the $q$-components, $q\in\scrS$, omitted. We have $\dbA_f=(\prod_{q\in\scrS} \dbQ_q)\times\dbA_f^{\scrS}$. Let $H^{\scrS}$ be a compact open subgroup of $G(\dbA_f^{\scrS})$ which is a product of its $q$-components (for primes $q\not\in\scrS$) and such that every $q$-component of it is a hyperspecial subgroup $H^q$ of $G(\dbQ_q)$. We call such a subgroup of $G(\dbA_f^{\scrS})$ hyperspecial. It is defined by the property that it is a compact subgroup of $G(\dbA_f^{\scrS})$ of maximal volume (with respect to any Haar measure on $G(\dbA_f^{\scrS})$): this is a consequence of [Ti, p. 55].
\finishproclaim

\Proclaim{6.4.4.* Theorem.}
Let $H_{\scrS}$ be an open subgroup of $G($$\prod_{q\in\scrS} \dbQ_p)$ such that $H_{\scrS}\times H^{\scrS}$ is smooth for $(G,X)$. We assume that $H_{\scrS}\times H^{\scrS}$ is $\scrS_1$-smooth for $(G,X)$, where $\scrS_1$ is the set of rational primes not belonging to $\scrS$ and dividing $t(G^{\ad})$. Then, there is a quasi-projective smooth scheme $\scrM(H_{\scrS})$  over the normalization $O_{(\scrS)}$ of $\dbZ\fracwithdelims[]1{\prod_{q\in\scrS} q}$ in $E(G,X)$, whose generic fibre is $\Sh_{H_{\scrS}\times H^{\scrS}}$$(G,X)$ and such that the normalization ${\Shtil}(G,X)$ of $\scrM(H_{\scrS})$ in the ring of fractions of $\Sh(G,X)$ has the properties:

\smallskip
a) It admits a $G(\prod_{q\in\scrS} \dbQ_q)\times H^{\scrS}$-continuous action;

b) For every prime $q\notin\scrS$, the group $G(\dbQ_q)$ acts continuously on ${\Shtil}(G,X)\times {O_{(\scrS)}\fracwithdelims[]1q}$ and the quotient of ${\Shtil}(G,X)\times\dbZ_{(q)}$ by $H^q$ gets a $G(\dbA_f^q)$-continuous action, together with which it is the integral canonical model of the triple $(G,X,H^q)$. 
\finishproclaim

\proof
It is enough to show that
there is a finite set $\scrS_1$ of rational primes containing $\scrS$ and a quasi-projective  smooth scheme $\scrM_1$ over the normalization $O_{(\scrS_1)}$ of $\dbZ\fracwithdelims[]1{\prod_{q\in\scrS_1} q}$ in $E(G,X)$, whose generic fibre is $\Sh_{H_{\scrS}\times H^{\scrS}}$$(G,X)$, and such that for any prime $p\notin\scrS_1$ the normalization of $\scrM_{1\dbZ_{(p)}}$ in $\Sh_{H^p}(G,X)$ is the integral canonical model of the triple $(G,X,H^p)$: if $q\in\scrS_1\setminus\scrS$, and if $\scrM^q$ is the integral canonical model of the triple $(G,X,H^q)$, then $\scrM^q/H_{\scrS}\times\prod_{p\notin {\scrS\cup\{q\}}} H^p$ is a quasi-projective smooth scheme over the normalization of $\dbZ_{(q)}$ in $E(G,X)$ (cf. 6.4.2.1); but now $\scrM_1$ and $\scrM^q/H_{\scrS}\times\prod_{p\notin {\scrS\cup\{q\}}} H^p$ (for $q\in\scrS_1\setminus\scrS$) can be glued together along their generic fibres. 

Part a) is trivial. We denote by $P(G,X)$ the statement of the existence of a set of rational primes $\scrS_1$ and of a scheme $\scrM_1$ as above for the Shimura pair $(G,X)$. Corollary 6.4.2.1 gives us the right to assume (for proving $P(G,X)$) that $H_{\scrS}$ is as small as desired. 
So the fact that $P(G,X)$ is true for $(G,X)$ of Hodge type  is a direct consequence of the proof of 3.4.7.

We treat now the case when $\Sh(G,X)$ is an arbitrary Shimura variety of preabelian type. Let $\Sh(G_1, X_1)$ be a Shimura variety of Hodge type having $\Sh(G^{\ad},X^{\ad})$ as its adjoint variety. Let $(G_2,X_2)\to (G^{\ad},X^{\ad})$ be a cover with $G_2^{\der}$ a simply connected semisimple group and with $E(G_2,X_2)=E(G^{\ad},X^{\ad})$ (cf. [MS, 3.4]). Let $(G_3,X_3)$ be the fibre product of $(G_1,X_1)$ and $(G_2,X_2)$ over $(G^{\ad},X^{\ad})$ (cf. 2.4.0). 

From 6.2.4.1 and the statement of 6.1.2 we deduce easily that $P(G_3,X_3)$ is true as $P(G_1,X_1)$ is true (i.e. the normalization of a scheme $\scrM_1$ as above, but for $(G_1,X_1)$,
in the ring of fractions of a quotient of $\Sh(G_3,X_3)$ by a subgroup of $G_3(\dbA_f)$ which is smooth for $(G_3,X_3)$, is a smooth scheme over the normalization 
$O_{(\scrS_1)}$ of $\dbZ\fracwithdelims[]1{\prod_{q\in\scrS_1} q}$ in $E(G_3,X_3)$, for $\scrS_1$ a large enough finite set of rational primes). 

We have $G_2^{\der}=G_3^{\der}$ (both are simply connected semisimple groups having the same adjoint group). From 3.2.14 and 3.2.15 (applied to the injective map $(G_3,X_3)\hookrightarrow (G_2,X_2)\times (G_3^{\ab},X_3^{\ab})$ defined by the natural projection of $(G_3,X_3)$ on $(G_2,X_2)$ and by the canonical map $(G_3,X_3)\to (G_3^{\ab},X_3^{\ab})$) we deduce easily that $P(G_2,X_2)$ is true as $P(G_3,X_3)$ is true. 

The proof of 6.2.2 implies that $P(G^{\ad},X^{\ad})$ is true as $P(G_2,X_2)$ is true. 

The same argument used in getting that $P(G_3,X_3)$ is true as $P(G_1,X_1)$ is true, we deduce from 6.4.5 below (applied to the canonical finite map $(G,X)\to (G^{\ad},X^{\ad})$), that $P(G,X)$ is true as $P(G^{\ad},X^{\ad})$ is true. This ends the proof of the Theorem.

\smallskip
As in the proof of 6.4.1, if $(G,X)$ is of abelian type, we do not need to use the statement of 6.1.2 (as we can use instead of it 6.2.3, 3.2.14 and 3.2.15).

\Proclaim{6.4.5. Lemma.}
Let $f^0\colon (G^0,X^0,H^0)\to (G^1,X^1,H^1)$ be a finite map of triples having integral canonical models $\scrM^0$ and respectively $\scrM^1$. We assume that the prime $p$ such that $H^1\subset G^1(\dbQ_p)$ is greater than 2 and that $\scrM^0$ and $\scrM^1$ have the EP. We also assume that either 

a) the order $q$ of the center of the simply connected semisimple group cover of $G^{0\der}$ is relatively prime to $p$ and $\scrM^0$ is a quasi-projective integral model, or

b)* $p\ge 5$ and $(G^0,X^0)$ is of preabelian type, or

c) $\scrM^0$ and $\scrM^1$ are pro-\'etale covers of proper smooth $\dbZ_{(p)}$-schemes. 

Then the natural morphism $\scrM^0\to\scrM^1$ makes $\scrM^0$ to be a pro-\'etale cover of an open closed subscheme of $\scrM^1$, and so $\scrM^0$ is the normalization of $\scrM^1$ in the ring of fractions of $\scrM^0$.
\finishproclaim

\proof
Let $V_0$ be the completion of the strict henselization of $\dbZ_p$. We can move over $V_0$ (i.e. we can shift from triples to quadruples). This is allowed as $\scrM^0$ is a scheme over the normalization of $\dbZ_{(p)}$ in $E(G^0,X^0)$ and as this normalization is an \'etale cover of the normalization of $\dbZ_{(p)}$ in $E(G^1,X^1)$ (cf. [Mi3, 4.7]) over which $\scrM^1$ is defined. Let $v^0$ be a prime of $E(G^0,X^0)$ dividing $p$ and let $v^1$ be the prime of $E(G^1,X^1)$ divided by $v^0$. For $i=\overline{0,1}$, let $\scrM^i_{V_0}$ be the extension to $V_0$ of the integral canonical model of the quadruple $(G^i,X^i,H^i,v^i)$. 

From 6.2.3 and rm. 10) of 3.2.7 we deduce that we can assume that $f^0$ is a cover. So case a) results from 6.2.2. To handle the other two cases we first remark that the normalization $\scrN$ of $\scrM^1_{V_0}$ in the ring of fractions of $\scrM^0_{V_0}$ has local rings of points of codimension 1 isomorphic to local rings of $\scrM^0_{V_0}$ of codimension 1. To see this it is enough (due to the EP enjoyed by $\scrN$ and $\scrM^0_{V_0}$) to check that any such  ring is a DVR. In case c) this is a consequence of [Mi4, 4.13], via the same argument used in $i_B)$ of 3.2.3.2 b). In case b) this is a consequence of 6.4.2.2 and 6.2.2: we can assume that $G^{0\der}$ is simply connected; so, based on 6.4.2.2, the proof of (Steps B), C) and D)) of 6.2.2 applies (it shows the existence of a natural morphism from the spectrum of such a ring into $\scrM^0_{V_0}$; using the natural morphism $\scrM^0_{V_0}\to\scrN$, we get the desired result).

From this and 3.4.5.2 we deduce that $\scrN$ is unramified over $\scrM^1_{V_0}$ in all these points. As $\scrM^0_{K_0}=\scrN_{K_0}$ is a pro-\'etale cover of $\scrM^1_{K_0}$, we deduce from the classical purity theorem that $\scrN$ is a pro-\'etale cover of $\scrM^1_{V_0}$. In particular $\scrN$ is a regular formally smooth scheme over $V_0$ having the EP (cf. C) of 3.2.2 4)). As $\scrM^0_{V_0}$ also has these two properties we get (cf. rm. 7) of 3.2.3.1) $\scrN=\scrM^0_{V_0}$.  
This ends the proof of the Lemma.

\smallskip
The proof of 6.8.1 shows that in fact we can handle the case a) as the other two cases, without reference to the involved 6.2.2, and so without assuming that $\scrM^0$ is a quasi-projective integral model.

\Proclaim{6.4.5.1.* Corollary.}
Let $f\colon (G_1,X_1,H_1,v_1)\to (G_2,X_2,H_2,v_2)$ be a finite map between two quadruples of preabelian type. We assume that $v_1$ is relatively prime to 6. Let $m\colon\scrM_1\to\scrM_2\times O_{(v_1)}$ be the natural morphism (cf. rm. 4) of 3.2.7) defined by $f$. Then $m$ is the composite of a pro-\'etale cover with an open closed embedding. A similar result is true if we work with triples.
\finishproclaim

\Proclaim{6.4.6. Remarks.} \rm
1)* If $(Y,U)$ is an extensible pair with $Y$ a healthy regular scheme over $\Spec(\dbZ\fracwithdelims[]1{\prod_{q\in\scrS} q}$, then any morphism $U\to\scrM(H_{\scrS})$  extends uniquely to a morphism $Y\to\scrM(H_{\scrS})$ (for a proof of this see 6.7). With the terminology to be introduced in [Va6] these schemes $\scrM(H_{\scrS})$ are integral canonical models of their generic fibres.  

2)* These smooth schemes $\scrM(H_{\scrS})$ are the analogue of the schemes attached to Siegel modular varieties which parameterize principally polarized abelian schemes (of a given dimension) and having a finite level symplectic similitude structure. Of course there are variants of 6.4.4 (and of 1)) with $\scrS$ replaced by a larger set of primes (not necessarily finite). But all these variants are a consequence of 6.4.4 (and resp. of 1)).

3)* We call ${\Shtil}(G,X)$ an extended integral canonical model of $\Sh(G,X)$ with respect to $H^{\scrS}$. The scheme ${\Shtil}(G,X)$ is also referred to as an unramified Shimura scheme defined by $(G,X)$. Let ${\text{Hyp}}(G;2)$ be the set of hyperspecial subgroup of $G(\dbA_f^{\scrS})$. 

We assume that $X=X^{\ad}$. Then, as schemes, ${\Shtil}(G,X)$ and $\scrM(H_{\scrS})$ do not depend on the hyperspecial subgroup $H^{\scrS}$ of $G(\dbA_f^{\scrS})$. 

To check this let $H_1^{\scrS}$ be another hyperspecial subgroup of $G(\dbA_f^{\scrS})$.  
It is enough to show the existence of cartesian squares of the form 
$$
\spreadmatrixlines{1\jot}
\CD
\Sh_{K^{\scrS}}(G,X) @>{i_{K^{\scrS}}}>> \Sh_{K_1^{\scrS}}(G,X)\\
@V{r}VV @VV{r_1}V\\
\Sh_{H_{\scrS}\times H^{\scrS}}(G,X) @>{i_H}>> \Sh_{H_{\scrS}\times H_1^{\scrS}}(G,X),
\endCD
$$
where $K^{\scrS}$ (resp. $K_1^{\scrS}$) stands for an arbitrary product of the factors of $H^{\scrS}$ (resp. of $H_1^{\scrS}$), where $r$ and $r_1$ are the natural quotient morphisms, and where $i_H$ and $i_{K^{\scrS}}$ are isomorphisms (cf. rm. 7) of  3.2.3.1). 

If $G$ is a torus then we have nothing to show. 
If $G$ is an adjoint group this is a consequence of 2.3 and of the fact that any two hyperspecial subgroups of $G(\dbA_f^{\scrS})$ are $G(\dbA_f^{\scrS})$-conjugate (cf. [Ti, p. 47]). The same argument works in the case when we have a cover $(G,X)\to (G^{\ad},X^{\ad})$ (as we have epimorphisms $G(\dbQ_l)\twoheadrightarrow G^{\ad}(\dbQ_l)$, for any prime  $l$). 

Using the strong approximation theorem for adjoint groups, we get that for any two hyperspecial subgroups $H^{\scrS}$ and $H_1^{\scrS}$ of $G(\dbA_f^{\scrS})$, there is $g\in G^{\ad}(\dbQ)$ normalizing $H_{\scrS}$ and such that $g(H^{\scrS})g^{-1}=H_1^{\scrS}$; as $X=X^{\ad}$, $g$ takes $X$ into itself and so we can take as $i_{K^{\scrS}}$ and $i_H$ the isomorphisms defined by the inner isomorphism of $G$ defined by $g$.

We now refer to the case when $X\neq X^{\ad}$. Let $N_{X,H_{\scrS}}$ be the subgroup of $G^{\ad}(\dbQ)$ normalizing (under inner conjugation) $X$ and $H_{\scrS}$. It acts (under inner conjugation) on ${\text{Hyp}}(G;2)$. The schemes ${\Shtil}(G,X)$ and $\scrM(H_{\scrS})$ depend only on the orbit $o$ of $H^{\scrS}$ under this action. We have at most ${{c(X^{\ad})}\over {c(X)}}$ such orbits, where $c(*)$ is the number of connected components of the complex manifold $*$.

Warning: (even subject to the restriction $X=X^{\ad}$) the association ${\Shtil}(G,X)$ to $(G,X)$ is not functorial. There are two obstructions to this: the first one is derived from  3.1.2.2 2), while the second one is derived from the fact that $\scrS$ depends on $(G,X)$. However 6.7.2 below is quite enough for many functorial purposes in the context of ${\Shtil}(G,X)$.

4) Lemma 6.4.5 has a variant for quadruples complementing 6.4.5.1: If $(G^0,X^0,H^0,v^0)\to (G^1,X^1,H^1,v^1)$ is  a finite map between two quadruples, with $(v^0,2)=1$, having integral canonical models $\scrM^0$ and respectively $\scrM^1$, and if either a) or c) of 6.4.5 is true, then the natural morphism $\scrM^0\to\scrM^1_{O_{(v^0)}}$ is the composite of a pro-\'etale cover with  an open closed embedding. 

5) In 6.2.4 a) it is enough to assume that $(v,6)=1$ and that $\scrM$ is a quasi-projective integral model. 

To argue this, as in 3.2.7 11), we consider two injective maps $(G_1,X_1,H_1)\hookrightarrow (G,X,H)$ and $(G_2,X_2,H_2)\hookrightarrow (G,X,H)$, such that all simple factors of $(G_1^{\ad},X_1^{\ad})$ and $(G_2^{\ad},X_2^{\ad})$ are of preabelian and respectively of special type, and we have a natural identification $G^{\ad}=G_1^{\ad}\times G_2^{\ad}$. Theorem 6.4.1 (resp. 3.4.1) points out that $(G_1,X_1,H_1)$ (resp. that $(G_2,X_2,H_2)$) has an integral canonical model $\scrN_1$ (resp. has a quasi-projective normal integral model $\scrN_2$ having the EP). As in 3.2.16 we get that $(G_1\times G_2,X_1\times X_2,H_1\times H_2)$ has a quasi-projective normal integral model $\scrN_{12}$ having the EP. Using the fact that the intersection $G_1^{\der}\cap G_2^{\der}$ (taken inside $G^{\der}$) is a finite group scheme of order relatively prime to $p$, from 6.2.3 and (the proof of) 6.2.2 b) we get that $(G,X,H)$ has a normal integral model $\scrN$ over the normalization of $\dbZ_{(p)}$ in $E(G_2\times G_1,X_2\times X_1)$, which is a quotient of $\scrN_{12}$ through a free (see 3.4.5.1) action. From 3.2.12 we get that: if $(G,X,H)$ has an integral canonical model, then $\scrN$ is smooth; so also $\scrN_{12}$ and $\scrN_2$ are smooth. So we can replace $(G,X,H)$ by $(G_1\times G_2,X_1\times X_2,H_1\times H_2)$. So our initial statement follows from 6.2.2, 6.2.3 and 6.4.1 or from 6.2.4 and 6.4.1 (here we need to add: the centers of simply connected  semisimple groups of $E_6$, $E_7$ or $D_l$ Lie type have orders a power of 2 or 3). If we exclude the $E_6$ Lie type, then we can replace $(v,6)=1$ by $(v,2)=1$.

The same applies to 6.2.2 a).  

6)* The philosophy of 6.4.4 is: to generalize Serre' s Lemma [Mu1, p. 207] to the context of Shimura varieties of preabelian type, most common we just have to check things in characteristic zero. 
\finishproclaim

\Proclaim{6.4.7. Remark.} \rm
If $(G,X)$ is the pair $(G_1,X_1)$ of 5.7.5 for $l=10$, then different open subschemes of the schemes $\scrM(H_{\scrS})\times O_{(\scrS)}\fracwithdelims[]1N$ are moduli schemes of polarized  (or just pseudo-polarized) $K3$-surfaces having some finite level-structure (cf. [Va6]).

\Proclaim{6.4.8. Remark.} \rm
For the $p=2$ and $p=3$ theory of Shimura varieties of preabelian type see [Va5] and [Va2]. In [Va2] we prove that 6.4.1 and 6.4.2 remain true for $p=3$. So in 6.4.3 we have $3\in\scrS$ iff $G$ is ramified over $\dbQ_3$.
\finishproclaim

\Proclaim{6.4.9.* Remark.} \rm
We do not know if all  integral canonical models whose existence is guaranteed by 6.4.1 do have the EEP (cf. 3.2.2 4)). However they do have an extension property broader than the EP. This is with respect to any healthy normal scheme (over the required localizations of $\dbZ$) whose local rings of mixed characteristic and of codimension 1 are DVR's (this can be easily checked starting from 6.1, 6.2 and A) of 3.2.2 4)). In fact it is enough that these local rings are certain inductive limits of discrete valuation rings (cf. the proof of 6.2.2; for instance if they are inductive limits of discrete valuation rings whose transition homomorphisms, at the level of fields of fractions, are of degree dividing a fix number $M\in\dbN$).  Similarly for the schemes $\scrM(H_{\scrS})$ constructed in 6.4.4 we have a broader extension property than the one mentioned in rm. 1) of 6.4.6. 
\finishproclaim

From 3.2.12 and 6.4.1 we get directly:

\Proclaim{6.4.10.* Criterion.}
Let $(G,X,H,v)$ be a quadruple of preabelian type, with $(v,6)=1$. Let $\scrM$ be a normal integral model of it over $O_{(v)}$ having the SEP. Then $\scrM$ is the integral canonical model of $(G,X,H,v)$ (in particular $\scrM$ is a smooth integral model and has the EP).
\finishproclaim

\Proclaim{6.4.11. The compact case.} \rm
In this section we assume that the pair $(G,X)$ of 6.4.3 is of compact type and that the expectation of 6.4.1.1 2) has been accomplished; for instance this is the case if each simple factor of $G^{\ad}$ has over $\dbR$ compact, simple factors. 
So $\Sh_{H_{\scrS}\times H^{\scrS}}(G,X)$ is a smooth projective scheme over $E(G,X)$. From 6.4.1.1 2) and 6.4.4 we get directly:
\finishproclaim

\Proclaim{A. Corollary.} 
$\Sh_{H_{\scrS}\times H^{\scrS}}(G,X)$ has good reduction with respect to any prime $v$ of $E(G,X)$ not dividing a prime of $\scrS$. 
\finishproclaim

\noindent
A similar thing can be stated for any connected component $\scrC^g$ of $\Sh_{H_{\scrS}\times H^{\scrS}}(G,X)_{\dbC}$: 

\Proclaim{B. Corollary.} The scheme $\scrC^g$ is naturally defined over a finite field extension $E(\scrC^g)$ of $E(G,X)$ unramified outside $\scrS$, and its canonical model over $E^g$ has good reduction with respect to any prime of $E(\scrC^g)$ not dividing a prime of $\scrS$.
\finishproclaim

\noindent
{\bf C.} Moreover, $\scrM(H_{\scrS})$ is the unique proper smooth scheme over $O_{(\scrS)}$ having $\Sh_{H_{\scrS}\times H^{\scrS}}(G,X)$ as its generic fibre. To see this let $\scrN(H_{\scrS})$ be a proper smooth scheme over $O_{(\scrS)}$ having $\Sh_{H_{\scrS}\times H^{\scrS}}(G,X)$ as its generic fibre. Using the extension type property enjoyed by $\scrM(H_{\scrS})$ (cf. 6.4.6 1)) we deduce the existence of a morphism $l\colon\scrN(H_{\scrS})\to\scrM(H_{\scrS})$ which is the identity on generic fibres. From [Hart, 11.3, p. 279] we deduce immediately that $l$ is an isomorphism. The same thing remains true if instead of $O_{(\scrS)}$ we work with any regular flat $O_{(\scrS)}$-scheme $D$ of dimension 1 such that any smooth $D$-scheme is healthy (see 3.2.2 1)), and if $\scrM(H_{\scrS})$ is replaced by its extension  to $D$: the same proof applies. 

\smallskip
{\bf D.} We can use C to give a second definition of an integral canonical model of a quadruple $(G,X,H,v)$ with $(v,6)=1$: 

\Proclaim{Theorem.} 
An integral model of $(G,X,H,v)$ over $O_{(v)}$ is the integral canonical model of $(G,X,H,v)$ iff it is a smooth proper integral model. 
\finishproclaim

This Theorem answers (slightly restricted) a question of M. Flach.

\smallskip
\Proclaim{6.5. A proof of 6.4.2 in the case when $p$ does not divide $B(G)$.} \rm
\finishproclaim

\Proclaim{6.5.1.} \rm
First we show that to prove 6.4.2 it is enough to treat the case when $G$ is a $\dbQ$--simple group. To check this let $(G,X)$ be a product of two Shimura pairs $(G^i,X^i)$ of adjoint type, $i=\overline{1,2}$, for which 6.4.2 is true. As $G$ is unramified over $\dbQ_p$ we deduce that the group $G^i$ is also unramified over $\dbQ_p$, $i=\overline{1,2}$. Let $(G^i_1,X^i_1)\hookrightarrow (GSp(W^i,\psi^i),S^i)$ be an embedding good with respect to $p$, with $(G^{i\ad}_1,X^{i\ad})=(G^i,X^i)$, and such that for any other Shimura pair $(G^i_2,X^i_2)$ of preabelian type having $(G^i,X^i)$ as its adjoint variety, there is an isogeny $G^{i\der}_1\to G^{i\der}_2$ ($i=\overline{1,2}$). Let $(G^3_1,X^3_1)$ be a Hodge quasi product of the two Shimura pairs $(G^1_1,X^1_1)$ and $(G^2_1,X^2_1)$ of Hodge type (cf. Example 3 of 2.5). Now the Segre embedding $(G^3_1,X^3_1)\to (GSp(W^1\oplus W^2,\psi^1\oplus\psi^2),S^0)$ is a good embedding with respect to $p$ (cf. 4.3.17). Moreover $G^{3\der}_1=G^{1\der}_1\times G^{2\der}_1$. So for any Shimura variety $(G_2^3,X_2^3)$ of abelian type such that its adjoint variety is the adjoint variety of $(G_1^3,X_1^3)$, there is an isogeny $G_2^{3\der}\to G_1^{3\der}$ (cf. [De2, 2.3.8]).

So we can assume that $G$ is a simple $\dbQ$--group.
We deduce the existence of a totally real number field $F$ and of an absolutely simple adjoint group $G^s$ over $F$ such that $G=\Res_{F/\dbQ}G^s$ [De2, 2.3.4]. As before $V_0:=W(\overline{\dbZ/p\dbZ})$. For any number field $E$ we denote by $E_{(p)}$ the normalization of $\dbZ_{(p)}$ in $E$. Let $G_{\dbZ_{(p)}}$ be an adjoint  group over $\dbZ_{(p)}$ having $G$ as its fibre over $\dbQ$ (cf. 3.1.3) and let $\Gtil_{\dbZ_{(p)}}$ be the simply connected semisimple group cover of it. We have:

\smallskip
a)\enspace
The group $G_{V_0}$ is a product of $[F:\dbQ]$ copies of a split adjoint group of the same Lie type as $G$ (this is obvious).

\smallskip
b)\enspace
As $G$ is unramified over $\dbQ_p$, $F$ is unramified over $p$ and $G_{F_i}^s$ is unramified over $F_i$, where $F\otimes\dbQ_p=\prod_{i\in I_p} F_i$, with $F_i$ local fields (we have $G_{\dbQ_p}=\prod_{i\in I_p} \Res_{F_i/\dbQ_p}G_{F_i}^s$).

\smallskip
Proposition [De2, 2.3.10] admits a $\dbZ_{(p)}$-version:

\Proclaim{6.5.1.1. Theorem.}
Let $K$ be a quadratic totally imaginary extension of $F$, unramified over $p$. Then there is a Shimura variety $\Sh(G_1,X_1)$ of Hodge type such that:

a) $\Sh(G,X)$ is its adjoint Shimura variety;

b) for any Shimura variety $\Sh(\Gtil_1,\Xtil_1)$ of abelian type with $(\Gtil_1^{\ad},\Xtil_1^{\ad})=(G,X)$, there is an isogeny $G_1^{\der}\to\Gtil_1^{\der}$;

c) its reflex field is the composite field of $E(G,X)$ and $E(\Res_{K/\dbQ} \dbG_m,h_T)$ (where $(\Res_{K/\dbQ} \dbG_m,h_T)$ is a 0 dimensional Shimura pair defined as in [De2, 2.3.9]); 

d)  it has a good embedding in a Siegel modular variety with respect to $p$. 
\finishproclaim

\proof
The proof is divided in two parts. First we treat the case when $p$ does not divide $B(G)$, then we continue in 6.6.5 with the general case.
In this section 6.5, the symbols $S$, $K$, $K_S$, $(G_2,X_2)$ and $(G_3,X_3)$ will have the same significance as in [De2, 2.3]. So $S$ is a set of extremal nodes of the Dynkin diagram of $G_{\dbC}$, $(G_2,X_2)$ and $(G_3,X_3)$ are Shimura pairs, while $K_S$ is a product of finite field extensions of $\dbQ$. If $(G,X)$ is of $B_l$, $C_l$ or $D_l^{\dbH}$ type (resp. of $A_l$ or $D_l^{\dbR}$ type), then to each simple factor of $G^{\ad}_{\dbR}$ corresponds one (resp. two) elements of $S$. We itemize the things we need.

\smallskip
i) We start with a self dual representation $W_{(p)}$ of $\Gtil_{\dbZ_{(p)}}$ over $\dbZ_{(p)}$ which over $V_0$ is isomorphic to $\oplus_{s\in S} V_p(s)^n$ for a convenient number $n\in\dbN$ (to be compared with [De2, 2.3.10]). Here $V_{p}(s)$ is the $V_0$-representation of $\Gtil_{V_0}$ given by the fundamental weight corresponding to $s\in S$ (cf. [De2, 2.3]).

\smallskip
ii) The totally imaginary quadratic extension $K$ of $F$ is assumed to be unramified above $p$ (i.e. $\Spec(K_{(p)})$ is an \'etale cover of $\Spec(\dbZ_{(p)}$)).

\smallskip
iii) The \'etale $\dbQ$-algebra $K_S$ is unramified above $p$ as $\Gtil_{\dbZ_{(p)}}$  splits over $V_0$.

\smallskip
iv) The Zariski closure of $G_3$ in $GL(W_{\dbZ_{(p)}})$, with $W_{\dbZ_{(p)}}:=K_{(p)}\otimes_{F_{(p)}} W_{(p)}$, is a reductive group ${G_3}_{\dbZ_{(p)}}$ over $\dbZ_{(p)}$ (cf. [De2, 2.3.10] for the meaning of $G_3$) (moving over $V_0$ this becomes obvious). Let $\Gtil_3^{c+{\der}}$ be the subgroup of $G_3$ generated by $G_3^{\der}$ and by the maximal subtorus of $Z(G_3)$ which over $\dbR$ is compact (cf. [De2, 2.3.3 and end of 2.3.10]). Let $\Gtil_3$ be the subgroup of $G_3$ generated by $\Gtil_3^{c+{\der}}$ and by the one dimensional split torus acting as scalar multiplication on 
$$
W:=W_{\dbZ_{(p)}}\otimes\dbQ.
$$ 
So any homomorphism $\dbS\to G_{3\dbR}$ defined by some $x\in X_3$ factors through $\Gtil_{3\dbR}$ (of course instead of $\Gtil_3$ we can work equally well with the smallest subgroup of $G_3$ satisfying this property). We get a Shimura pair $(\Gtil_3,\Xtil_3)$; here $\Xtil_3$ is a disjoint union of connected components of $X_3$ defined by a $\Gtil_3(\dbR)$-conjugacy class of an arbitrary $x\in X_3$. This is a slight restatement of [De2, 2.3.3]: we do not always have $\Xtil_3=X_3$, as it can be seen easily (to be compared with 2.5.1) through examples in which $F$ is a totally real quadratic extension of $\dbQ$.  

Let $\Gtil_{3\dbZ_{(p)}}$ (resp. $\Gtil_{3\dbZ_{(p)}}^{c+{\der}}$) be the Zariski closure of $\Gtil_3$ (resp. of $\Gtil_3^{c+{\der}}$) in ${G_3}_{\dbZ_{(p)}}$. 

From loc. cit. we get that $\Gtil_3$ is included in the group of symplectic similitude isomorphisms defined by a non-degenerate alternating form on $W$.

\smallskip
v) There is a perfect alternating form $\psi\colon W_{\dbZ_{(p)}}\otimes W_{\dbZ_{(p)}}\to\dbZ_{(p)}$ such that we get an injective map $f\colon (\Gtil_3,\Xtil_3)\hookrightarrow (G\Sp(W,\psi),S^0)$ (here we write as an exception $S^0$ for what we have always denoted by $S$, not to create confusion with the meaning of $S$ in [De2, 2.3]). 

This is so due to the fact that [De2, 1.1.18 b)] admits a $\dbZ_{(p)}$-version. To see this we first remark that the alternating bilinear forms $W_{\dbZ_{(p)}}\otimes W_{\dbZ_{(p)}}\to\dbZ_{(p)}$ fixed by $\Gtil_{3\dbZ_{(p)}}^{c+{\der}}$ form a free module $M$ over $\dbZ_{(p)}$. Choosing $n$ big enough (see 6.6.5 d) for an explicit presentation) we can assume that we have such bilinear forms which are as well perfect. 

In fact using the natural embedding $SL_m(\dbZ_{(p)})\hookrightarrow Sp_{2m}(\dbZ_{(p)})$ (as in 6.6.5 d1); here $Sp_{2m}(\dbZ_{(p)})$ is the group of symplectic isomorphisms defined by a perfect alternating form on $\dbZ_{(p)}^{2m}$, etc.), $m:=\dim_{\dbZ_{(p)}}(W_{\dbZ_{(p)}})$, we get the existence of such a perfect alternating bilinear form after we replace (if needed) $n$ by $2n$. This replacement corresponds to a replacement of $W_{(p)}$ by $W_{(p)}\oplus W_{(p)}$ and of $W_{\dbZ_{(p)}}$ by $W_{\dbZ_{(p)}}\oplus W_{\dbZ_{(p)}}$ (cf. the way we defined $W_{(p)}$ in i) and the definition of the connected component of the origin of $Z(\Gtil_3^{c+{\der}})$). We would like to point out that this fact is convenient for notations (and so used in what follows) but is irrelevant for what follows: we can work equally well (to be compared with 6.7.2) without having (or knowing) that the representation $\Gtil_{3\dbZ_{(p)}}^{c+{\der}}\to GL(\Wtil_{\dbZ_{(p)}})$ we get under the above natural embedding $SL_m(\dbZ_{(p)})\hookrightarrow Sp_{2m}(\dbZ_{(p)})$ is a sum of two copies of its representation on $W_{\dbZ_{(p)}}$. 

Now we look at $M$  as a group scheme over $\dbZ_{(p)}$. The intersection of a non-empty open (in the real topology) subset of $M(\dbR)$ with the set of $\dbZ_{(p)}$-valued points of the dense open subscheme $M(pa)$ of $M$ corresponding to perfect alternating bilinear forms is not void. Argument: $M(pa)$ has $\dbZ_{(p)}$-valued points; if $\tilde\psi\colon W_{\dbZ_{(p)}}\otimes W_{\dbZ_{(p)}}\to\dbZ_{(p)}$ corresponds to $\tilde z\in M(pa)(\dbZ_{(p)})$, then we can choose $\psi$ such that mod $p$ is $\tilde\psi$ mod $p$ (standard argument involving approximations with respect to non-equivalent valuations).

\smallskip
vi) Using 5.7.4 and 5.6.9 we get that if $p$ does not divide $B(G)$ (see 5.7.2 for the meaning of it), then $(\Gtil_3,\Xtil_3)\hookrightarrow (G\Sp(W,\psi),S^0)$ is a good embedding with respect to $p$. 

For checking this we first remark that we have 
$$W_{V_0}:=W_{\dbZ_{(p)}}\otimes V_0=\oplus_{(i,s)\in I\times S} V_{p}(s)^i$$ 
as ${G_3}_{V_0}^{\der}$-modules, with $I:=\{1,2,...,2n\}$, the upper indices $i$ just counting the numbers of copies of $V_{p}(s)$ we get. Moreover ${G_3}_{V_0}$ leaves invariant any summand of this direct sum decomposition. Let $\grg\grl(W_{V_0})=\grm_0\oplus\grm_1$, with $\grm_0$ the free $V_0$-submodule of $\text{End}(W_{V_0})$ leaving invariant any subspace $V_{p}(s)^i$ of $W_{V_0}$, and with $\grm_1$ the free $V_0$-submodule of $\text{End}(W_{V_0})$ taking, $\forall (i_0,s_0)\in I\times S$, the summand $V_{p}(s_0)^{i_0}$ of $W_{V_0}$ into $\oplus_{(i,s)\in I(i_0,s_0)}V_{p}(s)^i$ (here $I(i_0,s_0):=I\times S\setminus \{(i_0,s_0)\}$). Let $\pi_0$ be the projector of $\grg\grl(W_{V_0})$ on $\grm_0$ associated to the above direct sum decomposition. Now to get the first sentence of vi) we just have to apply 5.7.4 to the bilinear form $\overline{b}$ on $\grg\grl(W_{V_0})$ defined by 
$$\overline{b}(x,y):=\oplus_{(i,s)\in I\times S} {\gamma_{(i,s)}}\Tr_{(i,s)}(\pi_0(x),\pi_0(y)).$$ 
Here $x,y\in\grg\grl(W_{V_0})$, $\gamma_{(i,s)}$ are invertible elements of $V_0$ having all their partial sums still as invertible elements of $V_0$, and $\Tr_{(i,s)}$ is the trace form on $\text{End}(V_{p}(s)^i)$. The element $\Tr_{(i,s)}(\pi_0(x),\pi_0(y))$ of $V_0$ makes sense as $\grm_0=\oplus_{(i,s)\in I\times S} \text{End}(V_p(s)^i)$. Obviously $\overline{b}$ is fixed by ${G_3}_{V_0}$ and so by $\Gtil_{3V_0}$.

This ends the proof of 6.4.2 and 6.5.1.1 in the case when $p$ does not divide $B(G)$ (cf. [De2, 2.3.10 to 2.3.13] for the requirements on $E(\Gtil_3,\Xtil_3)=E(G_3,X_3)$ and on $\Gtil_3^{\der}=G_3^{\der}$ expressed in 6.5.1.1 b) and c)).

\smallskip
\Proclaim{6.6. The proof of 6.4.2 and 6.5.1.1 (the general case).} \rm
We continue to use the same notations as in 6.5. We present two proofs of the general case of 6.4.2: the first one (6.6.3), based on the (sophisticated) Proposition 6.6.2, and a second one (6.6.5) which is a simplified, down to earth, explicit version of the first one.

\Proclaim{6.6.1. Notation.} \rm
For any totally real number field $F_1\supset F$, we denote by $\Sh^{F_1}(G,X)$ the adjoint Shimura variety defined by the pair $(G^{F_1},X^{F_1})$, where $G^{F_1}:=\Res_{F_1/\dbQ}G_{F_1}^s$ and $X^{F_1}$ is the Hermitian symmetric domain obtained as the $G^{F_1}(\dbR)$-conjugacy class of homomorphisms $\dbS\to G^{F_1}_{\dbR}$ generated by the composite of any $x\in X$ with the natural inclusion $G_{\dbR}\hookrightarrow G^{F_1}_{\dbR}$.  So $X^{F_1}$ is a product of $[F_1:F]$ copies of $X$. We get a natural injective map $f_{F_1}\colon\Sh(G,X)\hookrightarrow\Sh^{F_1}(G,X)$. In particular $\Sh^F(G,X)=\Sh(G,X)$.

\Proclaim{6.6.2. Proposition.}
There are injective maps
$$(G_4,X_4)\operatornamewithlimits{\hookrightarrow}\limits^{f_0} (G^0,X^0) \operatornamewithlimits{\hookrightarrow}\limits^{f_1} (G^1,X^1)\operatornamewithlimits{\hookrightarrow}\limits^{f_2} (G\Sp(W,\psi),S^0)$$
having the properties:

\smallskip
\item{a)}
there is a $\dbZ_{(p)}$-lattice $L$ of $W$ such that $\psi$ induces a perfect bilinear form $\psi\colon L\otimes L\to\dbZ_{(p)}$ and  the Zariski closures of $G_4$, $G^0$ and $G^1$ in $G\Sp(L,\psi)$ are reductive groups over $\dbZ_{(p)}$ denoted respectively by ${G_4}_{\dbZ_{(p)}}$, ${G^0}_{\dbZ_{(p)}}$ and ${G^1}_{\dbZ_{(p)}}$;

\smallskip
\item{b)}
$(G_4^{\ad},X_4^{\ad})=(G,X)$ and there is a totally real number field $F_1\supset F$ such that $\Sh(G^{0\ad},X^{0\ad})=\Sh^{F_1}(G,X)$;

\smallskip
\item{c)}
the map $f_0$ induces the canonical homomorphism $f_{F_1}\colon G=G_4^{\ad}\to G^{0\ad}=G^{F_1}$;

\smallskip
\item{d)}
if $(G,X)$ is of $A_l$, $B_l$ or $D_l^{\dbR}$ (resp. of $C_l$ or $D_l^{\dbH}$) type, then $G_{\dbZ_{(p)}}^0$ is contained and has the same derived subgroup as the centralizer in $G_{\dbZ_{(p)}}^1$ of a torus of $G_{\dbZ_{(p)}}^1$ (resp. of a semisimple $\dbZ_{(p)}$-subalgebra $SA$ of $\End(L)$ such that we have a relative standard PEL situation $(G_{1\dbZ_{(p)}},SA)$ as defined in 4.3.16);

\smallskip
\item{e)}
the homomorphism $G^{0\der}\to G^{1\der}$ induced by $f_1$ is of the form $\Res_{F_1/\dbQ} f^{F_1}$ for $f^{F_1}\colon G_{F_1}^d\to\tilde G_{F_1}^1$ a group homomorphism between semisimple groups over $F_1$, with $\tilde G_{F_1}^{1\ad}$ a simple $F_1$-group and with $G_{F_1}^d$ a cover of $G_{F_1}^s$;

\smallskip
\item{f)}
$f_2$ is an injective map obtained through the $\dbZ_{(p)}$-version of [De2, 2.3.10] explained in 6.5.1.1 for the totally imaginary quadratic extension $K\otimes_F F_1$ of $F_1$, with $L=W_{\dbZ_{(p)}}$ and with the number $n$ (mentioned in i) of 6.5.1.1) a power of 2 (so the maximal torus of $Z(G^1)$ is naturally a subtorus of $\Res_{K/\dbQ} \dbG_m$, cf. [De2, 2.3.13]); 

\smallskip
\item{g)}
$p$ does not divide $B(G^{1\ad})$;

\smallskip
\item{h)}
$G_4^{\ab}=G^{0\ab}$;

\smallskip
\item{i)}
if $(G,X)$ is of $D_l^{\dbH}$ type, with $l\in\dbN$, $l\Ge 4$, then the embedding $G_4^{\der}\hookrightarrow G^{1\der}$ can be lifted to an embedding at the level of simply connected semisimple group covers.

\smallskip
Moreover if $(G,X)$ is of $A_l$, $B_l$ or $D_l^{\dbR}$ type, we can also get $E(G^1,X^1)=\dbQ$.
\finishproclaim

\proof
The proof of 6.6.2 presents no difficulty. The statement of the Proposition makes its proof  obvious (cf. also [Va5]). If $(G,X)$ is of $B_l$ (resp. $D_l^{\dbR}$) type, we can take $(G^1,X^1)$ of $B_{l+a}$ (resp. $D_{l+a}^{\dbR}$) type, with $a$ a non-negative integer; if $(G,X)$ is of $C_l$ (resp. $D_l^{\dbH}$) type, we can take $(G^1,X^1)$ of $C_{al}$ (resp. $D_{al}^{\dbH}$) type, with $a\in\dbN$; if $G$ is of $A_l$ Lie type we can take $G^1$ of $C_{l+1}$ Lie type (to be compared with 6.6.5 below). In practice we take the number $a$ to be 0 (when allowed), 1 or 2. We will just add that we need $F_1$ to be a totally real number field, containing $F$, unramified above $p$ and big enough.

For the last property (concerning the cases when we can take $E(G^1,X^1)=\dbQ$) needed for the proof of the Langlands--Rapoport conjecture (of 1.7) see [Va2]. We need 6.6.2 (presently) only for the $p=2$ and $p=3$ theory of Shimura varieties of preabelian type.

\Proclaim{6.6.3. Remark.} \rm
Property 6.6.2 a) implies that $G_4$ and $G^1$ are unramified over $\dbQ_p$. From 5.7.1 and 6.6.2 f) and g) we deduce that the injective map $(G^1,X^1)\hookrightarrow (G\Sp(W,\psi),S^0)$ is a (very) good embedding with respect to $p$. From this, 4.3.14 and 6.6.2 d) we deduce that $(G^0,X^0)\hookrightarrow (GSp(W,\psi,S^0)$ is a good embedding with respect to $p$. Now this, 4.3.16 and b), c), f) and h) of 6.6.2 imply that $(G_4,X_4)\hookrightarrow (G\Sp(W,\psi),S^0)$ is a good embedding with respect to p. This ends the first proof of the general case of 6.4.2.
\finishproclaim
We present now what 6.6.2 becomes in the case of classical Spin modular varieties of odd dimension (and rank 2).

\Proclaim{6.6.4. Example.} \rm
Let $l\Ge 3$ be an integer. Let $\Sh(G_i,X_i)$, $i=\overline{0,1}$, be two adjoint Shimura varieties showing up in 5.7.5, with $G_i=SO(2,2l-1+2i)$. The canonical inclusion $j_0\colon G_0\hookrightarrow G_1$ (corresponding to the identification of the group of invertible matrices of dimension $2l+1$ with the subgroup of invertible matrices of dimension $2l+3$ having on the last two lines and columns just two diagonal 1's) induces an injective map $j_0\colon (G_0,X_0)\hookrightarrow (G_1,X_1)$ and $G_0$ is contained and has the same derived subgroup as the centralizer in $G_1$ of a torus of $G_1$ of dimension 1. 

If $(G,X)=(G_0,X_0)=(G_0^{\ad},X_0^{\ad})$ and if $p$ is a prime not dividing $B(G)=6(2l-1)$, then in 6.6.2 we can take $G_4=G^0=G^1$ and for the map $f_2$ we can take the map associated to the Spin representation  described in 5.7.5. So $G_4^{\ab}=\dbG_m$. If $p\Ge 5$ divides $2l-1$ then in 6.6.2 we can take $G_4=G^0$, the adjoint of $f_1$ to be $j_0$, and as $f_2$ the map associated to the Spin representation of the simply connected group cover of $G_1$. So, regardless of how $p\Ge 5$ is, 4.3.14, 5.1 and 5.7.5 put together imply that $\Sh_p(G_0,X_0)$ exists; so (cf. 6.2.2) $\Sh_p(G_0,X_0)$ exists as well. 
\finishproclaim
     
\Proclaim{6.6.5. An  explicit proof of the above $\dbZ_{(p)}$-version of [De2, 2.3.10].} \rm
Here we present the second part of the proof of 6.5.1.1. Let $T$ be a maximal torus (cf. the argument in 3.1.4 based on [Ha, 5.5.3]) of a simply connected semisimple group $G_{F_{(p)}}^{\sc}$ (cf. 3.1.3) over $F_{(p)}$ having as its fibre over $F$ the simply connected semisimple group cover of $G^s$, such that for any embedding $F\hookrightarrow\dbR$, $T_{\dbR}$ is compact. Then $T_F$ splits over a Galois extension $E$ of $F$ unramified above $p$. Choosing the smallest such Galois extension, we get that $E$ is a CM-field (as $T_{\dbR}$ is a compact torus for any embedding $F\hookrightarrow\dbR$). We need $T$ (and $E$) just to fix a little bit the notations.

We consider homomorphisms (between reductive groups over $E_{(p)}$)
$$G_{E_{(p)}}^{\sc}\operatornamewithlimits{\to}\limits^{h_0} G_{E_{(p)}}^d\operatornamewithlimits{\hookrightarrow}\limits^{h_1}
\Gtil_{E_{(p)}}\operatornamewithlimits{\hookrightarrow}\limits^{h_2}
GL(W_{E_{(p)}})$$
such that the following hold:
 
\smallskip
\item{a)}
$W_{E_{(p)}}$ is a free $E_{(p)}$-module of finite rank.

\smallskip
\item{b)}
The group $\Gtil_{E_{(p)}}$ is semisimple and $\Gtil_{E_{(p)}}^{\ad}$ is a split simple group over $E_{(p)}$ such that $p$ does not divide $B(\Gtil_E^{\ad})$.

\smallskip
\item{c)}
The homomorphism $h_0$ is an isogeny. Here $G_{E_{(p)}}^{\sc}$ is the pull back of $G_{F_{(p)}}^{\sc}$ to $E_{(p)}$.

\smallskip
\item{d)}
The homomorphism $h_1$ is an $E_{(p)}$-version of the homomorphism $f^{F_1}$  mentioned in 6.6.2 e). Namely:

\smallskip
\item{d1)}
If $G^s$ is of $A_l$ Lie type, then we take $W_{E_{(p)}}$ of dimension $2(l+1)$ over
$E_{(p)}$ and we take $h_0$ to be an isomorphism. Let $\psi_0\colon W_{E_{(p)}}\otimes W_{E_{(p)}}\to E_{(p)}$ be a perfect alternating form. We choose  a basis $\{e_1,e_2,...,e_{2l+2}\}$ of $W_{E_{(p)}}$ with respect to which $\psi_0$ has the standard form, i.e. if $1\Le i\Le j\Le 2(l+1)$, then $\psi_0(e_i,e_j)=1$ if $j=i+l+1$ and $0$ otherwise. 
We identify $G_{E_{(p)}}^{\sc}$ with ${SL_{l+1}}_{E_{(p)}}$. We take $h_2\circ h_1$ such  that it takes $A\in {SL_{l+1}}_{E_{(p)}}(E_{(p)})$ into the element of $GL(W_{E_{(p)}})$  that acts as $A$ on the submodule of $W_{E_{(p)}}$ generated by the first $l+1$ elements of the chosen basis and as $(A^t)^{-1}$  on the submodule of $W_{E_{(p)}}$ generated by the last $l+1$ elements of the chosen basis. 
If $p$ does not divide $B(G^s)=6(l+1)$, we take $\Gtil_{E_{(p)}}=G^d_{E_{(p)}}$ (with $h_1$ as identity). If $p$ divides $6(l+1)$, we take $\Gtil_{E_{(p)}}=Sp(W_{E_{(p)}},\psi_0)$, and  $h_1$ and $h_2$ as the obvious monomorphisms (as $p$ does not divide  $B(\Gtil_E^{\ad})=6(l+2)$; we recall that $p\ge 5$).

\smallskip
\item{d2)}
Let now $(G,X)$ be of $D_l^{\dbR}$ type. We take $G_{E_{(p)}}^d=G_{E_{(p)}}^{\sc}=Spin(2l)_{E_{(p)}}$. We take $h_2\circ h_1$ to be the composition of the embedding $Spin(2l)_{E_{(p)}}\hookrightarrow Spin(2l+2)_{E_{(p)}}$ (which results by passage to simply connected group covers of the homomorphism  $SO(2l)_{E_{(p)}}\to SO(2l+2)_{E_{(p)}}$ described in terms of  matrices by the rule: $A\in SO(2l)_{E_{(p)}}({E_{(p)}})$ goes to the matrix having $A$ on the first $2l$ lines and columns and having on the last two lines and columns just two diagonal 1's) with the Spin representation of $Spin(2l+2)_{E_{(p)}}$. If $p$ divides $B(G^s)=6(2l-1)$, we take $\Gtil_{E_{(p)}}=Spin(2l+2)_{E_{(p)}}$ ($B(\Gtil_E^{\ad})=6(2l+1)$) and if $p$ does not divide $B(G^s)$ we take $\Gtil_{E_{(p)}}=G_{E_{(p)}}^d$ (and the obvious homomorphisms $h_1$ and $h_2$).

\smallskip
\item{d3)}
If $G^s$ is of $B_l$ Lie type, then the situation is entirely analogous to the situation described in d2) (to be compared with 5.7.5).

\smallskip
\item{d4)}
Let now $G^s$ be of $C_l$ Lie type. We take $G_{E_{(p)}}^{\sc}=G_{E_{(p)}}^d=Sp(W_{E_{(p)}}^1,\psi_1)$, with $W_{E_{(p)}}^1$ a free module over $E_{(p)}$ of dimension $2l$ and with $\psi_1\colon W_{E_{(p)}}^1\otimes W_{E_{(p)}}^1\to E_{(p)}$ a perfect alternating bilinear form. We take: $W_{E_{(p)}}=W_{E_{(p)}}^1\oplus W_{E_{(p)}}^1$ a direct sum of two copies of $W_{E_{(p)}}^1$. Let $\psi_0$ be an alternating form on it such that: $\psi_0(x,y)$ is $\psi_1(x,y)$ if $x,y$ belong to the same copy $W_{E_{(p)}}^1$ of $W_{E_{(p)}}$, and is equal to 0 otherwise. We take $h_2\circ h_1$ to be defined by: $A\in Sp(W_{E_{(p)}}^1,\psi_1)(E_{(p)})$ acts on $W_{E_{(p)}}$ as $A$ on each copy $W_{E_{(p)}}^1$. If $p$ does not divide $B(G^s)=6(l+1)$, then we take $\Gtil_{E_{(p)}}=G_{E_{(p)}}^d$, and if $p$ divides $B(G^s)$, then we take $\Gtil_{E_{(p)}}=Sp(W_{E_{(p)}},\psi_0)$ (as $p$ does not divide $B(\Gtil_E^{\ad})=6(l+2)$).

\smallskip
\item{d5)}
If $(G,X)$ is of $D_l^{\dbH}$ type, the situation is entirely analogous to the one described in d4) (we just have to replace alternating forms by symmetric bilinear forms), except that $h_0$ is not an isomorphism but an isogeny of degree 2. We have: $G_{E_{(p)}}^d$ is the split form of $SO(2l)_{E_{(p)}}$.

\smallskip
\item{e)}
If $(G,X)$ is of $A_l$, $B_l$ or $D_l^{\dbR}$ (resp. of $C_l$ or $D_l^{\dbH}$) type, then $G^d_{E_{(p)}}$ is the derived subgroup of the centralizer in $\tilde G_{E_{(p)}}$ of a torus $\tilde T$ of $\tilde G_{E_{(p)}}$ (resp. of a simple $E_{(p)}$-subalgebra $SA$ of $\End(W_{E_{(p)}})$ such that we have a relative PEL situation $(\tilde G_{E_{(p)}},SA)$), cf. d1) to d5).

\smallskip
The composition $h_2\circ h_1\circ h_0$ is the representation

\smallskip
-- in the case d1): direct sum of the representations associated to the fundamental weights corresponding to the roots $\alpha_1$ and $\alpha_l$ (see [De2] for the notations and the role of the roots; see also [Mi3, 1.21]);

\smallskip
-- in cases d2) and d3): direct sum of two copies of the Spin representation;

\smallskip
-- in  cases d4) and d5): direct sum of two copies of the representation associated to the fundamental weight corresponding to the root $\alpha_1$.
\finishproclaim
 
\Proclaim{6.6.5.1.} \rm
We now come back to i) to vi) of the proof of 6.5.1.1.
All the above part of 6.6.5 had just the role of making  6.5.1.1 i) well-fitted for the general case.  

We take $W_{(p)}:=W_{E_{(p)}}$. The group $\dbG_m(F)$ acts on $W_{(p)}\fracwithdelims[]1p$ by multiplication ($W_{E_{(p)}}$ is a module over $F_{(p)}$, cf. a)). We get the situation:
$$G_{3\dbZ_{(p)}}^{\der}\hookrightarrow \Res_{E_{(p)}/\dbZ_{(p)}}G_{E_{(p)}}^d\hookrightarrow \Gtil^0:=\Res_{E_{(p)}/\dbZ_{(p)}}\Gtil_{E_{(p)}}\hookrightarrow GL(W_{\dbZ_{(p)}}),$$
with $\Gtil^0(\dbZ_{(p)})=\Gtil_{E_{(p)}}(E_{(p)})$ acting on $W_{\dbZ_{(p)}}=K_{(p)}\otimes_{F_{(p)}} W_{(p)}$ through its canonical action on $W_{(p)}$. This is the explicit version of 6.5.1.1 i).

We keep ii) and iii) of 6.5.1.1. We have $n=2[E:F]$.

\smallskip
{\bf Case 1.}
We consider first the case when $\Sh(G,X)$ is a Shimura variety of $B_l$, $C_l$, $D_l^{\dbR}$ or $D_l^{\dbH}$ type, or of $A_l$ type but with trivial involution (cf. [De2, 2.3.12]). We choose $G_3$ as explained in [De2, 2.3.13] (i.e. we choose $G_2$ as small as allowed).
So the maximal torus of the center of ${G_3}_{\dbZ_{(p)}}$ commutes with $\Gtil^0$. This takes care of 6.5.1.1 iv). We keep  6.5.1.1 v).
The injective map $f\colon (\Gtil_3,\Xtil_3)\hookrightarrow (G\Sp(W,\psi),S^0)$ (we recall that $W=W_{\dbZ_{(p)}}\otimes\dbQ$) is a good embedding with respect to $p$, with $W_{\dbZ_{(p)}}$ a good $\dbZ_{(p)}$-lattice for the map $f$. This is a consequence of the fact that the family of tensors fixed by $\Gtil_3$ formed by the set of elements of the algebra $\Ltil$ of endomorphisms of $W_{\dbZ_{(p)}}$ fixed by $\Gtil_{3\dbZ_{(p)}}$ and by the family $\scrF$ of 3 tensors of degree 4 (as in 4.3.10 b) but for the embedding $\Gtil^0_{\dbQ}\hookrightarrow GL(W)$) is enveloped by $W_{\dbZ_{(p)}}$ and is $\dbZ_{(p)}$-very well positioned for $G_3$. To check this we use 4.3.6 2). Remark 4.3.13 takes care of the maximal torus of $Z(\Gtil_3)$, while 4.3.16 takes care of $\Gtil_3^{\der}$. To see this last part we just have to remark that (cf. 6.6.5 d) and e)):

\smallskip
-- if $(G,X)$ is of $A_l$, $B_l$ or $D_l^{\dbR}$ (resp. of $C_l$ or $D_l^{\dbH}$) type, then we have a relative PEL situation $(\Gtil^0,\Ltil,\Res_{E_{(p)}/\dbZ_{(p)}}\Ttil)$ (resp. relative PEL situations $(\Gtil^0,SAS)$ and $(\Gtil^0,\Ltil)$; here $SAS$ is $SA$ viewed as a $\dbZ_{(p)}$-algebra);

-- the family of tensors $\scrF$ is $\dbZ_{(p)}$-well positioned for the group $\Gtil^0_{\dbQ}$ and is enveloped by $W_{\dbZ_{(p)}}$, cf. 4.3.10 b) and 6.6.5 b) and d). 

\smallskip
This ends the explicit (second) proof of 6.4.2 as well as the proof of  6.5.1.1, in the case of the types listed above.

\smallskip
{\bf Case 2.}
We consider now the case when  $(G,X)$ is of $A_l$ type and has a non-trivial involution (as def. in [De1, 3.7]). We first remark that $\Res_{K_S/\dbQ} \dbG_m$ acts on $W_{(p)}\fracwithdelims[]1p$ (cf. the proof of [De2, 2.3.10]). We have to take some precautions: keeping 6.5.1.1 iv), the maximal torus $\scrG$ of the center of ${G_3}_{\dbZ_{(p)}}$ does not commute with $\Gtil^0$. However $\scrG_{\dbQ}$ is generated by two subtori: one is $\Res_{K/\dbQ} \dbG_m$ (it commutes with $\Gtil^0_{\dbQ}$), and another one which is a subtorus $T(K_S)$ of $\Res_{K_S/\dbQ} \dbG_m$ producing an isogeny $\Res_{F/\dbQ}\dbG_m\times T(K_S)\to\Res_{K_S/\dbQ} \dbG_m$ (cf. [De2, 2.3.10]). But $T(K_S)$ lies inside $\Gtil^0_{\dbQ}$ (cf. d1) above); in fact $T(K_S)$ is a subtorus of the generic fibre of $\Res_{E_{(p)}/\dbZ_{(p)}}\Ttil$ (cf. e) and d1) above). So keeping 6.5.1.1 v), we still get (the argument is the same as in Case 1 above) that the  map $f\colon (\Gtil_3,\Xtil_3)\hookrightarrow (G\Sp(W,\psi),S^0)$ is a good embedding with respect to $p$: again we have a relative PEL situation $(\Gtil^0,\Ltil,\Res_{E_{(p)}/\dbZ_{(p)}}\Ttil)$ (cf. 4.3.16). In other words the family of endomorphism of $W_{\dbZ_{(p)}}$ commuting with $\Gtil_{3\dbZ_{(p)}}$, together with the family $\scrF$ of three tensors (defined as in Case 1) is  $\dbZ_{(p)}$-very well positioned for $\Gtil_3$ and is enveloped by $W_{\dbZ_{(p)}}$ (cf. d1) above). This completes the explicit (second) proof of 6.4.2 as well as the proof of 6.5.1.1. 

\Proclaim{6.6.5.2. PEL type embeddings for the $A_l$ type.} \rm ${}^1$ $\vfootnote{1}{This section is added to this corrected version; it has been entirely incorporated in ``The Mumford--Tate Conjecture and Shimura Varieties, Part I,", math.NT/0212066.}$ 
We assume now that $(G,X)$ is of $A_l$ type. If $l=1$, it is easy to see that, replacing if needed $(\Gtil_3,\Xtil_3)$ by an enlargement (see def. 4.3.1) of it in $(GSp(W,\psi),S^0)$ (so we are not anymore interested to have $\Gtil_3$ as a subgroup of $G_3$), the injective map $(\Gtil_3,\Xtil_3)\hookrightarrow (GSp(W,\psi),S^0)$ is a PEL type embedding and the conditions of [Ko, Ch. 5] are satisfied for $p$ (i.e. we are in a situation as used in 4.3.11): we just need to choose $\tilde z\in M(pa)(\dbZ_{(p)})$ as mentioned in 6.5.1.1 v). We now check that a similar result holds for $l\Ge 2$.

As we took $W_{\dbZ_{(p)}}=K_{(p)}\otimes_{F_{(p)}} W_{(p)}$, the subgroup $G_{4\dbZ_{(p)}}$ of $GSp(W_{\dbZ_{(p)}},\psi)$ fixing all endomorphisms of $W_{\dbZ_{(p)}}$ fixed by $\tilde G_{3\dbZ_{(p)}}$, is reductive and has a derived subgroup which over $V_0$ is isomorphic (for $l\Ge 2$) to two copies of $G_{3V_0}^{\der}$ in such a way that the embedding of $G_{3V_0}^{\der}$ in $G^{\der}_{4\dbZ_{(p)}}$ is the diagonal embedding. Even if we replace $\tilde G_3$ by the smallest subgroup of it through which all homomorphisms $\Res_{\dbC/\dbR} \dbG_m\to G_{3\dbR}$ defining elements of $X_3$ factor, in general we can not ``get rid" of the second copy of $G_{3V_0}^{\der}$. 

There is a very simple way to adjust this so that we do get (with perhaps different notations) an injective map $\tilde f_3:(\tilde G_3,\tilde X_3)\hookrightarrow (GSp(W,\psi),S^0)$ such that $\tilde G_{3\dbZ_{(p)}}$ is the subgroup of $GSp(W_{\dbZ_{(p)}},\psi)$ fixing all endomorphisms of $W_{\dbZ_{(p)}}$ fixed by $\tilde G_{3\dbZ_{(p)}}$. It goes by: we entirely ``skip" the use of $K_{(p)}$ as follows. For the sake of uniformity, below we take $(G,X)$ of $A_l$ type, with $l\ge 1$.

We work with $W_{(p)}$ instead of $W_{\dbZ_{(p)}}$. If $l\Ge 2$ (resp. if $l=1$), then $K_S$ is a totally imaginary quadratic extension of $F$ (resp. is $F$) and so it makes sense to speak about ${K_S}_{(p)}$. If $l=1$, then $E$ is a totally imaginary quadratic extension of $F$. If $l\Ge 2$ (resp. $l=1$), let $GT$ be ${\Res}_{{K_S}_{(p)}/\dbZ_{(p)}} \dbG_m$ (resp. be ${\Res}_{E_{(p)}/\dbZ_{(p)}} \dbG_m$). It acts naturally in a  faithful way on $W_{(p)}$. Let $G_{4\dbZ_{(p)}}$ be the subgroup of $GL(W_{(p)})$ generated by $\tilde G_{\dbZ_{(p)}}$ and $GT$. Let $\tilde G_{4\dbZ_{(p)}}$ be the subgroup of $GL(W_{(p)})$ generated by $\tilde G_{\dbZ_{(p)}}^{\der}$, by $Z(GL(W_{(p)}))$ and by the maximal subtorus of $GT$ which over $\dbR$ is compact. It is reductive (cf. 3.1.6).

Let $RE$ be the set of real embeddings of $F$. For each $e_F\in RE$, let $V(e_F)$ be the maximal $\dbR$-vector subspace of $W_{(p)}\otimes_{\dbZ_{(p)}} \dbR$ on which the factor of $G_{4\dbR}^{\der}$ corresponding to $e_F$ acts non-trivially. So $GT_{\dbR}$ acts on it via its factor $F(e_F)$ which is a copy of $\Res_{\dbC/\dbR} \dbG_m$ and which is defined naturally by $e_F$; if $l\Ge 2$, the image $I(e_F)$ of $F(e_F)$ in $GL(V(e_F))$ is the center of the centralizer of the centralizer of the image of $\tilde G_{\dbR}^{\der}$ in $GL(V(e_F))$. We have a direct sum decomposition
$$W_{(p)}\otimes_{\dbZ_{(p)}} \dbR=\oplus_{e_F\in RE} V(e_F)$$
left invariant by $G_{4\dbR}$.
Let $x\in X_2$. We consider a monomorphism $h_x:{\Res}_{\dbC/\dbR} \dbG_m\hookrightarrow GL(W_{(p)}\otimes_{\dbZ_{(p)}} \dbR)$ having the properties:

\medskip
\item{{\bf 1)}} if $e_F$ is such that $G^s\times_{\Spec(F)} {}_{e_F}{\Spec}(\dbR)$ is non-compact, then the resulting homomorphism ${\Res}_{\dbC/\dbR} \dbG_m\to GL(V(e_F))$ is the one obtained by composing the homomorphism ${\Res}_{\dbC/\dbR} \dbG_m\to G_{2\dbR}$ defining $x$ with the natural homomorphism $G_{2\dbR}\to GL(V(e_F))$;

\item{{\bf 2)}} if $e_F$ is such that $G^s\times_{\Spec(F)} {}_{e_F}{\Spec}(\dbR)$ is compact, then the resulting homomorphism ${\Res}_{\dbC/\dbR} \dbG_m\to GL(V(e_F))$ is a monomorphism whose image is naturally identified with $I(e_F)$.

\medskip
So $h_x$ factors through $G_{4\dbR}$. The Hodge structure of $W_{(p)}\otimes_{\dbZ_{(p)}} \dbR$ it defines has type $\{(-1,0),(0,-1)\}$. So we can define two Shimura pairs $(G_4,X_4)$ and $(\tilde G_4,\tilde X_4)$ similar to $(G_3,X_3)$ and $(\tilde G_3,\tilde X_3)$. Moreover, taking perfect forms $\psi$ and $\tilde\psi$ on $W_{(p)}$ as mentioned in v) of 6.5.1.1 and above, the subgroup $G_{5\dbZ_{(p)}}$ of $GSp(W_{(p)},\psi)$ fixing all endomorphisms of $W_{(p)}$ fixed by $\tilde G_{4\dbZ_{(p)}}$, is reductive and has $G^{\der}_{4\dbZ_{(p)}}=\tilde G^{\der}_{4\dbZ_{(p)}}$ as its derived subgroup. Warning: we do not have to replace $W_{(p)}$ by a direct sum of two copies of itself (as this is implicitly done by 6.6.5 d1)). Let $G_5$ be the generic fibre of $G_{\dbZ_{(p)}}$. If $X_5$ is the $G_5(\dbR)$-conjugacy class of $h_x$ viewed as a homomorphism of $G_{5\dbR}$, then the pair $(G_5,X_5)$ is a Shimura pair whose adjoint is $(G,X)$. Moreover, we get a PEL type embedding
$$f_5:(G_5,X_5)\hookrightarrow (GSp(W_{(p)}[{1\over p}],\psi),S_0)$$
which has the desired property (i.e. we can take $\tilde f_3:=f_5$).
\finishproclaim

\Proclaim{6.6.5.2.1. Remark.} \rm
In 6.5.1.1 and Case 1 of 6.6.5.1, as above we can ``get rid" of $K$ for the $D_{2l+1}^{\dbR}$ type. But this is not true in general for the $B_l$, $C_l$, $D_{2l}^{\dbR}$ and $D_{l+2}^{\dbH}$ types. 
\finishproclaim

\Proclaim{6.6.6. Remark.} \rm
Except 6.5.1.1 vi), 6.6.3 and 6.6.4, everything else in 6.5 and 6.6 remains valid for $p=3$ (but working with $\frac{B(*)}3$ instead of $B(*)$; here $*$  substitutes a simple adjoint group of classical Lie type over a field). Even for $p=2$ some part of 6.5 and 6.6 remains valid. We will apply use remark in the construction of the $p=2$ and $p=3$ theories of Shimura varieties of preabelian type (cf. [Va5]).
\finishproclaim

\smallskip
\Proclaim{6.7. The proof of rm. 1) of 6.4.6.} \rm
See 2.11.1 for the meaning of $\scrU(*)$'s.
\finishproclaim

\Proclaim{6.7.1. Remark.} \rm
In 6.5.1.1 we can choose the number field $K$ and the Shimura pair $(G_1,X_1)$ such that $\scrU(G)\setminus\{2\}=\scrU(G_1)\setminus\{2\}$. This is a consequence of the proof of 6.5.1.1. Argument: $G_3^{\ab}$ is unramified over $\dbQ_l$ if $K$ and $K_S$ are unramified over $\dbQ_l$; moreover, if $G$ is unramified over $\dbQ_l$ then the number fields $K_S$ and $F$ are unramified over $l$ (cf. 6.5.1 b) and 6.5.1.1 iii)). So we just need  $K$ to be unramified over $l$ for all primes $l>2$ such that $G$ is unramified over $\dbQ_l$. For instance we can take $K=F(i)$. More generally: we can take $K=F(\sqrt{-d})$, where $d\in\dbN$ divides the discriminant of $F$. 

If there is a prime $l$ which mod $4$ is 2 or 3 or if there are two distinct primes $l$ such that $G_{\dbQ_l}$ is unramified, then we can choose $K$ and $(G_1,X_1)$ such that $\scrU(G)=\scrU(G_1)$. 

All these extend to the context of 6.4.2 (i.e. when $\Sh(G,X)$ is not a simple Shimura variety).
\finishproclaim

\Proclaim{6.7.2. Lemma.}
For any Shimura variety of Hodge type  $\Sh(G,X)$ there is an injective map $f\colon (G,X)\hookrightarrow (GSp(W,\psi),S)$ such that for any prime $l\in\scrU(G)$  there is a hyperspecial subgroup of $G(\dbQ_l)$ contained in a hyperspecial subgroup of $GSp(W,\psi)(\dbQ_l)$. 
\finishproclaim

\proof
We start with an arbitrary embedding $f\colon (G,X)\hookrightarrow (GSp(W,\psi),S)$. It takes care of all primes $l\in\scrU(G)\setminus\scrB(f)$, with $\scrB(f)\subset\scrU(G)$ a finite set. For any  $l\in\scrB(f)$ we choose arbitrarily a hyperspecial subgroup $H_l$ of $G(\dbQ_l)$. It is contained in a maximal compact open subgroup of $GSp(W,\psi)(\dbQ_l)$. But composing the natural map from $(G,X)$ to a Hodge quasi product (cf. Example 3 of 2.5) of $n$ copies of $(GSp(W,\psi),S)$, with $n\in\dbN$ big enough and suitable chosen, with the Segre embedding of this product into $(G_1,X_1):=(GSp(W^{\oplus n},\psi^{\oplus n}),S_n)$, we get that $H_l$ is contained in a hyperspecial subgroup of $G_1(\dbQ_l)$ (cf. the structure of maximal compact subgroups of $GSp(W,\psi)(\dbQ_l)$). 
The good values of $n$ depend only on the dimension of $W$ over $\dbQ$. So some $n\in\dbN$ works for all $l\in\scrB(f)$. 

In fact we can always take $n=2$: $H_l$ is contained in a hyperspecial subgroup of $GL(W)(\dbQ_l)$ (cf. 3.1.2.2 2)) and so 6.6.5 d1) applies. 

Now the injective map $(G,X)\hookrightarrow (G_1,X_1)$ has the desired property. This ends the proof of the Lemma.

\Proclaim{6.7.3.*} \rm
Now we are ready to prove rm. 1) of 6.4.6. We use the notations of 6.4.3 and 6.4.4. We assume that 6.4.1 and  6.4.2 are true for $p=3$ also (cf. 6.4.8) (otherwise we have to assume that $3\notin\scrS$). From 6.4.2.1 we deduce that we can assume that the open subgroup $H_{\scrS}$ of $G(\prod_{q\in\scrS}\dbQ_q)$ is as small as desired. This implies (cf. 6.4.5.1 and 3.2.3.1 5)) that we can assume that $(G,X)$ is of adjoint type. Remark 3.2.16 allows us to assume that $G$ is a simple $\dbQ$--group of adjoint type. From 6.7.1, 6.4.5.1, and C) of 3.2.2 4) (and 6.4.2.1) we deduce that we can assume that $(G,X)$ is of Hodge type. But this case is an easy consequence of 6.7.2 and 3.2.15: for $H_{\scrS}$ small enough we have a (special) (universal) principally polarized abelian scheme over $\scrM(H_{\scrS})$ (to be compared with 3.4.7 and 4.1). This ends the proof of rm. 1) of 6.4.6.
\finishproclaim

\Proclaim{6.8. About the proof of 6.1.2.} \rm
Here we present the proof of 6.1.2 as far as the tools presented in the present paper allow. For the last part of the non-compact case we have to refer either to [Va2] or to [Va3]. We keep the notations of 6.1. 
\finishproclaim

\Proclaim{6.8.0.} \rm
The part about triples implies and is implied by the part about quadruples. So we start using triples. For the case $p=3$ we refer to [Va2] or [Va3]. Here we consider $p>3$. From rm.  10) of 3.2.7 and 6.2.3 we deduce that we can assume that $f\colon (G_1,X_1,H_1)\to (G,X,H)$ is a cover. Moreover we can assume that $G_1^{\der}$ is a simply connected semisimple group. From rm. 11) of 3.2.7 we deduce that we can assume that $G^{\ad}$ is a simple $\dbQ$--group. 

We can assume that  $(G_1,X_1)$ is not of abelian type (cf. the proof of 6.4.1). So $(G_1,H_1)$ is of $D_l^{\dbH}$ type (cf. 6.4.2 and [De2, 2.3.10]). In particular the order of the center of $G_1^{\der}$ is a power of 2. From [De1, 2.4 and 2.5] and 3.2.8 we deduce that the connected components of $\Sh_{H_1}(G_1,X_1)_{\dbC}$ are defined over $K_0$ . As before $K_0$ is the field of fractions of $V_0=W(\dbF)$. 

Let $\scrN$ be the normalization  of $\scrM$ in the ring of fractions of $\Sh_{H_1}(G_1,X_1)$. It gets naturally a $G_1(\dbA_f^p)$-continuous action. So $\scrN$ is a quasi-projective  integral model of the triple $(G_1,X_1,H_1)$ (cf. 5.6.1 or 6.4.1 for the quasi-projectiveness part). Moreover it has the EEP. So we just need to show that 
it is a smooth integral model. For this it is enough to show that it is a pro-\'etale cover of the  open closed subscheme $\scrM^{\prime}$ of $\scrM$ defined as the image of $\scrN$ in $\scrM$. We can move over $V_0$, and so we come back to quadruples. From 6.2.3.1 we get:

\Proclaim{Fact.}
A connected component of $\Sh_H(G,X)_{K_0}$ is the quotient of a connected component of $\Sh_{H_1}(G_1,X_1)_{K_0}$ by a $16$-torsion pro-finite abelian group. 
\finishproclaim

\Proclaim{6.8.1. Lemma.} 
We assume that for any connected component $\scrC_{\dbF}$ of $\scrM^{\prime}_{\dbF}$ there is a $V_0$-valued point of $\scrN_{V_0}$ giving birth to an $\dbF$-valued point of $\scrN_{V_0}$ which is mapped into an $\dbF$-valued point of $\scrC_{\dbF}$. Then $\scrN$ is a pro-\'etale cover of $\scrM^{\prime}$.  
\finishproclaim

\proof 
Everything boils down (cf. the above Fact) in showing that: if $R=V_0[[x_1,...,x_d]]$ is a ring of formal power series in $d$  variables with coefficients in $V_0$, then there is no \'etale cover $Z$ of $\Spec(R\fracwithdelims[]1p)$ of degree 2, such that denoting by $R_1$ the normalization of $R$ in the field of fractions of $Z$, we do have a surjection $R_1\twoheadrightarrow V_0$ but $\Spec(R_1)$ is not an \'etale cover of $\Spec(R)$.

The proof of this is easy: $Z$ corresponds to a  field extension of the field of fractions of $R$ defined by an equation $x^2=z$, where $z$ is an invertible element of the unique factorization domain  $R\fracwithdelims[]1p$. As $R_1$ is not an \'etale cover of 
$R$, we deduce that we can assume that $z=pz_1$, with $z_1$ a unit of $R$. So we can not have surjections $\Spec(R_1)\twoheadrightarrow V_0$. This ends the proof of the Lemma.

\smallskip
In fact the result of the above proof remains true if we replace ``\'etale cover $Z$ of $\Spec(R\fracwithdelims[]1p)$ of degree 2" by: solvable Galois cover $Z$ of $\Spec(R\fracwithdelims[]1p)$ of degree relatively prime to $p$. Everything boils down to Kummer extensions, for which the above proof applies (to be compared with Step a) of 3.4.5.1). 

\Proclaim{6.8.2. Criteria.}
The hypothesis of 6.8.1 is satisfied if any one of the following conditions is satisfied:

a)  $\scrM$ admits smooth compactifications. 

b) The $\dbF$-valued points of $\scrM_{\dbF}$ obtained by specializing $K_0$-valued special points of $\scrM_{K_0}$ (cf. def. 2.10) are dense in $\scrM_{\dbF}$. 
\finishproclaim 

Criterion a) is a consequence  of 3.2.11 (which guarantees that $\scrN$ has plenty of $V_0$-valued points) and of 3.3.2. 
Criterion b) can be easily checked starting from [Mi4, 4.12], 2.7 and 2.8) (see [Va2]). 

In [Va3] we prove  a) (see 1.8), while in [Va2] we prove b) (cf. 1.6.1 and the density property referred to in 1.6.2). From 6.8.2 a) and 6.4.1.1 2) we get (without a reference to [Va3]) directly: 

\Proclaim{6.8.3. Corollary.}   
If $G_{\dbR}^{\ad}$ has compact factors, then 6.1.2 is true.
\finishproclaim

\Proclaim{6.8.4. Remark.} \rm
The condition 6.8.2 a) can be replaced by the condition that the connected components of $\scrM_{\overline{k(v)}}$ are permuted  transitively by $G(\dbA_f^p)$. This condition is satisfied (cf. 3.3.2) if there is an open subgroup $H_0\subset G(\dbA_f)$ such that $\scrM/H_0$ has smooth compactifications. 
\finishproclaim

From 6.8.2 a), 6.4.4, and the existence of smooth toroidal compactifications of $\Sh(G,X)$ (cf. [Har]), we get (without a reference to [Va3]):

\Proclaim{6.8.5. Fact.}
There is $\tilde N(G_1,X_1)\in\dbN$, depending only on the pair $(G_1,X_1)$, such that 6.1.2 is true if $p>\tilde N(G_1,X_1)$.
\finishproclaim

\Proclaim{6.8.6. The remaining cases.} \rm
From the above discussion we deduce that the cases of 6.1.2 which are not covered by 6.8.3 or by the abelian type situation and are needed for the complete proof of 6.1.2, can be summarized as follows. Keeping the notations of 6.1.2, we can assume (cf. also Example 5 of 2.5) that:

\smallskip
-- $(G^{\ad},X^{\ad})$ is a simple adjoint variety of $D_l^{\dbH}$ type ($l\in\dbN$, $l\ge 4$) such that $G^{\ad}_{\dbR}$ does not have compact factors.

\smallskip
We distinguish two cases: $(G^{\ad},X^{\ad})$ has a trivial or a non-trivial involution. If it has a trivial involution then $E(G^{\ad},X^{\ad})$ is a totally real number field, and we can assume that the embedding $f$ is a PEL type embedding (cf. Case 1 of 6.6.5 and [De2, 2.3.13]; the argument is the same as in Case 2 of 6.6.5). So we are reduced to the situation described in the case $D$ of [Ko, Ch. 5] (so $E(G,X)=E(G^{\ad},X^{\ad})$, cf. [De2, 2.3.13]; see also [Zi, p. 107]). If $E(G,X)=\dbQ$, it is an easy exercise to check that condition 6.8.2 b) is satisfied (Hint: use 1.6; in this case the results of the paragraph before 1.6.1 can be easily checked). 
\smallskip
If $(G^{\ad},X^{\ad})$ has a non-trivial involution, then $E(G^{\ad},X^{\ad})$ is a quadratic imaginary extension of a totally real number field and the situation still gets reduced to a PEL type situation. On the other hand, the ideas of 6.6.2 do not apply: with the notations of 6.6.2, if $(G,X)$ is of $D_l^{\dbH}$ type and has non-trivial involution, then $(G^1,X^1)$ is of $D_{al}^{\dbH}$ type and has as well non-trivial involution; here $a\in\dbN$. In particular 6.6.2 i) offers no simplification. So we do need, as mentioned above, either [Va2] or [Va3] to handle these two cases.  

\smallskip
\finishproclaim

\references{37}
{\nspace{

\Ref[BB]
W. Baily and A. Borel,
\sl Compactification of arithmetic quotients of bounded
symmetric domains,
\rm Ann. of Math. {\bf 84} (1966), p. 442--528.

\Ref[BHC]
A. Borel and Harish-Chandra,
\sl Arithmetic subgroups of algebraic groups,
\rm Ann. of Math. {\bf 75} (1962), 485--535.

\Ref[Bl]
D. Blasius,
\sl A $p$-adic property of Hodge cycles on abelian
varieties, 
\rm Proc. Symp. Pure Math. {\bf 55}, Part 2, p. 293--308.

\Ref[BLR]
S. Bosch, W. L\"utkebohmert, M. Raynaud,
\sl N\'eron models,
\rm Springer--Verlag, 1990.

\Ref[Bo]
A. Borel,
\sl Linear algebraic groups,
\rm Grad. Text Math. 126, Springer--Verlag, 1991.

\Ref[Bou1]
N. Bourbaki,
\sl Lie groups and Lie algebras, 
\rm Chapters 1-3, Springer--Verlag, 1989.

\Ref[Bou2]
N. Bourbaki,
\sl Groupes et alg\`ebres de Lie,
\rm Chapitre 7-8, Diffusion C.C.L.S., 1975.

\Ref[BT]
F. Bruhat and J. Tits,
\sl Groupes r\'eductifs sur un corps local,
\rm Publicationes Math., no. {\bf 60}, IHES, 1984.

\Ref[Ch]
C.-L. Chai,
\sl Every ordinary symplectic isogeny class in positive characteristic is dense in the moduli,
\rm Invent. Math. {\bf 121} (1995), no. 3, p. 439--479.

\Ref[De1]
P. Deligne,
\sl Travaux de Shimura,
\rm S\'eminaire  Bourbaki 389, LNM {\bf 244} (1971), p. 123--163.

\Ref[De2]
P. Deligne,
\sl Vari\'et\'es de Shimura: Interpr\'etation modulaire, et
techniques de construction de mod\`eles canoniques,
\rm Proc. Symp. Pure Math. {\bf 33}, Part 2 (1979), p. 247--290.

\Ref[De3]
P. Deligne,
\sl Hodge cycles on abelian varieties,
\rm Hodge cycles, motives, and Shimura varieties, LNM {\bf 900}, Springer--Verlag, 1982, p. 9--100.

\Ref[Fa1]
G. Faltings,
\sl Crystalline cohomology and $p$-adic  Galois
representations, Algebraic Analysis, Geometry, and Number
Theory,
\rm Johns Hopkins Univ. Press, 1990, p. 25--79.

\Ref[Fa2]
G. Faltings,
\sl The de Rham conjecture, Princeton notes correction to [Fa1], 4 pages.
\rm

\Ref[Fa3]
G. Faltings,
\sl Integral crystalline cohomology over very ramified
valuation rings,
\rm accepted for publication in  J. of A. M. S. (1998).

\Ref[Fa4]
G. Faltings,
Personal communication, June 1994.

\Ref[Fa5]
G. Faltings, 
\sl Hodge--Tate structures and Modular forms, 
\rm Math. Ann. {\bf 278}, 1987, p. 133--149.

\Ref[FC]
G. Faltings and C.L. Chai,
\sl Degeneration of abelian varieties,
\rm Springer--Verlag, 1990.

\Ref[Ha]
G. Harder,
\sl \"Uber die Galoiskohomologie halbeinfacher Matrizengruppen II,
\rm  Math. Z. {\bf 92} (1966), p. 396--415.

\Ref[Har]
M. Harris,
\sl Functorial properties of toroidal compactifications of locally symmetric varieties,
\rm Proc. London. Math. Soc. (3) {\bf 59} (1989) 1--22.

\Ref[Hart]
R. Hartshorne,
\sl Algebraic geometry,
\rm Grad. Text Math. 52, Springer--Verlag, 1977.

\Ref[He]
S. Helgason,
\sl Differential geometry, Lie groups, and symmetric spaces,
\rm Academic Press, New-York, 1978.

\Ref[Ja]
J. C. Jantzen,
\sl Representations of algebraic groups,
\rm Academic Press, 1987.

\Ref[dJO]
A. J. de Jong and F. Oort,
\sl On extending families of curves,
\rm J. Alg. Geom. {\bf 6} (1997), p. 545--562.

\Ref[Ko]
R. E. Kottwitz,
\sl Points on some Shimura Varieties over finite fields,
\rm Journal of the Am. Math, Soc., Vol. {\bf 5}, nr. 2, 1992, p. 373--444.

\Ref[La]
R. Langlands,
\sl Some contemporary problems with origin in the Jugendtraum,
\rm Mathematical developments arising from Hilbert's problems, Am. Math. Soc., Providence, RI, 1976, p. 401--418.

\Ref[LR]
R. Langlands and M. Rapoport,
\sl Shimuravarietaeten und Gerben, 
\rm J. Reine Angew. Math. {\bf 378} (1987), p. 113--220.

\Ref[Ma]
H. Matsumura,
\sl Commutative algebra,
\rm The Benjamin/Cummings Publishing Co., Inc., 1980.

\Ref[Me]
W. Messing,
\sl The crystals associated to Barsotti-Tate groups
with applicactions to abelian schemes,
\rm LNM {\bf 264}, Springer--Verlag, 1972.

\Ref[Mi1]
J. S. Milne,
\sl Canonical models of (mixed) Shimura varieties and automorphic vector bundles, 
\rm  Automorphic Forms, Shimura varieties and L-functions, vol I, Perspectives in Math., Vol. {\bf 10}, Acad. Press 1990.

\Ref[Mi2]
J. S. Milne,
\sl The action of an automorphism of $\dbC$ on a Shimura
variety and its special points,
\rm Prog. in Math., Vol. {\bf 35}, Birkh\"auser, Boston, 1983, p.
239--265.

\Ref[Mi3]
J. S. Milne,
\sl Shimura varieties and motives,
\rm Proc. Symp. Pure Math. {\bf 55} (1994), Part 2, p. 447--523.

\Ref[Mi4]
J. S. Milne,
\sl The points on a Shimura variety modulo a prime of good
reduction,
\rm The Zeta function of Picard modular surfaces, Les Publications CRM, Montreal 1992, p. 153--255.

\Ref[Mi5]
J. S. Milne,
\sl On the conjecture of Langlands and Rapoport,
\rm manuscript, september 1995.

\Ref[Mi6]
J. S. Milne,
\sl The conjecture of Langlands and Rapoport for Siegel modular varieties,
\rm Bulletin of the Am. Math. Soc., Vol. {\bf 24}, nr. 2, 1991, p. 335--341.

\Ref[MS]
J. S. Milne and K.-y. Shih,
\sl Conjugates of Shimura varieties,
\rm  Hodge cycles,
motives, and Shimura varieties, LNM {\bf 900}, Springer-Verlag, 1982, p. 280--356.

\Ref[Mu]
D. Mumford,
\sl Geometric invariant theory,
\rm Springer--Verlag, 1965.

\Ref[Mu1]
D. Mumford,
\sl Abelian varieties,
\rm Tata Inst. of Fund. Research, Oxford University Press, 1988.

\Ref[Na]
M. Nagata,
\sl Imbedding of an abstract variety in a complete variety,
\rm J. Math. Kyoto Univ., {\bf 2} (1962), p. 1--10.

\Ref[Pf]
M. Pfau,
\sl The reduction of connected Shimura varieties at a prime of good reduction,
\rm thesis, Univ. of Michigan, 1993.

\Ref[Ra]
M. Raynaud,
\sl Sch\'emas en groupes de type (p,...,p),
\rm Bull. Soc. Math. France, {\bf 102} (1974), p. 241--280.

\Ref[Sa]
I. Satake,
\sl Holomorphic imbeddings of symmetric domains into a Siegel space,
\rm  Am. J. Math. {\bf 87} (1965), p. 425--461.

\Ref[Se]
J.-P. Serre,
\sl Groupes alg\'ebriques associ\'es aux modules de Hodge-Tate,
\rm Journ\'ees de G\'eom\'etrie Alg\'ebrique de Rennes, part 3, J. Asterisques {\bf 65}, p. 155--188.

\Ref[SGA1]
\sl Rev\^etements \'etales et groupe fondamental,
\rm LNM {\bf 224}, 1971, Springer--Verlag.

\Ref[SGA3]
\sl Sch\'emas en groupes,
\rm LNM {\bf 152-3}, Vol. II-III, 1970, Springer--Verlag.

\Ref[Sh]
G. Shimura,
\sl On analytic families of polarized abelian varieties and automorphic functions,
\rm Ann. of Math. {\bf 78}, 1 (1963), p. 149--192.

\Ref[Ti]
J. Tits,
\sl Reductive groups over local fields, 
\rm Proc. Symp. Pure Math. {\bf 33}, part 1, p. 29--69.

\Ref[Va1]
A. Vasiu,
\sl Integral canonical models for Shimura varieties of Hodge type,
\rm thesis, Princeton University,  1994.

\Ref[Va2]
A. Vasiu,
\sl Points of the integral canonical models of Shimura varieties of preabelian type, $p$-divisible groups, and applications,
\rm to be published.

\Ref[Va3]
A. Vasiu,
\sl Toroidal compactifications of the integral canonical models of Shimura varieties of preabelian type, 
\rm to be published.

\Ref[Va4]
A. Vasiu,
\sl Examples of Shimura varieties and of their integral canonical models,
\rm preprint.

\Ref[Va5]
A. Vasiu,
\sl The $p=2$ and $p=3$ theory of Shimura varieties of preabelian type and the existence of the integral canonical models of the quotients of Shimura varieties of Hodge type with respect to non-hyperspecial subgroups,
\rm to be published.

\Ref[Va6]
A. Vasiu,
\sl Moduli schemes and the Shafarevich conjecture (the arithmetic case)  
for pseudo-polarized $K3$-surfaces,
\rm to be published.

\Ref[Vo]
P. Vojta,
\sl Nagata's embedding theorem,
\rm preprint (appendix to planned book), 1997.

\Ref[Za]
J. G. Zarhin,
\sl Endomorphisms of abelian varieties and points of finite order in characteristic $p$,
\rm Math. Notes of the Academy of Science of the USSR 21 (1977), 734-744.

\Ref[Zi]
T. Zink,
\sl Isogenieklassen von Punkten von Shimuramannigfaltigkeiten mit Werten in einem endlichen K\"orper,
\rm Math. Nachr. {\bf 112}, (1983), p. 103--124.

}}

\medskip
\hbox{address:}
\hbox{Adrian Vasiu}
\hbox{University of Arizona, Department of Mathematics}
\hbox{617 North Santa Rita, P.O. Box 210089}
\hbox{Tucson, AZ-85721, U.S.A.}
\hbox{login name: adrian\,\, e-mail address: math.arizona.edu}
\hbox{home page: www.arizona.edu/$\tilde {}$adrian}

\enddocument